\documentclass[3p, number, sort&compress, times, lefttitle]{elsarticle}

\journal{Computer Methods in Applied Mechanics and Engineering}

\usepackage[utf8]{inputenc}
\usepackage{gensymb}
\usepackage{hyphenat}
\usepackage{subcaption}
\usepackage{multirow}
\usepackage{graphicx}
\usepackage{bm}

\graphicspath{ {./figs/} {./figs/Meshes/} {./figs/RFunc/} {./figs/VennDiagram/} {./figs/GBC/} {./figs/Domains/} {./figs/eikonal/} {./figs/Heat2D/} {./figs/barycentric/} {./figs/kirchhoff/} {./figs/1D/} {./figs/Heat4D/} {./figs/Diagram_DNN/} }

\usepackage{tikz}
\usetikzlibrary{patterns,hobby}

\usepackage{amsmath,amsfonts}

\renewcommand{\Re}{\mathbb{R}}
\newcommand{\set}[1]{\mathbb{#1}}

\newcommand{\vm}[1]{{\bm{#1}}}
\newcommand{\vx}{\vm{x}}
\newcommand{\inn}{\vm{\nu}} 
\newcommand{\outn}{\vm{n}}  
\usepackage{mathtools}

\DeclareMathOperator*{\argmin}{arg\,min}
\newenvironment{proof}[1][Proof]{\begin{trivlist}
\item[\hskip \labelsep {\bfseries #1}.]}{\end{trivlist}}

\usepackage[colorlinks=true]{hyperref}

\usepackage{fancyvrb}
\usepackage{tgcursor}
\newcommand\boldverb[1]{\textbf{#1}}

\usepackage{xcolor}
\newcommand{\suku}[1]{{#1}}

\newcommand{\unn}{{\tilde{u}_{\textrm{nn}}}}
\newcommand{\unnx}{{\tilde{u}_{\textrm{nn}}(\vx)}}
\newcommand{\unnxt}{{\tilde{u}_{\textrm{nn}}(\vx;\vm{\theta})}}

\newcommand{\unnxs}{{\tilde{u}_{\textrm{nn}}(x)}}
\newcommand{\unnbcxs}{{\tilde{u}_{\textrm{nn}}^\textrm{bc}(x)}}
\newcommand{\unnxst}{{\tilde{u}_{\textrm{nn}}(x;\vm{\theta})}}
\newcommand{\unnbcxst}{{\tilde{u}_{\textrm{nn}}^\textrm{bc}(x ; \vm{\theta})}}

\newcommand{\unnbc}{{\tilde{u}_\textrm{nn}^\textrm{bc}}}
\newcommand{\unnbcx}{{\tilde{u}_\textrm{nn}^\textrm{bc}(\vx)}}
\newcommand{\unnbcxt}{{\tilde{u}_\textrm{nn}^\textrm{bc}(\vx;\vm{\theta})}}

\newcommand{\unnRxt}{{\tilde{u}_{\textrm{nn}}^R(\vx;\vm{\theta})}}

\newcommand{\unnRxst}{{\tilde{u}_{\textrm{nn}}^R(x;\vm{\theta})}}

\newcommand{\Lnnt}{{L_\textrm{nn}(\vm{\theta})}}
\newcommand{\Lnnbct}{{L_\textrm{nn}^\textrm{bc}(\vm{\theta})}}

\newcommand{\RELU}{{\texttt{ReLU}(x)}}
\newcommand{\REPU}{{\texttt{RePU}_n(x)}}

\newcommand{\calL}{{\cal L}}
\newcommand{\calN}{{\cal N}}

\newcommand{\fref}[1]{Fig.~\ref{#1}}
\newcommand{\sref}[1]{Section~\ref{#1}}

\begin{document}

\title{Exact imposition of boundary conditions with distance functions in physics-informed deep neural networks} 

\author[1]{N.\ Sukumar\corref{cor1}}
\ead{nsukumar@ucdavis.edu}

\author[2]{Ankit Srivastava}

\cortext[cor1]{Corresponding authors}

\address[1]{Department of Civil and Environmental Engineering, 
University of California, Davis, CA 95616, USA}

\address[2]{Department of Mechanical, Materials, and Aerospace Engineering,
            Illinois Institute of Technology, Chicago, IL 60616, USA}

\begin{abstract}
In this paper, we introduce a new approach based on 
distance fields to exactly impose
boundary conditions in physics-informed deep neural networks. 
The challenges in satisfying
Dirichlet boundary conditions in meshfree and particle methods are
well-known. This issue is also pertinent in the development 
of physics informed neural networks (PINN) for the solution of partial differential equations. We introduce {\em geometry-aware} trial functions in artifical neural networks to improve the training 
in deep learning for partial differential
equations. To this end, we use concepts from constructive solid geometry 
(R-functions) and generalized barycentric coordinates (mean value 
potential fields) to construct $\phi(\vx)$, an approximate  distance function to the boundary of a domain in $\Re^d$.  To exactly impose
homogeneous Dirichlet boundary conditions, the trial function is taken as $\phi(\vx)$ multiplied by the PINN approximation, and its generalization via transfinite interpolation is 
used to a priori satisfy inhomogeneous Dirichlet (essential), 
Neumann (natural), and Robin boundary conditions on complex geometries.
In doing so, we eliminate modeling error associated with the
satisfaction of boundary conditions 
in a collocation method and ensure that kinematic admissibility is met pointwise in a Ritz method.
With this new ansatz, the training for the neural network
is simplified: sole contribution to the loss function is from the residual error at interior 
collocation points where the governing equation is 
required to be satisfied. 
Numerical solutions are  computed using strong form collocation and
Ritz minimization.  To convey the main ideas and to assess
the accuracy
of the approach, we present numerical solutions for linear
and nonlinear boundary-value problems
over convex and nonconvex polygonal domains
as well as over domains with curved 
boundaries. Benchmark problems in one dimension for 
linear elasticity,
advection-diffusion, and beam bending; and in two dimensions 
for the steady-state
heat equation, Laplace equation, biharmonic equation
(Kirchhoff plate bending), and the nonlinear 
Eikonal equation are considered. 
The construction of approximate distance functions using R-functions extends to higher dimensions,
and we showcase its use by solving a Poisson problem with homogeneous Dirichlet boundary conditions over the four-dimensional 
hypercube. The proposed approach consistently outperforms a standard PINN-based collocation method, which underscores the importance of exactly (a priori) satisfying the boundary condition when constructing a loss function in PINN.  This study provides
a pathway for meshfree analysis to be conducted on the exact geometry without domain discretization. 
\end{abstract}


\maketitle

\section{Introduction}\label{sec:intro}
Machine learning algorithms based on supervised learning 
(deep neural networks) are relatively mature in fields such as 
computer vision, image processing, 
and speech recognition. Some of the earliest studies on physics-informed neural networks (PINN) to solve boundary-value 
problems can be traced to the contributions of
Lagaris et al.~\cite{Lagaris:1998:ANN,Lagaris:1997:ANN,Lagaris:2000:NNM}, and those of McFall~\cite{McFall:2006:ANN} and McFall and Mahan~\cite{McFall:2009:ANN}. These studies have provided the impetus for the recent interest and advancement of the approach. Over 
the past 3--4 years there have been many new
developments in a meshfree approach that is
based on PINN to solve
low-and high-order partial differential equations (PDEs). 
Some of the main contributions in this thread are based on: 
collocation~\cite{Raissi:2019:PIN,Berg:2018:UDA,Sirignano:2018:DGM},
variational principle (deep Ritz)~\cite{E:2018:DRM,Han:2018:SHD}, and Petrov-Galerkin
domain-decomposition~\cite{Kharazmi:2019:VPI,Kharazmi:2020:HPV}. Lu et al.~\cite{Lu:2020:DLL} present an overview of PINN for the solution of PDEs. 
These developments have been made possible due to the
advances in automatic differentiation tools that efficiently compute
the derivatives of nonlinear composite functions, and in stochastic gradient descent algorithms that can deliver accurate
solutions for nonlinear, nonconvex optimization problems. 
Furthermore, the availability of public-domain data analysis 
packages such as \texttt{Tensorflow}~\cite{Abadi:2016:TEF} and 
\texttt{PyTorch}~\cite{Paszke:2016:PYT} has also decreased the barriers to entry for newcomers to this field.
Within the purview of solving PDEs with neural networks, it is
instructive to highlight that unlike well-known computational 
methods (finite elements and its generalizations, finite volume, boundary
elements, meshfree, and others), where function approximation is based on
linear combinations of basis functions, which have to be chosen and
defined a priori, function approximation in deep neural networks is through nonlinear function composition of an activation function
$\sigma: \Re^d \to \Re$, which yields the {\em best} \suku{approximation function} 
via the solution procedure. 
A popular choice for $\sigma$
is: $\sigma(x) := \RELU = \max(0,x)$,
which is known as the Rectified Linear Unit (ReLU) 
activation function.
The collocation method using PINN is based on minimizing the least squares residual error (nonlinear and nonconvex mean squared error loss function) in order to satisfy the PDE and the boundary conditions at collocation points---whose solution yields the parameters that define the approximation function in the artificial neural network. 

In PINN-based \suku{deep collocation~\cite{Raissi:2019:PIN} and 
deep Ritz~\cite{E:2018:DRM}} methods, 
the inexact imposition of boundary conditions adversely affects the training of the neural network as well as the accuracy of the method~\cite{Wang:2020:UMG,Chen:2020:CSD,Lyu:2020:EEB}. On complex geometries, this shortcoming becomes more acute. This
issue has plagued meshfree Galerkin methods~\cite{Babuska:2003:SMG,Huerta:2018:MM} since their
inception, and a foolproof solution is still
unavailable for arbitrary two- and three-dimensional geometries. 
Meshfree 
basis functions are also nonpolynomial, and hence an additional
hurdle is a suitable \suku{cell-based numerical integration scheme 
that is consistent
(patch test is passed) and stable (absence of spurious
modes) for nonlinear
simulations.}
Since solving for the \suku{approximating function} 
is part of the minimization procedure in PINN, 
\suku{a simple (e.g., Monte Carlo) integration 
scheme for the entire domain can be used.}
Hence, if a reliable approach to 
impose boundary conditions on complex geometries in PINN is realized, it can lead to an 
accurate and robust meshfree method. In this paper, we tackle this problem by proposing 
a new approach in PINN to exactly impose boundary conditions, with
an eye towards enabling meshfree simulations over complex geometries in
$\Re^d$ ($d = 2,3$). Our main contributions are as follows:
\begin{enumerate}
\item We introduce a new {\em geometry-aware} method in 
      physics-informed neural networks
      that uses R-functions and transfinite interpolation to
      exactly impose boundary conditions over complex (affine,
      curved and multiply-connected) geometries. This geometry-aware approach, which is based on the construction of approximation
      distance functions (ADFs) to boundary sets,
      was first proposed by Kantorovich to satisfy Dirichlet boundary conditions~\cite{Kantorovich:1958:AMH}, and has been extended using the theory of R-functions 
      to exactly impose Dirichlet, Neumann and Robin boundary conditions in a meshfree method by initially Rvachev~\cite{Rvachev:1982:TRF}, and then subsequently by
      Rvachev and coworkers~\cite{Rvachev:1995:RBV,Rvachev:2000:OCR,Rvachev:2001:TII} and Shapiro and coworkers~\cite{Shapiro:1991:TRF,Shapiro:1999:MSD,Shapiro:2002:ASA,Biswas:2004:ADF,Shapiro:2007:SAG,Freytag:2011:FEA}. A related meshfree approach in the spirit of Kantorovich's method is that of 
      H{\"o}llig et al.~\cite{Hollig:2001:WEB}, who used web-splines that are multiplied by a weight function (approximates the distance function) to solve a Dirichlet problem.  We also draw from a previous study~\cite{Millan:2015:CBM}, where R-functions are used to define smooth approximate distance fields over polygonal domains.
\item Approximate distance fields
      that stem from mean value potentials~\cite{Floater:2003:MVC,Dyken:2009:TMV,Belyaev:2013:SLD} 
      are also used in PINN to solve PDEs.
\item On exactly satisfying boundary conditions in physics-informed
      neural networks, the training of the network is simplified, and
      this facilitates convergence and improved accuracy of the PINN
      approximation.
\item New application of R-functions in $\Re^4$: solution
      of the Poisson equation over the four-dimensional
      hypercube.
\item Solving PDEs on curved 
domains without discretizing the domain (mesh generation) is realized, which provides a pathway to conducting 
meshfree simulations on the exact geometry (isogeometric analysis)~\cite{Hughes:2005:IAC}.
\end{enumerate}

First, we present a few connections of prior work on finite elements and partition-of-unity meshfree methods to
better understand and place the PINN approximation and its use to 
solve PDEs. Contributions in $r$-adaptivity
with finite elements have a long history in
computational mechanics. For instance, in 
finite-deformation simulations using finite elements, 
the optimal nodal locations and the solution coefficients have both been
simultaneously treated as unknowns in the minimization of the potential energy functional~\cite{Thoutireddy:2004:VRA}. 
Since the PINN approximation that is composed by the ReLU 
activation function 
can exactly represent piecewise affine functions (Delaunay
basis functions)~\cite{He:2018:RDN},
one can view the ReLU network solution as a variational 
$r$-adaptive finite element solution procedure. 
Instead of refining elements in $h$-adaptive finite elements, 
adaptive solutions can be realized via a basis refinement strategy that has advantages (for example, ``hanging nodes'' are a nonissue), which was put forth by Grinspun~\cite{Grinspun:2003:BRM}, and a similar \suku{approximation}
refinement perspective can be associated with a multilayer
neural network~\cite{Cyr:2020:RTI}. The connections between
ReLU networks and $hp$-finite elements are studied in
Opschoor et al.~\cite{Opschoor:2020:DRN}.
Initiated by Kansa~\cite{Kansa:1990a:MUQ,Kansa:1990b:MUQ}, 
meshfree collocation schemes with positive-definite
radial basis functions (RBFs) such as the Gaussian and Hardy's inverse multiquadrics have been used to solve PDEs~\cite{Buhmann:2003:RBF,Fasshauer:2007:MAM}.
Schaback and Wendland~\cite{Schaback:2006:KTM} discuss the ties of kernel methods (for example, radial basis functions) to
machine learning and meshfree methods.
The choice on how to set the ``shape parameter'' in the Gaussian RBF (controls the support-width of the
Gaussian) is an unsettled issue since it is \suku{problem dependent} and
must be carefully selected for boundary and interior points if
exponential convergence is to be maintained without exacerbating the
condition number. A related approach is the {\em local}
maximum-entropy (max-ent) meshfree method~\cite{Arroyo:2006:LME} 
that yield compactly-supported 
basis functions of exponential form that constitute
a partition-of-unity~\cite{Babuska:1997:PUM}, possess 
linear completeness, and
provide a smooth transition from Delaunay basis functions~\cite{Rajan:1994:ODT} to
global maximum-entropy basis functions~\cite{Sukumar:2004:COP}. Consider
a nodal set with nodes (centers) that are located at $\{\vx_i\}_{i=1}^n$. When
viewed through the lens of Gaussian
weight functions~\cite{Sukumar:2005:MEA,Arroyo:2007:LME,Sukumar:2007:OAC}, a 
single parameter $\{\beta_i\}_{i=1}^n$ controls the support-width of 
each Gaussian weight function. 
Rosolen et al.~\cite{Rosolen:2010:OSS} proposed a variational
adaptivity formulation to find optimal values of $\{\beta_i\}_{i=1}^n$ 
that minimize
the potential energy functional for Poisson equation and nonlinear
elasticity. Since RBFs can be represented in a neural network
with a single hidden layer~\cite{Park:1991:UAU,Mhaskar:1996:NNO},
the neural network solution optimizes the location of the centers
$\{\vx_i\}_{i=1}^n$ as well as the 
support-widths $\{\beta_i\}_{i=1}^n$~\cite{Park:1991:UAU}. RBF-based partition-of-unity networks~\cite{Lee:2020:PUN} for $hp$-approximation have been introduced, and numerical experiments have been conducted using
sparse Gaussian networks to solve PDEs~\cite{Ramabathiran:2021b:SPI}.
Lastly, 
we mention the recent work 
of Greco and Arroyo~\cite{Greco:2020:HOM}, who presented a collocation scheme for PDEs that is 
based on high-order max-ent approximants, which
delivered accurate simulation results on
domains with affine and
curved boundaries.
In~\sref{subsubsec:RBFExample}, we present a one-dimensional
example using Gaussian neural networks.
Since PINN affords significant 
flexibility vis-{\`a}-vis existing
meshfree (basis-set) methods, solving PDEs over complex geometries
using collocation and Ritz methods with artificial neural
networks holds significant promise.

\suku{In a strong collocation PINN method,
the loss function consists of the residual errors from the
interior of the domain, which is known as {\em interior or PDE loss}, 
and from 
the boundaries of the domain, which is referred to as
{\em boundary (conditions) loss}~\cite{Raissi:2019:PIN}}. 
There 
are three distinct
contributions in the mean squared error loss function:
(1) residual error at interior collocation points where 
the PDE
    must be satisfied;
(2) residual error at boundary collocation points where the 
    essential (Dirichlet) boundary condition must be satisfied; and
(3) residual error at boundary collocation points where the 
    Robin or natural (Neumann) boundary condition must be satisfied.
Early approaches~\cite{Lagaris:1998:ANN,Lagaris:2000:NNM,McFall:2006:ANN,McFall:2009:ANN} had already recognized the importance of exact imposition of boundary conditions in artificial neural networks. Lagaris et al.~\cite{Lagaris:1998:ANN} considered two terms in the trial function, with the first term being an analytical function that exactly imposed the boundary conditions and the second term was chosen as the product of the PINN approximation and
a function that vanished on the boundary; for irregular boundaries, Lagaris et al.~\cite{Lagaris:2000:NNM} used a RBF network in the first term to satisfy the boundary conditions at a collection of discrete points on the boundary.  McFall et al.~\cite{McFall:2006:ANN,McFall:2009:ANN} and more
recently Sheng and Yang~\cite{Sheng:2021:PFN} introduced a length factor (measure of the distance to the boundary) associated with the boundary to impose boundary conditions, and 
\suku{Berg and Nystr{\"o}m~\cite{Berg:2018:UDA} approximated the distance 
function using a low-capacity neural
network to impose boundary conditions over complex geometries.}
In many recent studies~\cite{Dwivedi:2020:PIE,Wang:2020:UMG,Chen:2020:CSD,Lyu:2020:EEB}, the implications of
imposing essential boundary conditions via the
loss function in PINN have been studied, and numerical 
experiments have
affirmed that the presence of the boundary residual terms 
compromises the convergence of the 
\suku{stochastic gradient descent} algorithm and the accuracy
of the method. To address this problem, remedies have been 
introduced in the PINN literature, such as using two neural
networks, one 
for the PDE and the other to satisfy the essential 
boundary condition~\cite{Lagaris:2000:NNM,Berg:2018:UDA,Sheng:2021:PFN,Dwivedi:2020:PIE,Dwivedi:2020:SBE},
introduction of a penalty parameter via an augmented 
variational formulation
to weakly impose the essential boundary 
conditions~\cite{E:2018:DRM}, and Nitsche's method to impose the essential boundary condition~\cite{Liao:2019:DNM}. Some of these
approaches mirror those previously pursued in meshfree and
particle methods to satisfy essential boundary 
conditions~\cite{Babuska:2003:SMG,Li:2004:MPP,Huerta:2018:MM}. 
In meshfree Galerkin methods,
the choice of the space for Lagrange multipliers is delicate; penalty
approach leads to a saddle-point problem and the \suku{inf-sup} condition must be met; and though Nitsche's method is variationally consistent, the stabilization parameter in it must be judiciously chosen.
For low-dimensional problems over complex geometries, an accurate and robust meshfree approach remains elusive.
Since from the universal approximation theorem~\cite{Hornik:1989:MFN,Hornik:1991:ACM}
we know that a neural network with one hidden layer can represent 
any $L^2$ function to arbitrary accuracy, it stands to choose
an ansatz that satisfies the boundary conditions a priori, so that the
loss function is expressed solely in terms of the
residual error at only the interior collocation points where the PDE is
required to be satisfied. If the essential boundary conditions are exactly met, then this precludes ``variational crimes'' in a Ritz method~\cite{Strang:1973:AFE}. Lastly, and most importantly, 
\suku{in deep collocation~\cite{Raissi:2019:PIN},
multiple terms (interior loss and boundary losses) that have to be individually minimized are incorporated within a single objective (loss) function.} When this loss
function is minimized,
then the solution that is realized depends on the weight (equal weights is the \suku{unbiased choice}) that is assigned to each objective function, which reflects the importance of each residual error contribution. 
\suku{Rohrhofer et al.~\cite{Rohrhofer:2021:PFP} discuss network training in relation to the Pareto front that appears in multiobjective constrained optimization. In the NVIDIA SimNet\texttrademark\ toolkit~\cite{Hennigh:2021:NVI},
signed distance function weighting is used to dynamically assign the spatially varying weight functions for the PDE and boundary loss terms.
Since these weights are problem dependent, they should not be fixed a 
priori, 
since then the magnitude of the training
loss by itself is a misleading error measure.}
To establish the accuracy of the method, the error in $u$ as well as in $u'$ must be assessed.
We present numerical results in 
Sections~\ref{sec:1D}--\ref{sec:4D} that support this thesis, and which points to the merits of the new approach. 

In this paper, we solve this problem of competing loss terms in 
PINN formulations by 
constructing a trial function for the neural network that
a priori satisfies all boundary conditions in deep collocation, and meets
kinematic admissibility when used in a deep Ritz method.
This eliminates the boundary terms in the loss function.
Our approach is based 
on constructing distance functions (exact or approximate) to the boundary of the domain, and it can treat essential (Dirichlet)
as well as mixed (Dirichlet and Robin) boundary conditions over complex domains. We use the exact distance function whenever it is available and applicable.  However, in general, we 
construct approximate distance functions using two different techniques: the theory of 
R-functions~\cite{Rvachev:1995:RBV,Shapiro:2007:SAG}) 
and the theory of mean value potential fields~\cite{Floater:2003:MVC,Dyken:2009:TMV,Belyaev:2013:SLD}. These 
methods provide approximate distance functions that possess the desirable property of being zero on the boundary of the domain 
with unit (inward) normal directional derivative. In addition, they are
smooth in the interior of the domain, a property that 
the exact distance function does not always possess. Functions whose sign solely depend on the sign of its arguments encode
Boolean logic, and are known as R-functions.
R-functions provide an implicit function representation for line segments, curves, and solid regions, and are composed by
Boolean operations
(negation, conjunction, disjunction, equivalence).  Mean value
potential fields are specific forms of a singular
double-layer potential that yield
$L_p$-distance fields~\cite{Belyaev:2013:SLD}. For a 
domain in $\Re^2$, this singular
potential is defined as the integral of the reciprocal of the
$p$-th power of the distance from its boundary.
For a polygon with $p = 1$, closed-form expressions for the 
ADF are available~\cite{Floater:2003:MVC}, 
but for closed curves in $\Re^2$, numerical integration is required
to compute the ADF. Once the approximate
distance functions are formed, methods to
impose essential and Robin boundary conditions are available
that rely on transfinite interpolation~\cite{Rvachev:1995:RBV,Rvachev:2001:TII}.
R-functions with approximants such as B-splines~\cite{Shapiro:1999:MSD,Freytag:2011:FEA} and
RBFs~\cite{Tsukanov:2011:HME} have been used 
in a meshfree Galerkin method to solve boundary-value 
problems. 

The remainder of this paper is organized as follows. 
In~\sref{sec:distance}, we discuss the properties of 
distance functions, and in~\sref{sec:RFADF},
the essentials on R-functions and the construction of approximate distance functions are described. In particular, joining R-functions using R-equivalence composition is presented, which is used in this paper.  
The inverse of the normalizing weight function that appears in the expression for mean value coordinates (polygon) and transfinite mean value interpolant (closed curves), is an approximate distance field.
These are particular instances ($p = 1$) of $L_p$-distance fields,
and are discussed in~\sref{sec:mvp}. On using normalizing functions and
solution structures in the R-function method~\cite{Rvachev:1982:TRF,Rvachev:1995:RBV}, 
we describe in~\sref{sec:BCs} 
the use of ADFs to
construct an ansatz in PINN that exactly satisfies boundary
conditions for second- and fourth-order problems. The construction of 
the trial function in a deep neural network is presented 
in~\sref{sec:dnn}, along
with a summary of the feedforward neural network and 
\suku{the backpropagation (computation of the gradient of the loss function
and use of 
stochastic gradient descent) algorithm}. The loss function for  collocation and deep Ritz formulations are presented in~\sref{sec:formulations}. The numerical implementation is discussed 
in~\sref{sec:implementation}, where we also 
provide code snippets of some of the main functions.
One- and two-dimensional numerical simulations are presented 
in Sections~\ref{sec:1D} and~\ref{sec:2D}, where we apply this new approach to a broad suite of boundary-value problems (Poisson, harmonic coordinates, plate bending, Eikonal equation) on convex and nonconvex domains with affine and curved boundaries. In addition, the
Poisson equation over the four-dimensional hypercube is solved
in~\sref{sec:4D}. These numerical results clearly demonstrate the 
benefits of exactly imposing the boundary conditions in PINN---it simplifies the training of the network and enhances the \suku{accuracy and robustness}
of the method.
Finally, we conclude with~\sref{sec:conclusions}, where we summarize
the main developments in this paper and discuss some of the topics of future research.

\section{Distance Functions and their Properties} \label{sec:distance} 
The signed distance function is an
implicit representation for curves and surfaces, and also
provides fast evaluation of predicates for geometric objects. Let $S \subset \Re^d$ denote a \suku{domain} (open, bounded set) with
boundary $\partial S$. The exact distance function $d(\vx)$ gives the shortest distance between any point $\vx\in\Re^d$ to $\partial S$. It is clear that $d(\vx)$ is identically zero on $\partial S$. Computing the exact distance function requires solving the Eikonal 
equation~(see \sref{subsec:Eikonal}), which is computationally expensive. Therefore, it is desirable to construct an approximate distance 
function or ADF $\bigl($formally represented by $\phi(\vx) \bigr)$ 
that has a closed-form (non-iterative algorithm) expression. Furthermore, since the exact distance function may only be continuous and not continuously differentiable, it may not be suitable for use in a trial function to solve PDEs. Since our objective is to use ADFs in a collocation or Ritz method to solve boundary-value problems, their differential properties are important. If essential boundary conditions are imposed on the entire boundary $\partial S$, then the ADF must be zero on $\partial S$, positive in $S$, and its gradient must not vanish for any $\vx \in \partial S$.  In addition, since the exact distance function has derivative discontinuities on the medial axis of the domain, smooth approximations of the distance function must be used within the trial function for a collocation method with PINN. For a second-order problem, a $C^0$ distance function that has gradient discontinuities in the interior of the domain cannot be used in the collocation approach since the Laplacian of the distance function will be unbounded at a collocation point.
These considerations are crucial 
when used to solve
PDEs. For instance, positivity in $S$ precludes the presence of singularities within the 
domain, which in general is difficult to construct as noted in McFall~\cite{McFall:2006:ANN}.

For clarity of exposition in this paper, we use $\inn := \inn(\vx)$ to denote the unit inward normal vector (appears in the theory related to R-functions) on the boundary $\partial S$, and $\outn := \outn(\vx)$ as the unit outward normal vector (used when defining Neumann or Robin boundary condition) on $\partial S$. It is noted 
that $\outn = - \inn$. 
If $\partial S$ in $\Re^2$ 
is composed of piecewise line segments and curves, then
we use $\phi_i := \phi_i(\vx)$ to denote the ADF to each
curve or line segment.
For a point $\vx \in \Re^d$ on $\partial S$, it is essential that
any approximation to the distance
function satisfy $\phi =0$.
Furthermore, to mimic the exact distance function, the normal
derivative with respect to $\nu$ on the boundary should be unity, $\partial d / \partial \nu = \nabla d \cdot \vm{\nu} = 1$, and
it is desirable that all higher order normal derivatives vanish. 
An $m$-th
order approximate distance function requires that the second- to
$m$-th order normal derivatives
vanish on all regular points (unit normal is well-defined)
on $\partial S $~\cite{Rvachev:1995:RBV}:
\begin{equation}\label{eq:normalization}
\dfrac{\partial \phi}{\partial \nu} = 1, \quad
\dfrac{\partial^k \phi}{\partial \nu^k} = 0 \quad (k = 2,3,\ldots,m),
\end{equation}
and such a function is said to be \emph{normalized} to the $m$-th
order.  For finite $m$,
the normalized function matches the exact distance function only
in the vicinity of the boundary; for points that are away from the
boundary, it deviates from the
exact distance.  Apart from applications in solid modeling, 
mesh generation, real-time rendering and computer vision, where
distance functions
are used, normalized first-order distance functions are also
a suitable choice for the initialization and assignment of the extension velocity at points away from the interface in the level set method~\cite{Sethian:1999:LSM}. 
As noted in Biswas and Shapiro~\cite{Biswas:2004:ADF},
use of normalized distance functions mitigate the bulging
phenomenon in the vicinity of where the segments or curves are 
joined~\cite{Bloomenthal:1997:BEC}, since undulations (presence of local extrema) are
undesirable in the representation of the surface.

\section{R-functions and Approximate Distance Functions}
\label{sec:RFADF}
The theory of R-functions can be used to construct
a composite approximate distance function, $\phi(\vx)$, to any arbitrarily complex
boundary $\partial S$, when approximate
distance functions, $\phi_i(\vx)$, to the partitions
of $\partial S$ are known. Consider a real-valued function $F(\omega_1,\omega_2,\ldots,\omega_q)$, where 
$\omega_i(\vx): \Re^d \rightarrow \Re$
($i=1,\ldots,q$) are also real-valued functions.  If the sign of
$F(\cdot)$ is solely determined by the signs of its arguments
$\omega_i(\vx)$, then $F(\cdot)$ is known as an  
R-function~\cite{Shapiro:1991:TRF,Shapiro:2007:SAG}. R-functions were proposed
by T. L. Rvachev in 1963~\cite{Shapiro:1991:TRF}.
For example, $F_1(x,y) = 1 + x^2 + y^2$ and $F_2(x,y) = xy$ are R-functions in $\Re^2$, whereas 
$F_3(x,y) = \sqrt{x^2 + y^2} - 1 $ and $F_4(x,y) = \sin xy$ are not.
The important properties of R-functions are provided 
in Rvachev and Sheiko~\cite{Rvachev:1995:RBV} and 
Shapiro~\cite{Shapiro:1991:TRF}.  On combining set-theoretic Boolean operations with
such functions, the inverse problems of semi-analytic geometry (solid modeling) can be solved.

\begin{figure}
    \centering
    \begin{subfigure}{0.24\textwidth}
    \includegraphics[width=\textwidth]{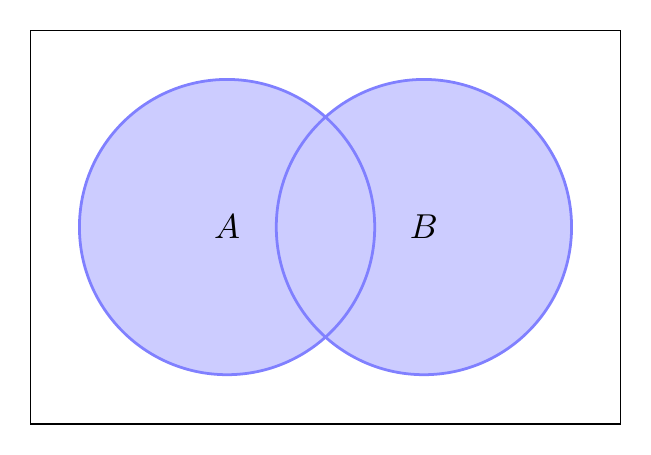}
    \caption{} 
    \end{subfigure}\hfill
    \begin{subfigure}{0.24\textwidth}
     \includegraphics[width=\textwidth]{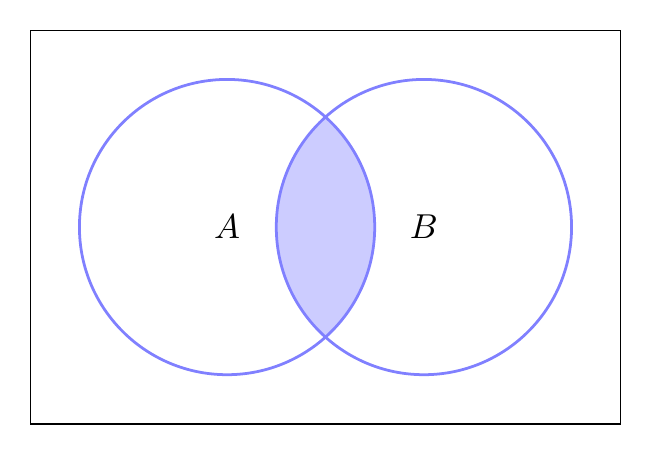}
    \caption{} 
    \end{subfigure}\hfill
    \begin{subfigure}{0.24\textwidth}
    \includegraphics[width=\textwidth]{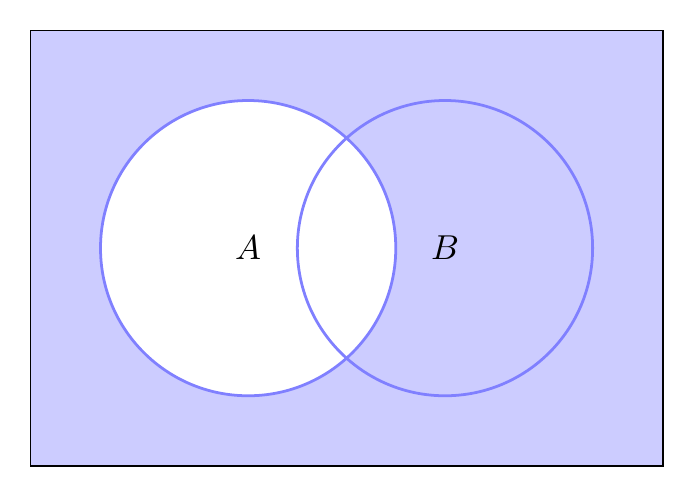}
    \caption{} 
    \end{subfigure}\hfill
    \begin{subfigure}{0.24\textwidth}
    \includegraphics[width=\textwidth]{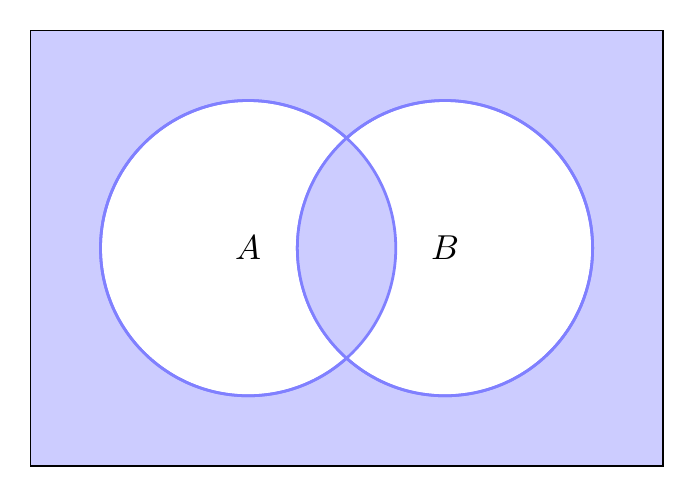}
    \caption{} 
    \end{subfigure}
    \caption{Venn diagram for union, intersection, 
    complement and equivalence in $\Re^2$. $\omega_A$ and
    $\omega_B$ are R-functions that are positive in
    $\Omega_A$ and $\Omega_B$ (open sets), respectively. 
    (a) $A \cup B \equiv \omega_A \vee \omega_B$; 
    (b) $A \cap B \equiv \omega_A \wedge \omega_B$; 
    (c) $\bar{A} \equiv - \omega_A$; and 
    (d) $(A \cap B) \cup (\bar{A} \cap \bar{B})
        \equiv \omega_A \sim \omega_B$. Examples of
        $\vee$ and $\wedge$ operations using R-functions
        are given in~\eqref{eq:R_alpha} and~\eqref{eq:R_s}.
     }
    \label{fig:venn}
\end{figure}

Consider a continuous 
function $\omega_i : \Re^d \to \Re$. Let
$\Omega \subset \Re^d$ be an open, bounded domain, 
$\bar{\Omega} = \Omega \cup \partial \Omega$ be the
closure of $\Omega$, and define
$\Omega^c$ to be the complement of 
$\Omega$ ($\Omega \cup \Omega^c = \Re^d$). 
If $\omega_i$ is strictly positive in $\Omega$,
identically equal to zero on
$\partial \Omega$, and strictly negative in
$\Omega^c \backslash \partial \Omega$, then it is evident that
$F(\omega_i) = \omega_i$ is an R-function.
Over
the region $\bar{\Omega}$, we
associate $\omega_i$ with the Boolean 1 (logical true)
and over the region $\Omega^c$ we associate it with
the Boolean 0 (logical false). 
Note 
that $\omega_i = 0$ ($0$ is assumed to be signed) is
included in both sets so that it can be assigned to either
the set of negative real values or the set of positive real
values~\cite{Shapiro:2007:SAG}. 
Hence, similar to Boolean functions, $\omega_i$ is
closed under composition. Furthermore, 
just as Boolean functions are written using the symbols $\neg$
$\vee$, and $\wedge$, which correspond to complement, union
and intersection in set theory, every 
R-function can be written as the composition of the corresponding
elementary R-functions: R-negation ($ -\omega$), R-disjunction ($\omega_1 \vee \omega_2$), and
R-conjunction ($\omega_1 \wedge \omega_2$).  On defining
R-functions for regions in $\Re^d$, a solid can then be composed using the set-theoretic operations of $\neg$, $\vee$, and $\wedge$.
For the universal set $\mathbb{U} = \Re^2$, Venn diagrams for some of the operations in set theory are shown in~\fref{fig:venn}, and the corresponding operations using R-functions are indicated.

The simplest examples of R-functions are the
R-disjunction (union) and the R-conjunction (intersection) functions. These are
\begin{equation*}
\omega_1 \vee \omega_2 = \frac{\omega_1 + \omega_2}{2}
+ \frac{\sqrt{(\omega_1 - \omega_2)^2}}{2} = \max(\omega_1,\omega_2), 
\quad
\omega_1 \wedge \omega_2 = \frac{\omega_1 + \omega_2}{2} -
\frac{\sqrt{(\omega_1 - \omega_2)^2}}{2} = \min(\omega_1,\omega_2) .
\end{equation*}
The generalization of the above R-functions
is~\cite{Shapiro:2007:SAG}: 
\begin{equation}\label{eq:R_alpha}
R_\alpha(\omega_1,\omega_2) := 
\dfrac{1}{1+\alpha} \left( \omega_1 + \omega_2 \pm
    \sqrt{\omega_1^2 + \omega_2^2 - 2\alpha \omega_1 \omega_2 } \right),
\end{equation}
with $(+)$ and $(-)$ signs defining R-disjunction and
R-conjunction, respectively. 
If $\omega_1$ and $\omega_2$ denote the sides of a triangle, then
the triangle inequality is expressed in~\eqref{eq:R_alpha} with $-1 < \alpha < 1$
being the cosine of the angle between the two sides. 
For $\alpha =  1$, the \texttt{max} and 
\texttt{min} R-functions are recovered.   If $\omega_1$ and
$\omega_2$ are positive, then so are $\omega_1 \vee \omega_2 $ and
$\omega_1 \wedge \omega_2 $.  The R-functions defined in~\eqref{eq:R_alpha}
are not analytic at points where $\omega_1 = \omega_2  =
0$. Smoothness can be obtained by defining the function ($\alpha =
0$ is selected)~\cite{Shapiro:2007:SAG}
\begin{equation}\label{eq:R_s}
R_s(\omega_1,\omega_2) := \left[ \omega_1 + \omega_2 \pm \sqrt{\omega_1^2 + \omega_2^2} \right] 
      \left(\omega_1^2 + \omega_2^2 \right)^{\frac{s}{2}},
\end{equation}
which renders these functions to be $C^s$-continuous at all points other than where $\omega_1 = \omega_2 = 0$.

\subsection{Normalized functions for line segments and curves}\label{subsec:norm_seg_curve}
Shapiro and Tsukanov~\cite{Shapiro:1999:IFG} describe the representation of line segments and curves using
R-functions and discuss their differential properties.
Let us consider one line segment that joins 
$\vx_1 := (x_1,y_1)$ and
$\vx_2 := (x_2,y_2)$. The center of this segment is denoted by
$\vx_c := (\vx_1 + \vx_2)/2$, and the length of the segment is: $L = ||\vx_2 -
\vx_1||$.  Now, we define~\cite{Rvachev:2001:TII}
\begin{equation}\label{eq:dline}
f := f(\vx) = \dfrac{(x - x_1)(y_2 - y_1) - (y-y_1)(x_2-x_1)}{L},
\end{equation}
which is the signed distance function from point $\vx$ to the line that
passes through $\vx_1$ and $\vx_2$.

Since the representation of the segment can be viewed as
the intersection of an infinite line with a disk of radius $L/2$, we
consider the following \emph{trimming} function that is normalized
to  first order~\cite{Rvachev:2001:TII}:
\begin{equation}\label{eq:trim}
t := t(\vx) = \dfrac{1}{L} 
\left[ \left(\frac{L}{2}\right)^2 - ||\vx - \vx_c||^2 \right],
\end{equation}
where $t \ge 0$ defines a disk with center at $\vx_c$. Now, with $f(\vx)$ and
$t(\vx)$ on-hand, we define a normalized function (up to first order)
$\phi(\vx)$ that is $C^2$ at all points 
away from the line segment~\cite{Rvachev:2001:TII,Biswas:2004:ADF}:
\begin{equation}\label{eq:phiconj}
\phi := \phi(\vx) = \sqrt{f^2 + \left(
\frac{\varphi - t}{2} \right)^2}, \quad
\varphi = \sqrt{t^2 + f^4},
\end{equation}
which is an approximation of the distance function to the
segment with end points $\vx_1$ and $\vx_2$.
The function $\phi$ in~\eqref{eq:phiconj} is a modification of the form $\varphi = |t|$~\cite{Shapiro:1999:IFG}, which 
has a derivative discontinuity at $t = 0$. 

\begin{figure}
\centering
\mbox{
\includegraphics[width=0.325\textwidth]{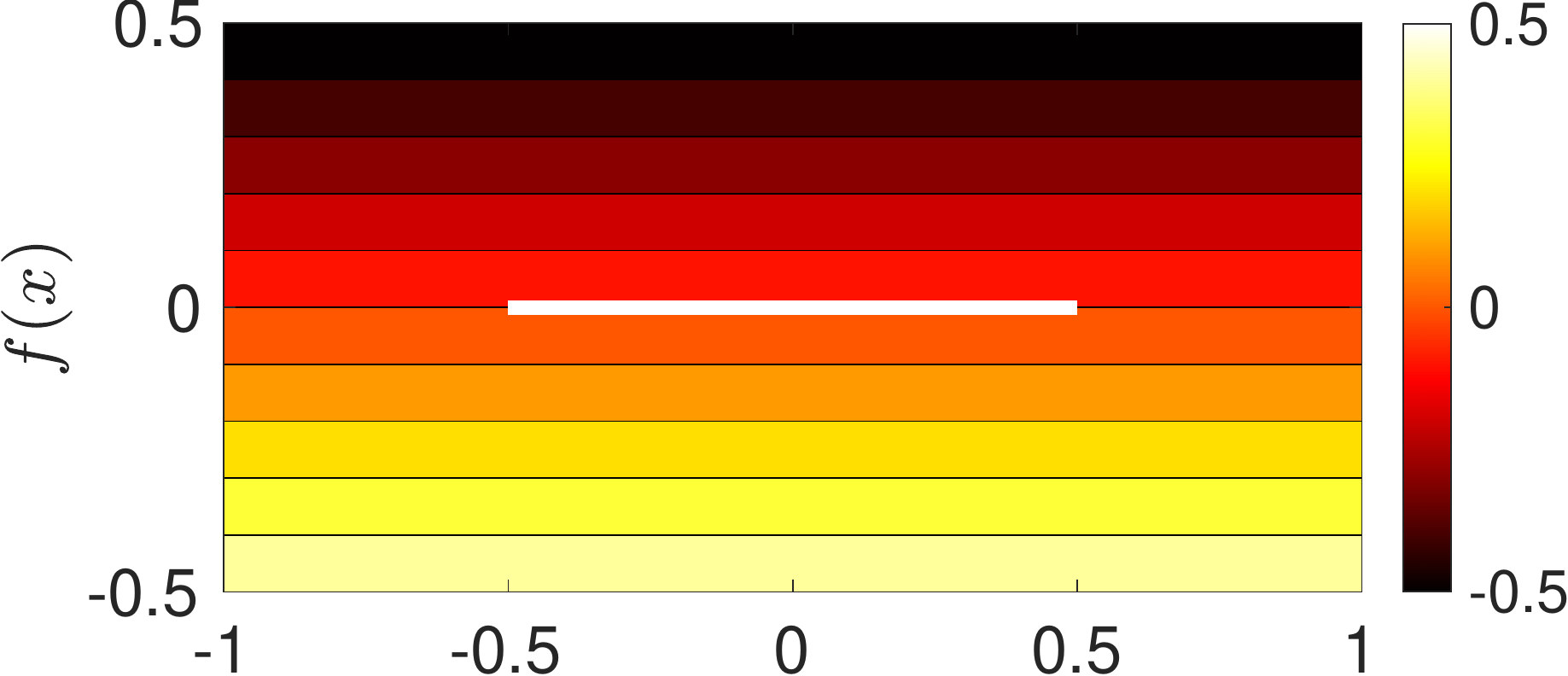}
\includegraphics[width=0.325\textwidth]{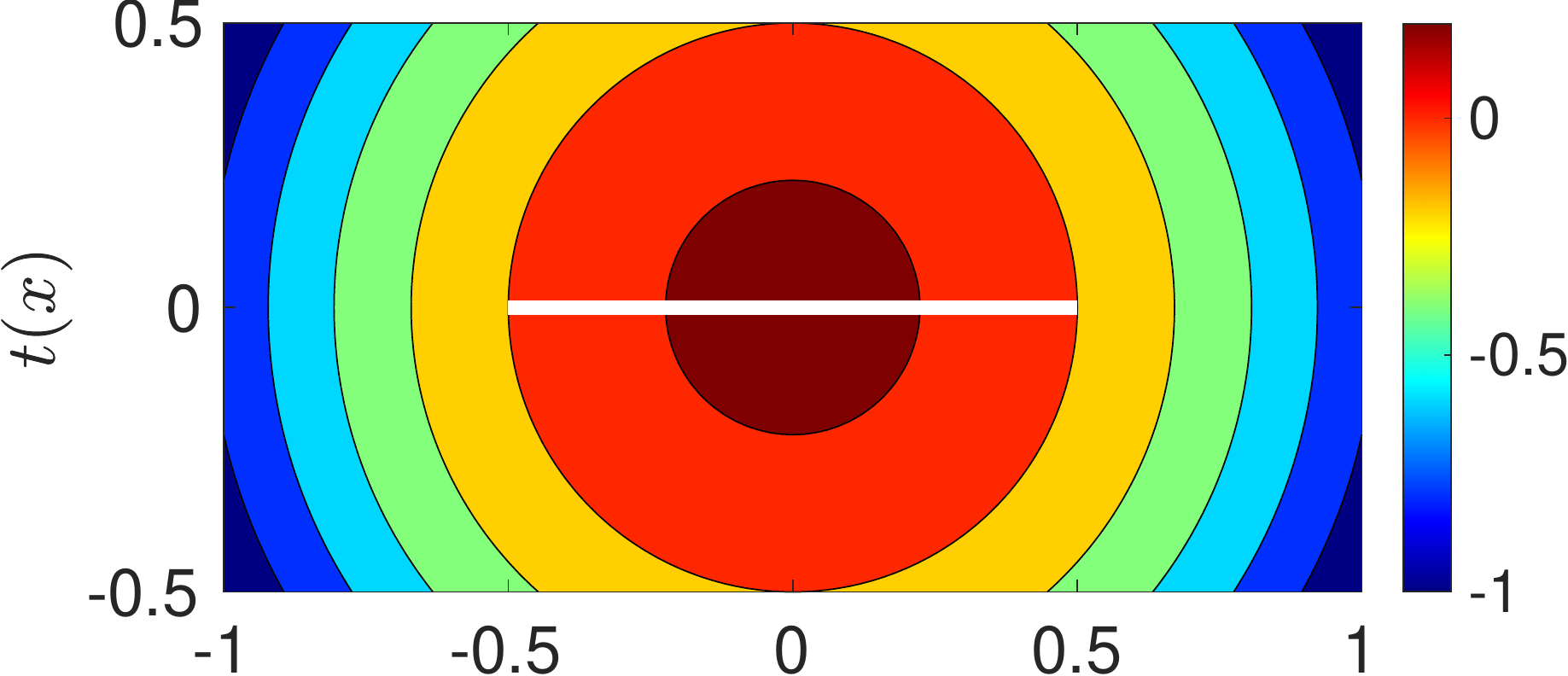}
\includegraphics[width=0.325\textwidth]{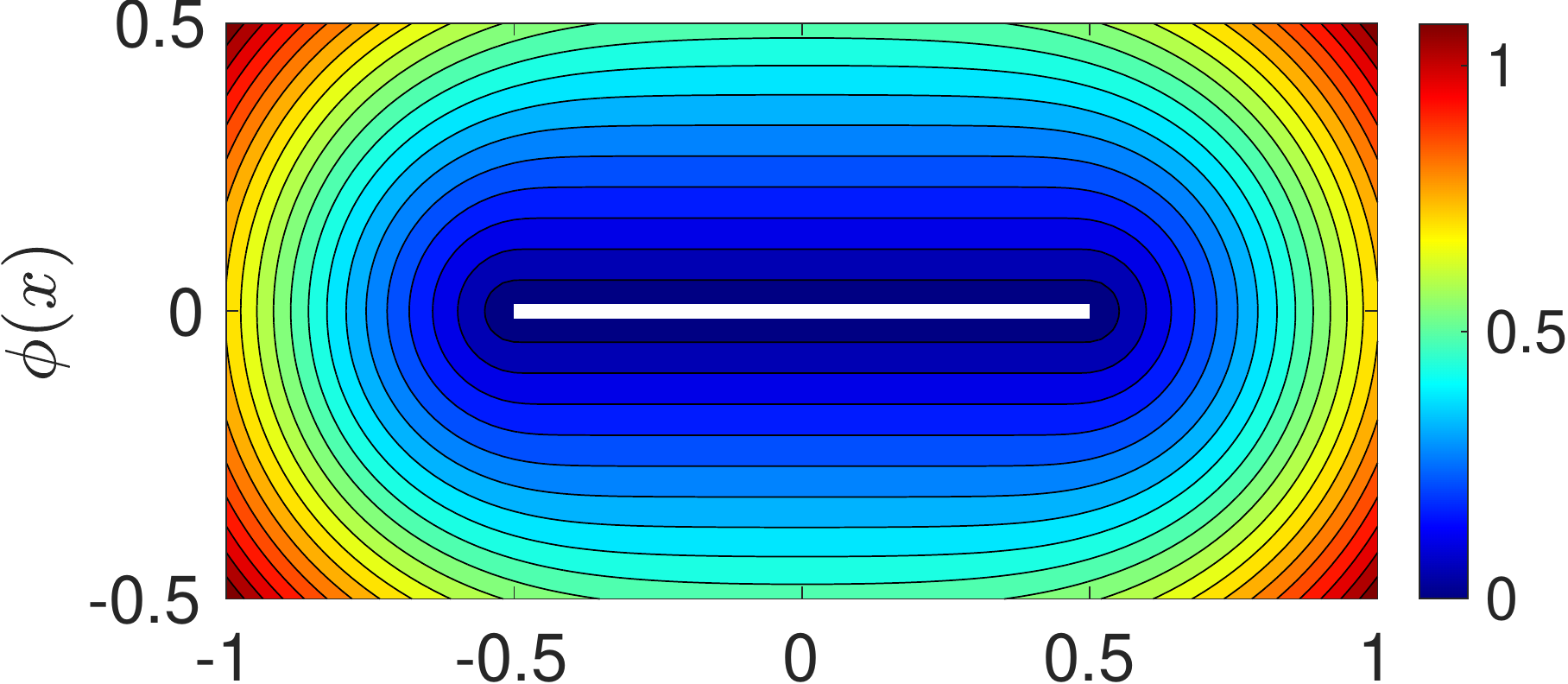}
}
\caption{
Construction of the approximate distance function to a line segment. The leftmost plot depicts the signed distance function~\eqref{eq:dline} to a line in $\Re^2$; the middle plot shows the trimming function~\eqref{eq:trim}; and the rightmost plot displays the approximate distance function~\eqref{eq:phiconj} to a line segment.}
\label{fig:phi_seg}
\end{figure}

\begin{figure}[!htb]
\centering
\mbox{
\includegraphics[width=0.325\textwidth]{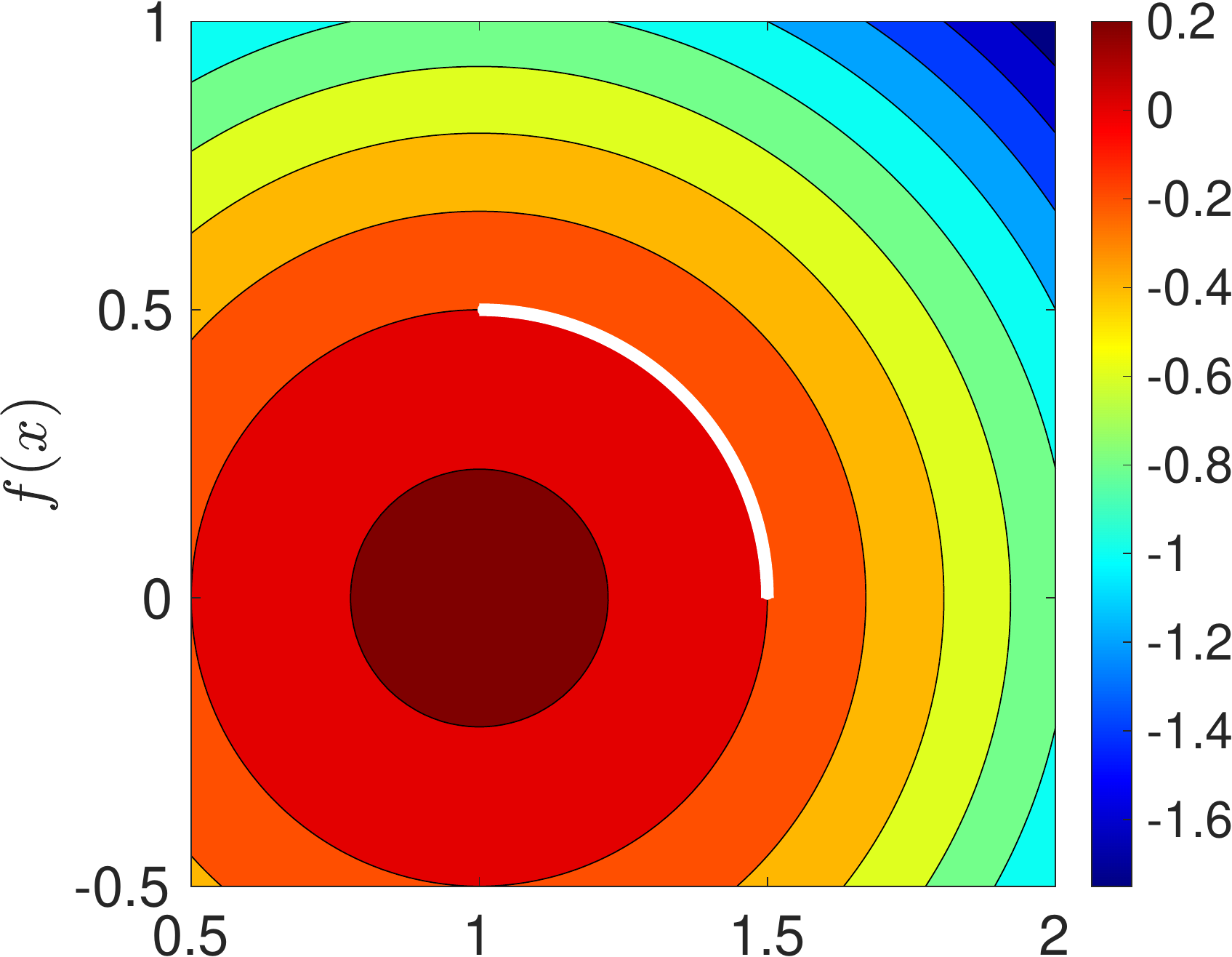}\hspace*{0.02in}
\includegraphics[width=0.325\textwidth]{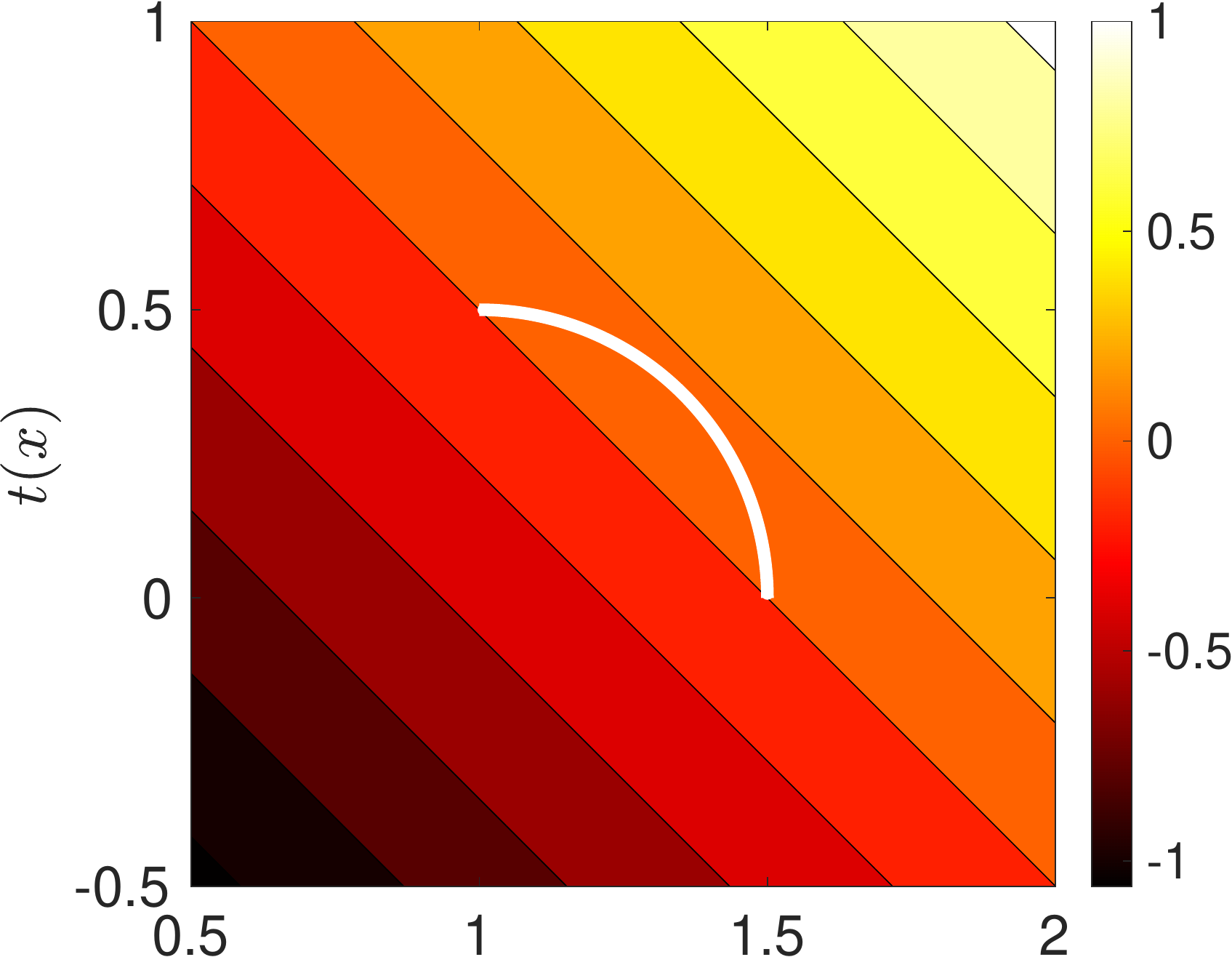}
\includegraphics[width=0.325\textwidth]{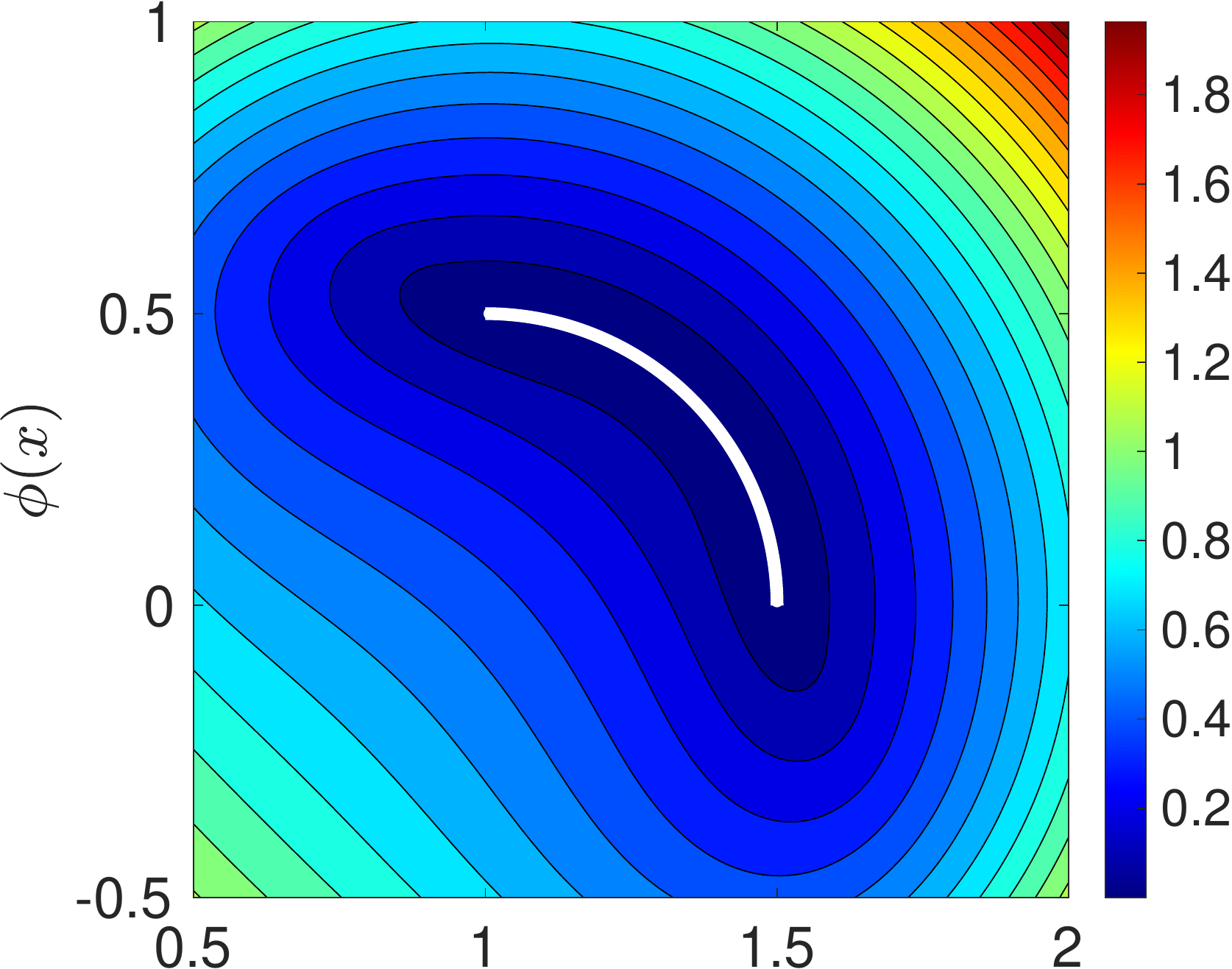}
}
\caption{Construction of the approximate distance function to a quarter-circular arc. The leftmost plot depicts $f$, the 
equation of the circular arc (normalized to first order); the trimming function
is shown in the middle plot; and the rightmost plot displays the approximate distance function given by~\eqref{eq:phiconj}.}
\label{fig:phi_arc}
\end{figure}

\begin{figure}[!htb]
\centering
\begin{subfigure}{0.49\textwidth}
\includegraphics[width=\textwidth]{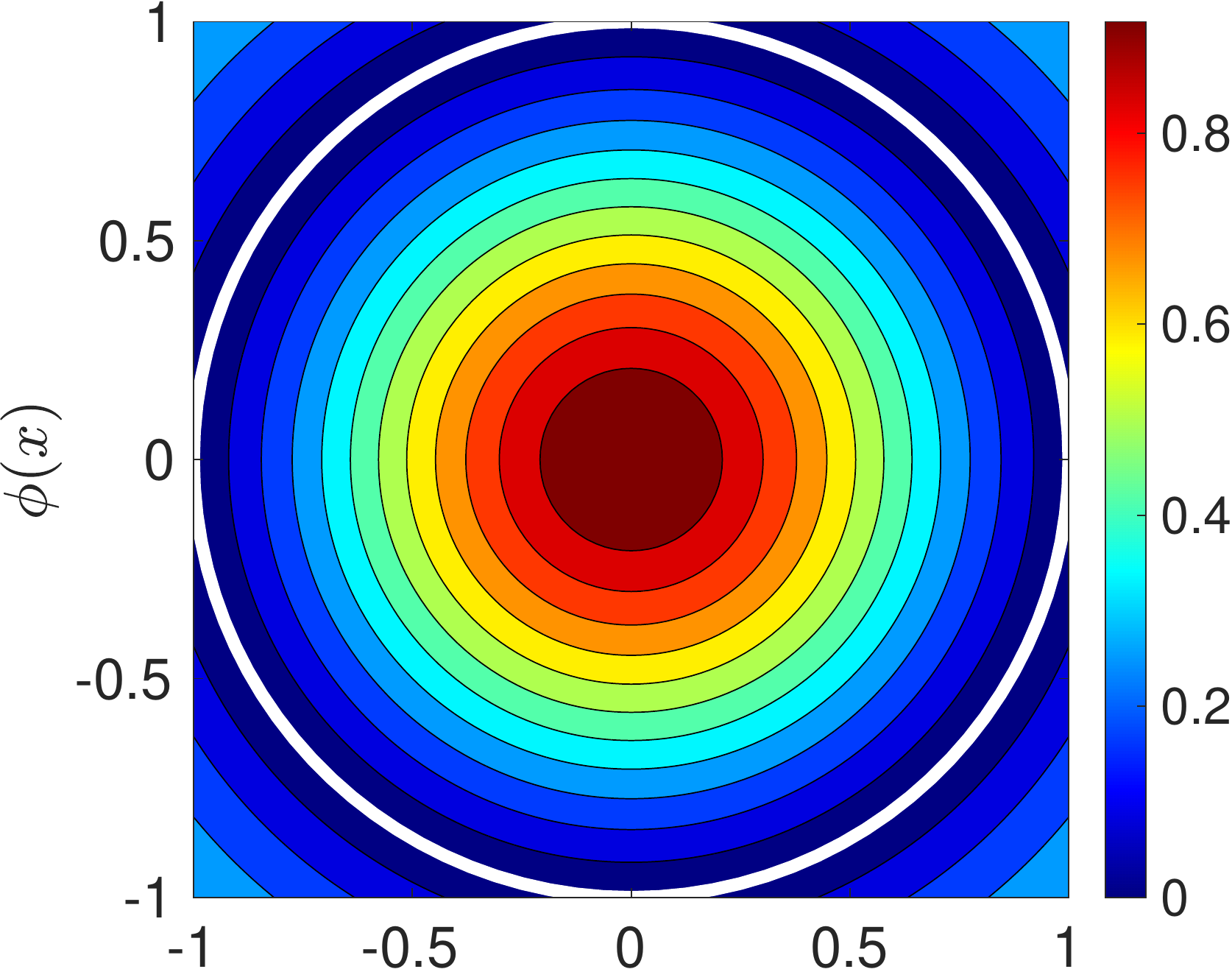}
\caption{}
\label{fig:phi_circle_ellipse-a}
\end{subfigure}
\begin{subfigure}{0.49\textwidth}
\includegraphics[width=\textwidth]{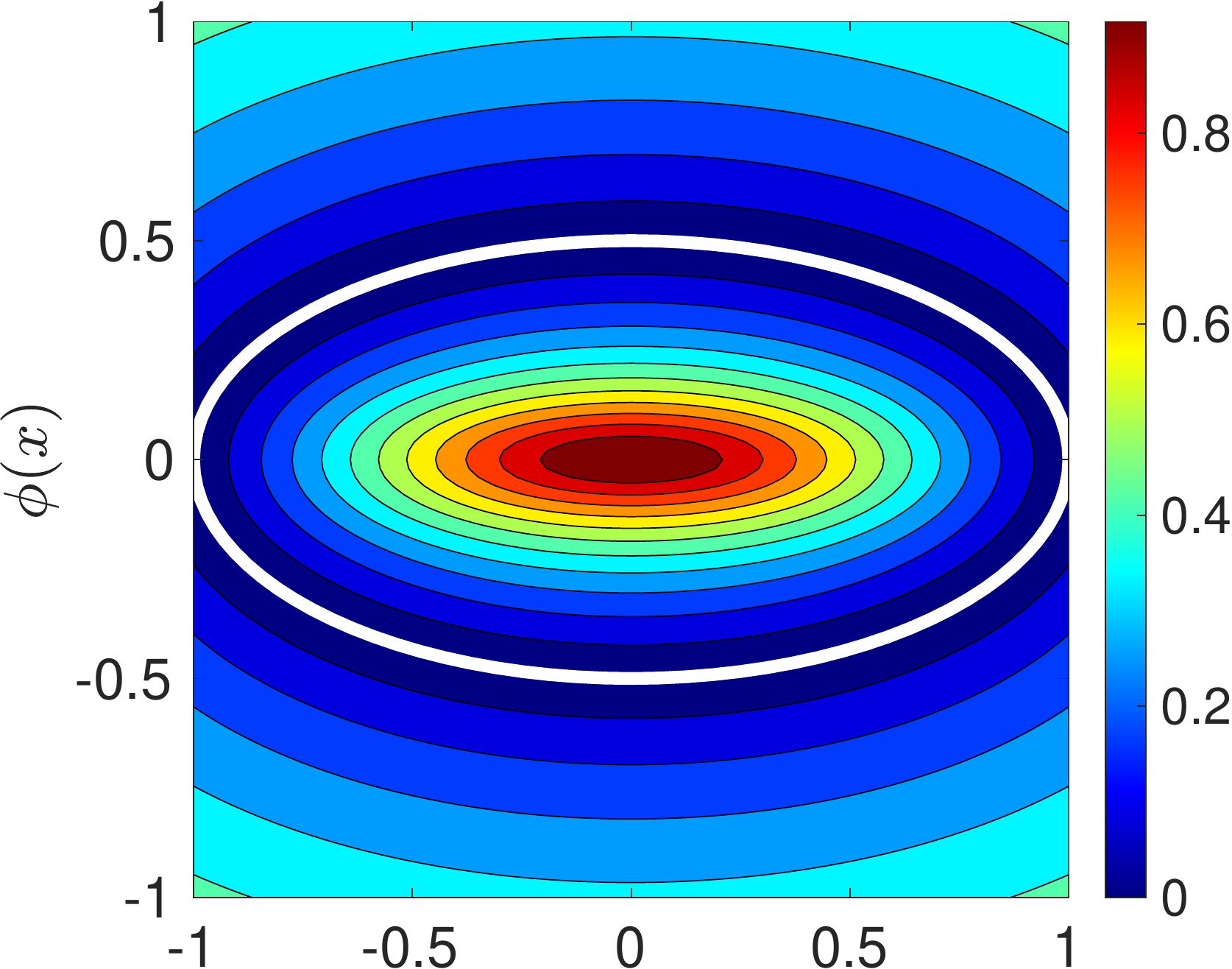}
\caption{}
\label{fig:phi_circle_ellipse-b}
\end{subfigure}
\caption{Approximate distance function (normalized to order 1) to a (a) circle and an (b) ellipse.} 
\label{fig:phi_circle_ellipse}
\end{figure}

Figure~\ref{fig:phi_seg} provides a graphical illustration of
$f$, $t$ and $\phi$ for a line segment. For a quarter-circular arc, the functions $f$, $t$ and $\phi$ are shown in~\fref{fig:phi_arc}.
In~\fref{fig:phi_circle_ellipse}, the approximate distance function (normalized to order 1) to a circle and an ellipse
are presented. The ADF to a circle of radius $R$ and center located at $\vm{x}_c := (x_c,y_c)$ is given by
\begin{equation}\label{eq:phi_circle}
\phi(\vx) = \frac{R^2 - (\vx - \vm{x}_c) \cdot (\vx - \vm{x}_c)}{2R},
\end{equation}
where $\phi(\vx)$ is a smooth (bivariate polynomial of degree 2) function.
For an elliptical disk whose closure (interior and boundary) is given by the R-function $ \omega(\vx) \ge 0$, we 
construct an ADF that is normalized to order 1 using~\cite{Rvachev:1995:RBV}
\begin{equation}\label{eq:phi_ellipse}
\phi(\vx) = \frac{\omega(\vx)} { \sqrt{ \omega^2(\vx) + 
|| \nabla \omega(\vx) ||^2 }} .
\end{equation}

In general for curves that are given in parametric form, such as B{\'ezier} and non-uniform rational B-spline (NURBS) curves, constructing ADFs require implicitization of the curve. These extensions are discussed and 
presented in Upreti et al.~\cite{Upreti:2014:ADE}. For some of the considerations and challenges in the representation (implicit and parametric) of curves in enriched computational methods, see Chin and Sukumar~\cite{Chin:2019:MCI}.

\subsection{R-equivalence operation}\label{subsec:Requiv}
Given the normalized distance functions $\phi_1$ and $\phi_2$ for two curves $c_1$ and $c_2$,  a distance field $\phi(\phi_1,\phi_2)$
for the union $c_1 \cup c_2$ must be zero when either
$\phi_1 = 0$ or $\phi_2 = 0$ and positive otherwise. When
$c_1$ and $c_2$ are line segments, the naive
formula
$\phi(\phi_1,\phi_2) = \phi_1 \phi_2$ is no longer normalized at the regular points of the segments. An R-equivalence
solution that preserves normalization up to order
$m$ of the distance function at all regular
points (nonvertices for polygonal curves) is given by~\cite{Biswas:2004:ADF}:
\begin{equation}\label{eq:phi_eq}
\phi(\phi_1,\phi_2) := \phi_1 \sim \phi_2 = 
\dfrac{\phi_1 \phi_2}{\sqrt[m]{\phi_1^m + \phi_2^m}} = 
\dfrac{1}{\sqrt[m]{\frac{1}{\phi_1^m}+\frac{1}{\phi_2^m}}} ,
\end{equation}
\suku{where $\lim_{m \to \infty} \phi(\phi_1,\phi_2) = \min(\phi_1,\phi_2)$.}
When $\partial S$ (closed  curve) is composed of $n$ pieces, 
then a $\phi$ that is normalized up to order $m$ is given
by (see the proof in~\cite{Millan:2015:CBM}):
\begin{equation}\label{eq:phin_eq}
\phi(\phi_1,\dots, \phi_n) := \phi_1 \sim \phi_2 \sim \cdots \sim 
\phi_n = 
\dfrac{1}{\sqrt[m]{\frac{1}{\phi_1^m}+\frac{1}{\phi_2^m}
          + \ldots + \frac{1}{\phi_n^m}}}.
\end{equation}
The $\phi$ that is formed in~\eqref{eq:phin_eq} can be viewed as
the reciprocal of the $L_m$-norm of inverse distance
measures, which bears similarity to $L_m$-distance fields~\cite{Belyaev:2013:SLD}.
An alternative joining procedure is to consider the R-conjunction given by 
\begin{equation}\label{eq:phi_Rconj}
\phi_s (\phi_1, \phi_2) := \phi_1 \wedge \phi_2 = 
\phi_1 + \phi_2 - \sqrt[s]{\phi_1^s + \phi_2^s} ,
\end{equation}
which is a function that is
normalized to the $(s-1)$-th order~\cite{Biswas:2004:ADF}.  However, the
joining operation is not associative, which makes this
choice less desirable. The R-equivalence joining relation in~\eqref{eq:phin_eq} is
associative. 

\begin{figure}[!htp]
\centering
\begin{subfigure}{0.49\textwidth}
\mbox{
\includegraphics[width=0.48\textwidth]{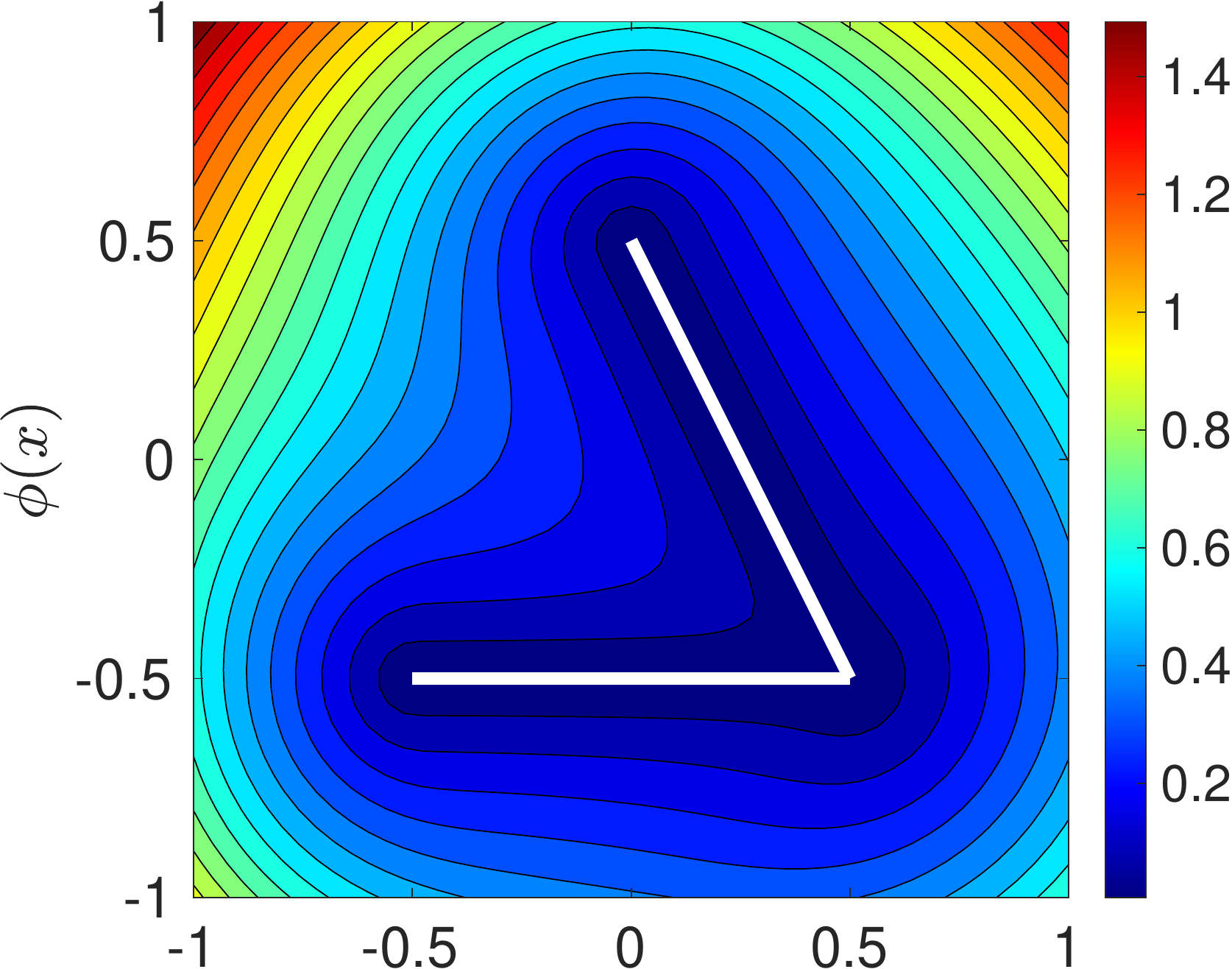}
\includegraphics[width=0.48\textwidth]{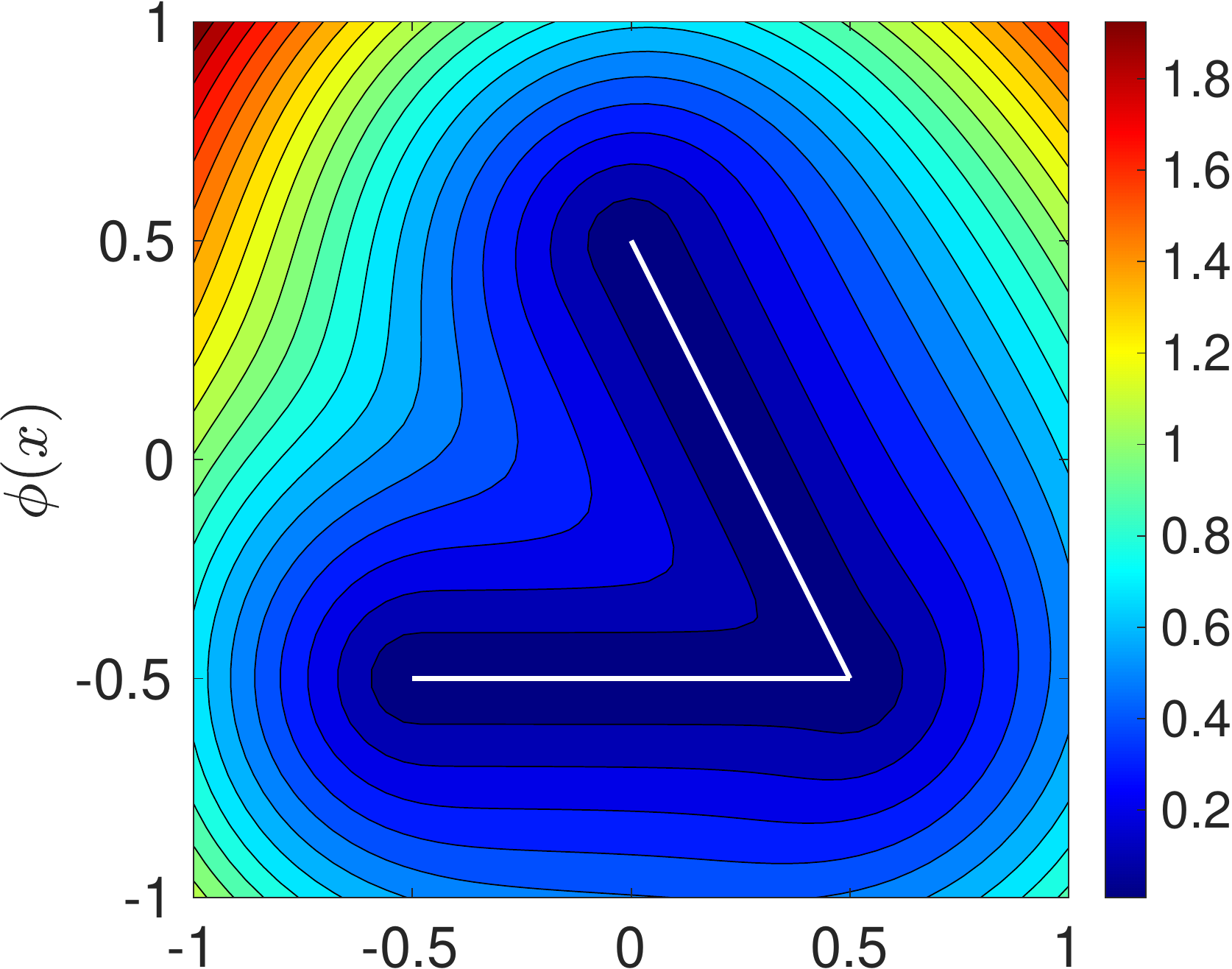}
}
\caption{}
\end{subfigure}
\begin{subfigure}{0.49\textwidth}
\mbox{
\includegraphics[width=0.48\textwidth]{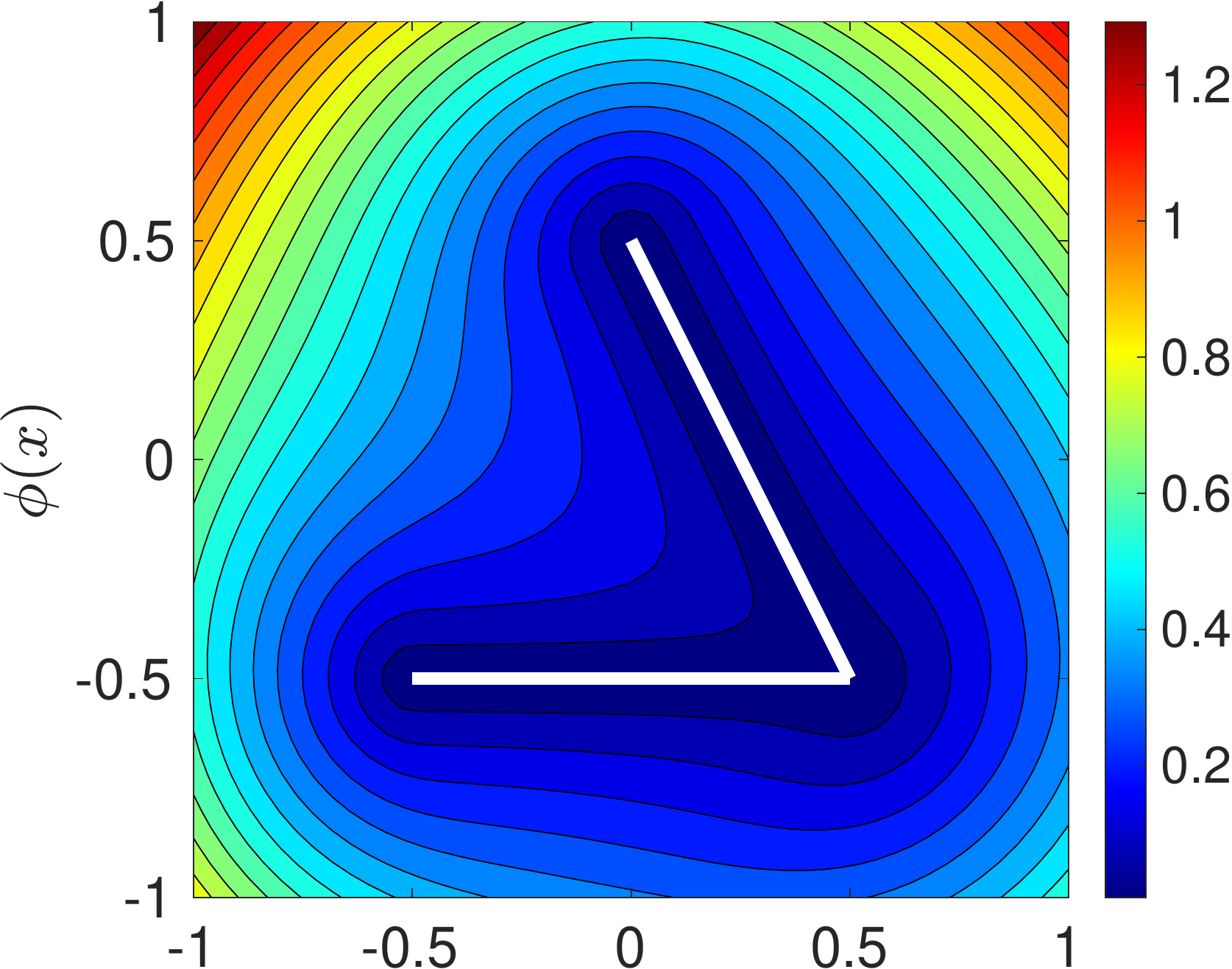}
\includegraphics[width=0.48\textwidth]{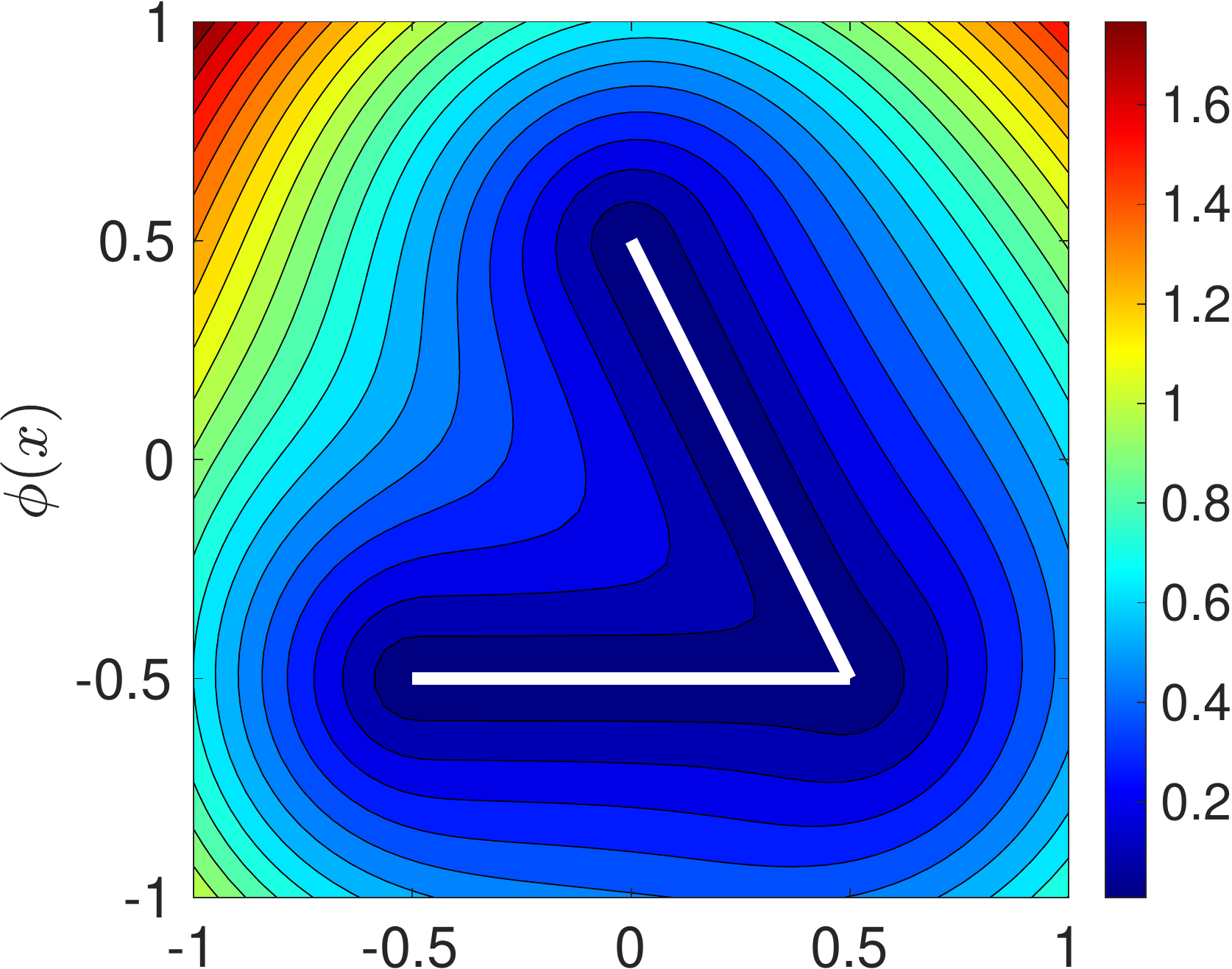}
}
\caption{}
\end{subfigure}
\caption{Approximation of the distance function to two line
         segments. (a) R-conjunction composition 
         with $s = 2,3$ in~\eqref{eq:phi_Rconj}, and (b) R-equivalence composition in~\eqref{eq:phi_eq} for the normalization parameter $m = 1,2$. The ADFs are
         normalized to order $s-1$ and $m$, respectively.}
         \label{fig:phi_twosegs}
\end{figure}

\begin{figure}[!htp]
\centering
\begin{subfigure}{0.24\textwidth}
\includegraphics[width=\textwidth]{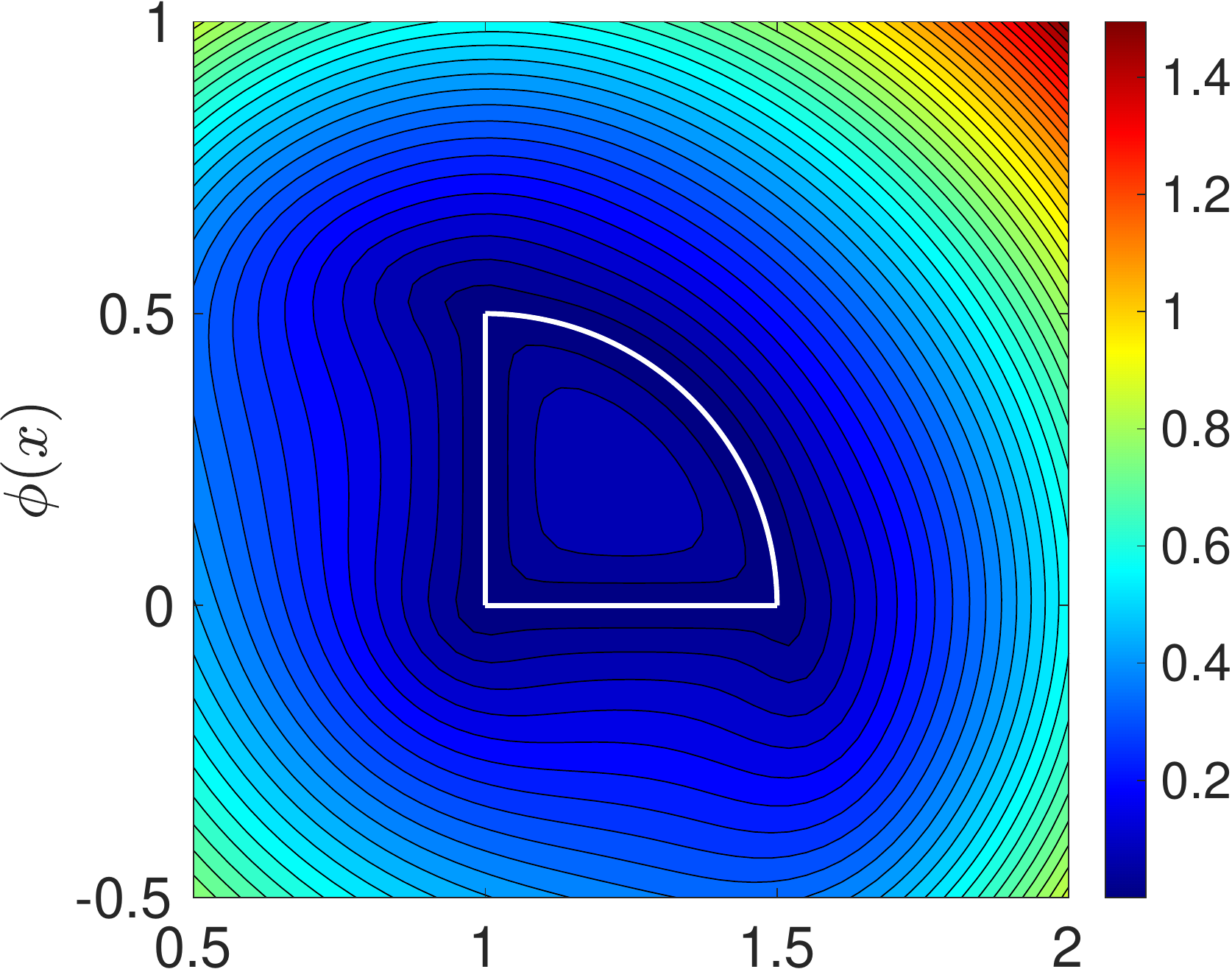}
\caption{$m = 2$}
\end{subfigure}
\begin{subfigure}{0.24\textwidth}
\includegraphics[width=\textwidth]{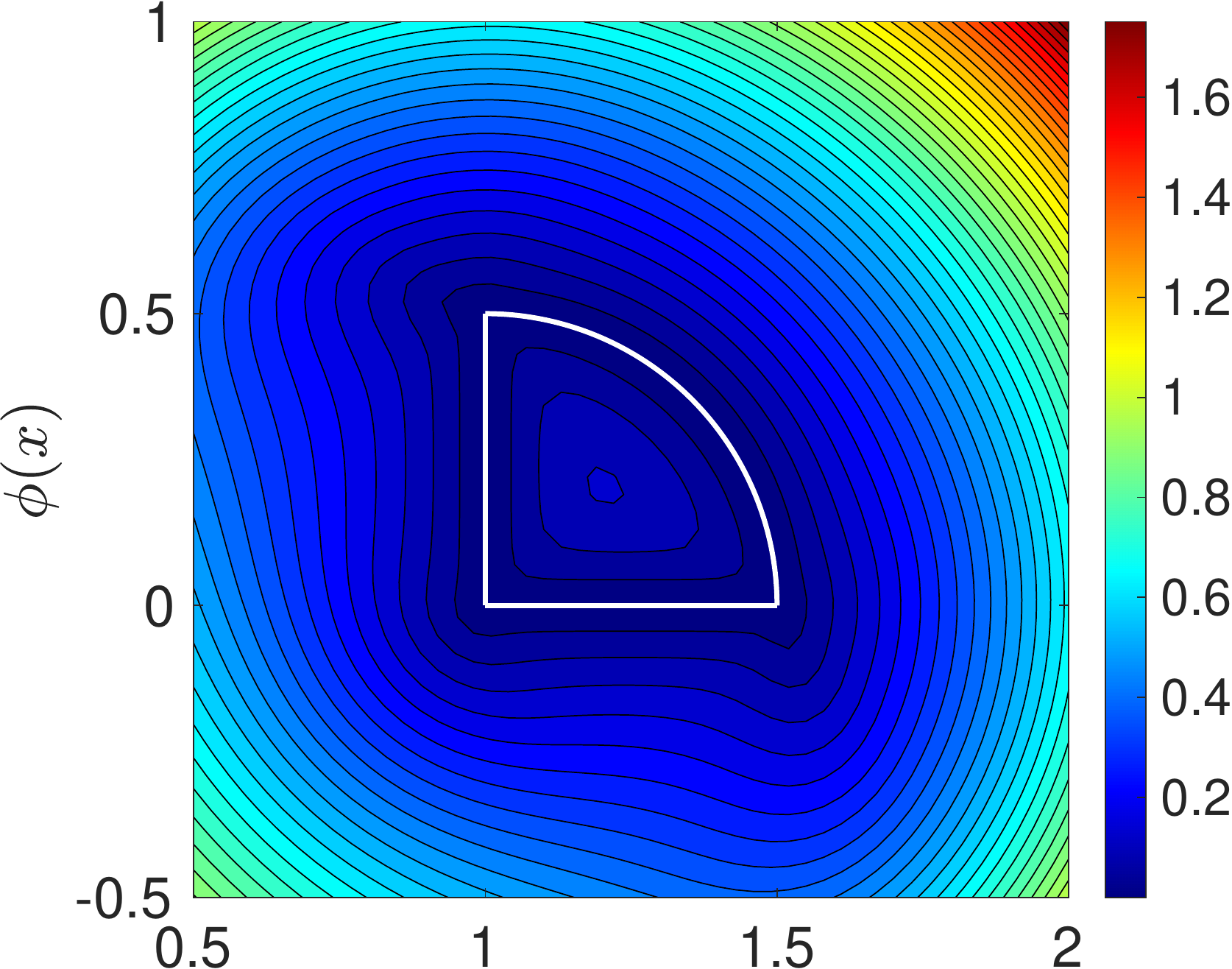}
\caption{$m = 3$}
\end{subfigure}
\begin{subfigure}{0.24\textwidth}
\includegraphics[width=\textwidth]{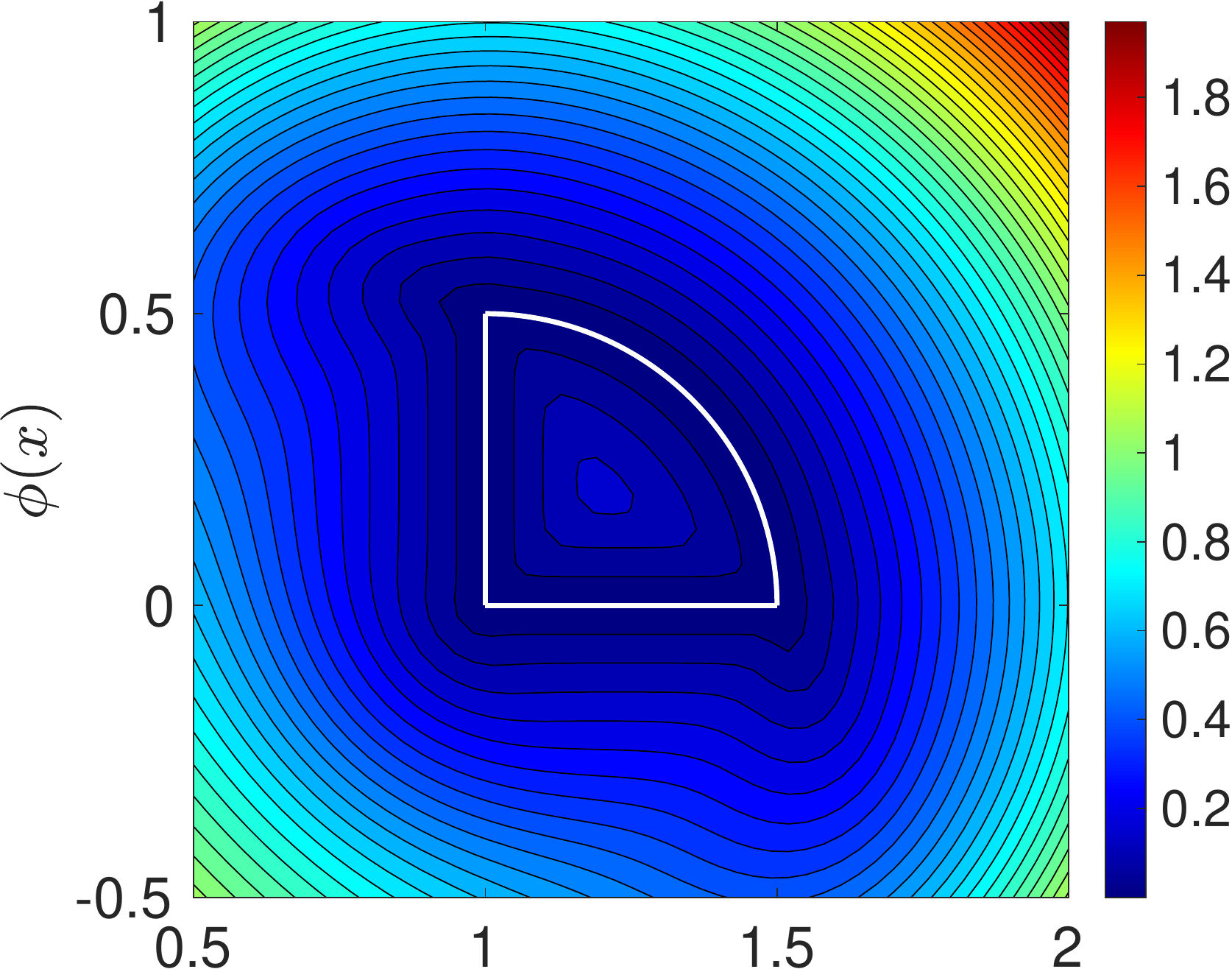}
\caption{$m = 6$}
\end{subfigure}
\begin{subfigure}{0.24\textwidth}
\includegraphics[width=\textwidth]{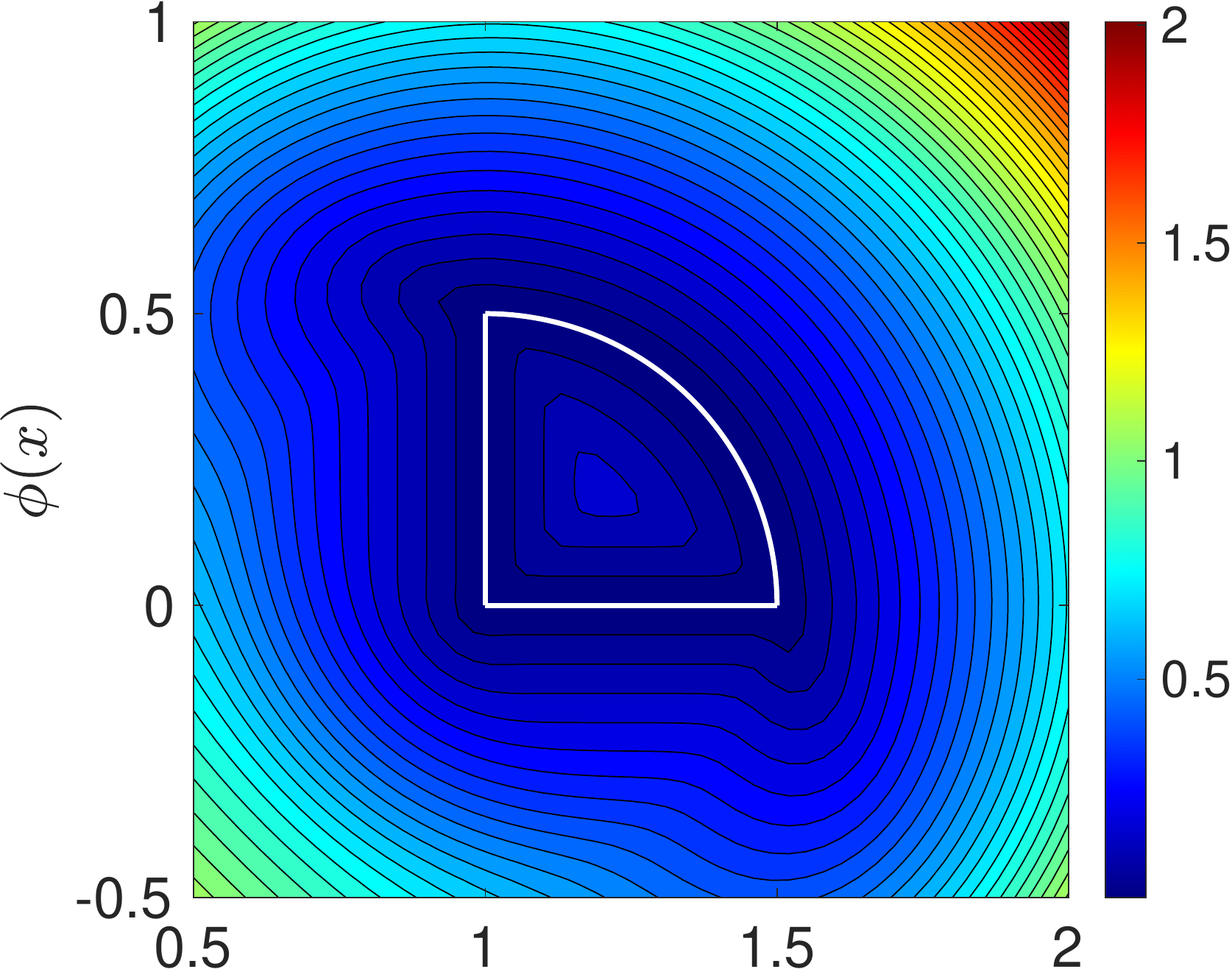}
\caption{$m = 10$}
\end{subfigure}
\caption{Plots of the approximate distance function to a curved
         triangle using R-equivalence for different choices of the normalizing parameter ($m = 2,3,6,10$).}
         \label{fig:phi_curvedtriangle}
\end{figure}

\begin{figure}[!htp]
\centering
\begin{subfigure}{0.24\textwidth}
\includegraphics[width=\textwidth]{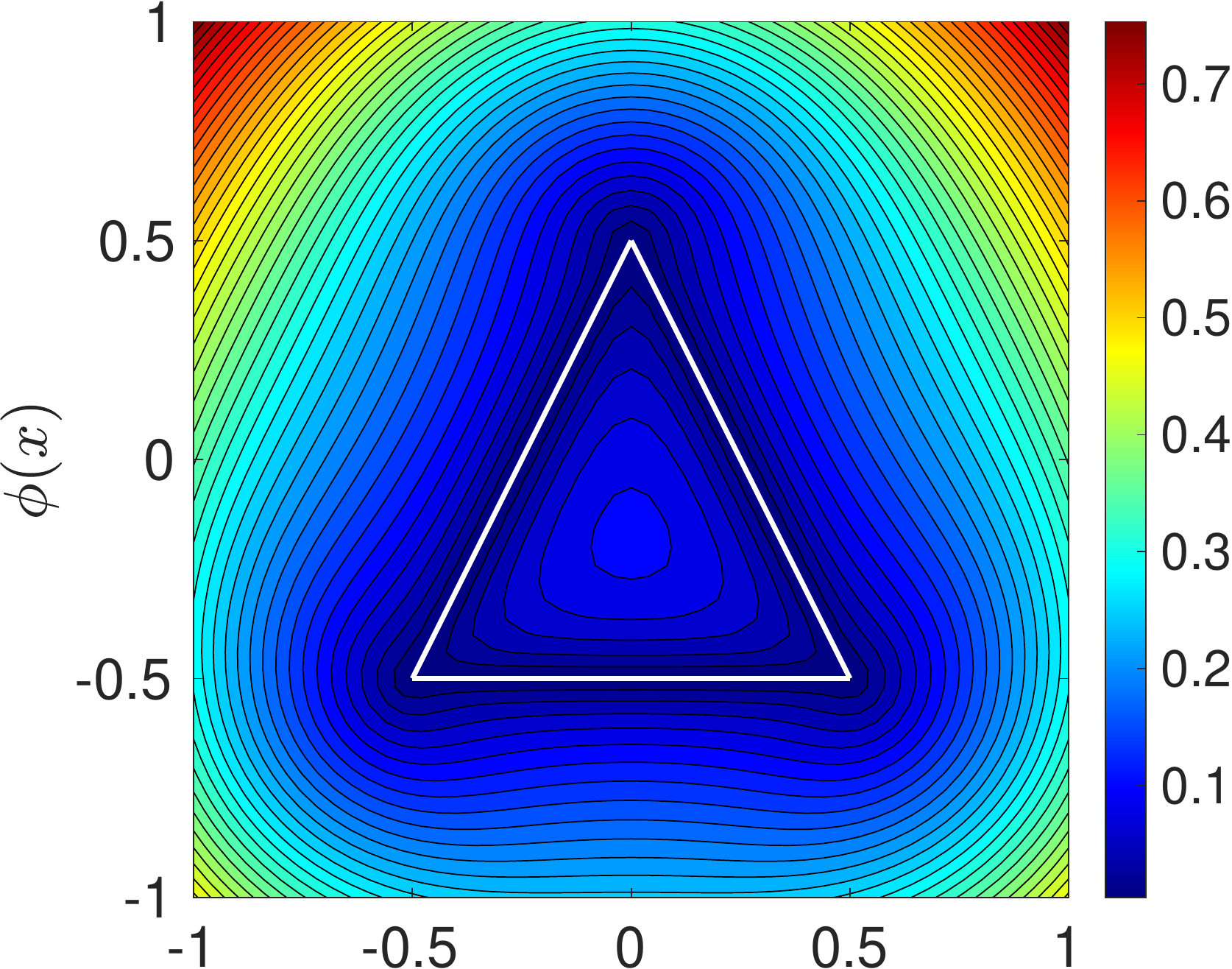}
\caption{}
\label{fig:phi_polygons-a}
\end{subfigure}
\begin{subfigure}{0.24\textwidth}
\includegraphics[width=\textwidth]{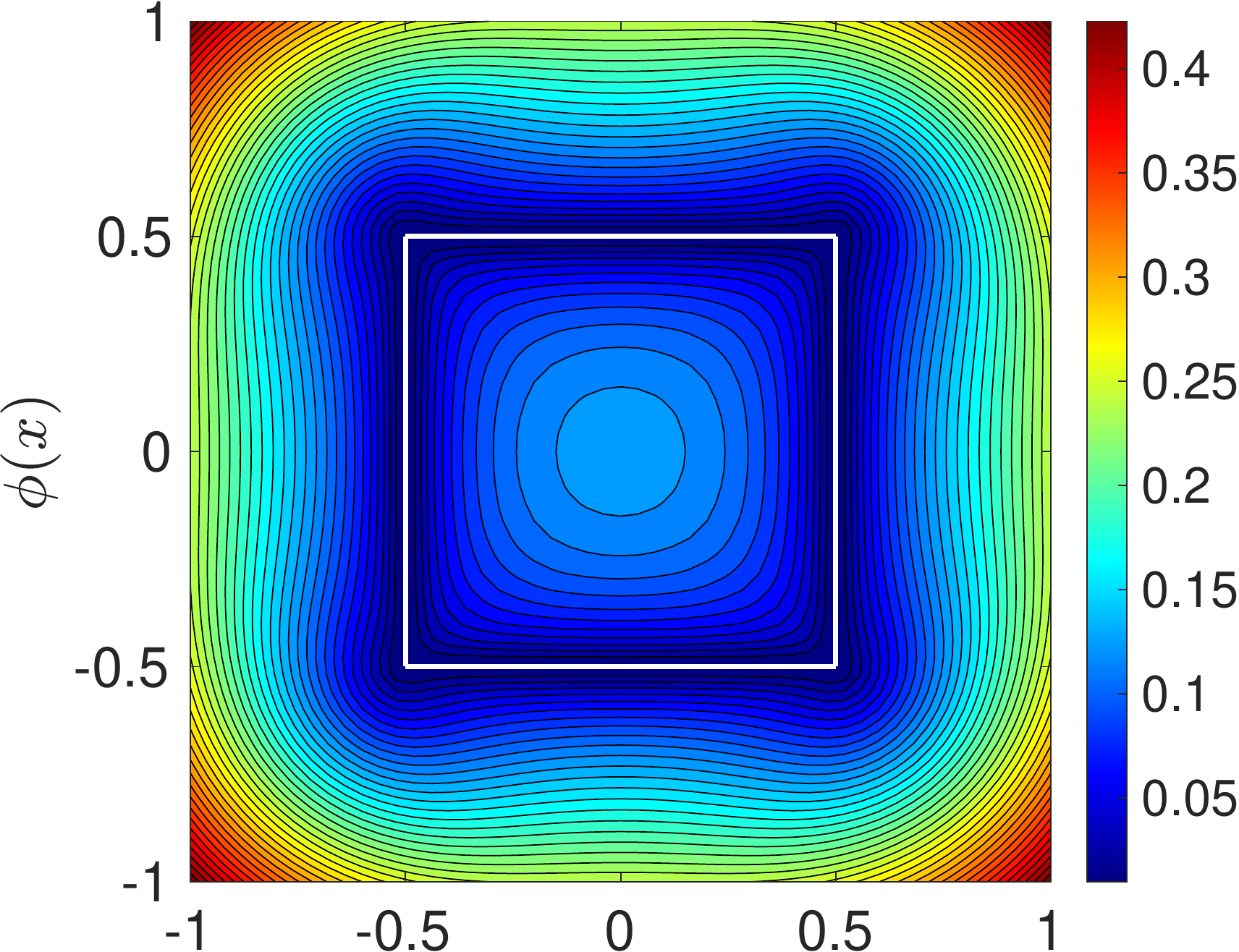}
\caption{}
\label{fig:phi_polygons-b}
\end{subfigure}
\begin{subfigure}{0.24\textwidth}
\includegraphics[width=\textwidth]{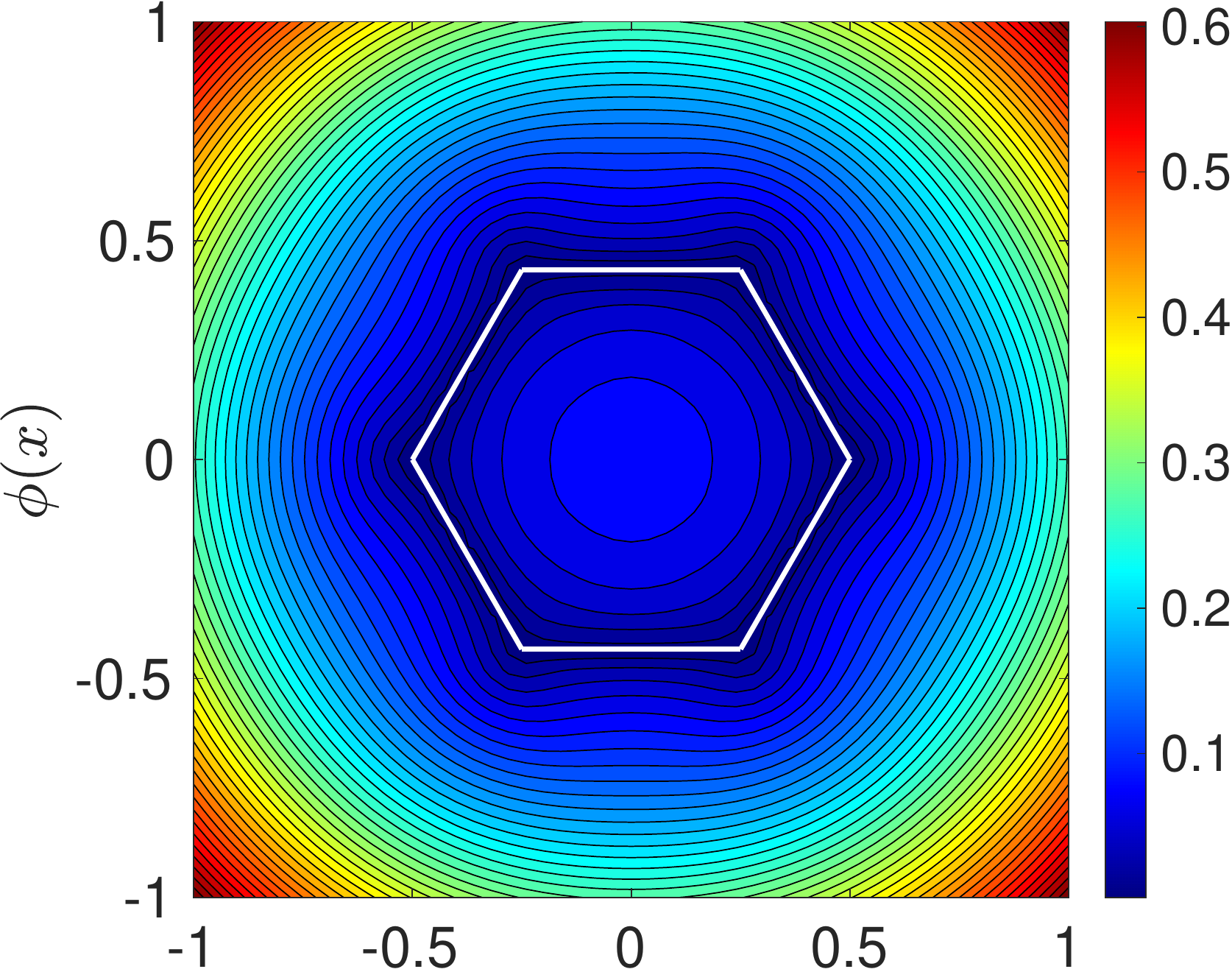}
\caption{}
\label{fig:phi_polygons-c}
\end{subfigure}
\begin{subfigure}{0.24\textwidth}
\includegraphics[width=\textwidth]{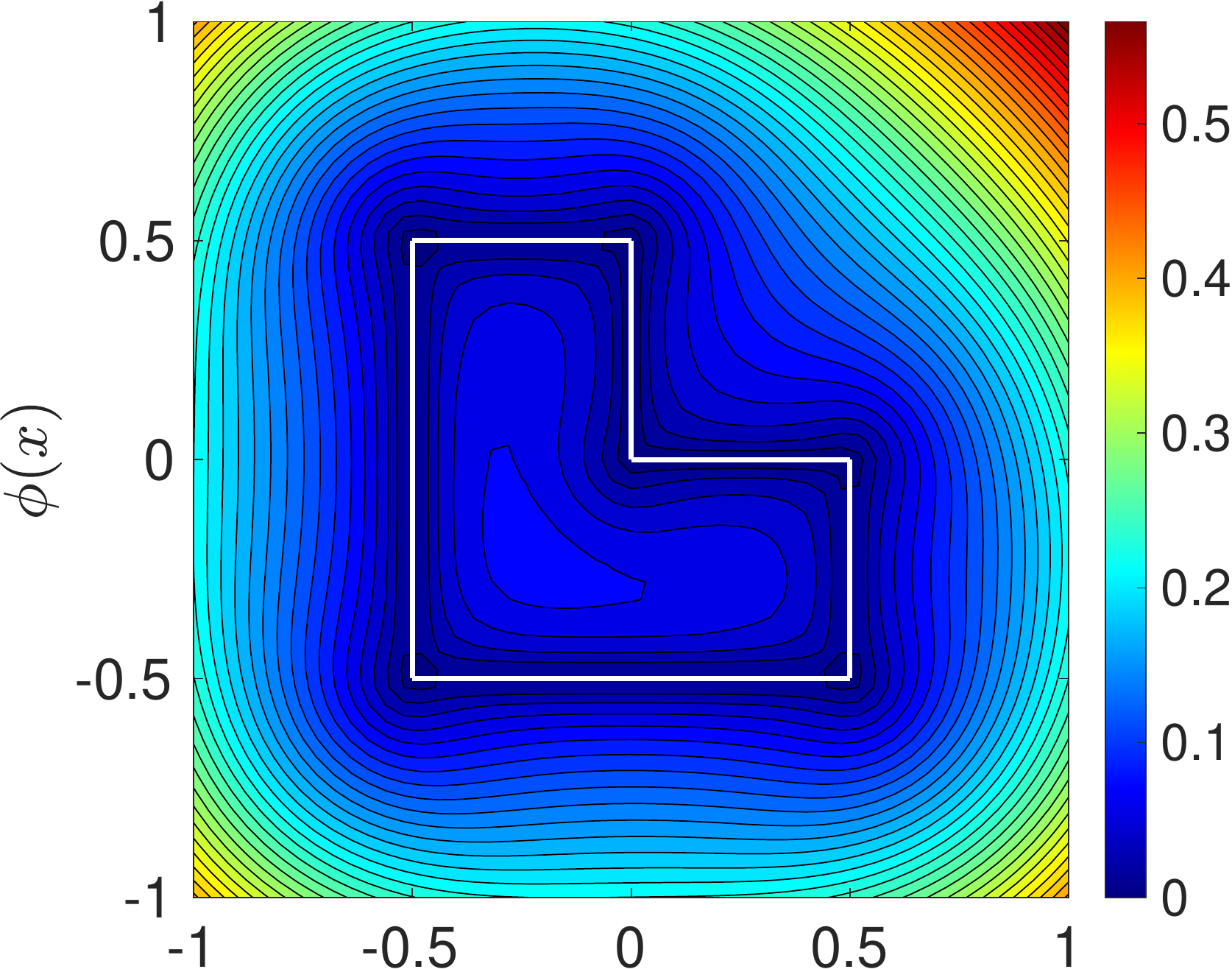}
\caption{}
\label{fig:phi_polygons-d}
\end{subfigure}
\caption{Plots of the approximate distance function using 
         R-equivalence ($m = 1$) for polygons. (a) triangle,  
         (b) square, (c) regular hexagon, and (d) L-shaped (nonconvex) polygon.
         } \label{fig:phi_polygons}
\end{figure}

\begin{figure}[!htp]
\centering
\includegraphics[width=0.5\textwidth]{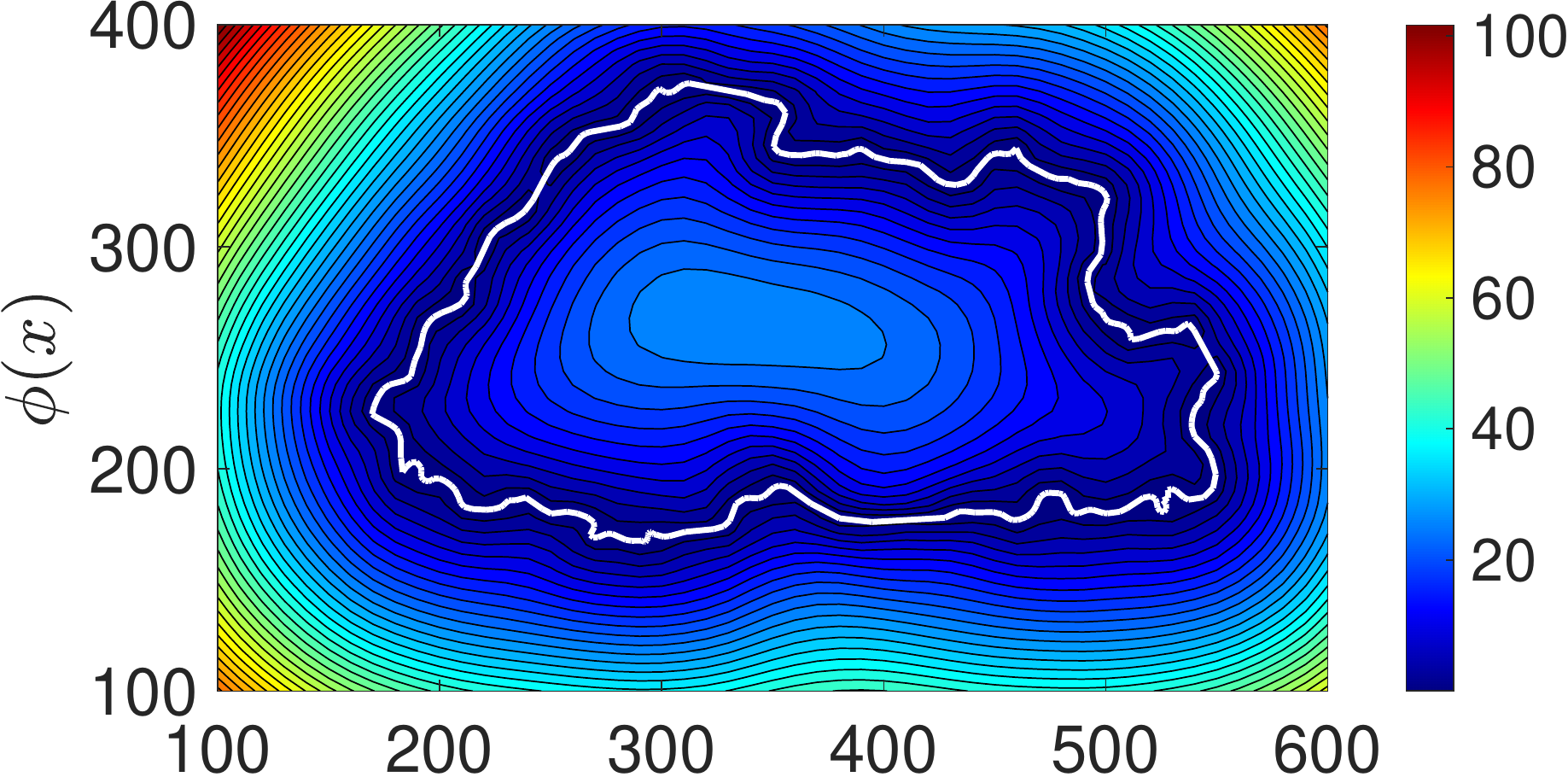}
\caption{Plot of the approximate distance function using R-equivalence            ($m = 1$) for polygonalized map of Bhutan. The polygon has 291 vertices.}
\label{fig:phi_bhutan_map}
\end{figure}

The approximate distance function to two line segments is shown
in~\fref{fig:phi_twosegs}, where the R-conjunction composition in~\eqref{eq:phi_Rconj} with $s = 2, 3$ and the R-equivalence relation in~\eqref{eq:phi_eq} with $m = 1,2$ are compared.
The bulging phenomenon~\cite{Bloomenthal:1997:BEC} 
is noticeable in the vicinity of the joining point.  In~\fref{fig:phi_curvedtriangle}, we present the approximate distance function to a curved triangular region using
R-equivalence for different orders of the normalizing parameter
$m$. As $m$ increases, the ADF approaches the exact distance, which is observed when inspecting the contours in the interior of the curved triangle. 
The ADFs for a triangle, square, hexagon, and an 
L-shaped polygon are shown in~\fref{fig:phi_polygons}. 
The ADF for a square in~\fref{fig:phi_polygons-b} bears similarities to a superellipse, $|x|^m + |y|^m = 1$, which
has rounded corners as $m$ increases.
As the last example, we present the ADF for a complex polygonal domain. 
We consider the polygonalized map of Bhutan,\footnote{Vectorized \texttt{eps} image obtained from \url{https://freevectormaps.com}} which has 291 boundary vertices. The contour plot for $\phi(\vx)$ is depicted in~\fref{fig:phi_bhutan_map}, and we
observe that the contours are smooth in the interior and
well-separated ($\phi$ is monotonic in $\Omega$). 
Finally,
we mention that a modified form of the R-equivalence relation~\eqref{eq:phin_eq} is also discussed in Biswas and Shapiro~\cite{Biswas:2004:ADF}. This modified form is constructed with an eye on better capturing the first-order normalization condition in the sector region between two line segments, where the normal is undefined and the closest point to either segment is the common vertex. With increase in $m$, the R-equivalence joining
operation provides a better approximation to the exact distance away from the segments and also improved normalization properties in
the vicinity of the segments; however, this comes at the expense of 
the function being higher order and hence its Laplacian will have
greater undulations.  Use of $m = 2$ has been 
adopted in prior computational studies~\cite{Biswas:2004:ADF,Upreti:2014:ADE,Millan:2015:CBM}, but
herein, we adopt $m = 1$ in most of the numerical simulations that are presented in Sections~\ref{sec:1D}--\ref{sec:4D}.

\section{Generalized Mean Value Potentials and \texorpdfstring{$L_p$}{Lp}-Distance Fields}\label{sec:mvp}
In addition to the theory of R-functions, another approach to construct approximate distance fields is
via the theory of mean value potential fields.  This has connections
to generalized barycentric coordinates and in particular to mean value interpolation over polygons and curved domains~\cite{Belyaev:2017:TBC}.
Generalized barycentric coordinates~\cite{Floater:2013:GBC,Anisimov:2017:BCP,Hormann:2017:GBC} are an
extension of barycentric coordinates
over simplices to polygons and polyhedra.  These
coordinates (shape functions) have linear precision and are nonnegative
over convex polygons.  Transfinite barycentric interpolation over domains bounded by curves is the
continuous counterpart of generalized barycentric coordinates over
polygons~\cite{Belyaev:2017:TBC}. Given a function 
$u :\Re^2 \to \Re$ that assumes the function $g(\vx)$ on the boundary (curved) of a
domain, a transfinite interpolant provides an approximation of $u(\vx)$ that matches $g(\vx)$ over the curved boundary of the domain. For a domain that is bounded by affine or curved boundaries, the reciprocal of the mean value normalization function
is a smoothed approximation to the exact distance function~\cite{Hormann:2006:MVC,Dyken:2009:TMV,Bruvoll:2010:TMV}, 
and is a specific instance ($p = 1$) of the reciprocal of 
a singular double-layer $L_p$-potential 
field~\cite{Belyaev:2013:SLD}. We refer
to the method that generates these smoothed distance 
functions by the acronym MVP, since they stem from {\em (generalized) 
mean value potential fields}~\cite{Belyaev:2017:TBC}. The
\suku{construction} of $\phi(\vx)$ over polygons and curved domains is
presented in Sections~\ref{subsec:mvp_polygon}
and~\ref{subsec:mvp_curved}, respectively.

\subsection{Approximate distance fields on arbitrary planar polygons}\label{subsec:mvp_polygon}
A popular generalized barycentric coordinate is
Floater's mean value coordinates~\cite{Floater:2003:MVC}, which
is derived using the circumferential mean value theorem for harmonic functions. This conception stemmed from the objective to approximate a harmonic map by a convex combination map (positive weights) over a triangulation, so that injectivity is
preserved. Mean value coordinates have many remarkable properties: for instance, they are 
valid on arbitrary planar polygons, including
nested polygons; $C^\infty$ smooth in $\Omega$ with derivative
discontinuities only at the vertices of the polygon;
they reduce to piecewise affine functions on the edges of a polygon; are 
nonnegative in the kernel of the polygon; reciprocal of the normalizing 
weight function is a smoothed ADF; and they  also have a smooth
extension outside the polygon~\cite{Hormann:2006:MVC}.

\begin{figure}[!htb]
\centering
\begin{subfigure}{0.49\textwidth}
\includegraphics[width=0.92\textwidth]{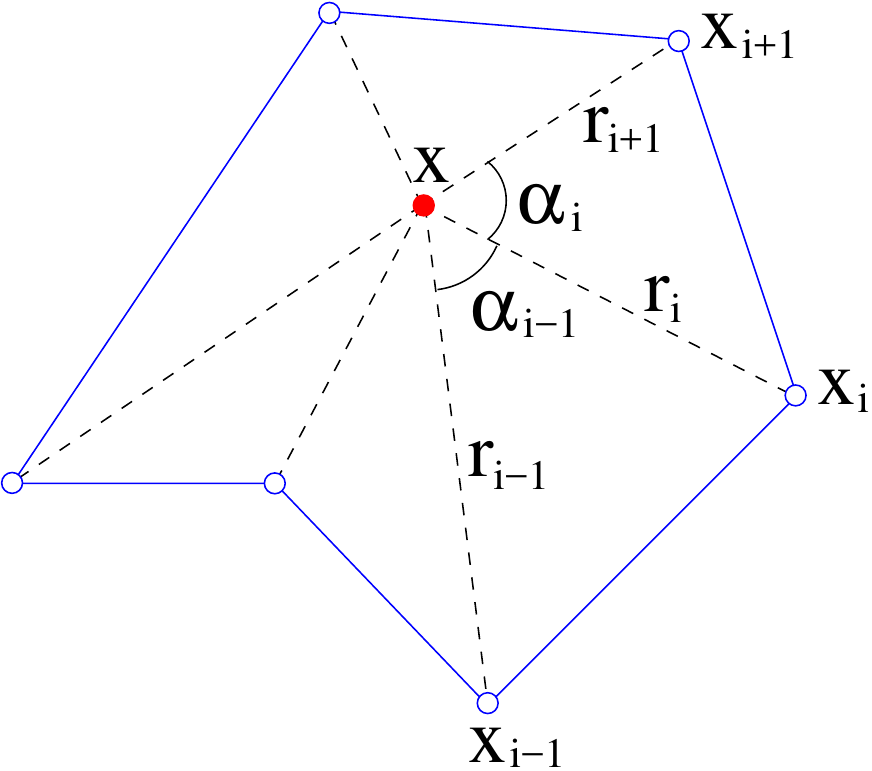}
\caption{}
\end{subfigure} 
\begin{subfigure}{0.42\textwidth}
\includegraphics[width=0.90\textwidth]{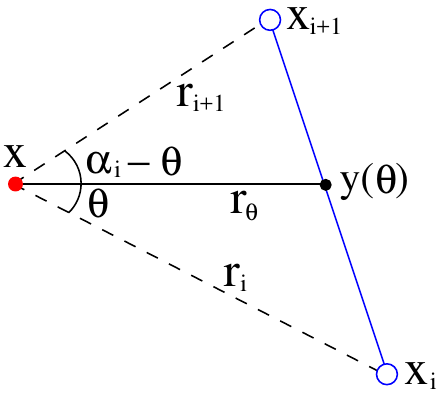}
\caption{}
\end{subfigure}
\caption{Notation used in the definition of (a) mean value coordinates~\protect\cite{Floater:2003:MVC} and (b) generalized mean value potentials~\protect\cite{Belyaev:2017:TBC}. In (b), the
parameters that are used to form $W(\vx)$ in~\eqref{eq:weight_mvc} are shown.}
\label{fig:mvc}
\end{figure}

Consider the nonconvex polygon ($n$-gon) shown in~\fref{fig:mvc}, 
whose $n$ vertices are defined in counterclockwise orientation. The coordinates of the vertices are
$\{\vx_i\}_{i=1}^n$, and $\vx$ is an arbitrary point in the interior of the polygon.  The mean value coordinates, $\{\varphi_i(\vx)\}_{i=1}^n$,
are defined as~\cite{Floater:2003:MVC}:
\begin{equation}\label{eq:mvc}
\varphi_i(\vx) = \dfrac{w_i(\vx)}{W(\vx)}, \quad
w_i(\vx) := \dfrac{\tan \left( \alpha_{i-1}/2 \right) + 
\tan \left( \alpha_i/2 \right) }{|| \vx_i - \vx||}, \quad
W(\vx) = \sum_{j=1}^n w_j(\vx), 
\end{equation}
where the angles $\alpha_{i-1}$ and $\alpha_i$ are shown in~\fref{fig:mvc}.
Let $\vm{r}_i := \vx_i - \vx$ with $r_i = ||\vx_i - \vx||$ represent
the Euclidean distance between $\vx$ and $\vx_i$.  On noting the half-angle
formula for $\tan (\cdot)$, we can define
\begin{equation}\label{eq:weight_mvc}
t_i := \tan \left( \dfrac{\alpha_i}{2} \right) = 
\dfrac{\sin \alpha_i}{1 + \cos \alpha_i} = 
\dfrac{r_i r_{i+1} \sin \alpha_i}{r_i r_{i+1} + \vm{r}_i \cdot \vm{r}_{i+1}} 
= \dfrac{ \textrm{det} \, ( \vm{r}_i, \vm{r}_{i+1}) }{r_i r_{i+1} + \vm{r}_i \cdot \vm{r}_{i+1}} , \quad
W(\vx) = \sum_{i=1}^n \left( \frac{1}{r_i} + \frac{1}{r_{i+1}}
\right) t_i  \quad (r_{n+1} := r_1),
\end{equation}
which is now valid for all points $\vx$ that are in the interior of a convex or nonconvex polygon.
The denominator vanishes
when $\alpha_i = \pi$, i.e., when $\vx$ lies on the boundary of the polygon,
but there $\varphi_i(\vx)$ are known.  The singularity of the
weight function on the boundary is a property shared by
nonnegative generalized barycentric coordinates.

\begin{figure}[!htb]
    \centering
    \begin{subfigure}{0.49\textwidth}
    \mbox{
    \includegraphics[width=0.53\textwidth]{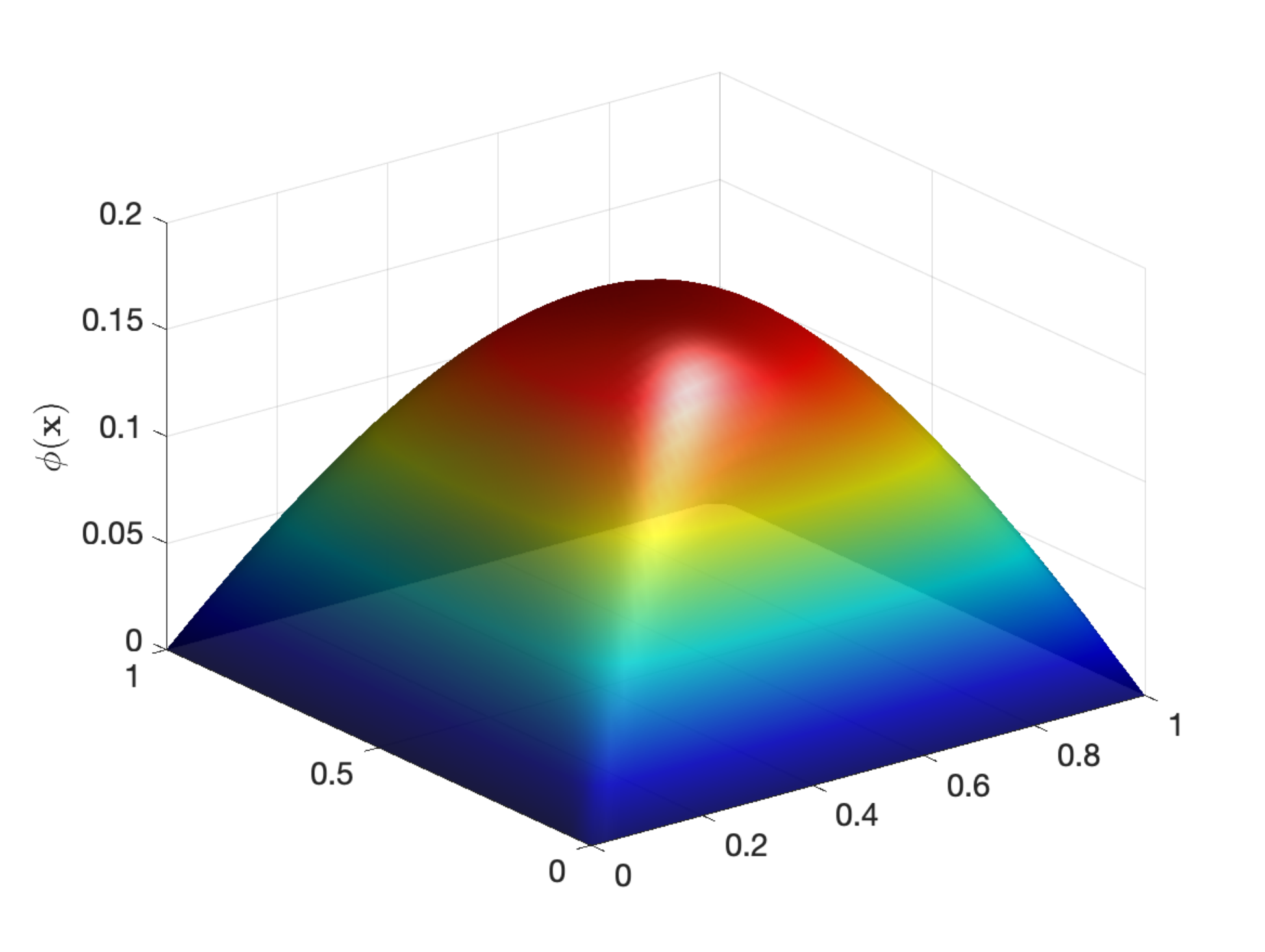}
     \includegraphics[width=0.4\textwidth]{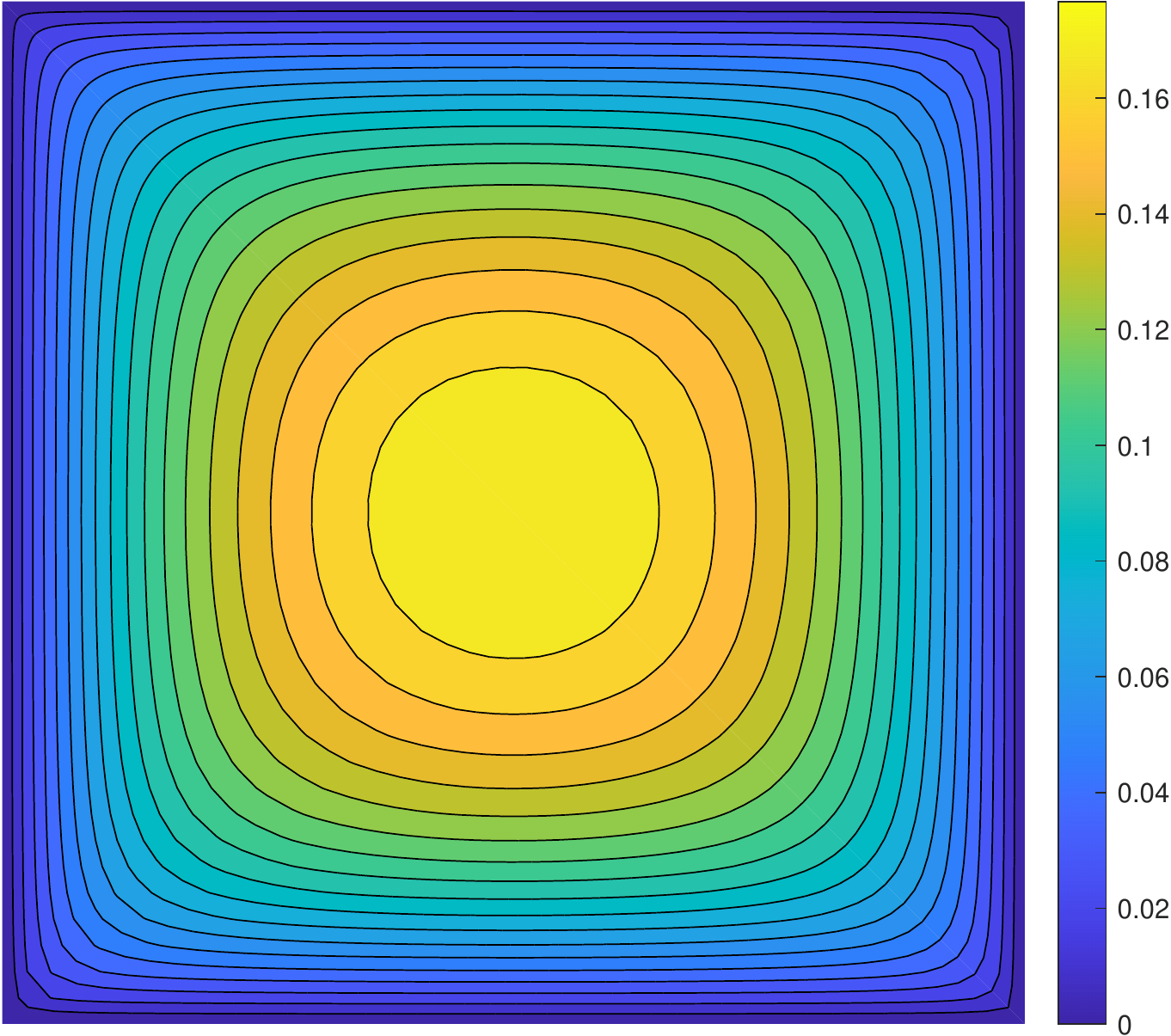}
     }
     \caption{}
     \label{fig:mvc_ADF-a}
    \end{subfigure}
    \begin{subfigure}{0.49\textwidth}
    \mbox{
    \includegraphics[width=0.53\textwidth]{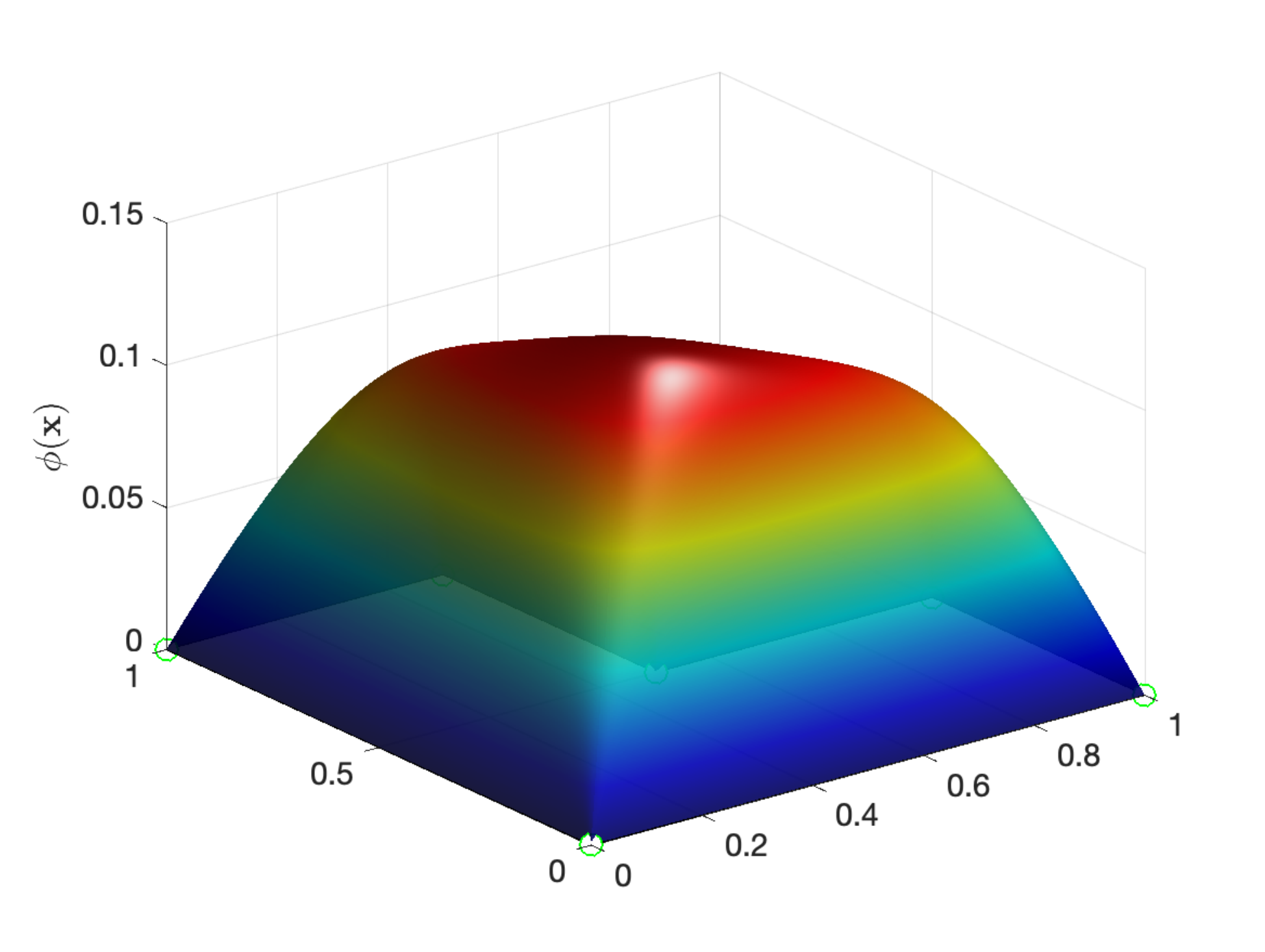}
    \includegraphics[width=0.43\textwidth]{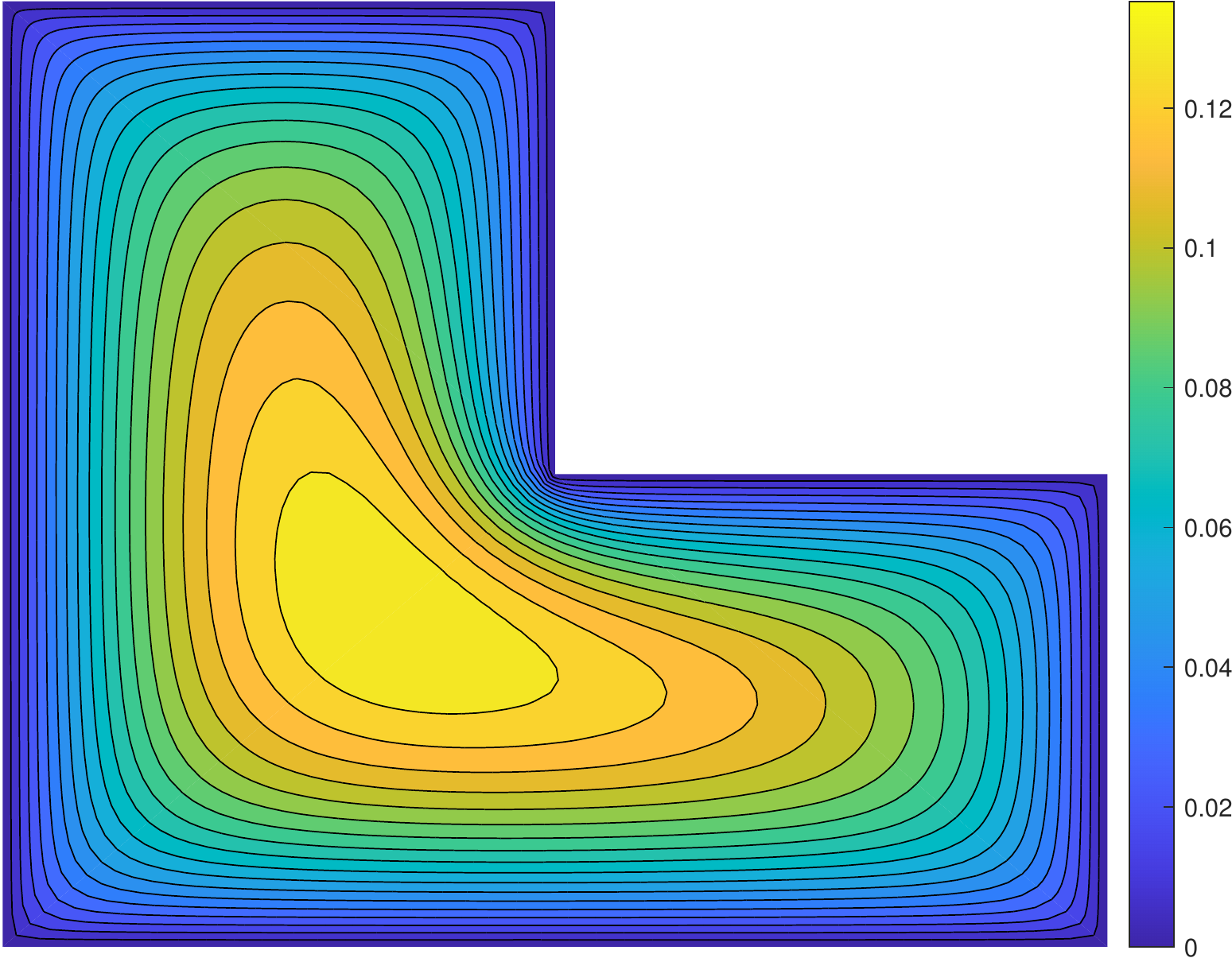}
    }
    \caption{}
    \label{fig:mvc_ADF-b}
    \end{subfigure}
    \begin{subfigure}{0.49\textwidth}
    \mbox{
    \includegraphics[width=0.53\textwidth]{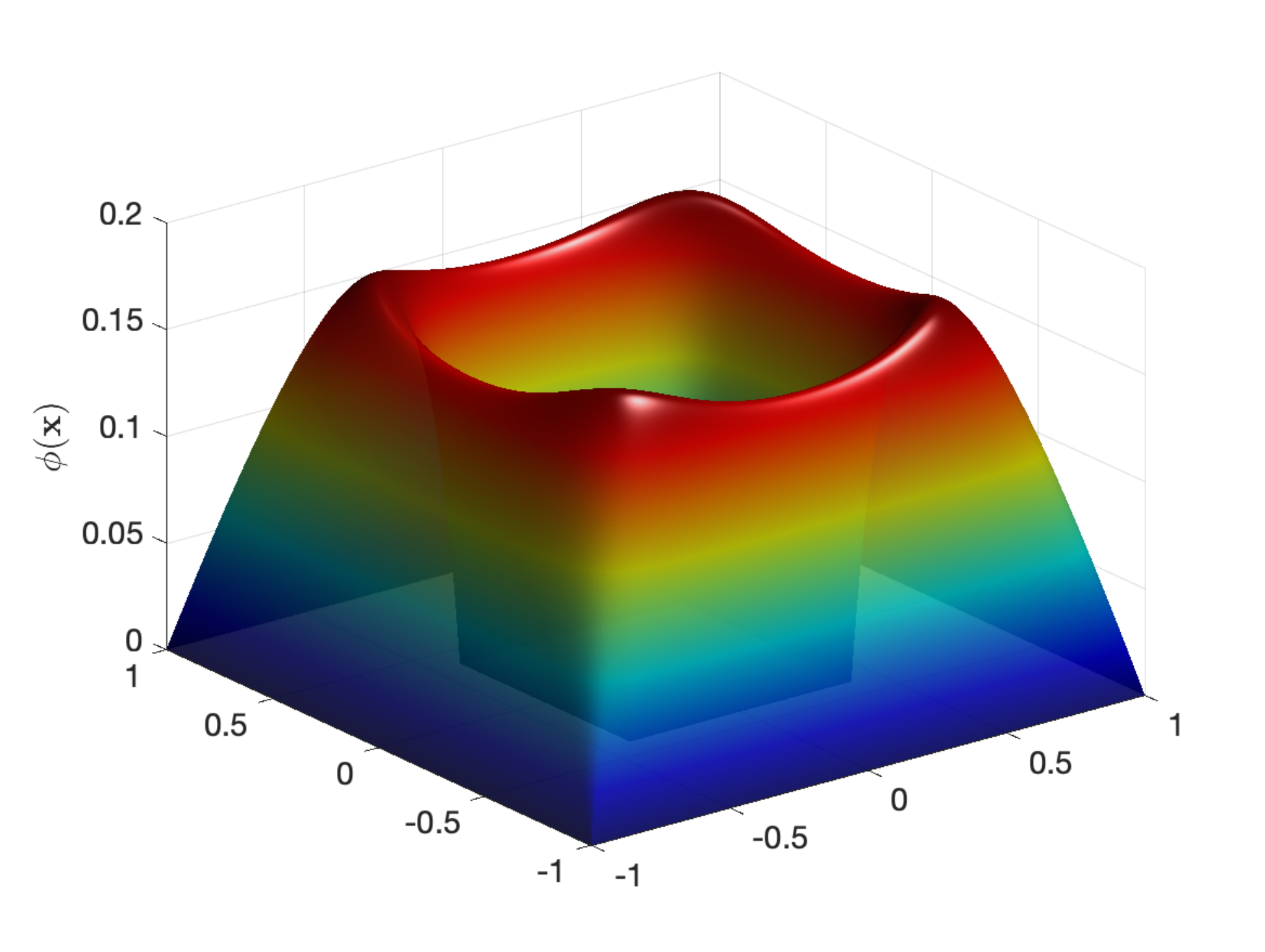}
    \includegraphics[width=0.43\textwidth]{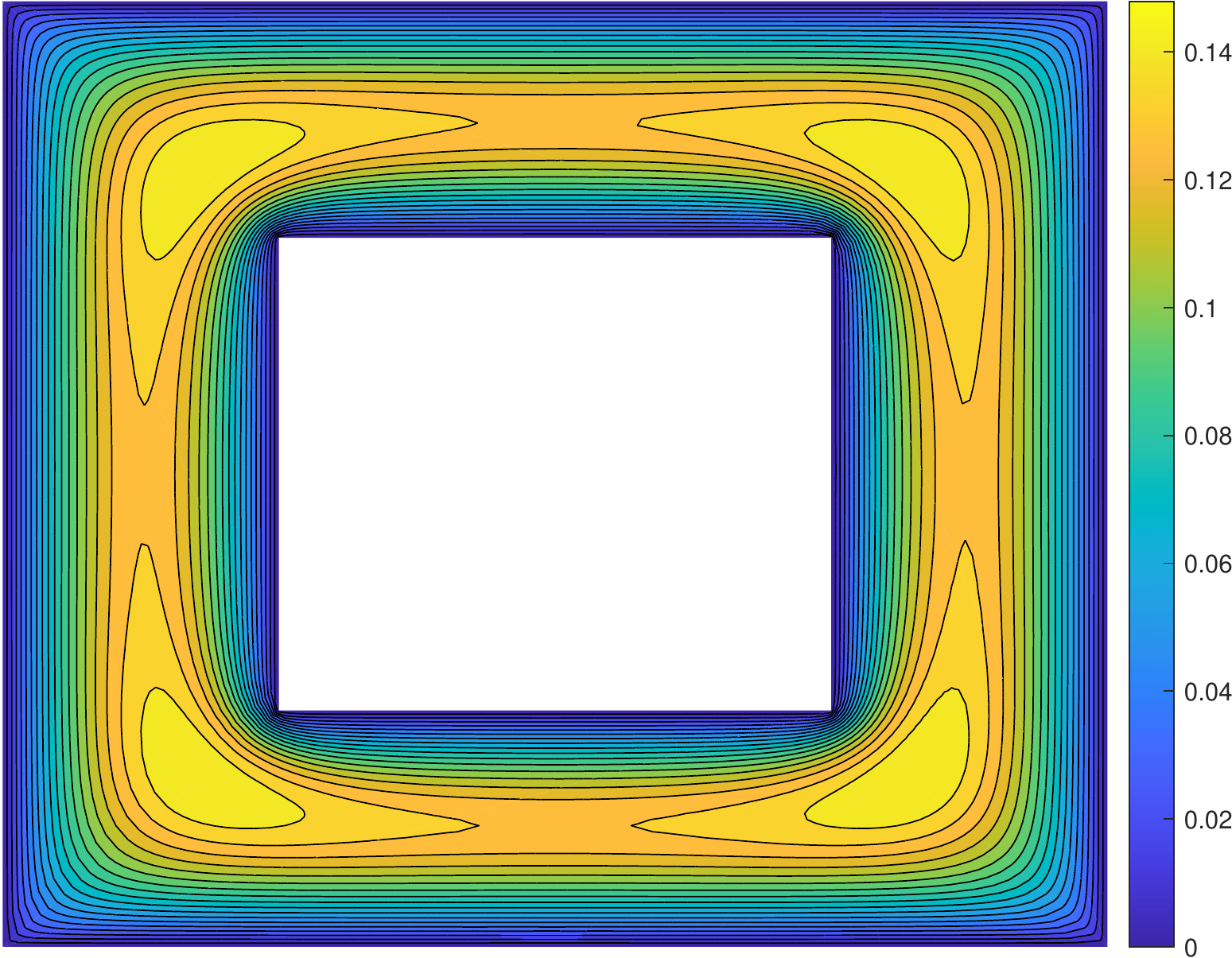}
    }
    \caption{}
    \label{fig:mvc_ADF-c}
    \end{subfigure}
    \begin{subfigure}{0.49\textwidth}
    \mbox{
    \includegraphics[width=0.53\textwidth]{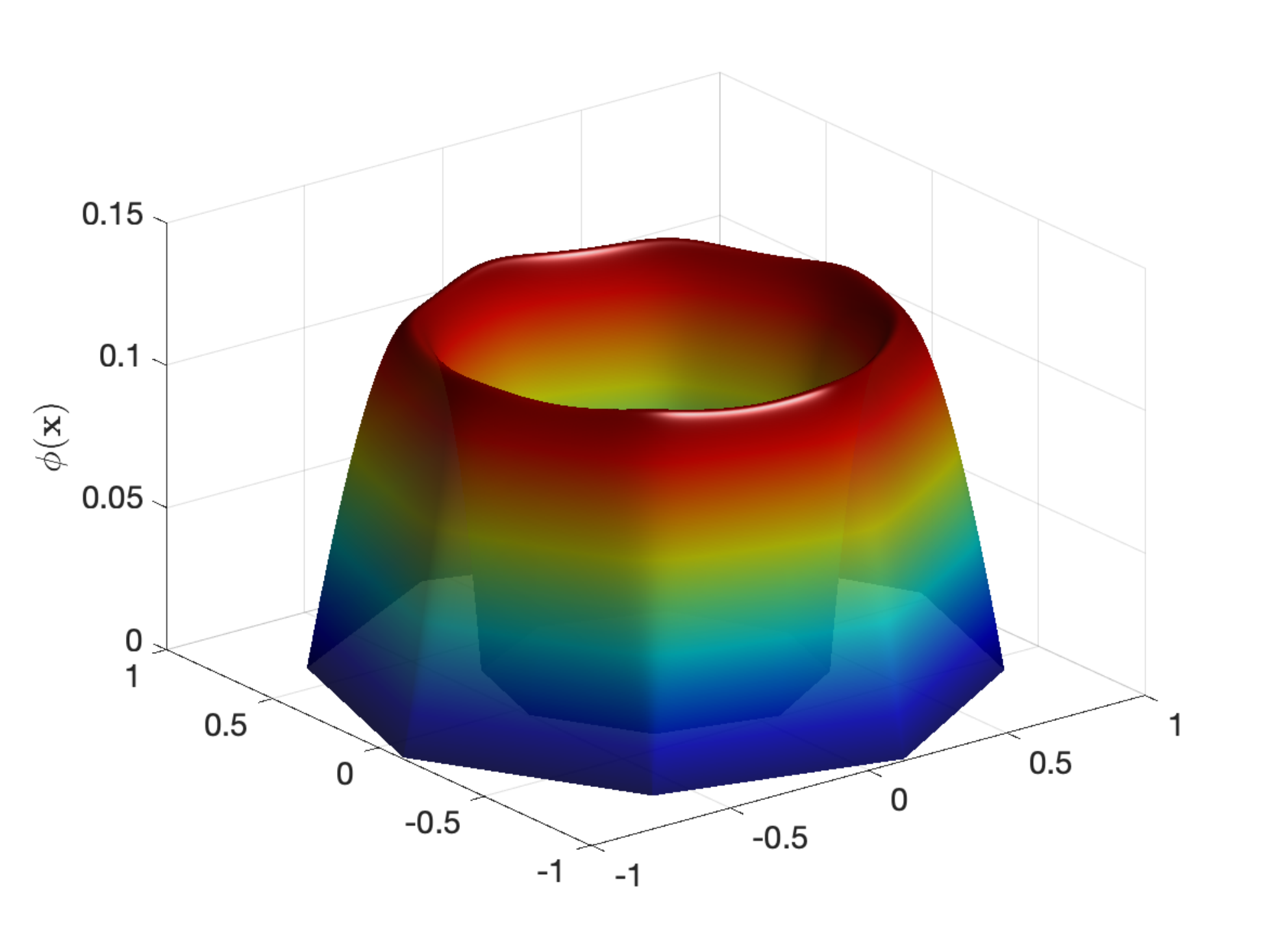}
    \includegraphics[width=0.43\textwidth]{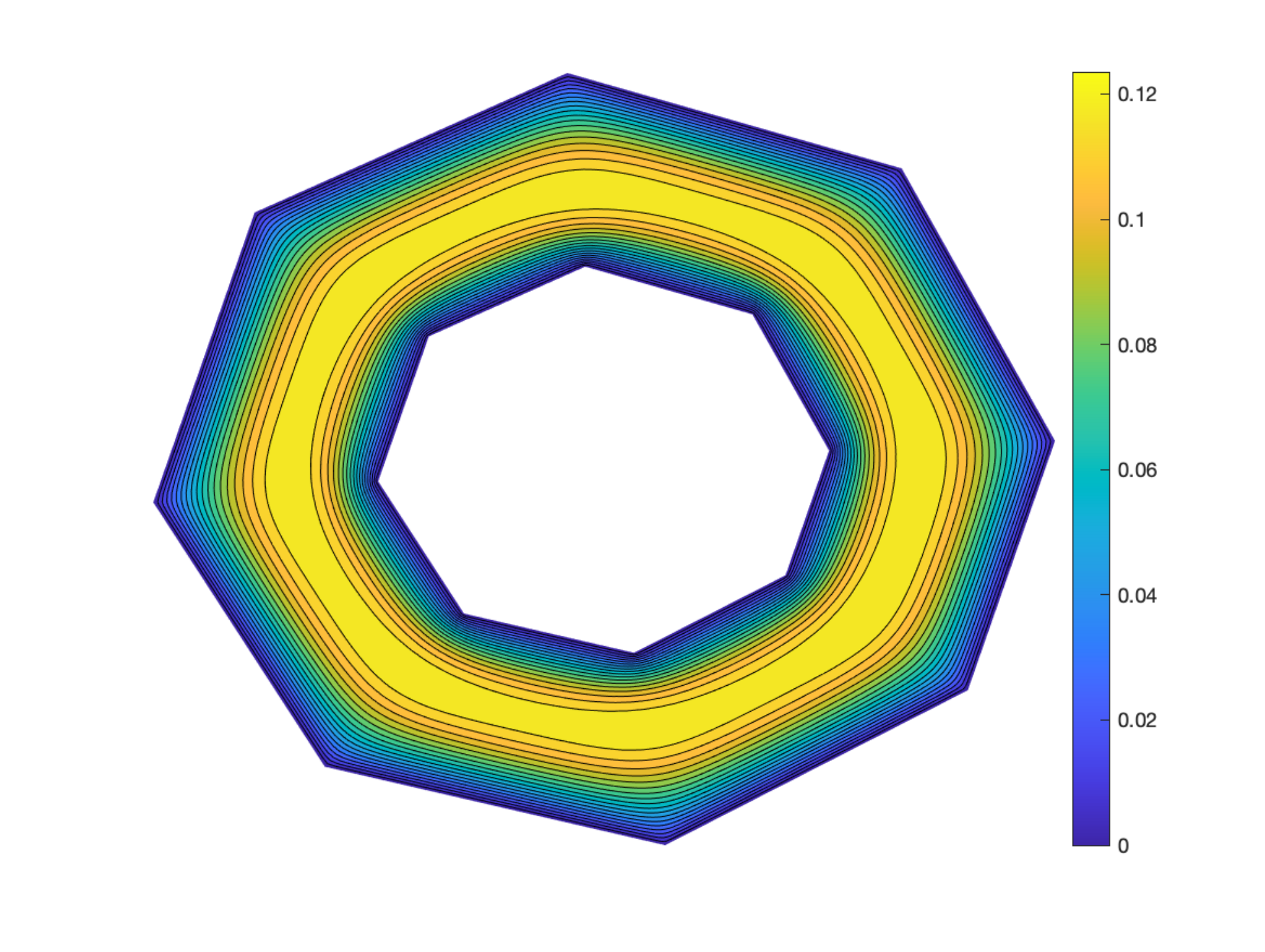}
    }
    \caption{}
    \label{fig:mvc_ADF-d}
    \end{subfigure}
 \caption{Approximate distance fields for polygonal domains.
          (a) square, (b) L-shaped, (c) nested squares, and (d) nested octagons.}
    \label{fig:mvc_ADF}
\end{figure}

\begin{figure}[!htb]
\centering
\begin{subfigure}{0.48\textwidth}
\includegraphics[width=\textwidth]{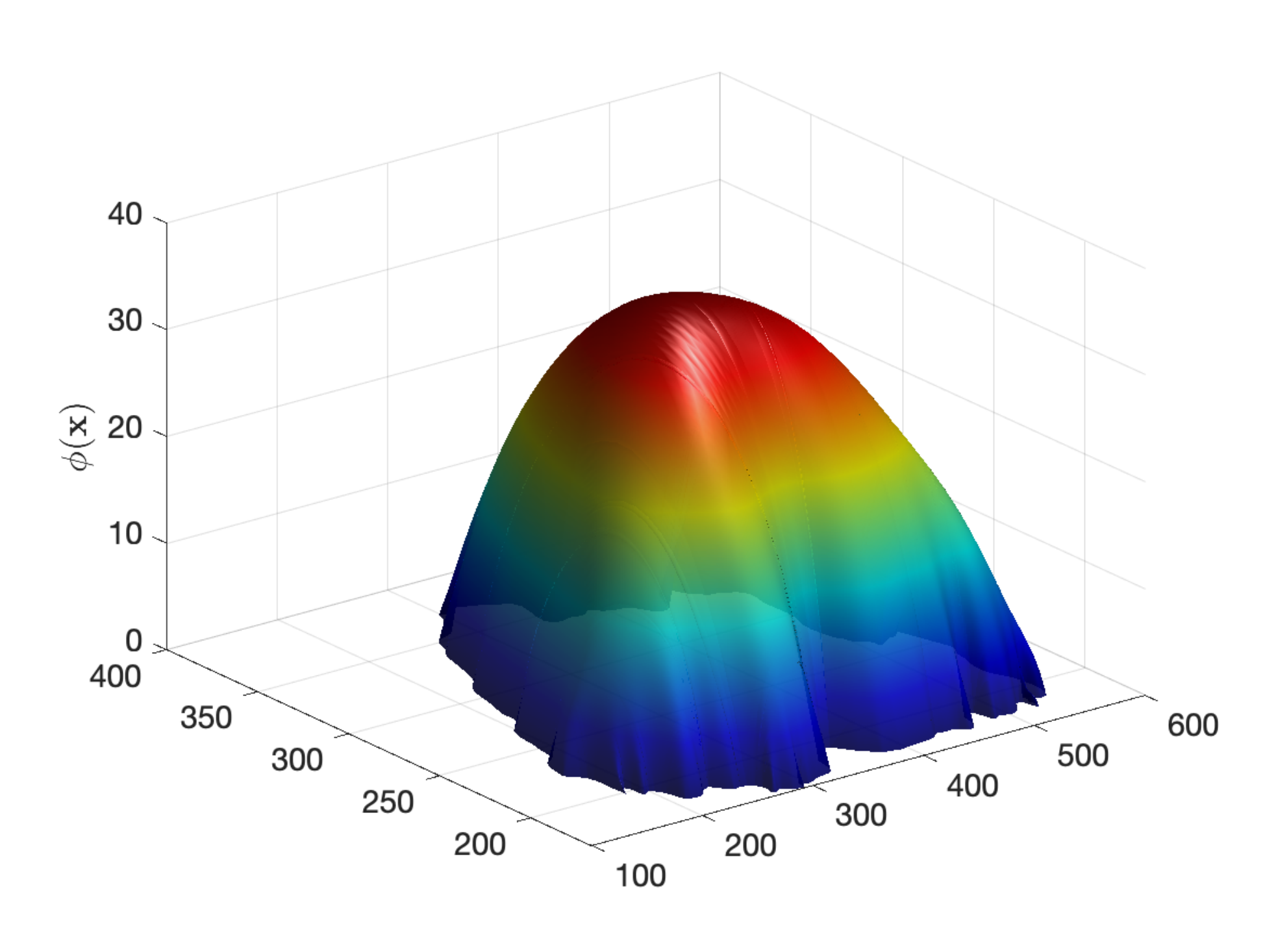}
 \caption{} 
\end{subfigure}
\begin{subfigure}{0.48\textwidth}
\includegraphics[width=\textwidth]{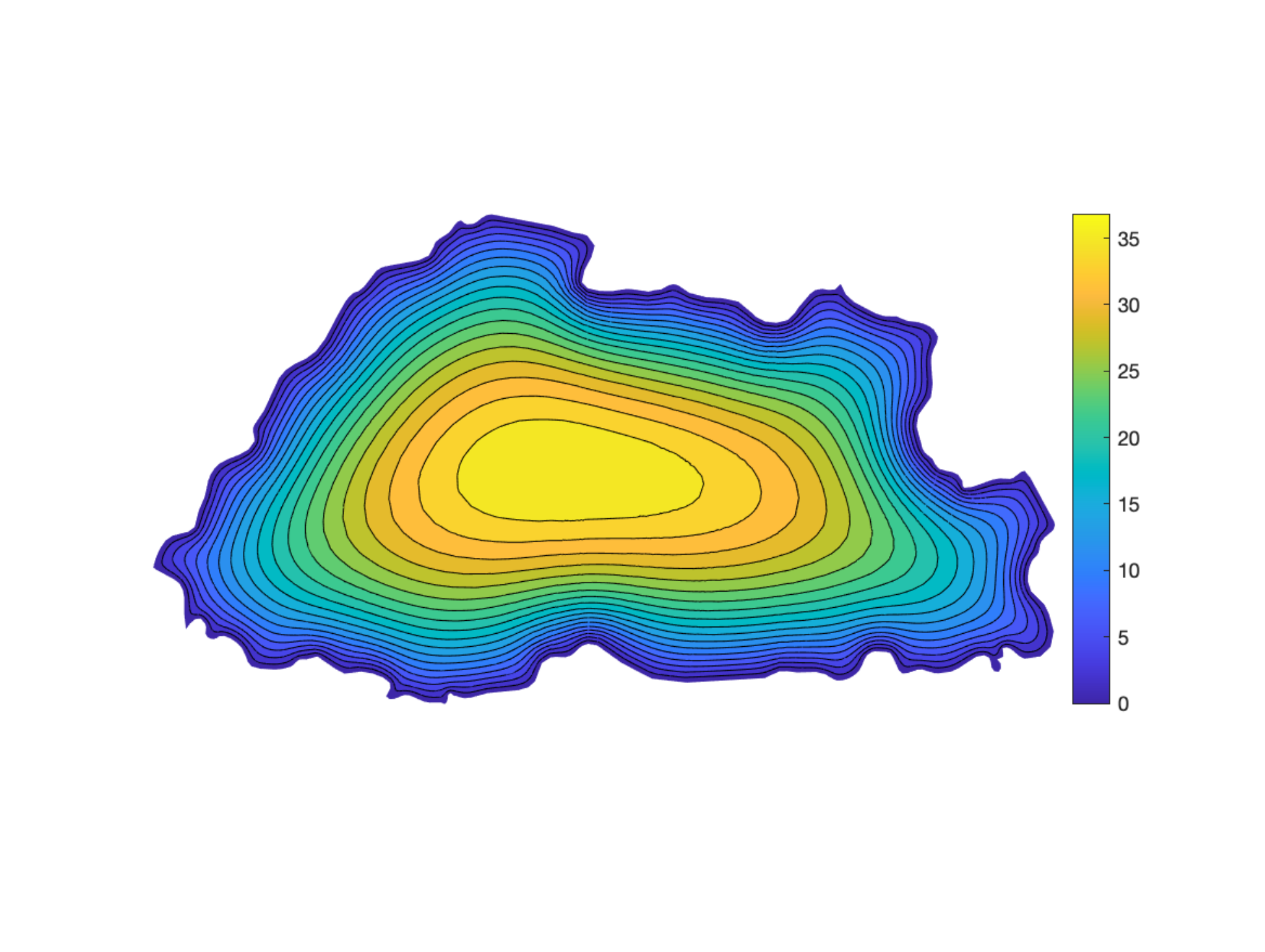}
\caption{} 
\end{subfigure}
\caption{Approximate distance field for the polygonalized 
         map of Bhutan. Surface plot is shown in (a) and the contour plot in (b). The polygon has 291 vertices.}
\label{fig:mvc_bhutan_map}
\end{figure}

For a polygon with one
interior (nested) $m$-gon,
the vertices of the inner polygon are defined in
clockwise orientation~\cite{Hormann:2006:MVC}.  The
contributions of $\{w_i(\vx)\}_{i=1}^{n+m}$ are used to form
$W(\vx)$ in~\eqref{eq:mvc}. Hormann and Floater~\cite{Hormann:2006:MVC}
showed that $\phi(\vx) = 1/W(\vx)$ is an ADF to the boundary
of the polygon, where its normal derivative is 1/2.  On taking
\begin{equation}\label{eq:phi_mvc}
\phi(\vx) = \frac{2}{ W (\vx)},
\end{equation}
where the
scaling factor is the volume of the unit sphere
in $\Re^{d-1}$~\cite{Bruvoll:2010:TMV} (equal to 2 when $d = 2$), the normal
derivative becomes $\partial \phi/\partial \nu = 1$ on $\partial \Omega$. As we \suku{discuss} in~\sref{subsec:mvp_curved}, the $W(\vx)$ that appears in~\eqref{eq:weight_mvc} and~\eqref{eq:phi_mvc}
is a particular instance ($p=1$) of the mean value potential field
$W_p(\vx)$.
We now have a smooth ADF for a polygon that is normalized to order 1. In~\fref{fig:mvc_ADF}, the surface and contour plots of $\phi(\vx)$ are shown for square and
L-shaped domains, as well as nested squares and octagons. The ADF using R-equivalence for the polygonalized map of Bhutan is presented in~\fref{fig:phi_bhutan_map}. In~\fref{fig:mvc_bhutan_map}, we show 
the approximate distance field (surface and contour plots) for the same polygonalized map.

\subsection{Approximate distance fields over curved domains}\label{subsec:mvp_curved}
Consider an open, bounded nonconvex domain
$\Omega$ with boundary $\Gamma = \partial \Omega$.  Given a function
$g : \Gamma \rightarrow \Re$ that is
prescribed on a curved boundary, the transfinite mean value interpolant $u :
\Omega \rightarrow \Re$ is defined as~\cite{Ju:2005:MVC,Dyken:2009:TMV}
\begin{equation}\label{eq:tmvi}
  u(\vx) = \frac{\int_{S_v} g \bigl( \vm{y} (\vx, \vm{v}) \bigr) K(\vx, \vm{y}) \, dS_v}
    {W(\vx)} , \quad
  W(\vx) = \int_{S_v} K(\vx, \vm{v}) \, dS_v , \quad
  K(\vx, \vm{v}) = \frac{1}{\|\vm{v} - \vx\|} ,
\end{equation}
where $\vx \in \Omega \backslash S_v$ and $\vm{v} \in S_v$, $S_v$ is the unit
circle that is centered at $\vx$, the ray from $\vx$ that passes through
$\vm{v}$ intersects the boundary $\Gamma$ at $\vm{y}$, and $K(\vx, \vm{v})$ is a
singular kernel function~\cite{Ju:2005:MVC,Dyken:2009:TMV}.

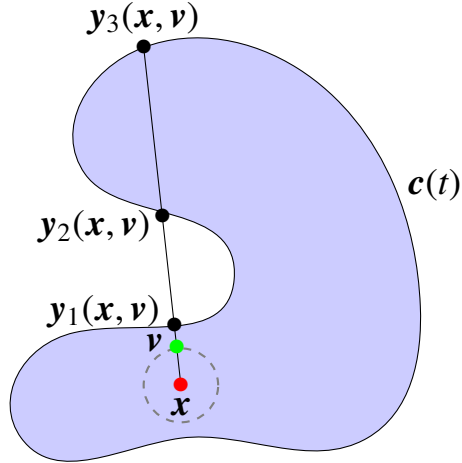
\begin{figure}[!htb]
  \centering
\begin{tikzpicture}[scale=0.7,every node/.style={scale=1.4},use Hobby shortcut]
  \path
  (0,0) coordinate (z0)
  (3,0) coordinate (z1)
  (6.5,0) coordinate (z2)
  (7.5,2.5) coordinate (z3)
  (7,4.8) coordinate (z4)
  (5.2,7) coordinate (z5)
  (2.3,7.4) coordinate (z6)
  (1.2,5) coordinate (z7)
  (4,3) coordinate (z8)
  (1,2) coordinate (z9);

\draw[closed,fill=blue!20] (z0) .. (z1) .. (z2) .. (z3) .. (z4) .. (z5) .. (z6) .. (z7) .. (z8) .. (z9);

\def\mycircle{(3,1) circle (0.7 cm)}
\draw[color=black!50,thick,dashed] \mycircle;

\coordinate (O) at (3,1);
\coordinate (Y3) at (z6);

\draw (O) -- (Y3)
      node[pos=0,color=red] {$\bullet$}
      node[pos=0,below] {$\bm{x}$}
      node[pos=0.11,color=green] {$\bullet$}
      node[pos=0.125,left] {$\bm{v}$}
      node[pos=0.175] {$\bullet$}
      node[pos=0.22,left] {$\bm{y}_1(\bm{x},\bm{v})$}
      node[pos=0.50] {$\bullet$}
      node[pos=0.47,left] {$\bm{y}_2(\bm{x},\bm{v})$}
      node[pos=1] {$\bullet$}
      node[pos=1,above] {$\bm{y}_3(\bm{x},\bm{v})$};

\coordinate (Q) at (z4);
\draw (Q) node[right] {$\bm{c}(t)$};

\end{tikzpicture}
  \caption{Nonconvex domain bounded by the curve $\vm{c}(t)$. The 
           variables that appear in~\eqref{eq:tmvi} and~\eqref{eq:phi_tmvi} are shown.}
  \label{fig:nonconvexdomain-mvp}
\end{figure}

\begin{figure}[!htb]
\centering
\begin{subfigure}{0.49\textwidth}
\mbox{
\includegraphics[width=0.56\textwidth]{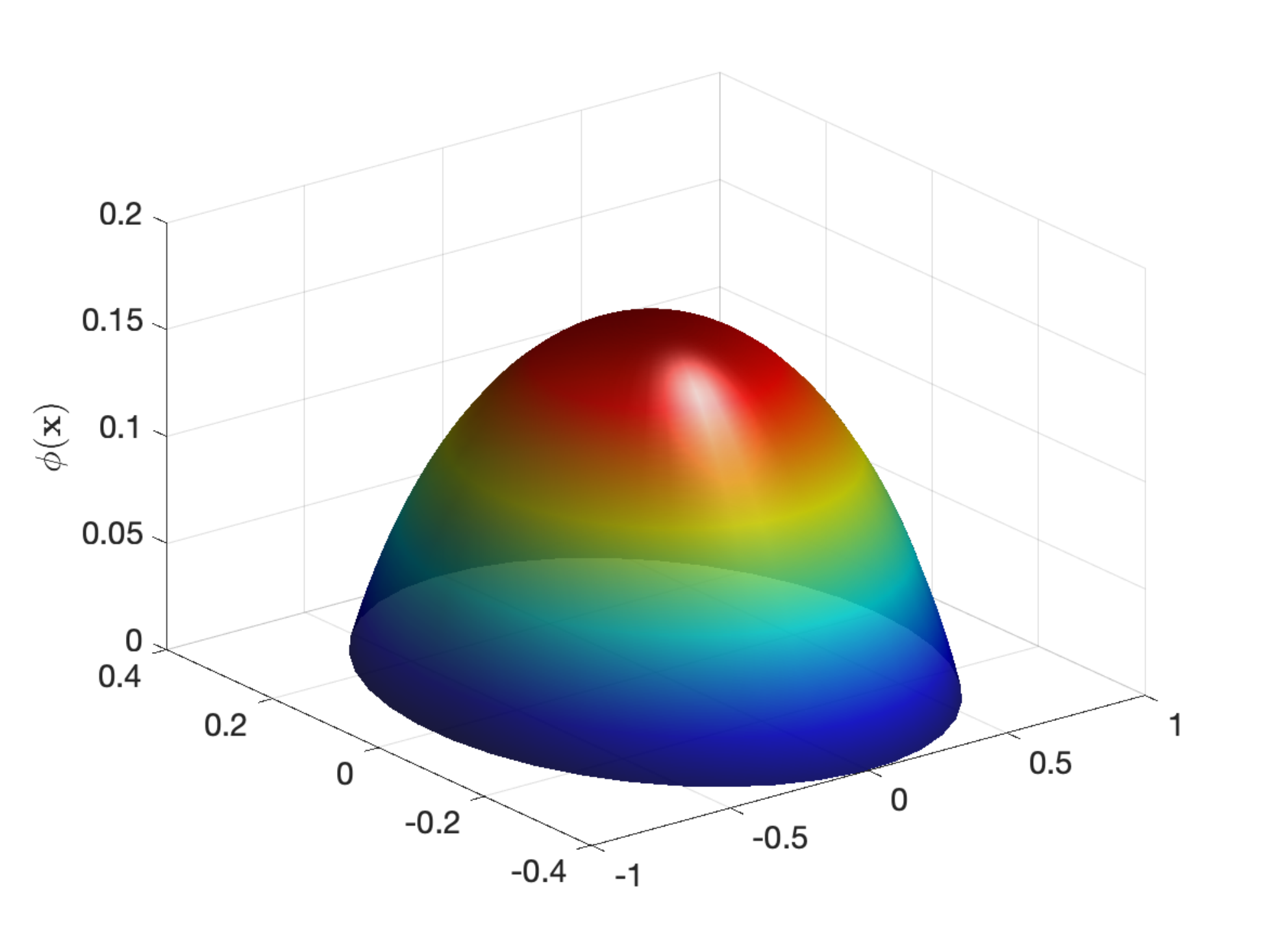}
\includegraphics[width=0.40\textwidth]{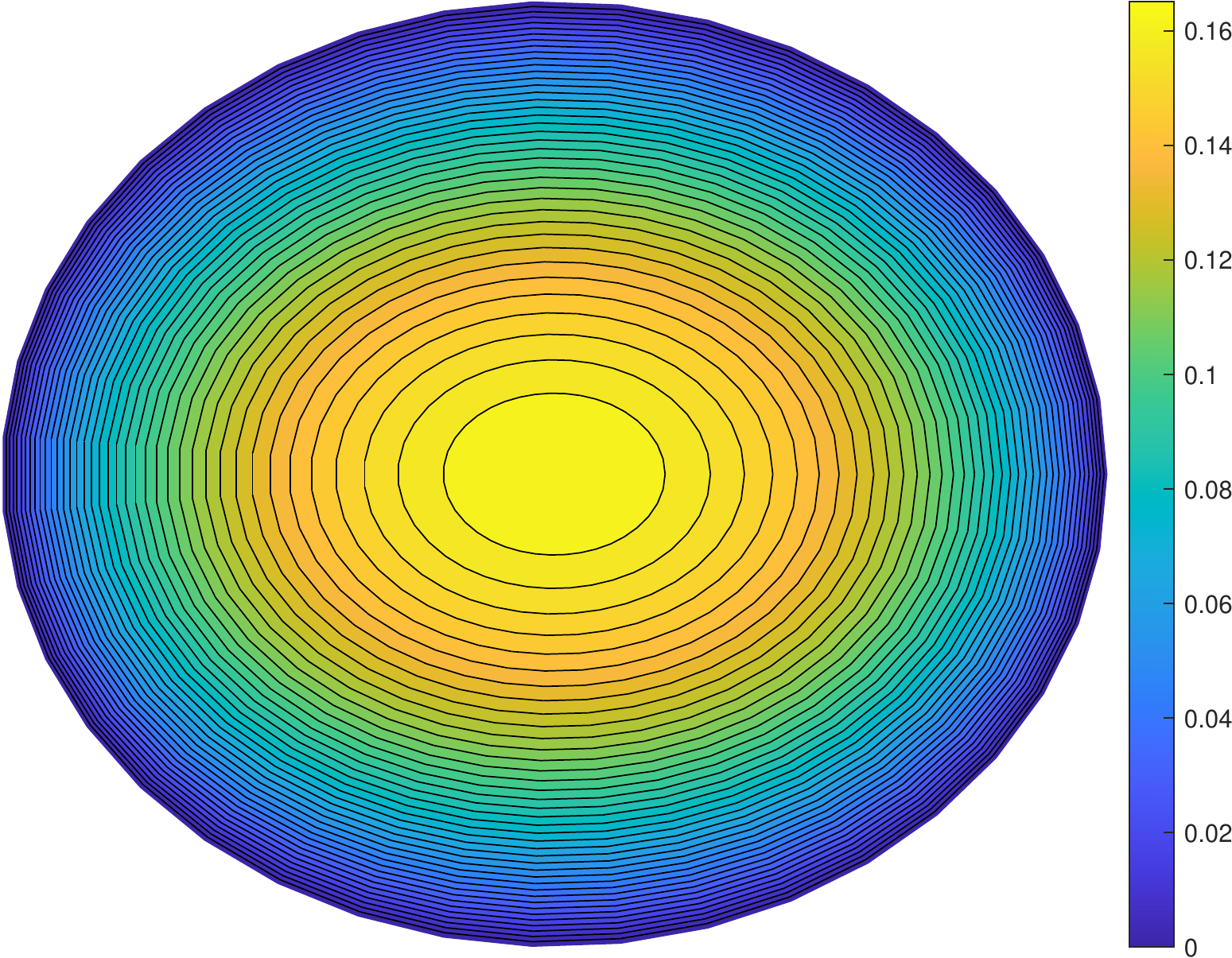}
}
\caption{}
\label{fig:tmvi_ADF-a}
\end{subfigure}
\begin{subfigure}{0.49\textwidth}
\mbox{
\includegraphics[width=0.56\textwidth]{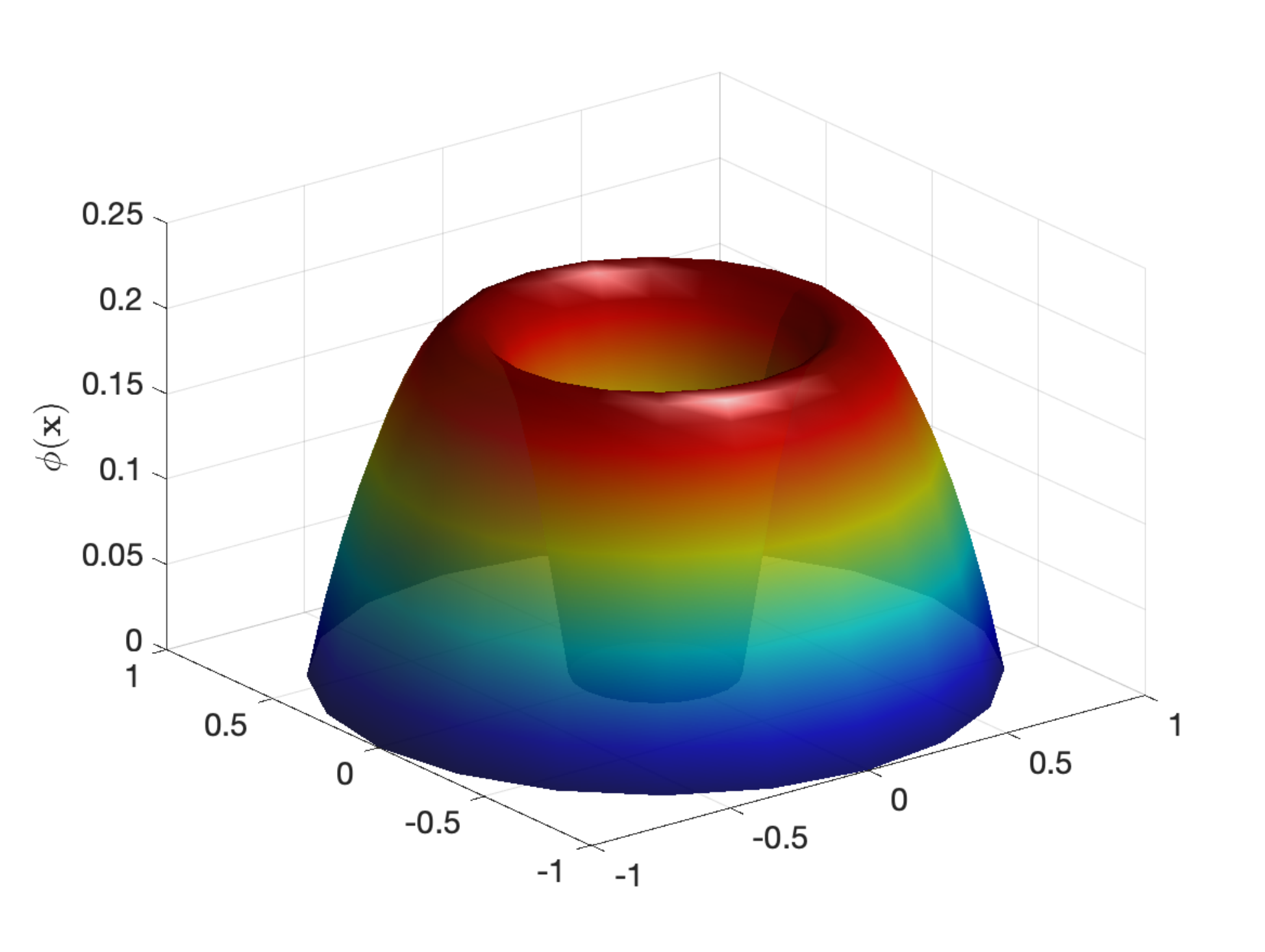}
\includegraphics[width=0.40\textwidth]{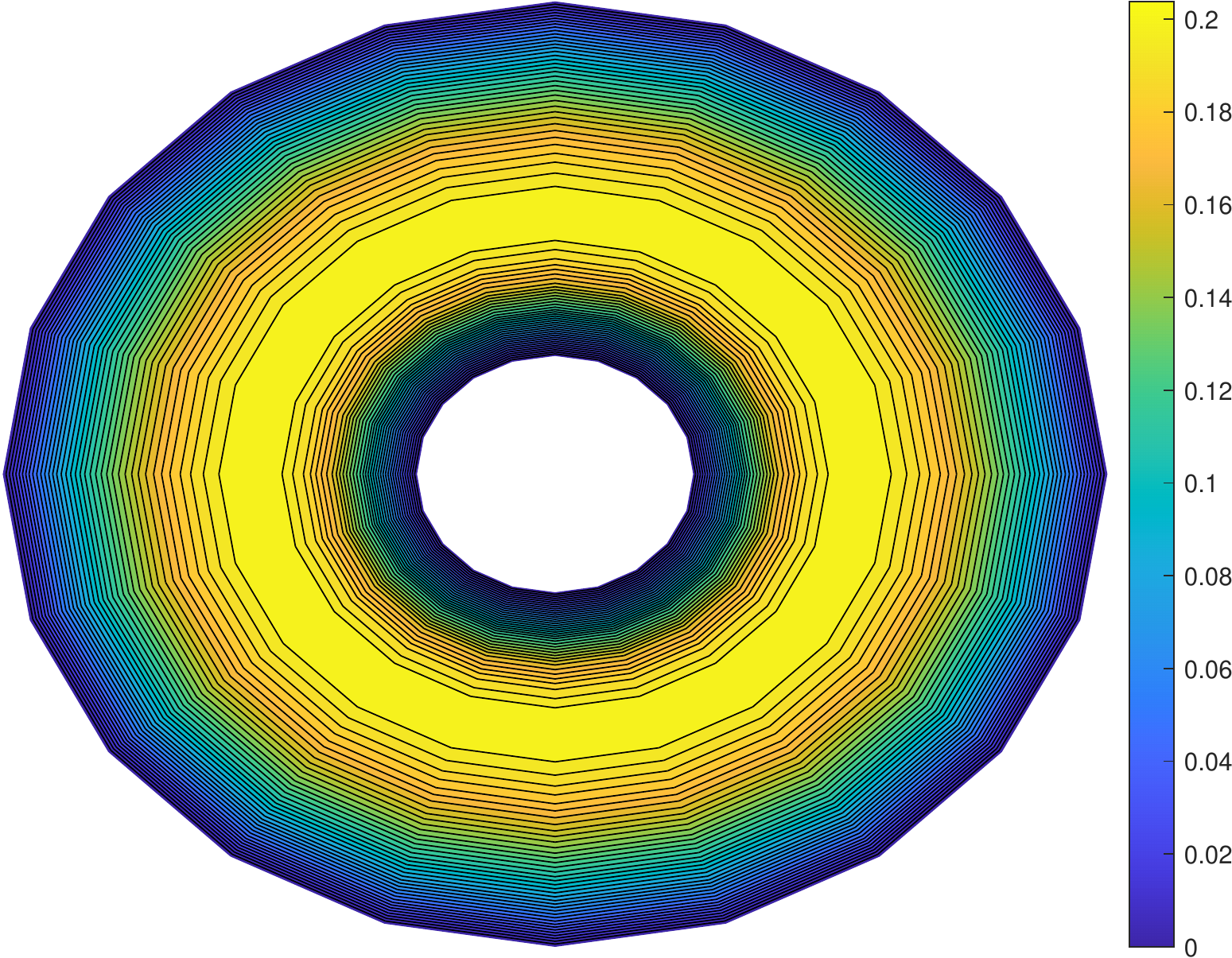}
}
\caption{}
\end{subfigure}
\begin{subfigure}{0.49\textwidth}
\mbox{
\includegraphics[width=0.53\textwidth]{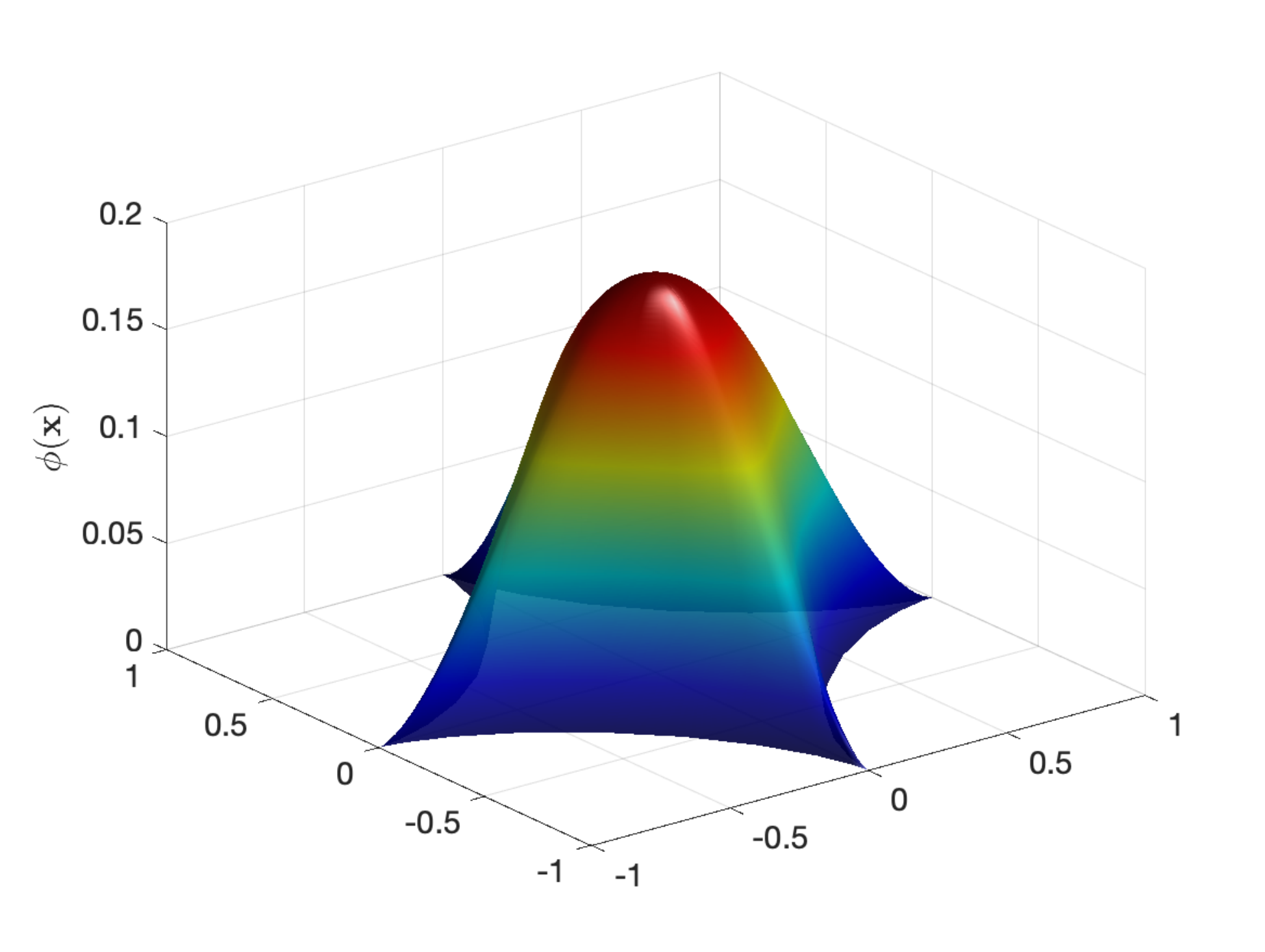}
\includegraphics[width=0.43\textwidth]{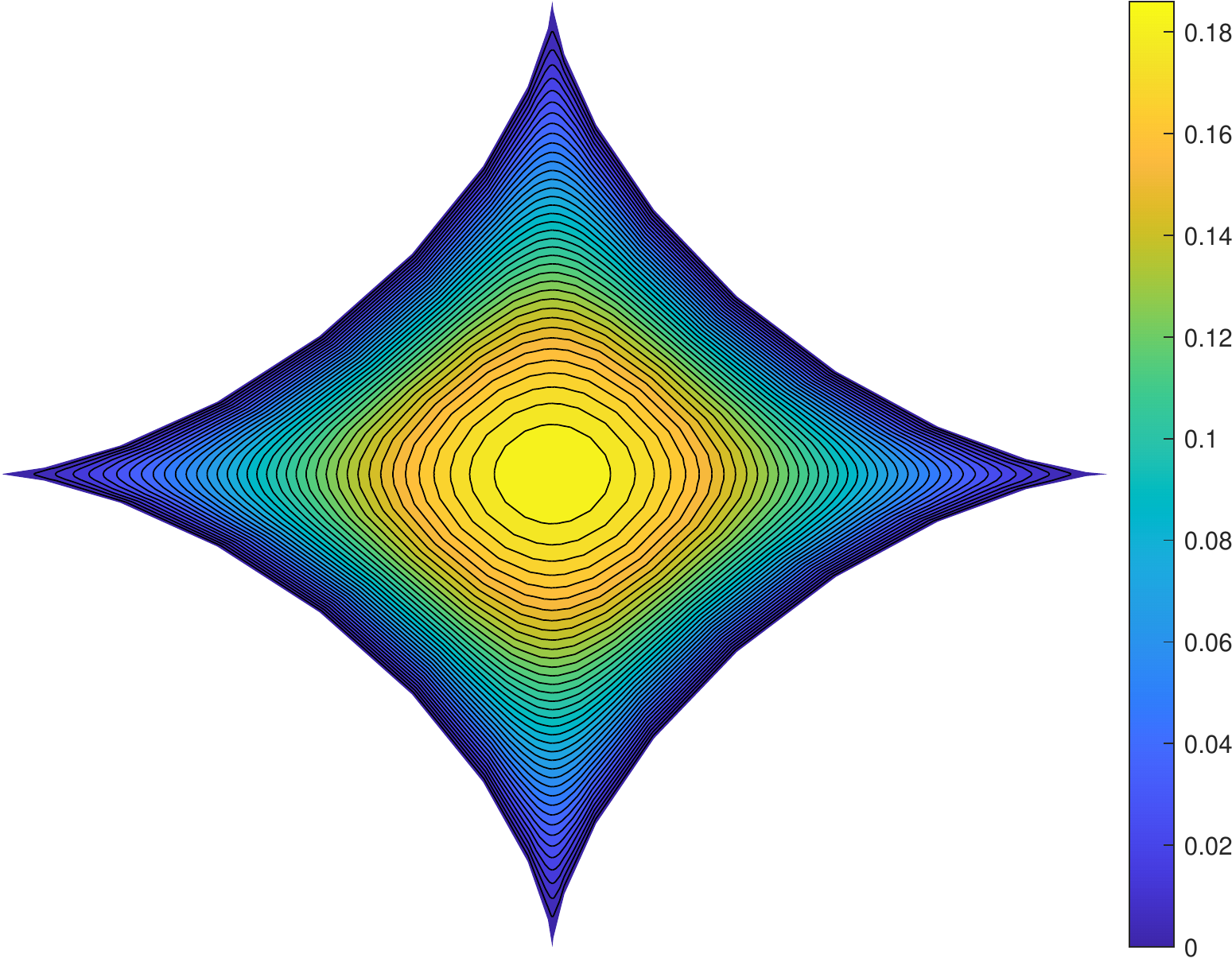}
}
\caption{}
\end{subfigure}
\begin{subfigure}{0.49\textwidth}
\mbox{
\includegraphics[width=0.53\textwidth]{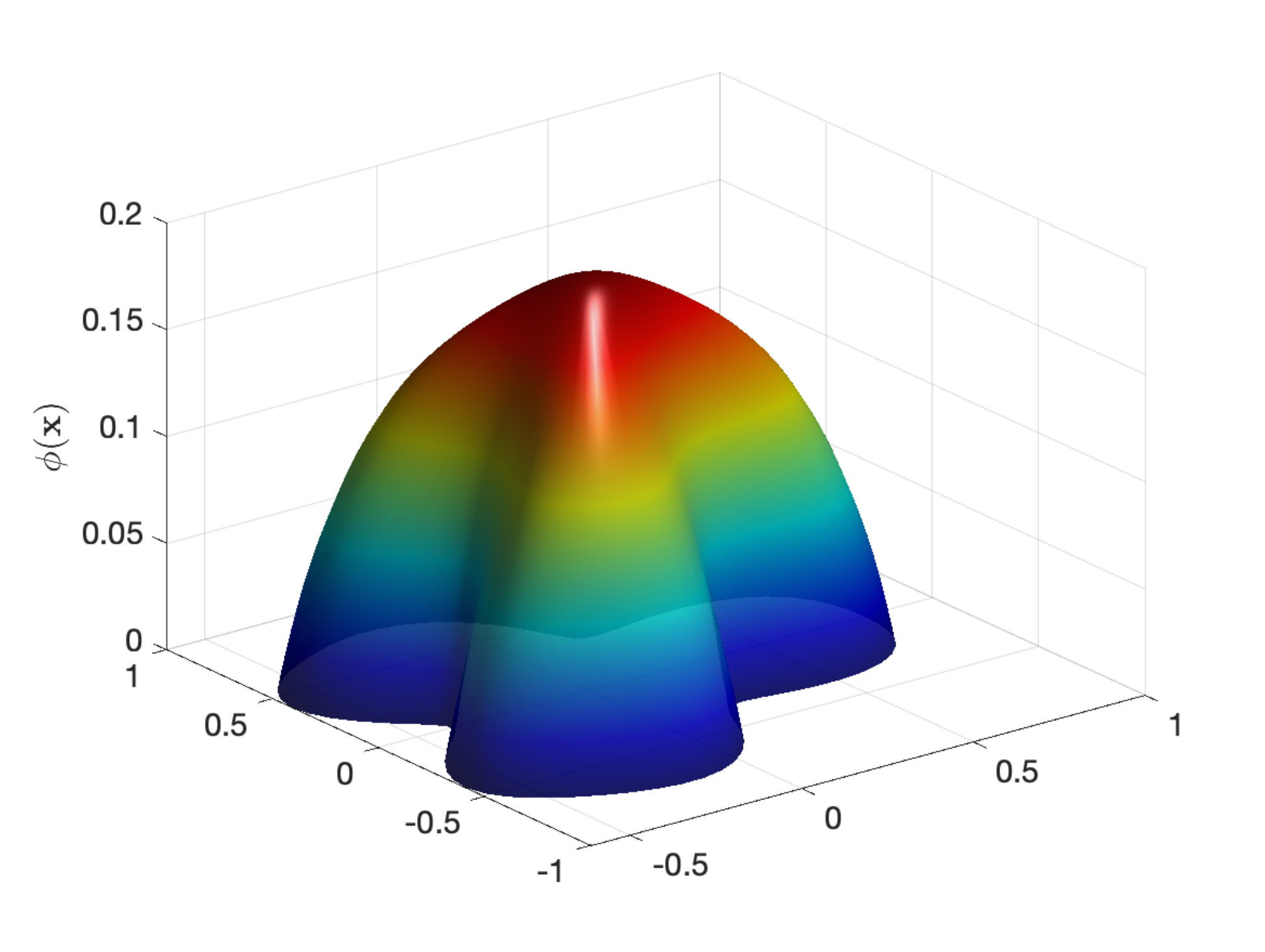}
\includegraphics[width=0.43\textwidth]{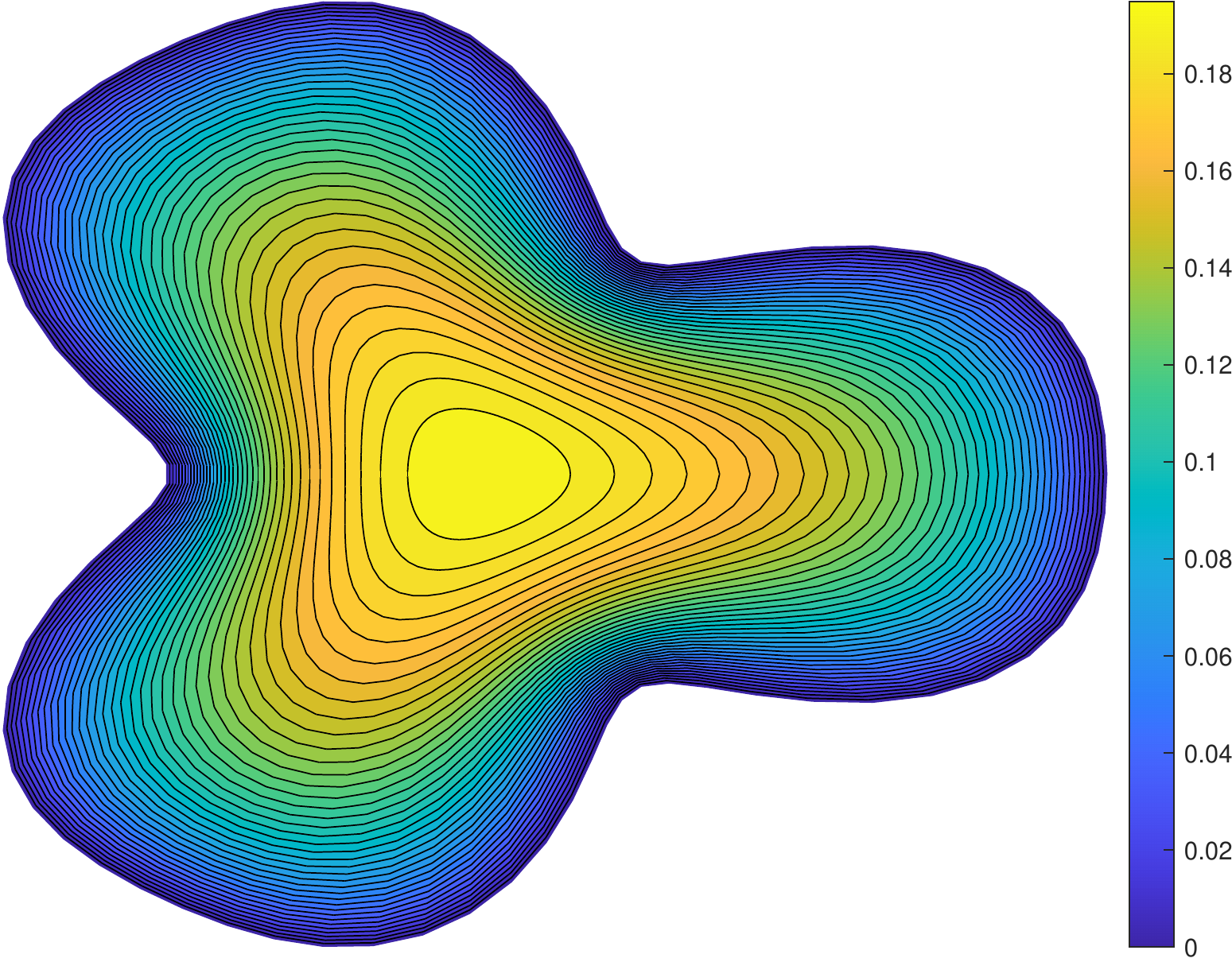}
}
\caption{}
\end{subfigure}
\caption{Approximate distance fields on curved domains using $\phi(\vx) = 2 /W_1(\vx)$, with $W_1(\vx)$ given in~\eqref{eq:phi_tmvi}. Surface and contour plots are shown for an (a) elliptical disk, (b) annulus, (c) hypocycloid, and (d) propeller~\cite{Chin:2019:MCI}.}
\label{fig:tmvi_ADF}
\end{figure}

Similar to the behavior of the inverse weight in mean value coordinates on polygons,
the function $\phi(\vx) = 1 / W(\vx)$ behaves like an approximate
distance function to the boundary and its normal derivative on the boundary
$\Gamma$ is $1/2$~\cite{Dyken:2009:TMV}.
Belyaev et al.~\cite{Belyaev:2013:SLD} introduced $L_p$-distance fields ($p \geq
1$), which approximates the exact distance function. These distance fields stem from a particular form of a singular double-layer potential, and hence the reference to them as generalized mean value potential fields. 
Consider the nonconvex domain shown in~\fref{fig:nonconvexdomain-mvp}. 
The parametrization of the curved boundary $\Gamma : [0,1] \to \Re$ is $\vm{c}(t)$, and its tangent is 
$\vm{c}^\prime (t)$. For a nonconvex domain, a ray from $\vx$ intersects the boundary at multiple points $\vm{c}(t_i) :=
\vm{y}_i(\vx,\vm{v})$. 
On projecting the boundary curve
onto the unit circle, an expression for $\phi(\vx)$ that is valid for convex as well as nonconvex domains is obtained
in terms of the
curve parameter $t \in [0,1]$~\cite{Ju:2005:MVC}:
\begin{equation}\label{eq:phi_tmvi}
  \phi (\vx) = \left( \frac{1}{W_p (\vx)} \right)^{1/p}, \qquad
  W_p(\vx) = \int_0^1 \frac{ \bigl( \vm{c}(t) - \vx \bigr) \cdot 
    \vm{c}^{\prime\perp} (t) }{ \| \vm{c}(t) - \vx \|^{2+p} } dt,
\end{equation}
where
$\vm{c}^{\prime\perp}(t) := \text{rot} \bigl( \vm{c}^\prime (t) \bigr)$ is
obtained by rotating $\vm{c}^\prime (t)$ through $90\degree$ in the clockwise direction.  For $\vx \in \Omega$, the
integral in~\eqref{eq:phi_tmvi} is numerically integrated; if $\vx \in \partial \Omega$ (integral is singular), we set $\phi(\vx) = 0$. 
In~\eqref{eq:phi_tmvi}, $W_p(\vx)$ is the generalized
mean value potential field, which is used to form the approximate distance function $\phi(\vx)$.
Equation~\eqref{eq:phi_tmvi} is also applicable for polygons: on
choosing $p = 1$ in~\eqref{eq:phi_tmvi}, we recover the $W(\vx)$ that appears in~\eqref{eq:weight_mvc}. The approximate distance fields ($p = 1$)
for an elliptical disk, annulus, hypocycloid, and a propeller-shaped domain are shown in~\fref{fig:tmvi_ADF}. The distance function is smooth in the interior of the domain, $\phi \in C^\infty(\Omega)$, and it is $C^k$ on the boundary for a $C^k$ curve (derivative discontinuities occur at the vertices for a polygonal curve). Over curved two-dimensional domains, Dyken and Floater~\cite{Dyken:2009:TMV} assessed the approximation properties of the transfinite mean value interpolant as well as its use to solve the Poisson equation with web-splines~\cite{Hollig:2001:WEB}, and Chin and Sukumar~\cite{Chin:2021:SBC} have used it in verification tests of a cubature rule for numerical integration  over curved regions.

\section{Imposing Boundary Conditions in Deep Neural Networks}\label{sec:BCs}
\suku{Let $\unnbcxt$ denote the PINN trial function. We
present the construction of $\unnbcxt$ so that it exactly satisfies all essential and Robin (natural boundary condition is a particular case) boundary conditions.}
The trial function includes the contribution from the neural network approximation, $\unnRxt$, where $\vm{\theta}$ contains the unknown parameters of the network. We defer the presentation of $\unnRxt$ to~\sref{sec:dnn}. 

\subsection{Normalizing functions and solution structures}
\label{subsec:solution_structures}
Let $\Omega \subset \Re^2$ be an open, bounded domain with boundary
$\partial \Omega$. Consider a smooth 
function $u(\vx): \Re^2 \to \Re$ and let
$\phi(\vx): \Re^2 \to \Re$ 
be an approximate distance function to $\partial \Omega$ 
that is normalized
to the $m$-th order (see~\sref{sec:distance}). Rvachev defined the normalizers to $u$ of the $m$-th order with respect to $\phi(\vx)$ via the transformation~\cite{Rvachev:1995:RBV,Rvachev:2001:TII}
\begin{equation}\label{eq:u_transformation}
u^*(\vx) = u(\vx - \phi \nabla \phi) ,
\end{equation}
which leads to
$u^*(\vx) = u(\vx)$ on $\partial \Omega$ and
$\partial^k u^* / \partial \nu^k = 0$ $(k = 1,2,\ldots,m)$
on $\partial \Omega$ .
Since we are treating Dirichlet and Robin boundary conditions in this paper, we proceed to show that~\eqref{eq:u_transformation} is normalized up to order 1.  If $u(\vx)$
is specified on $\partial \Omega$, then since $\phi = 0$ on $\partial \Omega$ it implies that $u^*(\vx) = u(\vx)$ on
$\partial \Omega$, which establishes zeroth-order normalization. The proof for first-order normalization follows.
\begin{proof}
Let $\vm{t} := \vx  - \phi \nabla \phi$.  Then, by the chain rule
we have
\begin{equation*}
\nabla u^* = \nabla_{\vm{t}}  u^* \cdot 
                      \nabla \otimes \vm{t} =
              \nabla_{\vm{t}} u^* \cdot [ \vm{I} - \phi \nabla \otimes
               \nabla \phi - \nabla \phi \otimes \nabla \phi ],
\end{equation*}
where $\otimes$ is the dyadic (tensor) product and $\nabla_{\vm{t}} (\cdot)$ is the gradient with respect to $\vm{t}$.  On $\partial \Omega$, $\phi = 0$, $\vm{t} = \vm{x}$ and $\nabla \phi = \inn$, since $\phi$ is normalized to the first order. Hence, we can write
\begin{equation*}
\left. \frac{\partial u^*}{\partial \nu} \right|_{\partial \Omega} =
    \left[ \nabla u^* \cdot \inn \right]_{\partial \Omega} = 
    \left[ \nabla u^* \right]_{\partial
    \Omega} \cdot 
    \biggl( \bigl[ \vm{I} - \phi \nabla \otimes \inn - 
        \inn \otimes \inn \bigr] \cdot \inn \biggr)_{\partial \Omega} =
    \left[ \nabla u^* \right]_{\partial
    \Omega} \cdot ( \inn - \inn) = 0,
\end{equation*}
which is the desired result. \qed
\end{proof}
Note that this result also holds if we consider the unit outward
normal vector, $\outn = -\inn$. 
By extension, if $\phi$ is normalized to the $m$-th order, one can establish that $u^*$ in~\eqref{eq:u_transformation} is normalized to the $m$-th order, i.e.,
\begin{equation}
    u^*(\vx) = u(\vx) \ \ \textrm{on } \partial \Omega, \quad
    \frac{\partial^k u^*}{\partial \nu^k} = 0 \ \
    \textrm{on } \partial \Omega \ \ (k = 1,2,\ldots,m).
\end{equation}
Let $u$ and its higher order normal derivatives along
$\inn$ be prescribed on $\partial \Omega$, i.e.,
\begin{equation}\label{eq:gen_bcs}
    u(\vx) = u_0(\vx) \ \ \textrm{on } \partial \Omega, \quad
    \frac{\partial^k u}{\partial \nu^k} = u_k(\vx)  \ \
    \textrm{on } \partial \Omega \ \ (k = 1,2,\ldots,m).
\end{equation}
Then, one can represent $u$ in the vicinity of $\partial \Omega$ and in the direction (inward) normal to the boundary using a polynomial Taylor series expansion of $u$ in terms of
$\phi$.  Rvachev et al.~\cite{Rvachev:1995:RBV,Rvachev:2000:OCR} referred to this as a generalized Taylor series expansion, which takes the form:
\begin{equation}\label{eq:gen_Taylor}
u(\vx) = u_0^*(\vx)
    + \sum_{k=1}^m \frac{u_k^*(\vx) }{k!} \, \phi^k(\vx)
    +  \phi^{m+1}(\vx) \, \Psi(\vx) ,
\end{equation}
where $u_k^*(\vx) = u_k(\vx - \phi \nabla \phi)$
and $\Psi(\vx)$ is an unknown function (approximation) in the remainder term. This equation resembles univariate Taylor series expansion about
$x = 0$, where instead of evaluating derivatives that are scalar
constants,
$u_k^*(\vx)$ in~\eqref{eq:gen_Taylor} are scalar-valued functions. Equation~\eqref{eq:gen_Taylor} is the general form
of the {\em solution
structure}  for $u$ that Rvachev introduced~\cite{Rvachev:1995:RBV,Rvachev:2000:OCR}. If a function
$u \in C^k(\Omega)$, and its first $m$ derivatives vanish on $\partial \Omega$, then the solution structure $\phi^{m+1} \Psi$
is sufficiently complete in the Hilbert space
$H^k(\Omega)$ to approximate $u$ and all its derivatives up to order $k$~\cite{Rvachev:2000:OCR}. 
On using~\eqref{eq:gen_Taylor}, the boundary conditions on
$\partial \Omega$ that are given in~\eqref{eq:gen_bcs} are exactly met.  
We now present the solution structure for three distinct sets of 
boundary conditions (Dirichlet, Neumann and Robin) that are imposed
on $\partial \Omega$.

\subsubsection{Solution structure for Dirichlet boundary condition}
\label{subsubsec:Dirichlet}
If $u = g$ is prescribed on $\partial \Omega$, then using 
$u_0^*(\vx) = g(\vx - \phi \nabla \phi)$
and $m = 0$ in~\eqref{eq:gen_Taylor}, we can write:
\begin{equation}
\begin{split}
    u(\vx) &= g(\vx - \phi \nabla \phi) 
    + \phi(\vx) \Psi(\vx) 
    = g(\vx) - \phi(\vx) 
    \left. \frac{\partial g(\vx) }{\partial \nu} \right|_{\partial \Omega}
    + \phi(\vx) \Psi(\vx) = g(\vx) + \phi(\vx) 
    \left. \frac{\partial g(\vx) }{\partial n} \right|_{\partial \Omega}
    + \phi(\vx) \Psi(\vx) \\
    &= g(\vx) + \phi(\vx) \Psi(\vx),
    \end{split}
\end{equation}
since $\outn = - \inn$ on $\partial \Omega$, and
a first order linearization of $g(\cdot)$ is used. Therefore, the solution structure is:
\begin{equation}\label{eq:Dirichlet}
    u = g + \phi \tilde{u} ,
\end{equation}
where  
$\tilde{u}$ is any suitable numerical approximation. So any
trial function of the form given in~\eqref{eq:Dirichlet} will
exactly satisfy the essential boundary condition $u = g$ on
$\partial \Omega$. In this paper,
we use deep neural networks to construct $\tilde{u}$.
\suku{For the homogeneous Dirichlet problem ($g = 0$) using
Kantorovich's method, 
Babu\v{s}ka et al.~\cite{Babuska:2003:SMG} provide an a priori error estimate in the $H^1$ norm.}

\subsubsection{Solution structure for Neumann boundary condition}
Since Neumann and higher order boundary conditions for PDEs are imposed in the direction of $\outn$, the unit outward normal vector to
$\partial \Omega$, we consider solution structures that are defined with respect to $\outn$. This is a departure from the literature on the R-function method to solve PDEs, but is aligned with how boundary-value problems are posed. 

For first order derivatives, linearizing $u^*(\vx) := u(\vx - \phi \nabla \phi)$ in the neighborhood of the boundary ($\phi = 0$) amounts to subtracting the variation in the normal direction
$\inn$, which leads to~\cite{Rvachev:1995:RBV}
\begin{equation}\label{eq:linearization}
\begin{split}
    u^*(\vx) &= 
    u(\vx) - \left[ \phi(\vx) \nabla u(\vx) \cdot \inn
    \right]_{\partial \Omega}  + \phi^2(\vx) \Psi(\vx) =
    u(\vx) - \left[ \phi(\vx) \nabla \phi(\vx) \cdot \nabla u (\vx)
    \right] _{\partial \Omega} + \phi^2(\vx) \Psi (\vx) \\
    &= (1 + \phi D_1^\phi) (u) + \phi^2 \Psi , \quad 
        D_1^\phi (\cdot) := \left[ -
     \nabla \phi \cdot \nabla (\cdot) \right]_{\partial \Omega}
     =  \left. \frac{\partial (\cdot)}{\partial n} 
         \right|_{\partial \Omega},
\end{split}
\end{equation}
where $\inn = - \outn = \nabla \phi$ on $\partial \Omega$
and $D_1^\phi(\cdot)$ is a differential operator that acts in the
outward normal direction to the boundary.

If
$\partial u / \partial n = h$ is prescribed on $\partial \Omega$, then using~\eqref{eq:gen_Taylor} and~\eqref{eq:linearization}, we can write
\begin{equation*}
\begin{split}
    u(\vx) &= u_0\bigl(\vx - \phi(\vx) \nabla \phi(\vx) \bigr) +
              \phi(\vx) u_1\big( \vx - \phi(\vx) \nabla \phi(\vx)
              \bigr) + \phi^2(\vx) \Psi_1(\vx)  \\
              &= [1 + \phi D_1^\phi] (u_0)  +
                 \phi u_1 +
                 \phi^2 \Psi(\vx)
              = [1 + \phi D_1^\phi] (u_0) - \phi h
               + \phi^2 \Psi ,
\end{split}
\end{equation*}
since $\partial u / \partial \nu = u_1 = -h$ on $\partial \Omega$, 
and therefore the most general
solution structure can be written as:
\begin{equation}\label{eq:Neumann}
    u = [1 + \phi D_1^\phi ] (\tilde{u}_1) - \phi h + \phi^2
        \tilde{u}_2,
\end{equation}
where $\tilde{u}_1$ and $\tilde{u}_2$ are arbitrary approximation functions.

\subsubsection{Solution structure for Robin boundary condition} 
If the Robin boundary condition,
$\partial u / \partial n + cu = h$, is prescribed on $\partial \Omega$ $\bigl(c := c(\vx)$, $h := h(\vx) \bigr)$, then following similar steps to that taken to obtain~\eqref{eq:Dirichlet} and~\eqref{eq:Neumann}, we can write the most general
solution structure for this boundary condition as~\cite{Rvachev:1995:RBV}:
\begin{equation}\label{eq:Robin}
  u = \left[ 1 + \phi \bigl(c + D_1^\phi \bigr) \right]
      (\tilde{u}_1) 
      - \phi h + \phi^2 \tilde{u}_2 .
\end{equation}
It is readily verified that the ansatz in~\eqref{eq:Robin}
satisfies the boundary condition
$\partial u / \partial n  + cu = h$ on $\partial \Omega$.

\subsection{Imposing inhomogeneous essential boundary conditions}
\label{subsec:EBCs}
Let us extend our analysis to the case when different
inhomogeneous essential boundary conditions are imposed on distinct subsets of $\partial \Omega$. Now, let the boundary
$\partial \Omega :=
\Gamma = \cup_{i=1}^N \Gamma_i$. The inhomogeneous essential boundary condition $u = g_i$ is imposed on
$\Gamma_i$ ($i = 1,2,\dots,M$), and on $\Gamma_i$ ($i = M+1,M+2,\ldots,N$) we assume natural boundary conditions are 
imposed through the potential energy functional in the
variational principle. 
Let $\phi_i$ be the ADF that is associated
with $\Gamma_i$ ($i = 1,2,\dots,M$), and let $\phi$ be
the approximate distance field that is composed either via R-equivalence using
$\phi_1 \sim \phi_2 \dots \phi_M$ in~\eqref{eq:phin_eq} or using the mean value potential field, $W(\vx)$, given in~\eqref{eq:phi_mvc} and~\eqref{eq:phi_tmvi}. Transfinite interpolation is the generalization of scattered data interpolation to interpolation of functions over curves and surfaces. On using the singular inverse-distance based Shepard weight
function~\cite{Shepard:1968:ATD}, we can write the transfinite 
interpolant as~\cite{Rvachev:1995:RBV,Rvachev:2001:TII}
\begin{equation}\label{eq:ti}
g(\vx) = \sum_{i=1}^M w_i(\vx) \, g_i(\vx), \quad
w_i(\vx) = \dfrac{ \phi_i^{- \mu_i} } { \sum_{j=1}^M \phi_j^{-\mu_i} }
= \frac{ \prod_{j=1; j\ne i}^M \phi_j^{\mu_j} } 
                { \sum_{k=1}^M \prod_{j = 1; j\ne k}^M  \phi_j^{\mu_j} }
,
\end{equation}
where the weights $w_i$ form a partition-of-unity, and~\eqref{eq:ti}
interpolates $g_i$ on the set $\Gamma_i$. In~\eqref{eq:ti}, $\mu_i \ge 1$
is a constant that controls the nature of interpolation that accrues. For $\mu_i = 1$, the
function $g_i$ is interpolated on $\Gamma_i$, whereas
if $\mu_i = 2$, both $g_i$ and $\partial g_i / \partial n$ are interpolated
on $\Gamma_i$.  Now, on using the solution
structure for the Dirichlet problem given in~\eqref{eq:Dirichlet} and following the work
of Rvachev et al.~\cite{Rvachev:1995:RBV,Rvachev:2000:OCR}, we can write the trial function in the deep Ritz method as:
\begin{equation}\label{eq:trial_EBC}
\tilde{u}_\textrm{nn}^\textrm{bc} (\vx;\vm{\theta}) = g(\vx) + 
\phi(\vx) \, \unnRxt ,
\end{equation} 
where $\phi$ can also be replaced by the product
$\prod_{i=1}^M \phi_i$. 
Since $g(\vx) = g_i(\vx)$ on $\Gamma_i$ and $\phi(\vx)$ vanishes on $\cup_{i=1}^M \Gamma_i$, kinematic admissibility of~\eqref{eq:trial_EBC} is verified.

\subsection{Imposing inhomogeneous essential and Robin 
            boundary conditions}\label{subsec:EBCRobin}
Let us consider the case that the boundary $\partial \Omega := \Gamma = \Gamma_1 \cup \Gamma_2$ with $\Gamma_1 \cap \Gamma_2 = \emptyset$. The boundary
conditions on $\Gamma_1$ and $\Gamma_2$ are:
\begin{equation}
    u = g \ \ \textrm{on } \Gamma_1, \quad
    \frac{\partial u}{\partial n} + cu = h \ \ \textrm{on } \Gamma_2.
\end{equation}
Let $\phi_1(\vx)$ and $\phi_2(\vx)$ be the approximate distance functions to the boundaries $\Gamma_1$ and $\Gamma_2$, respectively. We use
the R-equivalence ($m = 1$) relation in~\eqref{eq:phi_eq} to form
\begin{equation}\label{eq:phi1phi2}
\phi(\vx) = \frac{ \phi_1(\vx) \phi_2(\vx)}{\phi_1(\vx) +
\phi_2(\vx)}, 
\end{equation}
which is the
ADF to the boundary $\Gamma$. For the case when the boundary is partitioned into two disjoint sets, we consider
two approaches to form a trial function. The first approach uses superposition of two solution structures and has a simple form. The second method,
which is based on transfinite interpolation, is applicable
in general when there are multiple boundaries on which essential and Robin 
boundary conditions are imposed. In this paper, we apply the first approach to solve a one-dimensional problem
with mixed (Dirichlet and
Neumann) boundary conditions, 
and adopt the second approach to solve 
a problem with Dirichlet and
Robin boundary conditions in two dimensions.

\smallskip
\noindent {\em Approach I}:  We form solutions structures
$u_1$ and $u_2$ such that
\begin{equation} \label{eq:u1u2}
  u_1 = g \ \ \textrm{on } \Gamma_1, \ \
  \frac{\partial u_1}{\partial n} + u_1 = 0 \ \ \textrm{on } \Gamma_2, \quad
  u_2 = 0 \ \ \textrm{on } \Gamma_1, \ \
  \frac{\partial u_2}{\partial n} + u_2 = h \ \ \textrm{on } \Gamma_2, 
\end{equation}
and therefore the desired trial function is: $u(\vx) = u_1(\vx) + u_2(\vx)$.  Now, on using~\eqref{eq:Dirichlet}, we know that
the function $g + \phi_1 \tilde{u}_1$ satisfies the
essential boundary condition on $\Gamma_1$ but does not meet the
Robin boundary condition on $\Gamma_2$. To satisfy the Robin
boundary condition on $\Gamma_2$, we normalize $u_1 + \phi \tilde{u}_1$
to the first order using $h = 0$ in~\eqref{eq:Robin} to
obtain
\begin{subequations}\label{eq:u1}
\begin{align}
u_1 &= \left[ 1 + \phi \bigl(c + D_1^{\phi_2} \bigr) \right]
       (g + \phi_1 \tilde{u}_1)  =
       \left[ 1 + \phi D_1^{\phi_2} \bigr) \right]
       (g + \phi_1 \tilde{u}_1)  + c \phi g + c \phi_1 \phi
       \tilde{u}_1, \\
       \intertext{where}
D_1^{\phi_2} (\cdot) &= \left[ - \nabla \phi_2 \cdot \nabla (\cdot)
                        \right]_{\Gamma_2}
                      =  \left. \frac{\partial (\cdot)}{\partial n} 
         \right|_{\Gamma_2}
         \end{align}
\end{subequations}
is the differential operator that acts in the outward normal direction on the boundary $\Gamma_2$.
Similarly, using~\eqref{eq:Robin}, the minimal structure for
$u_2$ is:
\begin{equation}\label{eq:u2}
u_2 = \phi  ( \phi_2 \tilde{u}_2 - h ) .
\end{equation}
Since $\phi = 0$ on $\Gamma$ and $\phi_2 = 0$ on $\Gamma_2$, the conditions on $u_2$  
in~\eqref{eq:u1u2} are satisfied.
On choosing $\tilde{u}_1 = \tilde{u}_2 := \unnRxt$ 
in~\eqref{eq:u1} and~\eqref{eq:u2} and adding them up,
the
ansatz $\unnbcxt$ is:
\begin{equation}\label{eq:trial_mixedBC-I}
\begin{split}
\unnbcxt &= \phi_1(\vx) \, \unnRxt \\ 
& \quad + \phi(\vx)  \Biggl[ \biggl\{ \phi_2(\vx) + c(\vx) \phi_1(\vx) \biggr\} \unnRxt  +             D_1^{\phi_2} \biggl( \phi_1 (\vx) \, \unnRxt \biggr)
           + D_1^{\phi_2} \bigl (g(\vx) \bigr) + 
           c(\vx) g(\vx) - h(\vx) \Biggl] \ +  \ g(\vx) ,
\end{split}
\end{equation}
where $\phi(\vx)$ is given in~\eqref{eq:phi1phi2}.
For mixed (Dirichlet and Robin) boundary conditions, the form of the
trial function in~\eqref{eq:trial_mixedBC-I} appears
in Rvachev and Sheiko~\cite{Rvachev:1995:RBV}.

\smallskip
\noindent {\em Approach II}:
On using~\eqref{eq:Dirichlet} and~\eqref{eq:Robin}, we 
select the boundary functions 
$u_1$ and $u_2$ on $\Gamma_1$ and $\Gamma_2$, respectively, as:
\begin{subequations}\label{eq:solutionstructure:mixedBC:II}
\begin{align}
u_1 &= g , \\
u_2 &= \left[ 1 + \phi_2 \bigl(c + D_1^{\phi_2} \bigr) \right]
       (\tilde{u}) - \phi_2 h ,
\end{align}
\end{subequations}
with $\phi_1 \phi_2^2 \tilde{u}$ as the composite remainder term.
From~\eqref{eq:ti}, we recall the transfinite interpolant, where we now choose $\mu_1 = 1$ and $\mu_2 = 2$. On carrying out a few algebraic simplifications and
using~\eqref{eq:solutionstructure:mixedBC:II} with
$\tilde{u} := \unnRxt$, we can write the trial function $\unnbcxt$ that exactly imposes the mixed boundary conditions as:
\begin{subequations}\label{eq:trial_mixedBC-II}
\begin{align}
\label{eq:trial_mixedBC-II-a}
\unnbcxt &= w_1(\vx) u_2(\vx;\vm{\theta}) + w_2(
\vx)u_1(\vx) + \phi_1(\vx) \phi_2^2(\vx) \unnRxt, 
\intertext{where}
\label{eq:trial_mixedBC-II-b}
w_1(\vx) &= \frac{ \phi_1(\vx) } { \phi_1(\vx) + \phi_2^2(\vx) },
\quad w_2(
\vx)= \frac{ \phi_2^2 (\vx) } { \phi_1(\vx) + \phi_2^2(\vx) } , \\
\intertext{and}
\label{eq:trial_mixedBC-II-c}
u_1(\vx) &= g(\vx) , \quad
u_2(\vx;\vm{\theta}) = \left[ 1 + \phi_2(\vx) \bigl(c(\vx) + D_1^{\phi_2} \bigr) \right]
       \bigl( \unnRxt \bigr) - \phi_2(\vx) h(\vx) .
\end{align}
\end{subequations}

\section{Approximation of Trial Functions in a Deep Neural Network}\label{sec:dnn}

In this paper, we exclusively use the densely connected neural network architecture, also known as the multi-layer perceptron (MLP), which has its origin in the early works of Rosenblatt~\cite{rosenblatt1958perceptron}. MLPs consist of multiple layers of neurons, where each neuron has the task of converting its input to an output by generally passing it through a nonlinear function called activation. MLPs are characterized by an architecture where neurons in a given layer are connected densely to the neurons in the neighboring layers (\fref{fig:diagram_dnn}). We note in passing, however, that the latest revolution in deep learning began with a different neural network 
architecture---convolutional neural networks (CNN)---applied to 
image classification tasks~\cite{lecun1998mnist}. Independent
of research in PINNs, some of these modern architectures have also been applied to mechanics problems~\cite{finol2019deep}. In keeping
with the broadly accepted definition, we consider any 
\suku{deep network to have two or more hidden layers}.

In a standard collocation or Ritz method, the trial
function is expanded as a linear combination of known basis functions.
The point of departure in using deep neural networks is that the ansatz
herein is represented by a nonlinear map that consists of unknown
parameters. These parameters are obtained via
the solution of a minimization problem, which in general but not necessarily, is a least squares optimization problem. Once the parameters are determined, one obtains the numerical solution to the boundary-value problem.

\subsection{Feedforward neural network}\label{subsec:feedforward}
\begin{figure}
\centering
\begin{subfigure}{0.2\textwidth}
\includegraphics[width=\textwidth]{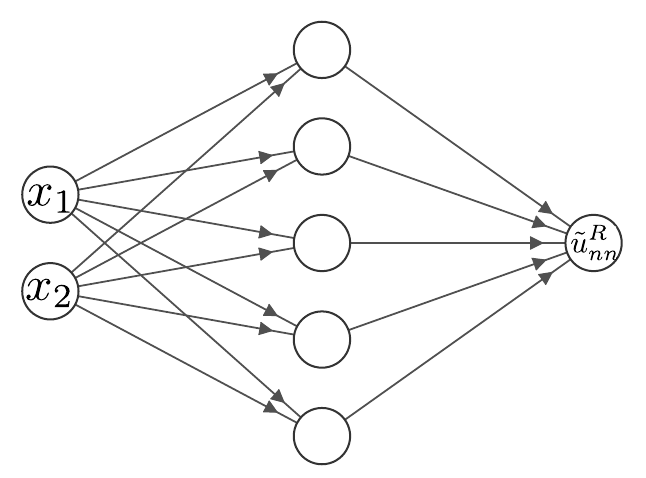}
\caption{}
\end{subfigure}\hfill
\begin{subfigure}{0.3\textwidth}
\includegraphics[width=\textwidth]{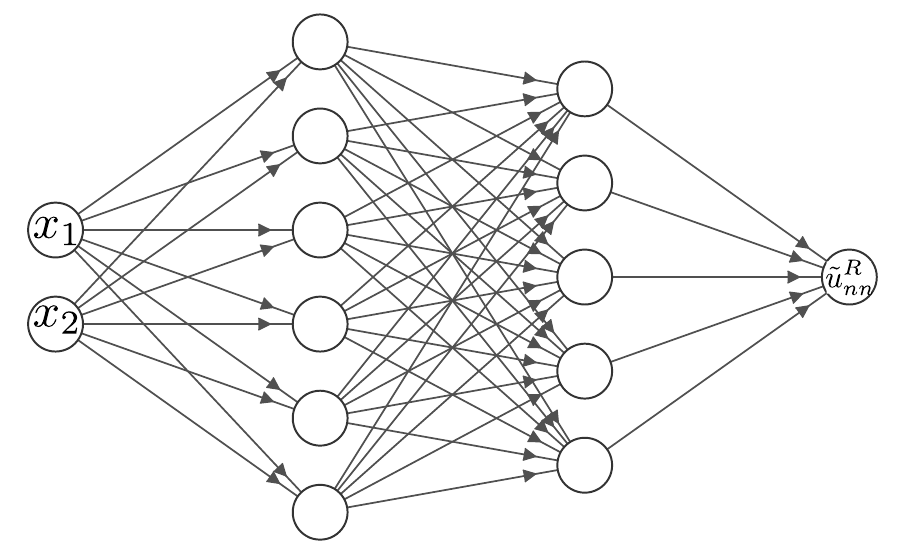}
\caption{}
\end{subfigure}\hfill
\begin{subfigure}{0.5\textwidth}
\includegraphics[width=\textwidth]{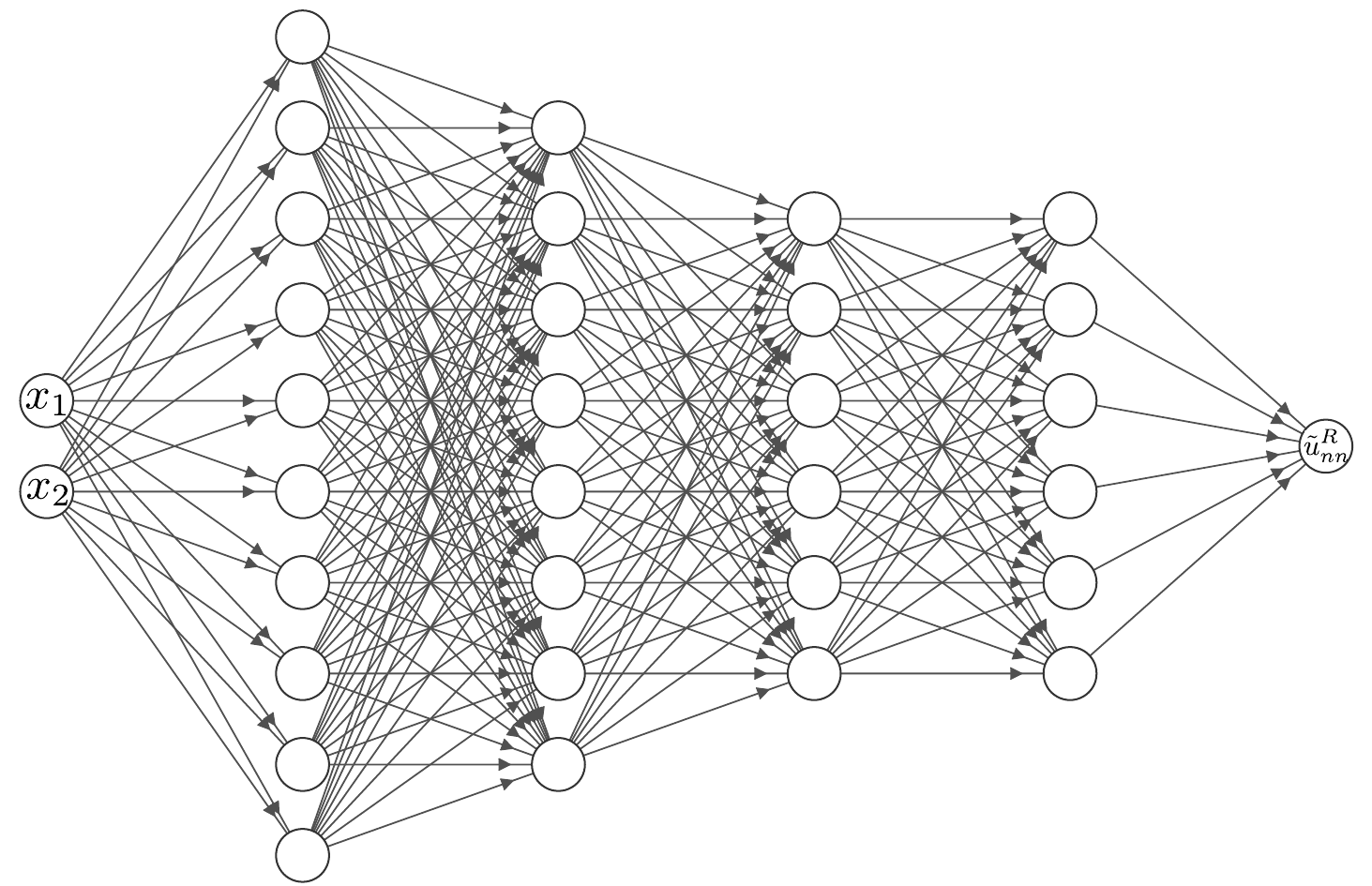}
\caption{}
\end{subfigure}
\caption{Deep neural network with (a) one, (b) two, and (c) four hidden  
         layers~\cite{LeNail:2019:NNS}. Input layer
         is a point $\vx \in \Re^2$, and 
         the output layer is the
         PINN approximation, $\unnRxt$. 
         }\label{fig:diagram_dnn}
\end{figure}

Given a point $\vx \in \Re^d$, we use a multilayer feedforward 
deep neural network to construct $\unnRxt$, which is then used to 
build $\unnbcx$, the approximation to $u(\vx):\Re^d \rightarrow \Re$.  
The layers in between
the input and output layer are known as hidden layers. Each hidden layer consists of
neurons (hidden units), and each neuron in a hidden layer takes its
input from the neurons in the preceding layer and computes its own activation.
A network diagram of neural networks with one, two, and four 
hidden layers
is shown in \fref{fig:diagram_dnn}.

In this paper, we consider the following boundary-value problem:
\begin{subequations}\label{eq:model_bvp}
\begin{align}
\label{eq:model_bvp-a}
\calL u(\vx) &= f(\vx) \ \ \textrm{in } \Omega \subset \Re^d , \\
\label{eq:model_bvp-b}
{\cal B}_u u(\vx) &= g(\vx) \ \ \textrm{on } \Gamma_u, \\
\label{eq:model_bvp-c}
{\cal B}_n u(\vx) &= h(\vx) \ \ \textrm{on } \Gamma_n ,
\end{align}
\end{subequations}
where $\calL$ is in general a  differential operator plus the
identity and $d$ is the spatial dimension. 
In~\eqref{eq:model_bvp},
$\Omega$ is an open bounded domain with boundary
$\partial \Omega = \Gamma_u \cup \Gamma_n$ and
$\Gamma_u \cap \Gamma_n = \emptyset$.  Equation~\eqref{eq:model_bvp-b}
represents essential boundary conditions for second- and fourth-order
problems, and~\eqref{eq:model_bvp-c} imposes natural and higher order
boundary conditions.
For a second-order problem, \eqref{eq:model_bvp-c} is a Robin 
boundary condition.
As noted in~\sref{sec:BCs}, we form the trial function to
approximate $u(\vx)$ as a combination of terms that involve the 
approximate distance functions $\bigl(\phi_k(\vx)$, $\phi(\vx)\bigr)$ and the neural network approximation. Consider
a neural network that consists of $L$ hidden layers with $\calN_\ell$
neurons in the hidden layer $\ell$ and activation function
$\sigma : \Re \rightarrow \Re$. The size of the neural network is:
$\calN = \sum_{\ell=1}^L \calN_\ell$. 
 Let $\unnRxt$ be the PINN
approximation, where
$\vm{\theta} := \{ \vm{W}, \ \vm{b} \}$ is the unknown parameter vector, with 
weights $\vm{W}_\ell \in \Re^{\calN_\ell \times \calN_{\ell-1}}$ and 
biases $\vm{b}_\ell \in \Re^{ \calN_\ell }$. 
We write $\unnRxt$ via the composition of $T^{(\ell)}$
($\ell=1,2,\ldots,L)$ and a linear map ${\cal G}$ as:
\begin{equation}\label{eq:composition}
\unnRxt = {\cal G} \circ T^{(L)} \circ T^{(L-1)} \circ \ldots \circ 
T^{(1)}(\vx),
\end{equation}
where ${\cal G}: \Re^{\calN_L} \rightarrow \Re$ is the linear mapping
for the output layer and in each hidden layer ($\ell = 1,2,\ldots,L$), the
nonlinear mapping is:
\begin{equation}\label{eq:Tmap}
T^{(\ell)}(\vm{z}) = \sigma(\vm{W}_\ell \cdot \vm{z} + \vm{b}_\ell ),
\end{equation}
where $\vm{z} \in \Re^{\calN_{\ell-1}}$.
For a neural network with activation function $\sigma$ and
 a single hidden layer that consists of 
$\calN$ neurons, the PINN approximation is:
\begin{equation}
\unnRxt = \sum_{i=1}^\mathcal{N} c_i \, \sigma(\vm{w}_i \cdot \vx + b_i) ,
\end{equation}
where $\vm{w}_i\in \Re^d$ and $b_i, c_i \in \Re$.
It is known that a multilayer neural network is sufficiently rich to be 
able to approximate any $L^2$ function to arbitrary accuracy~\cite{Hornik:1989:MFN,Hornik:1991:ACM,Hornik:1993:SNR}. However, realization of this in
practice hinges on the
width and depth of the network, choice of the activation function and the
computational algorithm to solve the optimization problem.

\subsection{Backpropagation algorithm}\label{subsec:backprop}
Determination of the optimal parameters of the network is done through an iterative optimization process called backpropagation and a popular algorithm for backpropagation is stochastic gradient descent, which is the stochastic version of the gradient
descent algorithm. \suku{An important step in this procedure is the efficient computation of the gradient of the loss function using 
automatic differentiation.}
In this paper, we exclusively use the \texttt{Adam} backpropagation algorithm whose details are given in~\cite{kingma2014adam}. \texttt{Adam} is an extension of the stochastic gradient descent algorithm and differs from it primarily in its implementation of per-parameter learning rates that are continuously tweaked during learning through both the first and second moments of the gradients.

\section{Formulations}\label{sec:formulations}
We now present the formulations for deep collocation and deep Ritz for
second- and fourth-order problems. In the collocation approach,
all boundary conditions are exactly satisfied; 
for the Ritz method, the essential boundary conditions are met.

\subsection{Deep collocation}\label{subsec:collocation}
Let us consider a second-order boundary-value problem with mixed boundary conditions. We require the trial function given in~\eqref{eq:trial_mixedBC-II} (meets all boundary conditions) to 
also satisfy the governing equation in~\eqref{eq:model_bvp-a} at $N_I$ interior collocation points. We label
these collocation points as 
$\{\vx_k\}_{k=1}^{N_I}$.  When substituted in~\eqref{eq:model_bvp-a}, this defines a residual error at
each point $\vx_k$ and to determine the parameter $\vm{\theta}$,
we minimize the mean squared residual error:
\begin{equation}\label{eq:collocation_loss}
\vm{\theta}^* = \argmin_{\vm{\theta}} \Lnnbct, \quad
\Lnnbct := || \calL \unnbcxt - f(\vx) ||_{\Omega,N_I}^2
= \frac{1}{N_I} \sum_{k=1}^{N_I}   \left[ 
\calL \unnbc (\vx_k;\vm{\theta}) - f(\vx_k) \right]^2  ,
\end{equation}
where $\Lnnbct$ is known as the loss function, and
$||\cdot||_{\Omega,N_I}$ denotes the mean discrete $L_2$ norm
of its argument over the domain $\Omega$ that is discretized using
$N_I$ collocation points. 
The \texttt{Adam} algorithm~\cite{kingma2014adam} is used to solve~\eqref{eq:collocation_loss}.

When a standard PINN trial function, $\unnxt$, is used, additional residual error
contributions from the boundary conditions are present in the loss function. If only essential boundary conditions are imposed, with
$u = g$ on
$\partial \Omega$, then the solution for the
parameter $\vm{\theta}$ is given by
\begin{equation}\label{eq:collocation_loss_unn}
\vm{\theta}^* = \argmin_{\vm{\theta}} \Lnnt, \quad
\Lnnt := || \calL \unnxt - f(\vx) ||_{\Omega,N_I}^2
+ || \unnxt - g(\vx) ||_{\partial \Omega,N_B}^2 ,
\end{equation}
where $||\cdot||_{\partial \Omega, N_B}$ is the mean discrete $L_2$
norm of its argument over the boundary $\partial \Omega$, and
$N_B$ is the number of collocation points on
$\partial \Omega$. \suku{As noted in~\sref{sec:intro}, the}
first and second terms in $\Lnnt$ are referred to as the interior (PDE) loss
and the boundary loss, respectively.

\subsection{Deep Ritz}\label{subsec:Ritz}
E and Yu~\cite{E:2018:DRM} introduced the deep Ritz method to solve low- and high-order boundary-value problems that have a variational structure. Samaniego et al.~\cite{Samaniego:2020:EAS} 
have applied the Ritz approach to solve problems in computational solid mechanics.
We consider second-order (Poisson) and fourth-order (plate bending) boundary-value problems.  Essential and mixed boundary conditions are considered for the Poisson problem and clamped boundary conditions are imposed for the plate bending problem.
We use a variational principle (minimization of the potential energy functional) to solve both problems. A trial function, $\unnbcxt$, from a finite-dimensional space is used in the
variational principle, which forms the basis of the deep Ritz method.

\subsubsection{Second-order problems}\label{subsubsec:2nd_order}
Referring to the model boundary-value problem in~\eqref{eq:model_bvp}, we first consider a 
Poisson problem with Dirichlet boundary conditions:
\begin{subequations}\label{eq:poisson-ritz}
\begin{align}
\label{eq:poisson-ritz-a}
-\nabla^2 u &= f \ \ \textrm{in } \Omega, \\
\label{eq:poisson-ritz-b}
u &= g \ \ \textrm{on } \partial \Omega.
\end{align}
\end{subequations}
The variational principle for this problem is: 
\begin{subequations}\label{eq:poisson_var}
\begin{align}
\label{eq:poisson_var-a}
\min_{u \in {\cal S}} \ \Biggl[
\Pi[u] &:= 
\underset{W_\textrm{int}[u]}
 {\underbrace{\frac{1}{2} a(u,u) }}
 -
 \underset{W_\textrm{ext}[u]}
 {\underbrace{\ell(u) }}, \  \ 
 {\cal S} = \biggl\{ u : u \in H^1(\Omega), \ u = g \ \textrm{on }
\partial \Omega \biggr\} \ \Biggr], \\
 \intertext{where}
 \label{eq:poisson_var-b}
 a(u,w) &=  \int_\Omega \nabla u \cdot \nabla w \, d\vx,
 \quad
 \ell(w) = \int_\Omega f \, w d\vx
\end{align}
\end{subequations}
are a symmetric bilinear functional and a linear functional, respectively,
and $u$ and $w$ are the trial and test functions, respectively.
In~\eqref{eq:poisson_var-a},
$W_\textrm{int}[u]$ is the internal work (strain energy)
and $W_\textrm{ext}[u]$ is the external work, and
$H^k(\Omega)$ is the Sobolev space that consists of 
functions that have square integrable derivatives up to 
order $k$ in $\Omega$. 

As a second problem, we consider a 
Poisson problem with mixed (Dirichlet and Robin) boundary conditions:
\begin{subequations}\label{eq:poisson2-ritz}
\begin{align}
\label{eq:poisson2-ritz-a}
-\nabla^2 u &= f \ \ \textrm{in } \Omega, \\
\label{eq:poisson2-ritz-b}
u = g \ \ \textrm{on } \Gamma_u, \quad 
& \frac{\partial u}{\partial n}  + c u = h \ \ \textrm{on }
\Gamma_n, 
\end{align}
\end{subequations}
where $\vm{n}$ is the unit outward normal on
$\Gamma_n$, $c := c(\vx)$ and $h := h(\vx)$ are boundary
data, and  
\suku{$\Gamma = \partial \Omega$ with} $\Gamma = \Gamma_u \cup
\Gamma_n$ and $\Gamma_u \cap \Gamma_n = \emptyset$.
The variational principle for this problem is: 
\begin{subequations}\label{eq:poisson2_var}
\begin{align}
\label{eq:poisson2_var-a}
\min_{u \in {\cal S}} \ \Biggl[
\Pi[u] &:= \frac{1}{2} a(u,u) - \ell(u), \ \ 
 {\cal S} = \biggl\{ u : u \in H^1(\Omega), \ u = g \ \textrm{on }
\partial \Omega \biggr\} \ \Biggr], \\
 \intertext{where}
 \label{eq:poisson2_var-b}
 a(u,w) &=  \int_\Omega 
 \nabla u \cdot \nabla w \, dx + 
 \int_{\Gamma_n} c u w  \, dS,
 \quad
 \ell(w) = \int_\Omega f w \, d\vx 
 + \int_{\Gamma_n} h w  \, dS.
\end{align}
\end{subequations}

As indicated in~\eqref{eq:trial_EBC}, we choose a kinematically admissible trial function so that the essential boundary condition is satisfied. On
substituting the PINN trial function 
in~\eqref{eq:poisson_var} or~\eqref{eq:poisson2_var}, the
finite-dimensional minimization problem becomes
\begin{equation}
\min_{\vm{\theta} } \left[
\Pi[\unnbc] := \frac{1}{2} \, a \, \Bigl( \unnbcxt , \unnbcxt \Bigr) 
- \ell \, \Bigl( \unnbcxt \Bigr) \right]. 
\end{equation}
To compute the potential energy functional, we use a 
\suku{Monte Carlo integration rule with points} that are distributed 
(randomly or
quasi-uniformly) in the domain with a constant weight that is attached to each point. To obtain the unknown parameters,
we solve a discrete nonlinear minimization problem 
that is posed as:
\begin{equation}\label{eq:poisson_var_Lt}
\vm{\theta}^* = \argmin_{\vm{\theta}} \Lnnbct, \quad
\Lnnbct  = \frac{1}{N_I} \sum_{k=1}^{N_I}
\left[ 
\frac{1}{2} \, a \, \Bigl( \unnbc (\vx_k; \vm{\theta}) , 
\unnbc (\vx_k; \vm{\theta}) \Bigr) 
- \ell \, \Bigl( \unnbc (\vx_k ; \vm{\theta} ) \Bigr)
\right] ,
\end{equation}
where $\Lnnbct$ is the loss function, and
$a(\cdot,\cdot)$ and $\ell(\cdot)$ are given
in~\eqref{eq:poisson_var-b} and~\eqref{eq:poisson2_var-b}
for the Poisson problems with Dirichlet and mixed boundary conditions,
respectively.
In general,
the stochastic gradient descent algorithm or a variant of it
is used to solve~\eqref{eq:poisson_var_Lt}.

\subsubsection{Fourth-order problem}\label{subsubsec:4th_order}
We consider the fourth-order problem of Kirchhoff plate bending with
clamped boundary conditions. The strong form of the boundary-value
problem is:
\begin{subequations}\label{eq:plate_bvp}
\begin{align}
\label{eq:plate_bvp-a}
\nabla^4 u &= f \ \ \textrm{in } \Omega \subset \Re^2, \\
\label{eq:plate_bvp-b}
u = 0 \ \ \textrm{on } \partial \Omega, & \quad 
u_{,n} :=  \frac{\partial u}{\partial n} = 
0 \ \ \textrm{on } 
\partial \Omega,
\end{align}
\end{subequations}
where $u := u(x,y)$ is the out-of-plane plate deflection, $f$
is the transverse load per unit area, and 
$\vm{n}$ is the unit outward normal on the boundary 
$\partial \Omega$. The 
flexural rigidity of the plate is
assumed to be unity. 

The variational principle that is associated with the strong form in~\eqref{eq:plate_bvp} is: 
\begin{subequations}\label{eq:plate_var}
\begin{align}
\label{eq:plate_var-a}
\min_{u \in {\cal S}} \
\Biggl[ 
\Pi[u] &:=
\frac{1}{2} a(u,u) - \ell(u), \ \ 
{\cal S} = 
\biggl\{ u : u \in H^2(\Omega), \ u = 0 \ \textrm{on }
\partial \Omega, \ u_{,n} = 0 \ 
\textrm{on }  \partial \Omega \biggr\} \Biggr]
, \\
\intertext{where}
\label{eq:plate_var-b}
& a(u,w) =  \int_\Omega 
\bigl( \nabla^2 u \bigr) \bigl( \nabla^2 w \bigr) \, d\vx, \quad
\ell(w) = 
\int_\Omega f w \, d\vx ,
\end{align}
\end{subequations}
and now both essential boundary conditions in~\eqref{eq:plate_bvp-b} 
must be met to satisfy kinematic admissibility.  To meet this
objective, the solution structure for the PINN trial
function is chosen as:
\begin{equation}\label{eq:trial_plate}
    \unnbcxt =
     [ \phi(\vx) ]^2 \, \unnRxt,
\end{equation}
where $\phi(\vx)$ is an approximate distance function to the boundary $\partial \Omega$.

On substituting the trial function in~\eqref{eq:plate_var}, the
finite-dimensional minimization problem is of the same
form as~\eqref{eq:poisson_var}, where $a(\cdot,\cdot)$ and $\ell(\cdot)$ are now given by~\eqref{eq:plate_var-b}. 
Similar to the \suku{deep} Ritz approach for the Poisson problem, we use
a \suku{Monte Carlo integration rule with equally-weighted points}. Following
the same steps, the loss function is of the same form as in~\eqref{eq:poisson_var_Lt} with the
bilinear and linear functionals defined in~\eqref{eq:plate_var-b}. 

\section{Numerical Implementation}\label{sec:implementation}
Numerical solutions using deep neural networks are presented for one-, two- and four-dimensional boundary-value problems. Both, second- and fourth-order problems are considered. 
The trial function, $\unnbcxt$, contains $\unnRxt$, the
approximation that is composed by the neural network. The trial function
in the standard PINN~\cite{Raissi:2019:PIN} is denoted by
$\unnxt$.
Likewise, $\Lnnbct$ is used to denote the loss function in our approach and $\Lnnt$ is that for standard PINN. 
We refer to the numerical solution that is obtained (after training) by our approach as $\unnbcx$ and that obtained (after training) using standard PINN~\cite{Raissi:2019:PIN} as $\unnx$. 
\suku{Deep collocation} and deep Ritz methods are used to solve  boundary-value problems.  For the Ritz method, either ReLU or cubic ReLU 
activation function is used, with the problem that is solved in~\sref{subsubsec:RBFExample} being the 
sole exception where the Gaussian activation function is used. 
Unless otherwise stated, R-equivalence (REQ)
composition with $m = 1$ and mean value potential (MVP) with $p = 1$
are used to form approximate distance functions in $\unnbcxt$. 
For all problems that are considered in this paper, whenever REQ is indicated, we 
employ~\eqref{eq:phin_eq} to form $\phi(\vx)$; when MVP is mentioned, 
the expression for $\phi$ given 
in~\eqref{eq:phi_mvc} or~\eqref{eq:phi_tmvi} applies. 
For collocation, 
all boundary conditions are exactly satisfied, and essential boundary conditions are exactly met in the Ritz method,
which ensures that the trial function 
is kinematically admissible. All collocation points are considered in a single batch in the network; points are not sorted in bins and passed in batches, which in general is a more efficient approach.

The formulation described in~\sref{sec:formulations} has been implemented using Google's \texttt{JAX} library in Python~\cite{jax2018github}, which can automatically differentiate native Python and NumPy functions. As an example, the code listing below sets up the calculation for the approximate distance function over an arbitrary polygon through the R-equivalence operation 
in~\eqref{eq:phin_eq}.

\begin{Verbatim}[commandchars=\\\{\}]
\boldverb{def} dist(x1,y1,x2,y2):
  \boldverb{return} jax.numpy.sqrt((x2-x1)**2+(y2-y1)**2)

\boldverb{def} linseg(x,y,x1,y1,x2,y2):
  L   = dist(x1,y1,x2,y2)
  xc  = (x1+x2)/2.
  yc  = (y1+y2)/2.
  f   = (1/L)*((x-x1)*(y2-y1) - (y-y1)*(x2-x1))
  t   = (1/L)*((L/2.)**2-dist(x,y,xc,yc)**2)
  varphi = jax.numpy.sqrt(t**2+f**4)
  phi = jax.numpy.sqrt(f**2 + (1/4.)*(varphi - t)**2)
  \boldverb{return} phi
 
\boldverb{def} phi(x,y,segments):
  m = 1.
  R = 0.
  for i in range(len(segments[:,0])):
    phi = linseg(x,y,segments[i,0],segments[i,1],segments[i,2],segments[i,3])
    R = R + 1./phi**m
  R = 1/R**(1/m)
  \boldverb{return} R
\end{Verbatim}
Here \verb!segments! is a NumPy array, which
contains the coordinates of the line segments that make up the polygon and \verb!phi(x,y,segments)! returns an approximate distance function from the location \verb!(x,y)! to the polygon. A multilayer perceptron neural network is created using the following set of functions:

\begin{Verbatim}[commandchars=\\\{\}]
\boldverb{def} RePU3(x):
    \boldverb{return} (jax.numpy.maximum(0, x**3))

\boldverb{def} repu3_layer(params, x):
    \boldverb{return} RePU3(jax.numpy.dot(params[0], x) + params[1])

\boldverb{def} NN(params, x, y):
    """ Compute the forward pass for each example individually """
    activations = jax.numpy.array([x,y])
    # Loop over the RePU3 hidden layers
    for w, b in params[:-1]:
        activations = repu3_layer([w, b], activations)

    final_w, final_b = params[-1]
    final = jax.numpy.sum(jax.numpy.dot(final_w, activations)) + final_b
    \boldverb{return} (final[0])
\end{Verbatim}
In the above, the cubic ReLU function is used as the activation in the hidden layers, which is also an instance of the Rectified Power Unit (RePU) activation function $\REPU := [\max(0,x)]^n$ 
($n = 3$ for cubic ReLU). A linear activation 
function is used in the output layer. Here, \verb!params! is a NumPy array consisting of the optimizable parameters, \verb!(w,b)!, of the network, and its shape is determined by the architecture of the network. A network architecture that has been used frequently in this paper for problems in $\Re^d$ is a 2 hidden layer network $d$--$\calN$--$\calN$--$1$, where $\calN$ is the number of neurons in each hidden layer and is typically $50$. Finally, the following functions and their derivatives are used to construct the ansatz $\unnbcxt$, which satisfies homogeneous essential boundary conditions that are
imposed on the boundary of a polygonal domain:

\begin{Verbatim}[commandchars=\\\{\}]
\boldverb{def} u(params, x, y):
    \boldverb{return} phi(x,y,segments)*NN(params,x,y)
    
#Examples of first-order partial derivatives    
gradx = grad(u,1)
grady = grad(u,2)

#Examples of second-order partial derivatives    
gradxx = grad(grad(u,1),1)
gradyy = grad(grad(u,2),2)
gradxy = grad(grad(u,1),2)
\end{Verbatim}

\begin{figure}[htp]
\centering
\includegraphics[width=\textwidth]{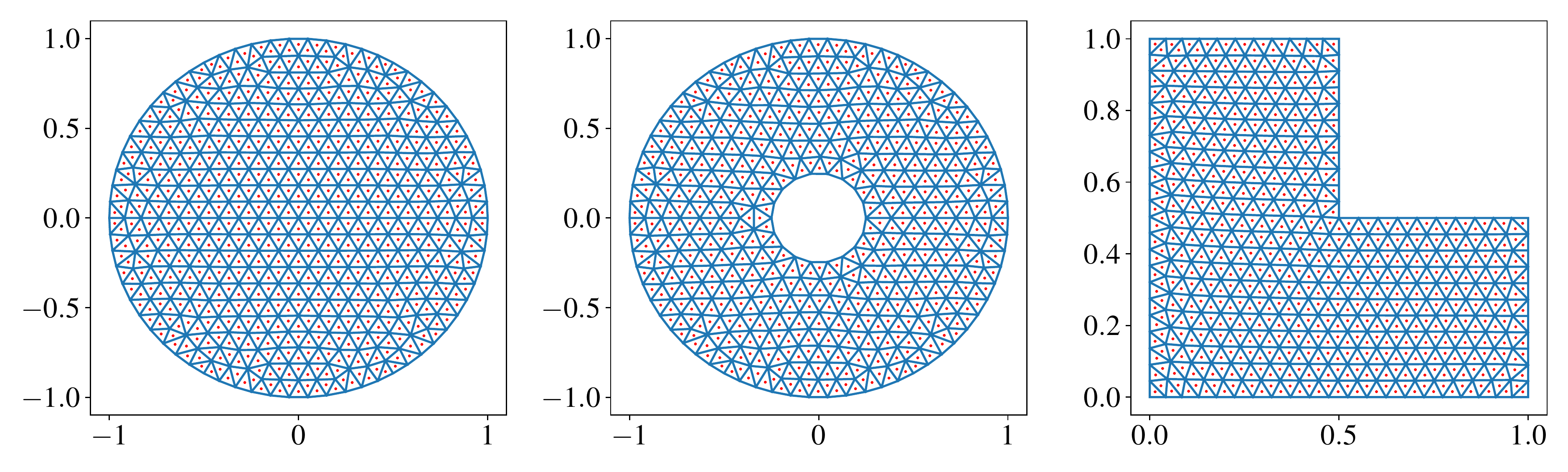} 
\caption{Representative meshes generated from \texttt{dmsh} and the corresponding collocation points.} 
\label{fig:FEM-mesh}
\end{figure}

Once the derivatives are formed, the appropriate
loss functions are constructed in the interior of the
domain to solve the problem through Ritz or collocation. For all
problems that are solved in this paper, network training is done using Google's Colaboratory cloud platform~\cite{bisong2019google}. 
\suku{Single-precision arithmetic is used in the computations}. 
An important consideration is the generation of collocation points in the interior and on the boundary of the domain to evaluate the loss
terms. When the domain is simple---for example, a square or a hypercube---then these points are generated on a uniform grid in $\Re^d$. This is the case for the solution of the heat equation on square domains, Eikonal problems in~\sref{subsec:Eikonal}, and for the Poisson problem over the 4-dimensional hypercube in~\sref{sec:4D}. In other examples that involve more complicated domains, we use the Python library \texttt{dmsh}~\cite{Schlomer:2020}, which draws 
inspiration from 
\texttt{distmesh}~\cite{Persson:2004:SMG},
to create triangular meshes and we use the centroids of the generated triangles as the 
interior collocation points. Figure~\ref{fig:FEM-mesh} shows a few representative meshes generated by \texttt{dmsh} that are used in this paper. The corresponding interior collocations points are shown as  dots. 
When needed, \texttt{dmsh} is also used to create collocation points on the boundaries of the domain.

\section{Numerical Examples in One Dimension}\label{sec:1D}
We consider several second-order problems that involve essential
and mixed boundary conditions. 
As a prototypical fourth-order problem, we solve
the deflection for a clamped Euler-Bernoulli beam, 
\suku{and also investigate a single hidden layer meshfree RBF-network solution in~\sref{subsubsec:RBFExample}.
Our objective is to demonstrate the benefits of the new formulation vis-{\`a}-vis standard PINN~\cite{Raissi:2019:PIN} in which equal weights are chosen for the PDE and
boundary loss terms. Hence, our emphasize is not on obtaining
the most accurate PINN solution by finding the optimal
hyperparameters nor comparing its performance versus the finite element method for forward problems.
We compare our results with those obtained using standard PINN~\cite{Raissi:2019:PIN} and to the exact solution.}

\subsection{Deformation of a homogeneous elastic rod}\label{subsec:rod}
Consider the boundary-value problem for the deformation of an elastic rod
(Youngs's modulus and cross-sectional area are taken as unity):
\begin{subequations}\label{eq:elasticRod}
\begin{align}
u'' + b &= 0 \ \ \textrm{in } \Omega = (x_1,x_2) \\
u(x_1) = g &,  \quad u'(x_2) + c u(x_2) = h,
\end{align}
\end{subequations}
where $u := u(x)$ is the axial displacement field, $u'(x)$ is the strain field,
$b := b(x)$ is the axial
body force per unit length, and $c$, $g$ and
$h$ are constants.  The second boundary condition is a Robin boundary condition; if the bar is connected to a spring that is attached to a fixed end, then $h = 0$. 
\suku{We select test problems with different boundary conditions,  and vary the regularity of $b(x)$ from it being a smooth function to a $\delta$-function, and even choose $b(x)$ that has a singularity at 
the origin.}

\subsubsection{Example 1}
As the first example, a Dirichlet problem in 
$\Omega = (0,1)$ is selected with body force
$b(x)=1-2x+10x^2$, and boundary conditions 
$u(0)=1/2$ and $u(1)=-1/2$. The exact solution is:
\suku{$u(x) = 1/2 - x^2/2 + x^3/3 - 10x^4/12$}.
The exact signed distance functions to $x = 0$ and $x = 1$ are $\phi_1(x) = x$ and 
$\phi_2(x) = 1 - x$, \suku{respectively}.  Now, we join these to form a smooth approximate distance function to the boundary $\partial
\Omega = \{0,1\}$. On using the product of $\phi_1$
and $\phi_2$, and the R-conjunction ($\alpha = 0$) and R-equivalence ($m = 2$) relations
in~\eqref{eq:R_alpha} and~\eqref{eq:phi_eq}, respectively,
we obtain the following combined distance functions:
\begin{equation}\label{eq:phiABC}
    \phi_A(x) = \phi_1(x) \phi_2(x), \ \ 
    \phi_B(x) = \phi_1(x) + \phi_2(x) 
                          - \sqrt{ \phi_1^2 (x) +
                                   \phi_2^2 (x) }, \ \
    \phi_C(x) = \frac{\phi_1(x) \phi_2(x)}
                          {\sqrt{\phi_1^2(x) + \phi_2^2(x)}}.
\end{equation}
Note that in one dimension, the product formula is normalized to order 1, but 
\suku{this} does not generalize to higher
dimensions. Coincidentally, for any domain $\Omega = (x_1,x_2)$ in one dimension, the product formula scaled by $L:=x_2-x_1$ coincides with $m = 1$ in the R-equivalence 
relation.

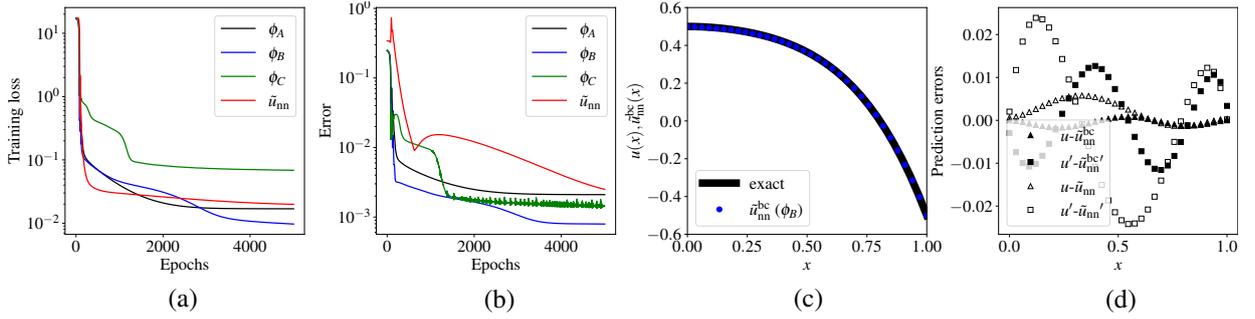
\begin{figure}[hbt]
\centering
\begin{tikzpicture}
\node at (0,0) {\includegraphics[width=0.98\textwidth]{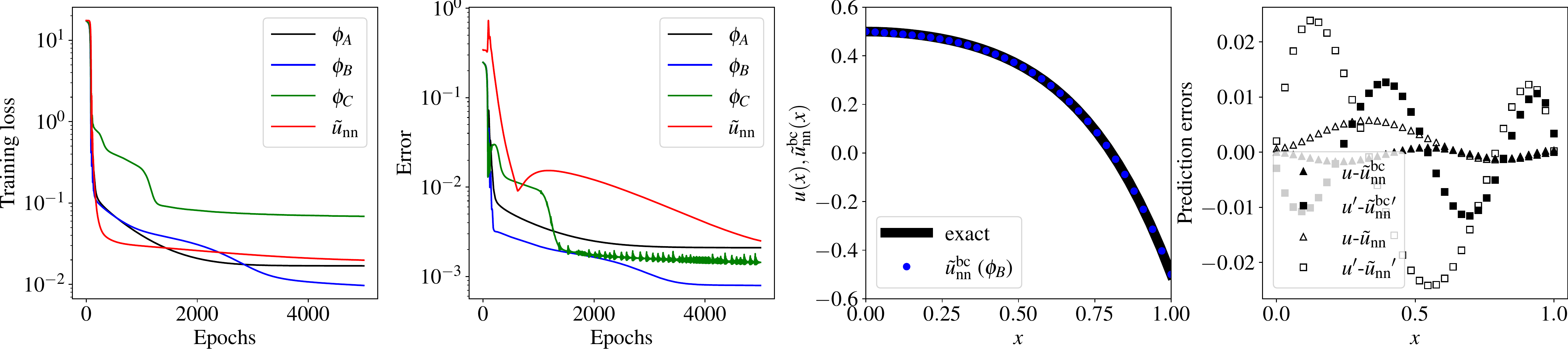}};
\node at (-2.28in,-0.85in) { (a)};
\node at (-0.64in,-0.85in) { (b)};
\node at ( 0.96in,-0.85in) { (c)};
\node at (2.58in,-0.85in) { (d)};
\end{tikzpicture}
\caption{Collocation solutions for a Dirichlet problem. The body force $b(x)=1-2x+10x^2$, and the essential boundary conditions are $u(0)=1/2$ and $u(1)=-1/2$. The network architecture is 1--30--30--1. Numerical solutions, $\unnbcxs$, are computed using different ADFs 
($\phi_A$, $\phi_B$, $\phi_C$), and are compared to the solution obtained using standard PINN, $\unnxs$. (a), (b) Training loss 
and normalized absolute error in the displacement field as a function of epochs. (c) 
Exact solution $u(x)$ and $\unnbcxs$ (using $\phi_B$ ) are shown, and
(d) Errors in the displacement and strain fields are compared.}
\label{fig:elasticRod_i} 
\end{figure}

In the numerical computations, the trial
function is formed using~\eqref{eq:Dirichlet}: 
\begin{equation}\label{eq:uhelasticRod_i}
    \unnbcxst = g(x) + \phi(x) \,
               \unnRxst, \quad g(x) = \frac{1-2x}{2},
\end{equation}
where $g(x)$ is formed using the transfinite interpolant 
in~\eqref{eq:ti}, $\unnRxst$ is the neural network approximation, and $\phi(x)$
is chosen to be either $\phi_A(x)$, $\phi_B(x)$ or $\phi_C(x)$. 
\suku{Note that} we did not solve the patch test ($b = 0$) since the exact
solution, $u(x) = g(x)$, is already captured by the presence 
of $g(x)$ in~\eqref{eq:uhelasticRod_i}.
The network architecture 1--30--30--1 is used.  
We compute collocation solutions using $\phi_A$, $\phi_B$ and
$\phi_C$ in~\eqref{eq:uhelasticRod_i}, and also for the standard
PINN approximation, $\unnxst$~\cite{Raissi:2019:PIN}. In~\fref{fig:elasticRod_i}, the numerical results are presented. All approaches are able to reach losses of the same order, and there is a good match between $\unnbcxs$ and the exact displacement field. From~\fref{fig:elasticRod_i}d, we observe that the errors in the
displacement and strain fields
using $\unnbc$ are uniformly smaller than those obtained using $\unn$. This difference stems from the
fact that in our approach the trial function is constructed with the  exact satisfaction of the boundary conditions.

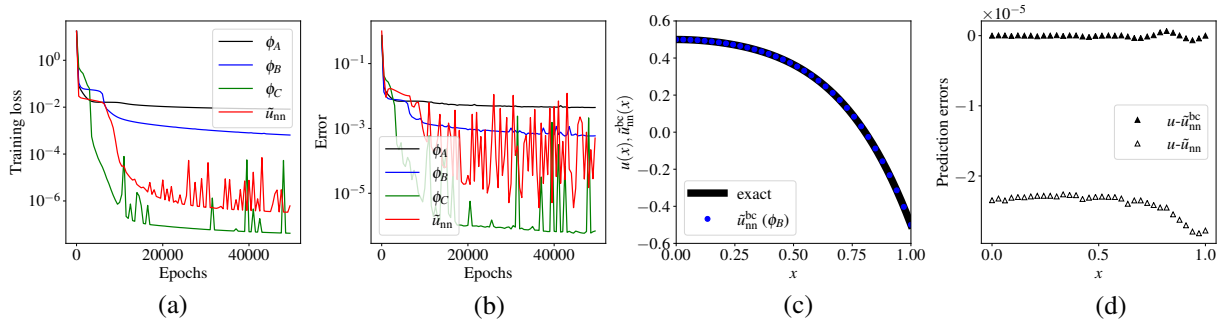
\begin{figure}[hbt]
\centering
\begin{tikzpicture}
\node at (0,0) {\includegraphics[width=0.98\textwidth]{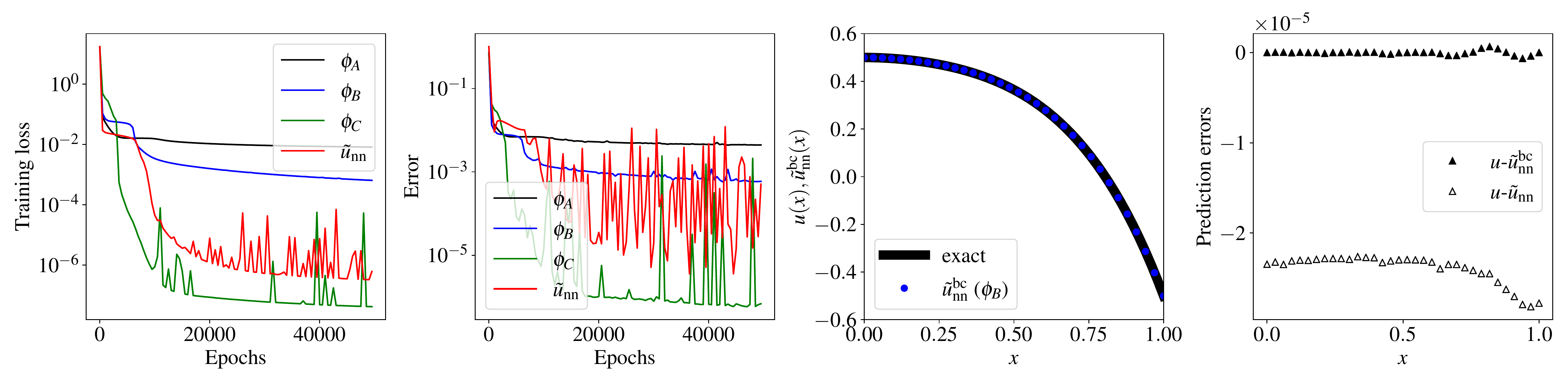}};
\node at (-2.28in,-0.85in) { (a)};
\node at (-0.64in,-0.85in) { (b)};
\node at ( 0.96in,-0.85in) { (c)};
\node at (2.58in,-0.85in) { (d)};
\end{tikzpicture}
\caption{\suku{Collocation solutions for the Dirichlet problem that is
presented in Example 1.  The network architecture is 1--50--50--1 with
300 interior collocation points. The captions for (a), (b), (c), (d)
mirror those shown in~\protect\fref{fig:elasticRod_i}.}}
\label{fig:elasticRod_50k} 
\end{figure}

\suku{We emphasize that it is possible for both standard PINN and our formulation to deliver better accuracy if a larger network, more interior collocation points, and network training for a longer duration are chosen. This is realized at the expense of more computing time. So for just this example, we demonstrate the same. We now use a 1--50--50--1 network architecture with 300 interior collocation points. The training is conducted until 50,000 epochs.  The solutions are presented in~\fref{fig:elasticRod_50k}, which reveal that both approaches are now much more accurate than the results that are shown in~\fref{fig:elasticRod_i}. From~\fref{fig:elasticRod_50k}, we observe that the solution obtained using standard PINN and our approach have relative errors of
${\cal O}(10^{-5})$ and ${\cal O}(10^{-6})$, respectively.  So even at smaller relative errors, we note that our approach is more accurate than standard PINN for the same values of the hyperparameters.}

\subsubsection{Example 2}
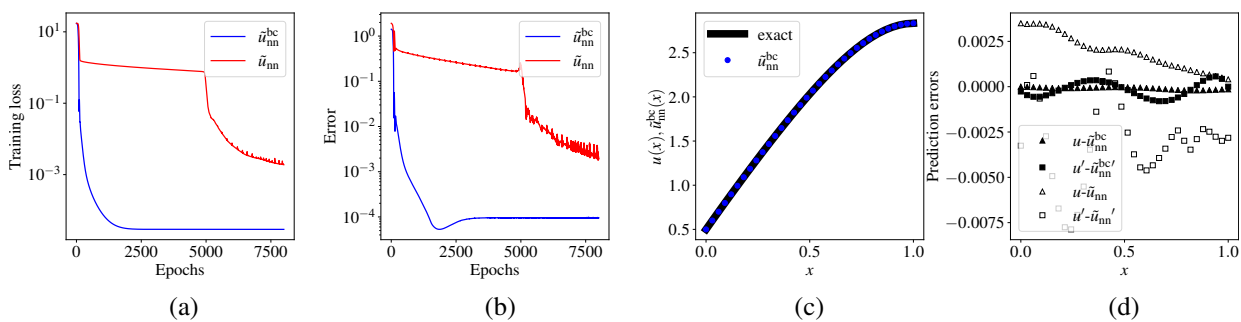
\begin{figure}[!htp]
\centering
\begin{tikzpicture}
\node at (0,0) {\includegraphics[width=0.98\textwidth]{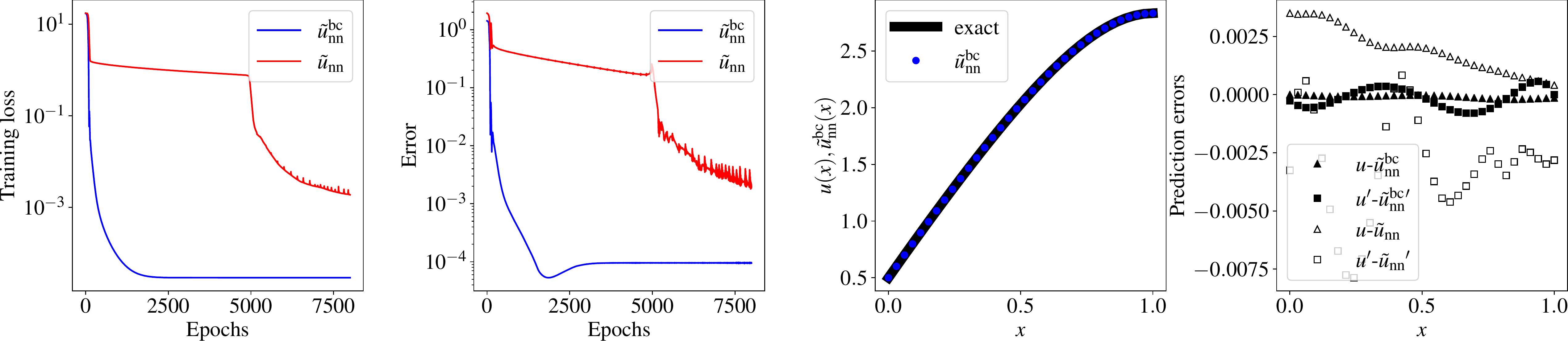}};
\node at (-2.28in,-0.85in) { (a)};
\node at (-0.64in,-0.85in) { (b)};
\node at ( 0.96in,-0.85in) { (c)};
\node at (2.59in,-0.85in) { (d)};
\end{tikzpicture}
\caption{Collocation solutions for a Neumann problem. The body force $b(x)=1-2x+10x^2$, and $u(0)=1/2$ and $u'(1)= 0$. The network architecture is 1--50--50--1. Numerical solution using
$\unnbcxst$ and $\unnxst$ are compared.
 (a), (b) Training loss and normalized absolute error as a function of epochs. (c) $\unnbcxs$ is compared to the exact solution.
(d) Errors in the displacement and strain fields for
    $\unnbcxs$ and $\unnxs$.} 
\label{fig:elasticRod_ii}
\end{figure}

We reconsider the problem posed in Example 1,
 \suku{with Dirichlet boundary condition at $x = 0$ but
 homogeneous Neumann (traction-free) condition at $x = 1$.} 
 For the domain $\Omega = (0,1)$, body force $b(x)=1-2x+10x^2$, and boundary conditions
$u(0)=1/2$ and $u'(1)=0$,
the exact solution is: 
\suku{$u(x) = 1/2 + 10x/3 - x^2/2 + x^3/3 - 10x^4/12$}.
The trial
function, $\unnbcxst$, is formed 
using~\eqref{eq:trial_mixedBC-I}, 
with $\phi_1(x) = x$, $\phi_2(x) = 1 - x$, $\phi(x) = x(1-x)$,
$g = 1/2$, and $c = h = 0$:
\begin{equation}\label{eq:trial_Example2}
\unnbcxst = x \, \unnRxst + x(1-x)
\left[ (1-x) \, \unnRxst  +  \bigl\{ x \, \unnRxst \bigr\}_{x=1}^\prime \right] +  \frac{1}{2} .
\end{equation}
The network architecture 1--50--50--1 is used, 
since the network 1--30--30--1 did not converge for standard
PINN. 
The collocation solutions $\unnbcxs$ and
$\unnxs$ are compared to the exact solution $u(x)$ in~\fref{fig:elasticRod_ii}. We observe from~\fref{fig:elasticRod_ii}a that the
training loss for
$\unnbcxst$ is about two orders smaller than $\unnxst$, whereas
in~\fref{fig:elasticRod_ii}b, the corresponding normalized absolute 
error is one order smaller. The numerical solution
$\unnbcxs$ is in excellent agreement with the
exact solution in~\fref{fig:elasticRod_ii}c.
Over the entire interval $x \in [0,1]$, we find that the 
displacement and strain fields \suku{from} $\unnbcxs$ are markedly
more accurate than those from $\unnxs$ (see~\fref{fig:elasticRod_ii}d).

\subsubsection{Example 3}
\begin{figure}[!htb]
\centering
\begin{tikzpicture}
\node at (0,0) {\includegraphics[width=0.98\textwidth]{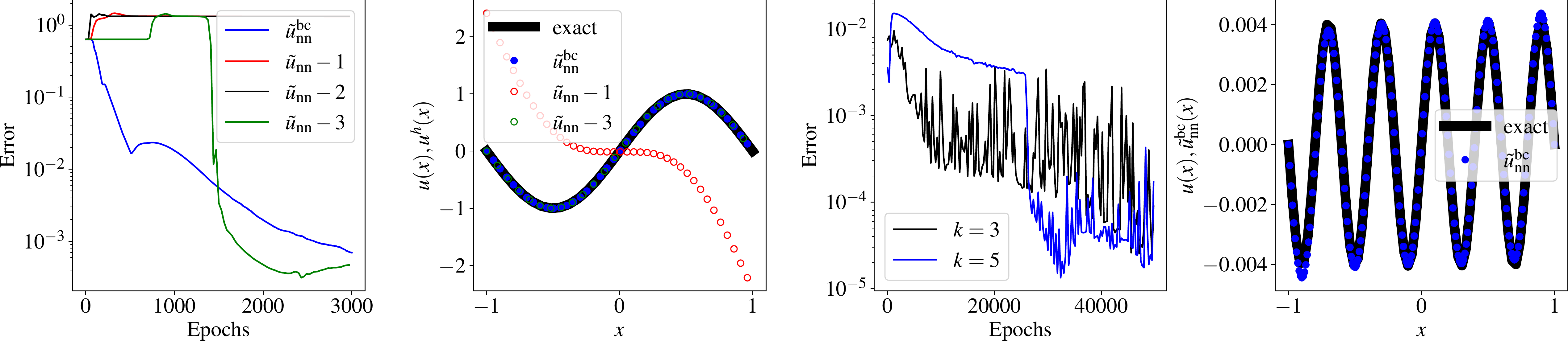}}; 
\node at (-2.30in,-0.85in) { (a) $k=1$};
\node at (-0.66in,-0.85in) { (b)};
\node at ( 0.96in,-0.85in) { (c)};
\node at (2.58in,-0.85in) { (d) $k=5$};
\end{tikzpicture}
\caption{
Collocation solutions for a homogeneous Dirichlet problem with a sinusoidal body force $b(x) = -\sin( k \pi x)$ with varying $k$.
(a) Normalized absolute error as a function of epochs for $k=1$. $\unnbcxst$ on a 1--30--30--1 architecture is compared to $\unnxst$. 
$\unn$--1, $\unn$--2 and $\unn$--3 are solutions that are obtained on network architectures 1--30--30--1, 1--100--100--1 and
1--100--100--100--1, respectively. (b) $\unnxs$, $\unnbcxs$, and $u(x)$ are plotted for $k = 1$.  (c) Normalized absolute error during
training 
using $\unnbcxst$ for $k = 3$ and $k =5$. Network architecture is 1--100--100--1. (d) For $k = 5$, error in the displacement field for $\unnxs$.} 
\label{fig:elasticRod-sin}
\end{figure}
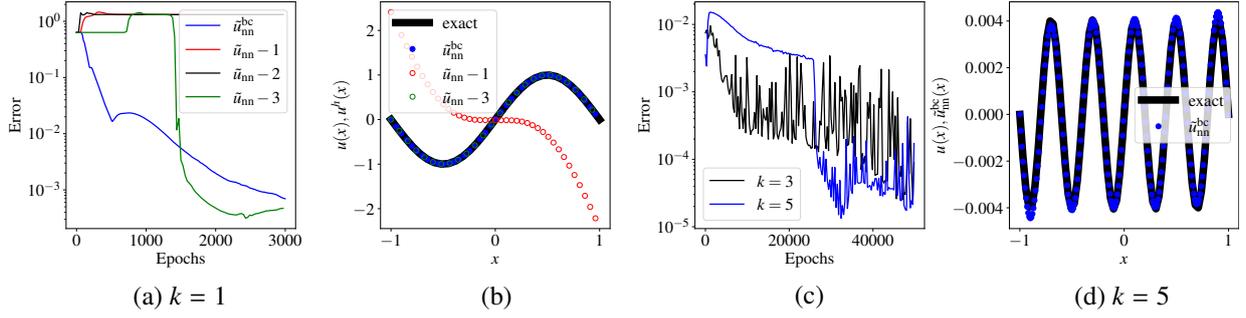

For this example, we choose $\Omega = (0,1)$, a sinusoidal body
force $b(x) = - \sin( k \pi x)$ with varying $k$, and essential boundary conditions $u(-1) = u(1) = 0$. The exact solution is:
$u(x) = - \sin(k \pi x)/\pi^2k^2$. 
\suku{This example serves to reveal the spectral (low-frequency) bias~\cite{Rahaman:2019:SBN,Wang:2021:EBF} of neural network approximations.}
In~\fref{fig:elasticRod-sin}, the numerical results are presented. For $k = 1$ and using standard PINN, we found that the normalized absolute error for the networks 1--30--30--1 and even 1--100--100--1 did not converge; it took a network architecture of 1--100-100-100-1 for the error to be comparable to that obtained using $\unnbcxst$ on a 1--30--30--1 
architecture. The normalized absolute errors as a function of epochs
is shown in~\fref{fig:elasticRod-sin}a. The numerical solutions are compared to the exact solution $u(x)$ in~\fref{fig:elasticRod-sin}b, and we notice the poor solution that is generated by $\unnxst$ on the 1--30--30--1 network. For 
$k = 3,\, 5$ on a
1--100--100-1 network, the normalized absolute error in the $\unnbcxs$ displacement field is shown in~\fref{fig:elasticRod-sin}c, and $\unnbcxs$ and $u(x)$ are compared in~\fref{fig:elasticRod-sin}d for $k=5$. Good agreement between $\unnbcxs$ and $u(x)$ is realized.  We attribute the poor performance of the standard PINN approach to the fact that the boundary conditions are poorly approximated, whereas with $\unnbcxst$ the boundary conditions are exactly satisfied. This observation is in broad agreement with the findings of Wang 
et al.~\cite{Wang:2020:UMG}, who further analyze the source of this 
discrepancy by drawing attention to the contributions from the boundary 
and interior terms in the loss function. 
For $k \ge 10$ with $\unnbcxst$, accurate PINN solutions 
\suku{using a single neural network for the entire domain becomes infeasible 
due to the spectral bias of the neural network approximation}; 
one \suku{can} adopt a domain-decomposition strategy to obtain 
\suku{accurate} numerical solutions for high-frequency problems.

\subsubsection{Example 4}
\begin{figure}[!htb]
\centering
\begin{tikzpicture}
\node at (0,0) {\includegraphics[width=0.98\textwidth]{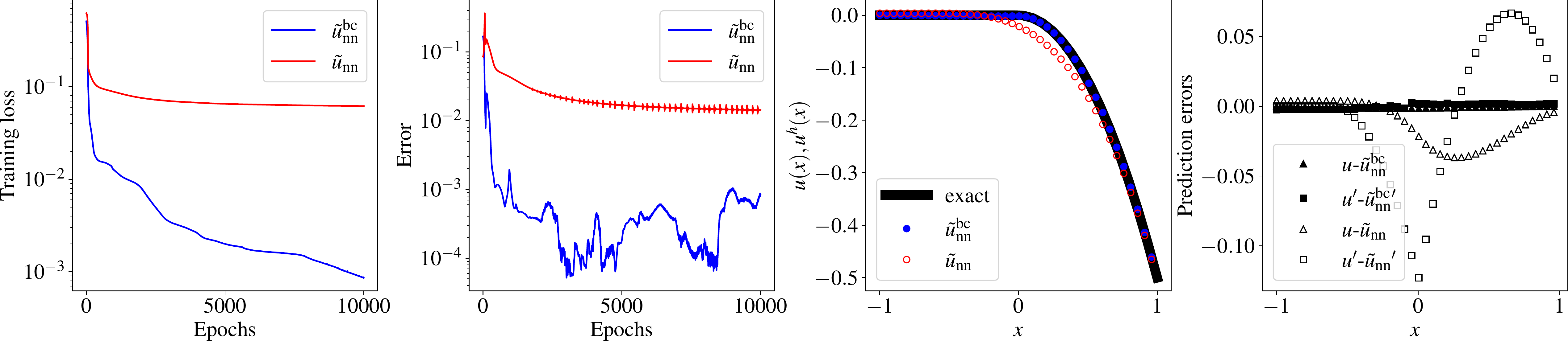}}; 
\node at (-2.28in,-0.85in) { (a)};
\node at (-0.65in,-0.85in) { (b)};
\node at ( 0.96in,-0.85in) { (c)};
\node at (2.58in,-0.85in) { (d)};
\end{tikzpicture}
\caption{Collocation solutions for a Dirichlet problem with a discontinuous body force, $b(x) = H(x)$. Essential boundary
conditions are $u(-1) = 0$ and $u(1)= -1/2$, and the network
architecture is 1--50--50--1. (a), (b) Training loss and normalized absolute error as a function of epochs for $\unnbcxst$ and $\unnxst$. 
(c) Comparisons of $\unnbcxs$ and $\unnxs$ with the exact solution.
(d) Error in the displacement and strain fields.
}
\label{fig:heaviside}
\end{figure}

Consider an elastic rod that occupies $\Omega = (-1,1)$ and is subjected to a discontinuous body force $b(x) = H(x)$, where
$H(x)$ is the Heaviside function. Essential boundary conditions are imposed at both ends: $u(-1) = 0$ and $u(1) = -1/2$.  The exact solution, $u(x) \in C^1(\Omega)$, is: 
\begin{equation}
u(x) = \begin{cases} 
        \ 0 &, \ - 1 \le x < 0 \\
        \ -\frac{x^2}{2}  &, \  \ 0 \le x \le 1
       \end{cases} .
\end{equation}
Numerical computations are performed on a 1--50--50--1 network. Numerical results for $\unn$ and $\unnbc$ are presented in~\fref{fig:heaviside}. The training loss of $\unnbcxst$ at 10,000 epochs is ${\cal O}(10^{-3})$ but $\unnxst$ stagnates to a value just below $0.1$. These losses correspond to normalized absolute errors of 
${\cal O}(10^{-5})$ and ${\cal O}(10^{-4})$ at the end of the training for $\unnbcxs$ and
$\unnxs$, respectively (see~\fref{fig:heaviside}b). Figures~\ref{fig:heaviside}c and~\ref{fig:heaviside}d reveal that
$\unnbcxs$ is in good agreement with $u(x)$, whereas the errors in $\unnxs$ are significant, and are especially pronounced in the vicinity of $x = 0$. The PDE loss is dominant over the boundary losses. We point out that if both boundary conditions are homogeneous then the accuracy of $\unnxs$ is comparable to that of $\unnbcxs$.

\subsubsection{Example 5}
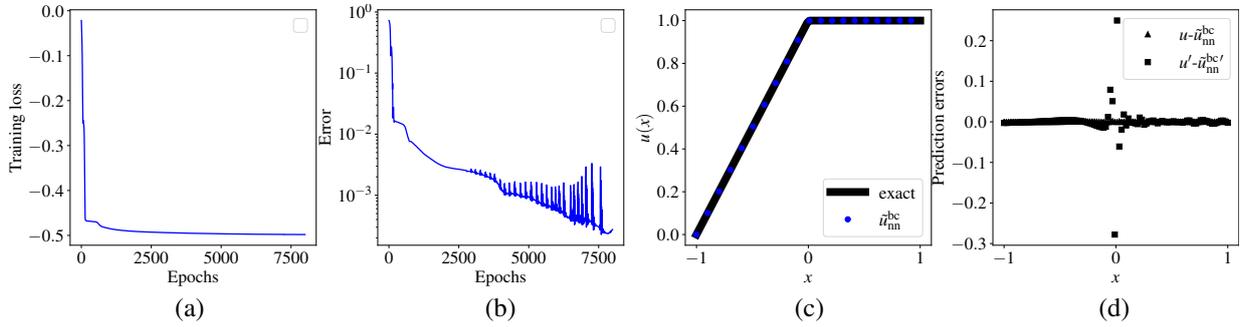
\begin{figure}[!htb]
\centering
\begin{tikzpicture}
\node at (0,0) {\includegraphics[width=0.98\textwidth]{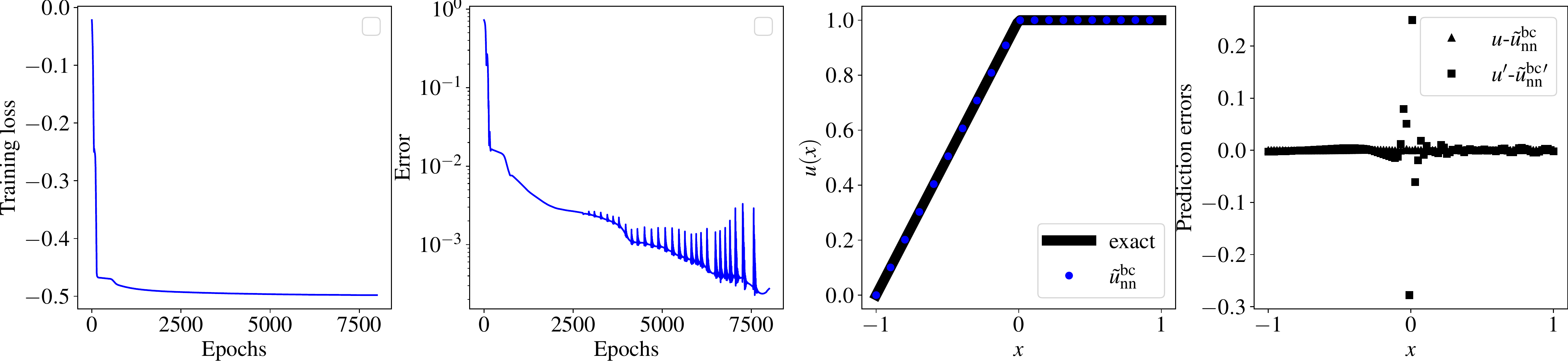}}; 
\node at (-2.25in,-0.86in) { (a) };
\node at (-0.64in,-0.86in) { (b)};
\node at ( 0.96in,-0.86in) { (c)};
\node at (2.54in,-0.86in) { (d)};
\end{tikzpicture}
\caption{Ritz solution for a point load, $b(x)=\delta(x)$. 
Homogeneous essential boundary condition is imposed at $x = -1$ and traction-free conditions prevail at $x = 1$. The network architecture is 1--50--50--1.
(a), (b)  Training (Ritz) loss and normalized absolute error as a function of epochs.
(c) Comparison of $\unnbcxs$ with the exact solution.
(d) Errors in displacement and strain fields.
} 
\label{fig:elasticRod-Ritz}
\end{figure}

Let us consider an elastic rod that occupies $\Omega = (-1,1)$ and is subjected to a unit point load at the origin, i.e., $b(x) = \delta(x)$, where $\delta(x)$ is the $\delta$-function.  Homogeneous essential boundary condition is prescribed at $x = -1$ and traction-free conditions prevail at $x = 1$. The exact solution, $u(x) \in C^0(\Omega)$, is:
\begin{equation}\label{eq:u_exact_Example5}
u(x) = \begin{cases} 
        \ 1 + x &, \  -1 \le x < 0 \\
        \ 1     &, \  0 \le x \le 1
       \end{cases} ,
\end{equation}
\suku{which has a kink at $x = 0$.}
We use the \suku{deep} Ritz method with $\unnbcxst$ to solve 
this problem since
it is not possible to solve this problem using the collocation 
method ($\delta(x)$ is a distribution that is defined over a zero measure). The
network
architecture is 1--50--50--1. Numerical results for $\unnbc$ are presented in~\fref{fig:elasticRod-Ritz}. 
From~\fref{fig:elasticRod-Ritz}a, we observe that the
training loss converges to a value close to $-0.5$ in a few epochs;
this corresponds to a small normalized absolute error of ${\cal O}(10^{-3})$ (see~\fref{fig:elasticRod-Ritz}b).  
On using $u$ given 
in~\eqref{eq:u_exact_Example5}, we find that the
potential energy of the exact solution is:
\begin{equation*}
    \Pi[u] = \frac{1}{2} \int_{-1}^1 u^{\prime 2} \, dx - u(0) =
    -\frac{1}{2} ,
\end{equation*}
which is the target loss that $\Pi[\unnbc]$
seeks to attain.
Figure~\ref{fig:elasticRod-Ritz}c shows 
excellent agreement between $\unnbcxs$ and the exact solution. The
errors in the displacement and strain fields are presented 
in~\fref{fig:elasticRod-Ritz}d. The errors in $u$ are uniformly small for all $x \in [0,1]$. The errors in the strain field 
follow the same trends, but have larger errors in the vicinity
of the origin.  This is
not surprising since $\unnbcxs$ is $C^2(\Omega)$ (cubic ReLU activation function), whereas the exact solution $u \in C^0(\Omega)$, with $u'(x)$ being discontinuous at $x = 0$.
This is also the source for the small discrepancy in the Ritz energy loss.

\subsubsection{Example 6}
\begin{figure}[!htb]
\centering
\begin{tikzpicture}
\node at (0,0) {\includegraphics[width=0.98\textwidth]{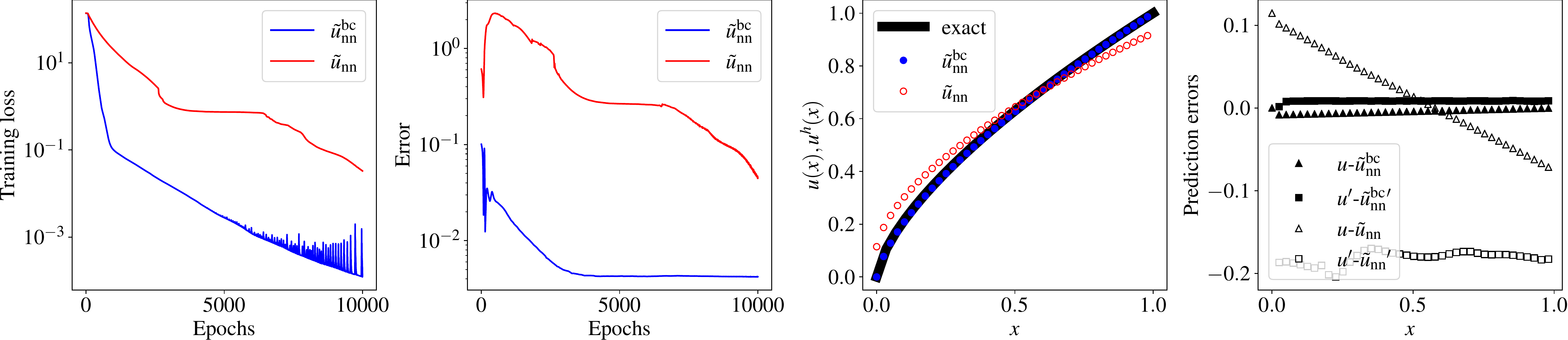}}; 
\node at (-2.28in,-0.84in) { (a) };
\node at (-0.67in,-0.84in) { (b)};
\node at ( 0.94in,-0.84in) { (c)};
\node at (2.55in,-0.84in) { (d)};
\end{tikzpicture}
\caption{Collocation solutions for a Dirichlet problem with a singular body force,  $b(x) = 2 x^{-4/3} / 9$. Essential boundary
conditions are $u(0) = 0$ and $u(1)= 1$, and the network
architecture is 1--50--50--1. (a), (b) Training loss and normalized absolute error as a function of epochs for $\unnbcxst$ and $\unnxst$. 
(c) Comparisons of $\unnbcxs$ and $\unnxs$ with the exact solution.
(d) Error in the displacement and strain fields.
}
\label{fig:elasticRod-viii}
\end{figure}
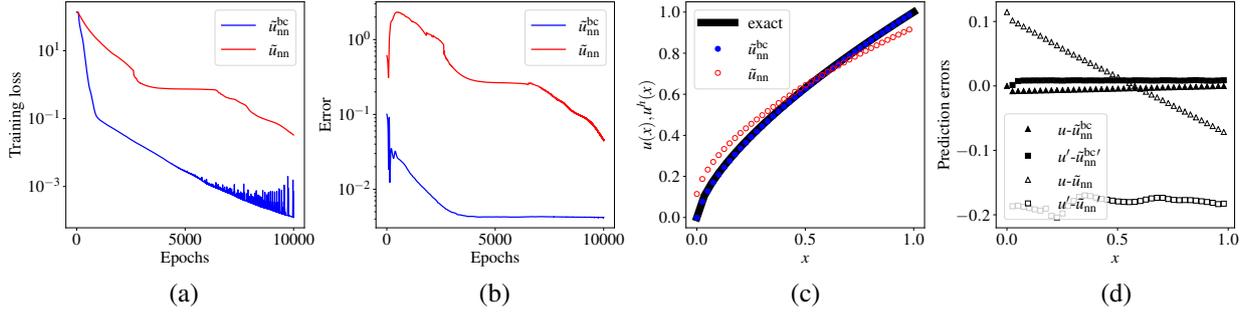

To obtain a weakly singular solution for the elastic rod problem, we consider a body force that has a singularity at the origin. We choose
$\Omega = (0,1)$ with $b(x) = 2 x^{-4/3} / 9$ with essential
boundary conditions $u(0) = 0$ and $u(1) = 1$. The exact solution is $u(x) = x^{2/3}$, and $u \in H^1(\Omega)$ is weakly singular. 
This problem is solved using the collocation approach on a 1--50--50--1 network architecture. Numerical results for $\unnbcxs$ and $\unnxs$ are presented in~\fref{fig:elasticRod-viii}.
The training loss of $\unnbcxst$ at 10,000 epochs is close to
${\cal O}(10^{-4})$ and $\unnxst$ is more than two orders larger. This same trend is observed in the normalized absolute error as a function of epochs (see~\fref{fig:elasticRod-viii}b).  Figures~\ref{fig:elasticRod-viii}c and~\ref{fig:elasticRod-viii}d show that
$\unnbcxs$ is in fairly good agreement with $u(x)$, whereas the error in the displacement and strain fields of $\unnxs$ are appreciable. It appears that $\unnxs$ and $u(x)$ differ by close to an affine function, which one can infer as being present within $\unnxs$ to meet the essential boundary conditions.

\subsubsection{Example 7}\label{subsubsec:RBFExample}
To draw connections to meshfree methods based 
on RBFs~\cite{Fasshauer:2007:MAM} and local maximum-entropy approximants~\cite{Arroyo:2006:LME}, which are discussed in~\sref{sec:intro}, we solve~\eqref{eq:elasticRod} with
inhomogeneous 
Dirichlet boundary conditions. The domain
$\Omega = (0,1)$, and we choose the exact solution as:
\begin{equation}\label{eq:exact_RBF}
    u(x) = \sum_{i=1}^2 \exp \left[ - \gamma_i ( x - a_i)^2 
    \right] ,
\end{equation}
which is a sum of two Gaussian functions,
and $\gamma_i$ and $a_i$ ($i=1,2$) are constants. 
The body force $b(x) = - u''(x)$.  Essential boundary conditions 
are imposed at $x = 0$ and $x = 1$ that are consistent with the exact solution in~\eqref{eq:exact_RBF}.

\begin{figure}[htp]
\centering
\begin{tikzpicture}
\node at (0,0) {\includegraphics[width=0.98\textwidth]{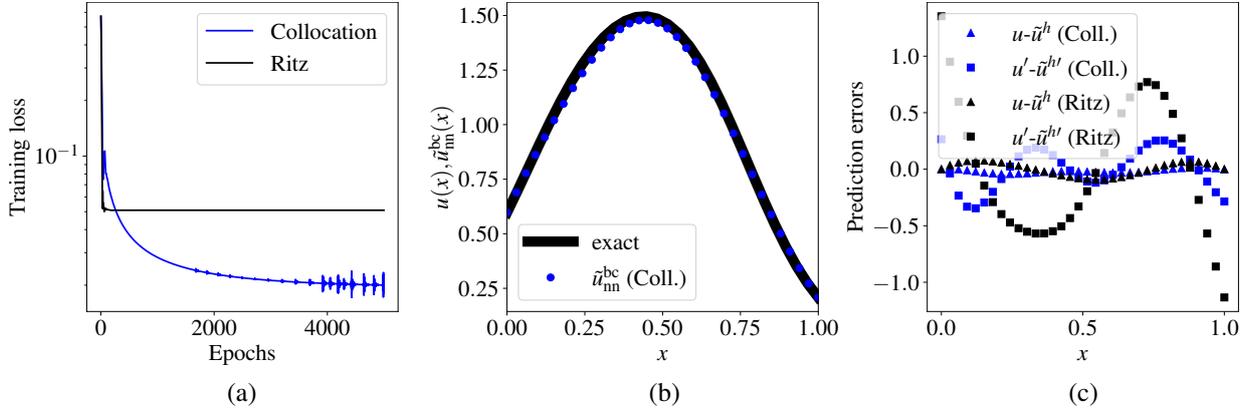}}; 
\node at (-1.98in,-1.10in) { (a) };
\node at (0.21in,-1.10in) { (b)};
\node at ( 2.38in,-1.10in) { (c)};
\end{tikzpicture}
\caption{Collocation and Ritz solutions using a Gaussian
activation function for a Dirichlet problem with the exact solution
as the sum of two Gaussians.  Network
architecture is 1--10--1. (a) Training loss as a function of epochs for $\unnbcxst$ (collocation and Ritz methods). 
(b) Comparisons of $\unnbcxs$ (collocation) with the exact solution.
(c) Error in the displacement and strain fields for collocation 
    and Ritz methods.
}
\label{fig:elasticRod-RBF}
\end{figure}

In the numerical
computations, we choose $a_1 = 1/4$,
$a_2 = 6/10$, $\gamma_1 = 9$ and $\gamma_2 = 10$. The
network architecture is 1--10--1 (1 hidden layer). For the hidden layer, we select
a Gaussian activation function, $\sigma(x) = \exp(-x^2)$, and
a linear activation function for the output layer.  The centers 
of the Gaussian are chosen to be fixed, and only the support-widths of the Gaussians and the weights in the output layer are the unknown parameters in the network.
The centers for the neurons are chosen as
$\vm{b} = \{0,\,1/9, \,2/9, \dots ,1\}$. In~\fref{fig:elasticRod-RBF}, the numerical
results for $\unnxst$ using collocation and Ritz are presented.
By the end of the training, the loss for the collocation and Ritz solutions are ${\cal O}(10^{-1})$ or better. Note that the loss measures
are distinct for the collocation and Ritz solutions; PDE loss is shown for the former whereas it is the (Ritz) energy loss that is presented for the latter.  The collocation solution is in
good agreement with the exact solution (\fref{fig:elasticRod-RBF}b). From~\fref{fig:elasticRod-RBF}c, we observe that the errors in the displacement field for both methods are small, but the errors in the derivative (strain fields) are appreciable. With more number of neurons in the hidden layer, it is expected that the accuracy in
the strain field will substantially improve.
This example reveals the flexibility that the PINN affords in that 
variational adaptive solutions can
be captured by a Gaussian neural network with a single hidden 
layer.  Realizing this is much more difficult using meshfree basis functions, since the
underlying Ritz formulation becomes a nonlinear, nonconvex minimization problem. 

\subsection{Longitudinal vibrations of a homogeneous elastic rod}
\label{subsec:rod_vibrations}
The eigenproblem for the longitudinal vibrations of an elastic bar that is fixed at both ends is:
\begin{subequations}
\begin{align}
u^{\prime \prime} + \omega^2 u &= 0 \quad \textrm{in } \Omega = (0,1),\\
u(0) &= u(1) = 0 .
\end{align}
\end{subequations}
The exact eigenfunctions are: $u_n(x) = \sin(\omega_n x)$,
where $\omega_n = n \pi$ ($n \in \set{N}$) are the natural frequencies. The eigenvalue $\lambda_n = \omega_n^2$ corresponds to the eigenfunction $u_n(x)$.  

We use the Ritz method to solve this problem using $\unnbcxst$ and
$\unnxst$.  The Rayleigh quotient minimization problem for the smallest eigenvalue (lowest mode) is~\suku{\cite{E:2018:DRM}}:
\begin{equation}\label{eq:rodvibration_var_problem}
    \min_{u \in {\cal S}} \frac{\int_0^1 u'^2 \, dx}
                               {\int_0^1 u^2 \, dx}, \quad
    \textrm{subject to} \ \int_0^1 u^2\, dx =1,
\end{equation}
where ${\cal S} = \{ u : u \in H^1(\Omega), \ u(0) = u(1) = 0\}$ and the normalization constraint on the eigenfunction appears 
in~\eqref{eq:rodvibration_var_problem}. For the trial function $\unnbcxst$, the loss function is:
\begin{equation}\label{eq:RQ}
   \Lnnbct = \frac{ \sum_{k=1}^{N_I} \bigl(  
    \unnbc {}^\prime  (x_k;\vm{\theta}) \bigl)^2 }
         {  \sum_{k=1}^{N_I}
            \bigl( \unnbc (x_k;\vm{\theta}) \bigr)^2    }
    + \left[ \left\{ \frac{1}{N_I} \sum_{k=1}^{N_I}  
      \big( \unnbc(x_k;\vm{\theta}) \bigr)^2 \right\} - 1 \right]^2 ,
\end{equation}
where $N_I$ is the number of interior integration points. 
\suku{Note that 
apart from the PDE loss term, we have an additional loss term due to the
normalization constraint in~\eqref{eq:rodvibration_var_problem}.}
The  loss function for $\unnxst$ consists of the 
\suku{two contributions}
that appears in~\eqref{eq:RQ}, and in addition it will include two boundary loss terms to impose the essential
boundary conditions.

\begin{figure}[htp]
\centering
\begin{tikzpicture}
\node at (0,0) {\includegraphics[width=0.98\textwidth]{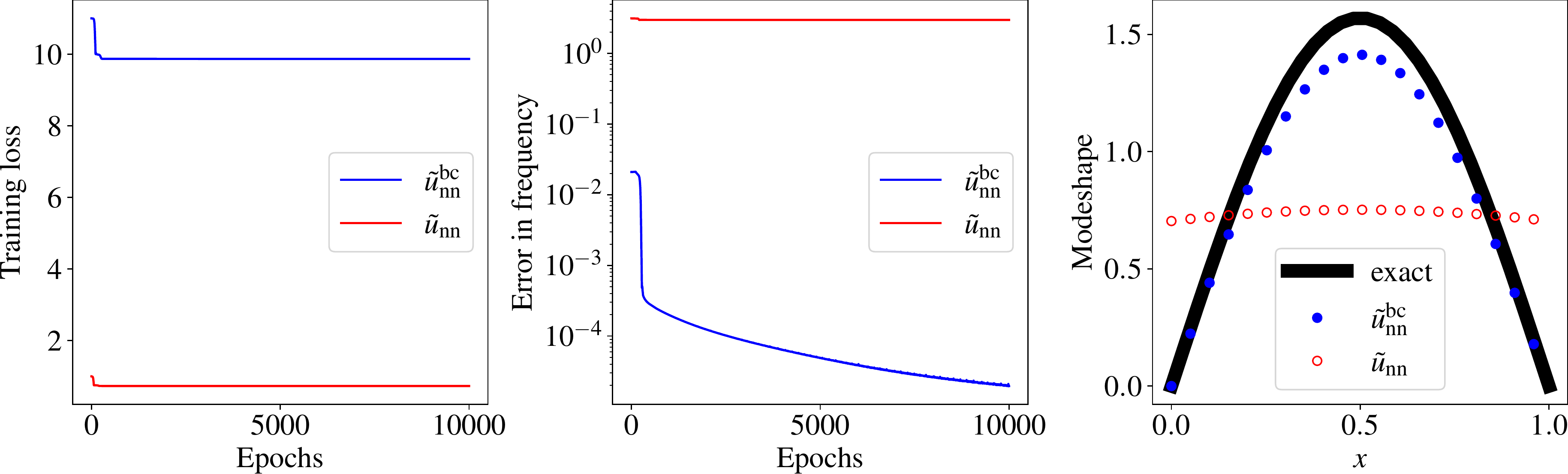}}; 
\node at (-2.04in,-1.12in) { (a) };
\node at (0.16in,-1.12in) { (b)};
\node at ( 2.35in,-1.12in) { (c)};
\end{tikzpicture}
\caption{Ritz solutions for the longitudinal vibrations (lowest mode) of
a homogeneous elastic rod. The network architecture is 1--50--50--50--1. (a), (b) Training loss and error in lowest natural frequency as a function of epochs for $\unnbcxst$ and $\unnxst$.
(c) Comparisons of $\unnbcxs$ and $\unnxs$ for the lowest mode (eigenfunction) that corresponds to the lowest natural frequency $\omega_1 = \pi$.} 
\label{fig:Ritz-eig}
\end{figure}

The network architecture 1--50--50--50--1 is used. In~\fref{fig:Ritz-eig}, the Ritz solutions for $\unnbcxs$ and $\unnxs$ are presented. The loss function for $\unnbc$ and
$\unn$ saturate to values of 1 and 10, respectively. The error in the natural frequency for $\unnbcxst$ and $\unnxst$ are ${\cal O}(10^{-4})$ and ${\cal O}(1)$ at 10,000 epochs (see~\fref{fig:Ritz-eig}b). In~\fref{fig:Ritz-eig}c, we compare the lowest mode (eigenfunction) from the numerical solutions to the exact solution: $\unnbcxs$ is in much better agreement with the exact solution than $\unnxs$. The exact mode shape is well-captured by 
$\unnbcxs$ but it is has not been exactly normalized, which leads to the observed discrepancy in the maximum amplitude.

\subsection{Advection-diffusion problem}\label{subsec:adv-diff}
\begin{subequations}\label{eq:ad-bvp}
\begin{align}
u'' &= \alpha u' \ \ \textrm{in } \Omega = (0,1) \\
& u(0) = 0,  \quad  u(1) = 1,
\end{align}
\end{subequations}
where $\alpha$ is the Peclet number, which measures the ratio of the advective rate
to the diffusion rate.  The exact solution of the problem posed in~\eqref{eq:ad-bvp} is:
\begin{equation}
u(x) = \frac{e^{\alpha x} - 1}{e^{\alpha} - 1}.
\end{equation}

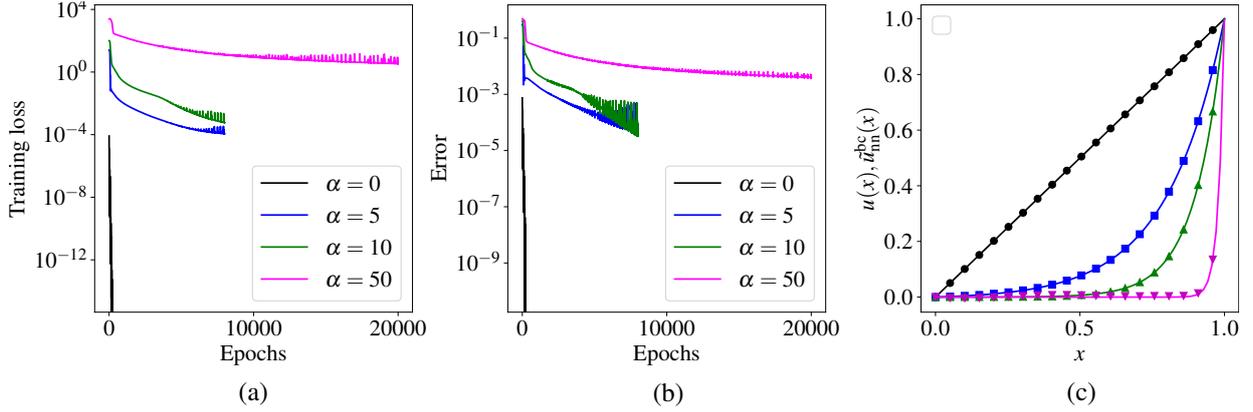
\begin{figure}[htp]
\centering
\begin{tikzpicture}
\node at (0,0) {\includegraphics[width=0.98\textwidth]{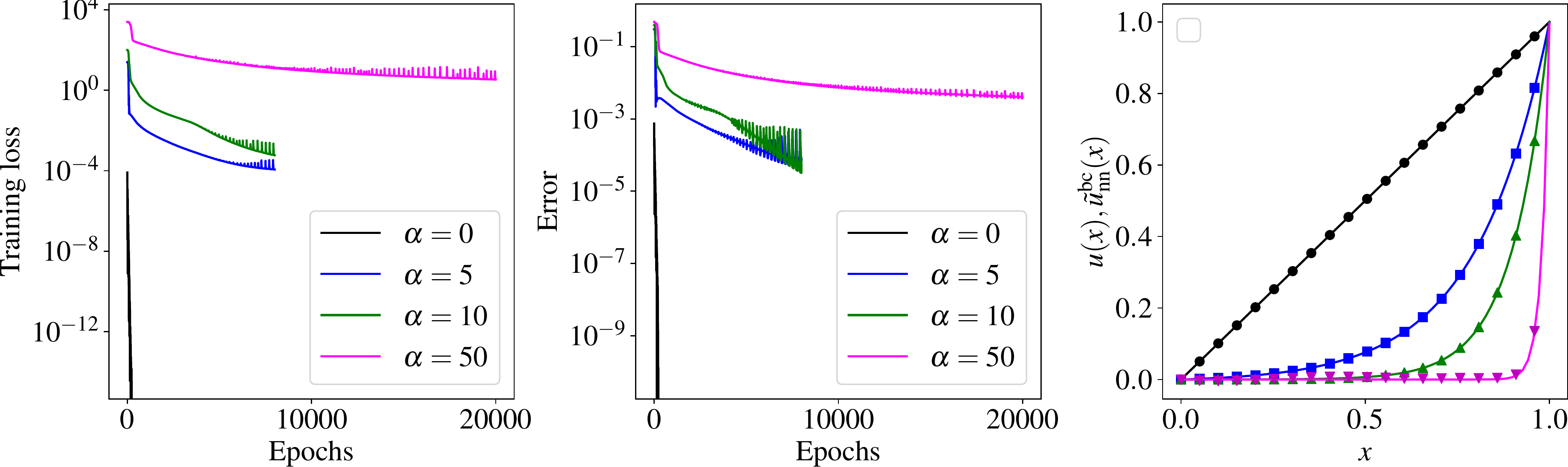}}; 
\node at (-1.92in,-1.10in) { (a) };
\node at (0.23in,-1.10in) { (b)};
\node at ( 2.37in,-1.10in) { (c)};
\end{tikzpicture}
\caption{Collocation solutions using $\unnbcxst$ for the advection-diffusion problem. For $\alpha = 1,5,10$, the network architecture 1--50--50--1 is used, and for $\alpha = 50$, the architecture is 1--50--50--50--1. 
(a), (b) Training loss and normalized absolute error as
a function of epochs. 
(c) Comparisons of $\unnbcx$ to the exact solutions for
different $\alpha$.
The solid lines are for the
exact solutions and the markers (colors are consistent with
those shown in (b)) represent the numerical solutions. 
}
\label{fig:AD}
\end{figure}

We choose $\alpha = 0,5,10,50$ in this study (pure diffusion for $\alpha = 0$ to strongly advective flow for $\alpha = 50$) and run collocation simulations using $\unnbcxst$. For $\alpha = 1,5,10$, the network architecture 1--50--50--1 is used, and for $\alpha = 50$, the architecture is 1--50--50--50--1. In~\fref{fig:AD},
the simulation results are presented. For all cases
shown in~\fref{fig:AD}b, we observe that the normalized absolute error during the training is ${\cal O}(10^{-3})$ or less. For all $\alpha$ that are selected, \fref{fig:AD}c shows an excellent match between  $\unnbcxs$ and the exact solutions.

\subsection{Euler-Bernoulli beam bending}\label{subsec:beam}
Consider the boundary-value problem for the deflection of a cantilever (Euler-Bernoulli) beam of unit length that is clamped at both ends
and is subjected to a distributed load $q(x)$: 
\begin{subequations}\label{eq:beam_bvp}
\begin{align}
\label{eq:beam_bvp-a}
EI v^{\prime \prime\prime \prime} &= q \ \ \textrm{in } \Omega = (0,1),\\
\label{eq:beam_bvp-b}
v(0) = v'(0) &= v(1) = v'(1) = 0 ,
\end{align}
\end{subequations}
where $EI$ is the flexural rigidity of the beam.

We apply a point moment (clockwise orientation) of magnitude $M_0$ at $x = 1/2$ so that $q(x) = M_0 \delta^\prime(x - 1/2)$. 
For this
$q(x)$, the variational
principle that is associated with the strong form in~\eqref{eq:beam_bvp} is:
\begin{equation}
\min_{v \in S} \left[ \Pi[v] 
= \frac{1}{2} \int_0^1 EI (v'')^2 dx + M_0 v'(1/2), \ \ 
S = \biggl\{ v : v \in H^1(0,1), \ v(0) = v'(0) = v(1) = v'(1) = 0
\biggr\} \right].
\end{equation}
\suku{
For this problem, all homogeneous
boundary conditions associated with $v$ and $v^\prime$ that appear in~\eqref{eq:beam_bvp-b} are essential boundary conditions, and hence a 
kinematically admissible PINN trial function is
given by~\eqref{eq:trial_plate}:}
\begin{equation}
    \tilde{v}_\textrm{nn}^\textrm{bc}(x;\vm{\theta}) =
    [ \phi (x) ]^2 \tilde{v}_\textrm{nn}^R (x;\vm{\theta}) ,
\end{equation}
where $\phi(x)$ is an ADF (normalized to order 1) that vanishes at $x = 0$ and $x = 1$, and its normal derivative has unit magnitude on the boundary.

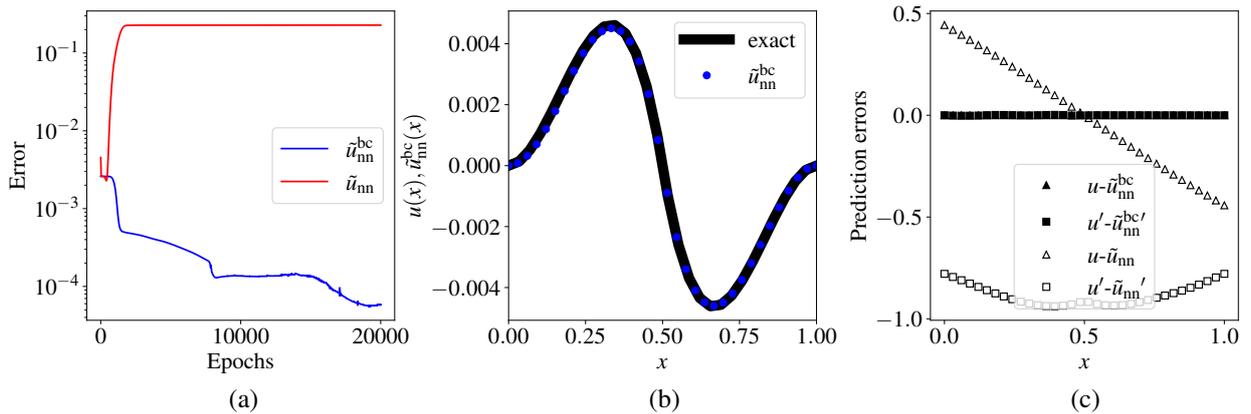
\begin{figure}[htp]
\centering
\begin{tikzpicture}
\node at (0,0) {\includegraphics[width=0.98\textwidth]{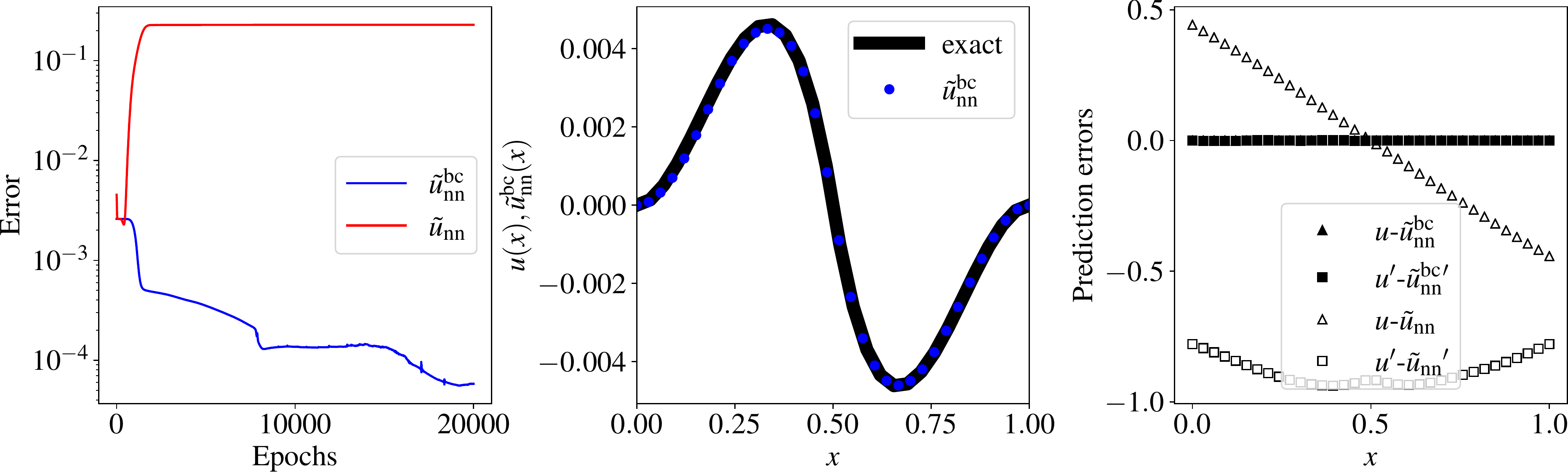}}; 
\node at (-1.97in,-1.10in) { (a) };
\node at (0.21in,-1.10in) { (b)};
\node at ( 2.39in,-1.10in) { (c)};
\end{tikzpicture}
\caption{Deep Ritz solutions using $\unnbcxst$ and $\unnxst$ for the Euler-Bernoulli beam problem. The network architecture 1--50--50--1 is used. 
(a) Normalized absolute eror as a function of epochs.
(b) Comparison of $\unnbcxs$ with the exact solution.
(c) Errors in the deflection and rotation fields.}
\label{fig:euler-b}
\end{figure}

For the computations, we choose $EI = 1$ and $M_0 = 1$.  The exact solution of~\eqref{eq:beam_bvp} is:
\begin{equation}\label{eq:beam_exact_v}
v(x) = \frac{ \big[\texttt{ReLU} (x-1/2) \bigr]^2}{2} + \frac{x^2}{8} 
       - \frac{x^3}{4}, \quad x \in \Omega,
\end{equation}
which is a $C^1(\Omega)$ (piecewise cubic) function with a moment discontinuity of 
unit magnitude at $x = 1/2$. 
Deep Ritz solutions $\unnbcxst$ and $\unnxst$ are computed on the network architecture 1--50--50--1 using the cubic ReLU activation
function, and the results are presented in~\fref{fig:euler-b}. From~\fref{fig:euler-b}a, we observe that $\unnbcxst$ converges to a loss of ${\cal O}(10^{-4})$ at 20,000 epochs but the loss for $\unnxst$ remains at about 0.1. Figure~\ref{fig:euler-b}b shows that the numerical solution $\unnbcxs$ is in excellent
agreement with the exact solution. In~\fref{fig:euler-b}c, 
the error fields are shown: 
the deflection and rotation fields using $\unnbcxs$ are accurate, whereas both fields have large errors for $\unnbcxs$. We attribute the poor performance of $\unnxst$ due to scaling issues of the interior and boundary terms in the loss function. If the ratio of the weights assigned to the boundary loss term and the interior loss terms is set to $10^3:1$, we find that the results improve but they are still worse than those obtained using $\unnbcxs$. 
It requires a weight ratio of $10^4$ for the two solutions to have comparable accuracy.

\section{Numerical Examples in Two \suku{Dimensions}}\label{sec:2D}
\suku{The promise of PINN is intriguing for inverse and
parameteric (design) problems. But this rests on its accuracy, robustness and reliability on solving the forward problem, which is the emphasis in this contribution. To this end, we focus on
the performance of the PINN formulation with exact imposition of boundary conditions versus the standard PINN~\cite{Raissi:2019:PIN} with equally-weighted loss terms.}
For the two-dimensional problems, we consider polygonal domains and also domains with curved boundaries.  We consider four distinct types of problems: steady-state heat conduction; computation of 
harmonic (Laplace equation) coordinates~\cite{Joshi:2007:HCF}, which is an instance of generalized barycentric coordinates~\cite{Hormann:2017:GBC}; clamped Kirchhoff plate bending (fourth-order PDE); and the Eikonal equation to compute the signed distance function to a boundary. 
\suku{For these problems, we present our solutions and 
compare them to either the exact solution (if available) or a reference finite element solution, and to deep collocation~\cite{Raissi:2019:PIN}.}
In addition, we identify key differential
properties of the approximate distance function and
bring to fore the issue of exact satisfaction of boundary conditions and its implications in the training of the network and the accuracy that the PINN approximation delivers.

Prior to presenting the numerical examples, it is instructive to
understand the properties and behavior of approximate distance functions that are obtained by either R-functions with R-equivalence as presented in~\sref{subsec:Requiv} or via mean value potential fields that are discussed in~\sref{sec:mvp}. Since these functions are used in the PINN ansatz $\unnbcxt$ that is presented in~\sref{sec:BCs}, one must consider the regularity of these functions when used in a 
deep collocation or a deep Ritz method.

\subsection{Laplacian of approximate distance fields
            }\label{subsubsec:Laplacian}
\begin{figure}[htp]
\centering
\begin{tikzpicture}
\node at (0,0) {\includegraphics[width=\textwidth]{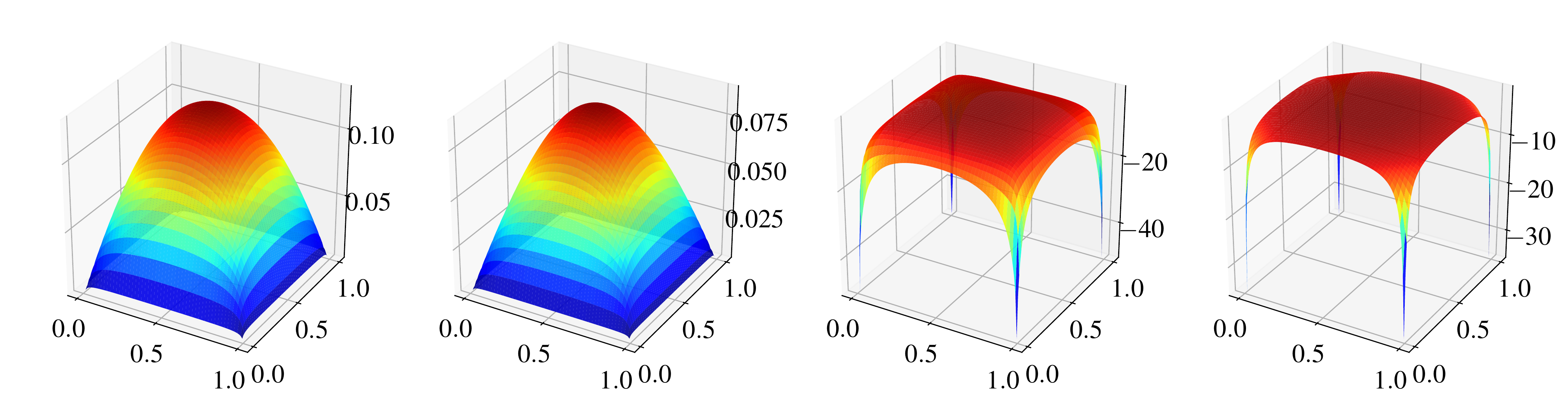}}; 
\node at (-2.35in,-0.95in) { (a) $\phi_R (\vx)$ };
\node at (-0.75in,-0.95in) { (b) $\phi_M (\vx)$};
\node at ( 0.85in,-0.95in) { (c) $\nabla^2 \phi_R (\vx)$};
\node at (2.45in,-0.95in) { (d) $\nabla^2 \phi_M (\vx)$};
\end{tikzpicture}
\caption{Computation of $\phi(\vx)$ and $\nabla^2 \phi(\vx)$ over the unit square for ADFs constructed from R-equivalence and mean value potential fields.
}
\label{fig:Heat-Laplacian}
\end{figure}
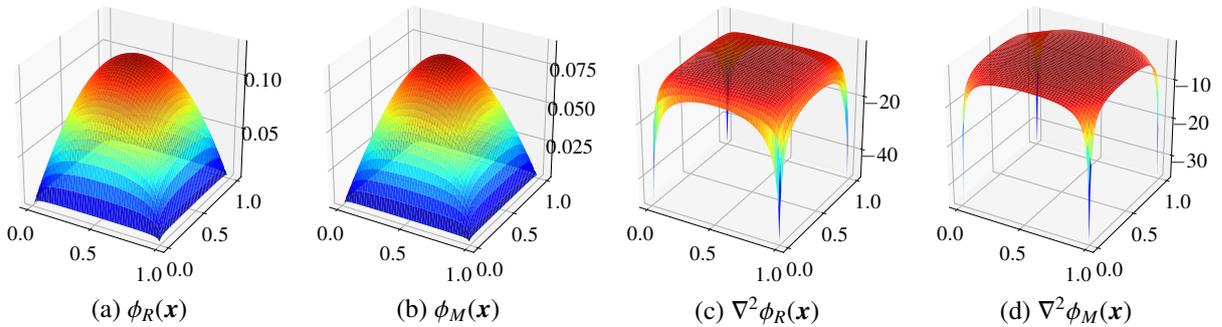

For Poisson or Laplace boundary-value problems that involve the Laplace operator, we must understand the behavior of the Laplacian of the ADFs that stem from R-functions and
mean value potentials. Let us consider the unit
square, $\Omega = (0,1)^2$. The boundary $\partial \Omega$ consists of four line segments. On using either~\eqref{eq:phin_eq} or~\eqref{eq:phi_mvc}, we can construct an approximate distance 
function to the boundary, $\phi(\vx)$, which is 
normalized to order $1$. Let us refer to these functions as 
$\phi_R(\vx)$ (REQ)
and $\phi_M(\vx)$ (MVP).
In~\fref{fig:Heat-Laplacian}, $\phi(\vx)$ and its Laplacian over the unit square are presented. We observe that both $\phi_R$ and 
$\phi_M$ are zero 
on the entire boundary and monotonic (concave) inside the domain. 
This property of these functions on $\partial \Omega$ is used to impose essential boundary conditions, as described in~\sref{subsec:EBCs}. From Figs.~\ref{fig:Heat-Laplacian}c and~\ref{fig:Heat-Laplacian}d, we observe that the
Laplacians, $\nabla^2\phi_R$
and $\nabla^2\phi_M$, blow up at the vertices of the square. In fact, it is known that both $\nabla^2\phi_R$ and 
$\nabla^2\phi_M$ are singular at the vertices of a polygon, and therefore very large in magnitude near any of its vertices. Therefore, in a collocation-based approach to solve the Poisson equation, which involves the Laplacian, the contributions to the total loss from regions near the vertices can be very large. This inference does not influence Ritz-based solutions of the Poisson equation since the highest derivative in the
variational principle is of order 1, and both $\phi_R$ and $\phi_M$ 
and its first-order derivatives are well-behaved (bounded) over the entire domain. 
There are two possible remedies to address this issue. The first involves 
modifying the $\phi_i$
that are obtained from R-equivalence and the ADFs that stem
from mean value potential fields, so that the corresponding
Laplacians are bounded in the domain. The second is to consider collocation points inside the domain that are not very close to the
vertices. We leave the first route as part of future work, and proceed in this paper with the second choice. For instance, we show in~\sref{subsubsec:heat:differentebcs} that if all interior collocation points are located in $\Omega_\delta = [\delta,\, 1-\delta]^2$ ($\delta = 0.01$), which is a subset of the unit square, then both methods
perform well. Finally, we reemphasize that it is imperative that in most instances $\phi$ be smooth in the interior of a computational domain; otherwise, $\nabla^2\phi$ will blow up at an interior collocation point and then one cannot use a trial function that uses $\phi$ in a collocation-based PINN method. So in most instances in 2D or 3D, this precludes the use of exact distance functions in the ansatz, and hence {\em approximate distance functions} should be used. 
\suku{To show this, let us} consider the exact signed distance to the unit disk, $\phi(\vx) = 1 - \sqrt{x^2 + y^2}$. We can write
\begin{equation*}
\frac{\partial^2 \phi}{\partial x^2} = 
\frac{x^2}{(x^2+y^2)^{3/2}} - \frac{1}{\sqrt{x^2+y^2}},
\end{equation*}
and since the second term is unbounded at the origin, 
the Laplacian of $\phi$ blows up at the origin. 
There are exceptions when the exact distance function can be used. It is a suitable choice when the medial axis of a domain is not part of the computational domain. For example, when solving a boundary-value problem over an annulus in 2D (see the problem solved in~\sref{subsubsec:heatM}) or a hollow cylinder in 3D, then the exact distance function can be used since the origin (where the exact distance function has derivative 
discontinuities) lies outside the computational domain. 

\subsection{Steady-state heat conduction}\label{subsec:heat}
Let us consider the following model problem for isotropic 
steady-state heat conduction:
\begin{subequations}\label{eq:heat_bvp}
\begin{align}
-\nabla^2 u &= f \ \ \textrm{in } \Omega \subset \Re^2 \\
u = g \ \  \textrm{on }  \Gamma_u &,  \quad
 \frac{\partial u}{\partial n} + c u = h \ \ \textrm{on } \Gamma_n
\end{align}
\end{subequations}
where $u(\vx)$ is the temperature field and
$f(\vx)$ is the heat source. The boundary $\partial \Omega
= \Gamma_u \cup \Gamma_n$ is partitioned into two parts, with
$\Gamma_n \cap \Gamma_n = \emptyset$. The temperature field
$g(\vx)$ is imposed on the essential boundary
$\Gamma_u$, and boundary data $h(\vx)$ that is associated with a Robin boundary condition is prescribed on
$\Gamma_n$ ($c$ is in general a spatially varying field).

\subsubsection{Essential boundary conditions} \label{subsubsec:heatE}
As the first example, we consider the
biunit square,
$\Omega = (-1,1)^2$, with $u = g = 0$ prescribed on $\partial 
\Omega$. If $k \in \mathbb{N}$ and
$f(x,y)=\sin(k\pi x)\sin(k\pi y)$ is the forcing function, then the
exact solution for this problem is:
\begin{equation*}
u(\vx) = \frac{  \sin(k\pi x) \sin(k\pi y) } {2k^2\pi^2}. 
\end{equation*}
In the numerical computations, we consider two distinct forms
of trial functions in the neural network. The first ansatz is 
the standard PINN that is given by $\unnxt$, which does not a priori satisfy the boundary condition. The second form consists of trial functions that are given by $\unnbcxt = \phi(\vx)\unnRxt$, where $\phi(\vx)$ is a function that is zero on $\partial\Omega$. This property of $\phi(\vx)$ ensures that $\unnbcxt$ automatically satisfies the essential boundary conditions. While an obvious choice for $\phi(\vx)$ is $(1-x^2)(1-y^2)$, here we consider $\phi(\vx)$ that are constructed using R-functions with
R-equivalence composition (see~\sref{subsec:Requiv}) and by mean value potential fields (see~\sref{sec:mvp}) as they readily
generalize to more complex domains. When needed, we
use the acronyms REQ
and MVP to distinguish the numerical solutions, $\unnbcx$, which are obtained using these two methods.
The plot of $\phi(\vx)$ using REQ and MVP over a square are shown in Figures~\ref{fig:phi_polygons-b} and~\ref{fig:mvc_ADF-a}, and it can be observed that $\phi$
is zero on the boundary of the domain. 
For the collocation scheme, we randomly sample $N_I$ number of points in $\Omega$ and $N_B$ number of points on $\partial\Omega$. To solve the problem using standard PINN, we minimize the loss $\Lnnt$  
given in~\eqref{eq:collocation_loss_unn}, which is reproduced below
(using $g = 0$):
\begin{equation*}
\Lnnt =  ||\nabla^2\unnxt+f(\vx)||^2_{\Omega,N_I} 
         +  ||\unnxt||^2_{\Omega,N_B} ,
\end{equation*}
where $|| \cdot ||_{\Omega,N_I}$ and 
$|| \cdot ||_{\partial \Omega,N_B}$ are defined in~\eqref{eq:collocation_loss} and~\eqref{eq:collocation_loss_unn}. 
Since $\unnbcxt$ automatically satisfies the boundary condition, the parameters in this ansatz are found by minimizing the loss given
in~\eqref{eq:collocation_loss}, which for this problem is:
\begin{equation*}
\Lnnbct  =  ||\nabla^2\unnbcxt+f(\vx)||^2_{\Omega,N_I} .
\end{equation*}

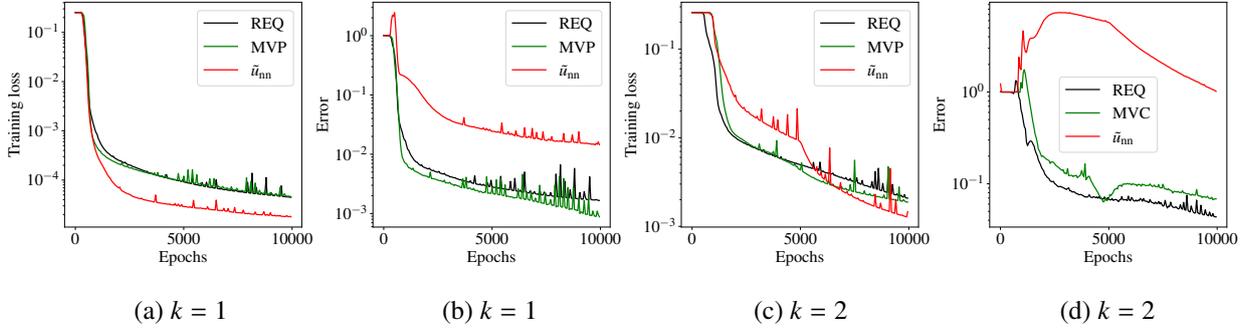
\begin{figure}[htp]
\centering
\begin{tikzpicture}
\node at (0,0) {\includegraphics[width=0.98\textwidth]{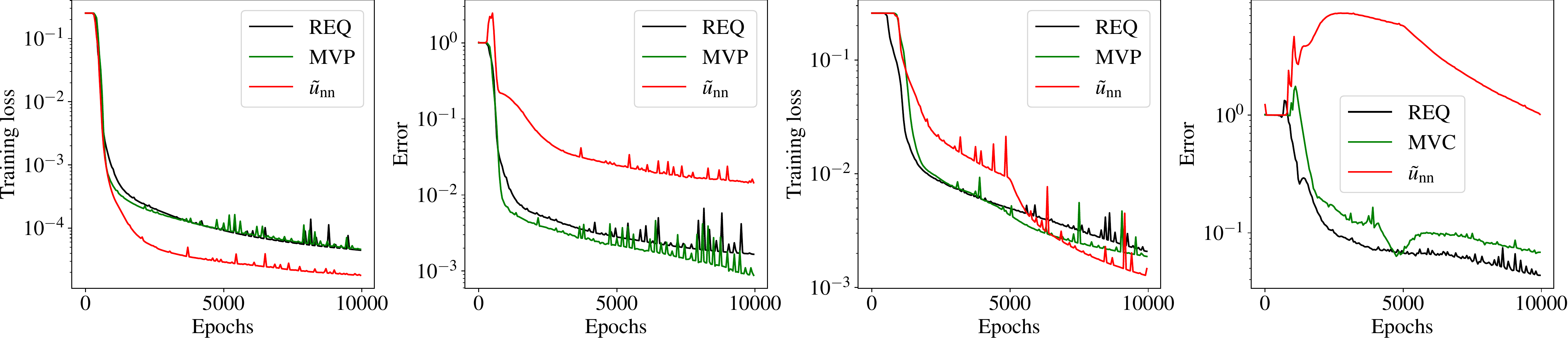}}; 
\node at (-2.29in,-0.93in) { (a) $k=1$ };
\node at (-0.69in,-0.93in) { (b) $k=1$};
\node at ( 0.92in,-0.93in) { (c) $k=2$};
\node at (2.51in,-0.93in) { (d) $k=2$};
\end{tikzpicture}
\caption{Training loss and normalized absolute errors for 
$\unnbcxt$ (REQ and MVP) and $\unnxt$
in the heat conduction problem with
homogeneous Dirichlet boundary conditions and forcing function
$\sin k\pi x$. (a), (b) $k = 1$, and (c), (d) $k = 2$.} 
\label{fig:HeatEq}
\end{figure}

Figure~\ref{fig:HeatEq} shows the training loss and 
normalized absolute error as functions of the training epochs for $k=1,\,2$ ($N_I = 5,000,\, N_B = 400$). We observe from~Figs.~\ref{fig:HeatEq}a 
and~\ref{fig:HeatEq}c that the training loss for
$\unnxt$ is either comparable to or less than the loss for $\unnbcxt$ over the same number of epochs. However, this does not translate into better prediction accuracy. Figures~\ref{fig:HeatEq}b and~\ref{fig:HeatEq}d show the 
normalized absolute errors for $\unnbcxt$ (REQ and MVP) as well as
$\unnxt$. These plots
reveal that both REQ- and MVP-based schemes deliver
an order of magnitude more accurate solutions compared to $\unnxt$. 

\begin{figure}[!htb]
\centering
\begin{tikzpicture}
\node at (0,0) {\includegraphics[width=0.98\textwidth]{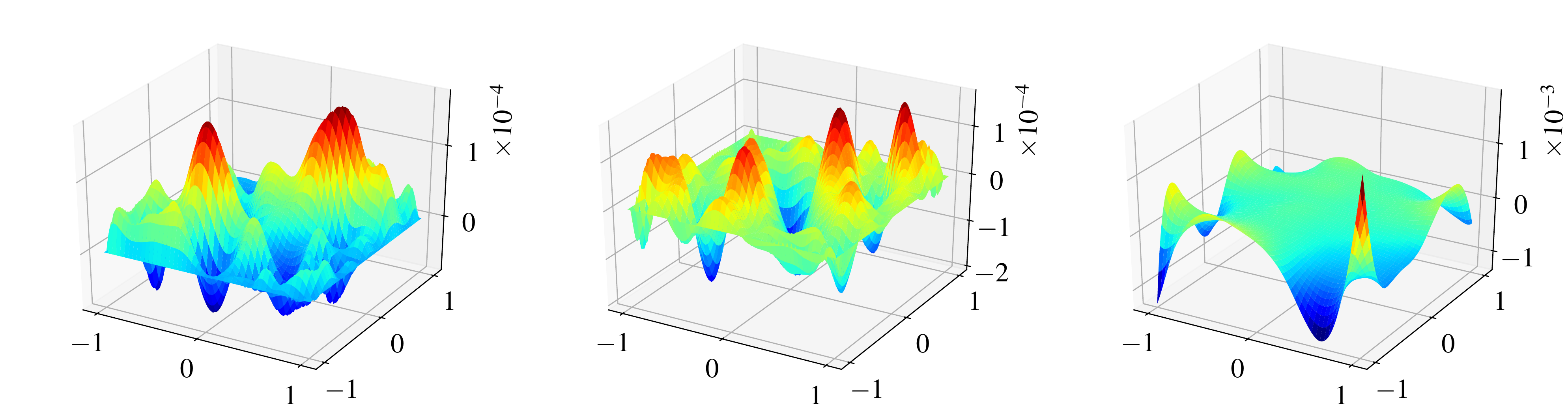}};
\node at (-2.06in,-1in) { (a)};
\node at (0.04in,-1in) { (b)};
\node at ( 2.18in,-1in) { (c)};
\end{tikzpicture}
\caption{Prediction errors for $k = 1$ in Example 1. 
(a) $\unnbcx$ (REQ), (b) $\unnbcx$ (MVP), and (c) $\unnx$.} 
\label{fig:HeatEq-errors}
\end{figure}

It is interesting to note that $\unnbcxt$ produces smaller 
normalized absolute errors during training than $\unnxt$ even though it has larger losses. This observation is noticed in almost all cases, and it deserves some comments here. We mention that there is no reason to assume that the losses for $\unnbcxt$ and 
$\unnxt$ should be comparable. For the problem under
consideration, $\Lnnbct$ comprises of terms that involve the derivatives of $\phi(\vx)$ and these terms are not present in $\Lnnt$. A further issue with $\unnxt$ is the relative scaling of the losses in $\Lnnt$, which comprises of a loss on $\nabla^2 \unnxt$ and a loss on $\unnxt$. In this problem, the norm of $u$ is much smaller than the norm of its Laplacian, and therefore the $\nabla^2 \unnxt$ term dominates in $\Lnnt$.  However, the optimizer can drive $\Lnnt$ to very
small values without adequately addressing the boundary loss term. We can see that this is indeed the issue if we compare the prediction errors for $\unnx$ over the domain (see~\fref{fig:HeatEq-errors}). We see that the prediction errors for REQ and MVP schemes are much smaller than $\unnx$, and furthermore, the prediction errors for $\unnx$ are large near the boundary of the domain. This suggests that $\unnxt$ is undervaluing the boundary loss term in $\Lnnt$. In principle, it is possible to improve the results for $\unnxt$ by assigning a larger weight to the boundary loss term in $\Lnnt$~\cite{Wang:2020:UMG}, but this is an ad-hoc remedy,  which is not needed in our approach since the boundary condition is exactly met. To see that this is indeed the case, one can assume the loss function for $\unnxt$ to be a
convex combination of the interior and boundary loss terms:
\begin{equation*}
\Lnnt =  \underset{\textrm{PDE loss}}{\underbrace{
w \, ||\nabla^2\unnxt+f(\vx)||_{\Omega,N_I}^2 }}
        +  
        \underset{\textrm{Boundary loss}}
        {\underbrace{(1-w) \, ||\unnxt||_{\partial \Omega, N_B}^2 }},
\end{equation*}
where $w \in [0,1]$ is a scalar that can be used to tune the relative weights of the two losses. 

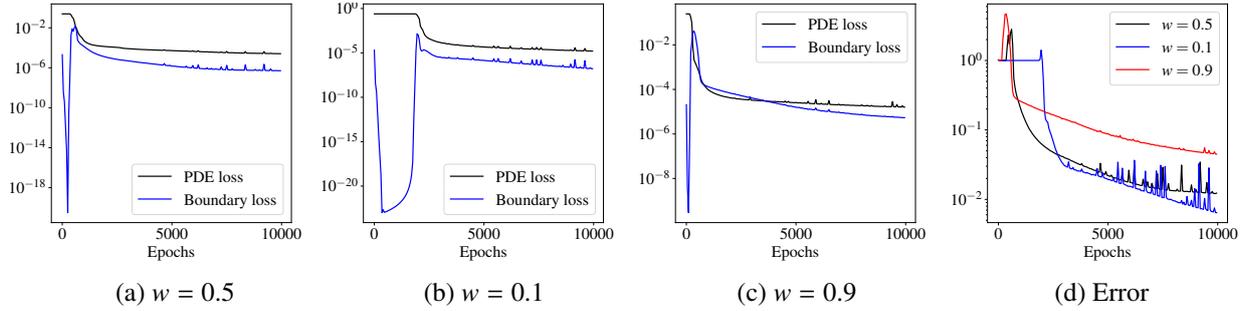
\begin{figure}[!htb]
\centering
\begin{tikzpicture}
\node at (0,0) {\includegraphics[width=0.98\textwidth]{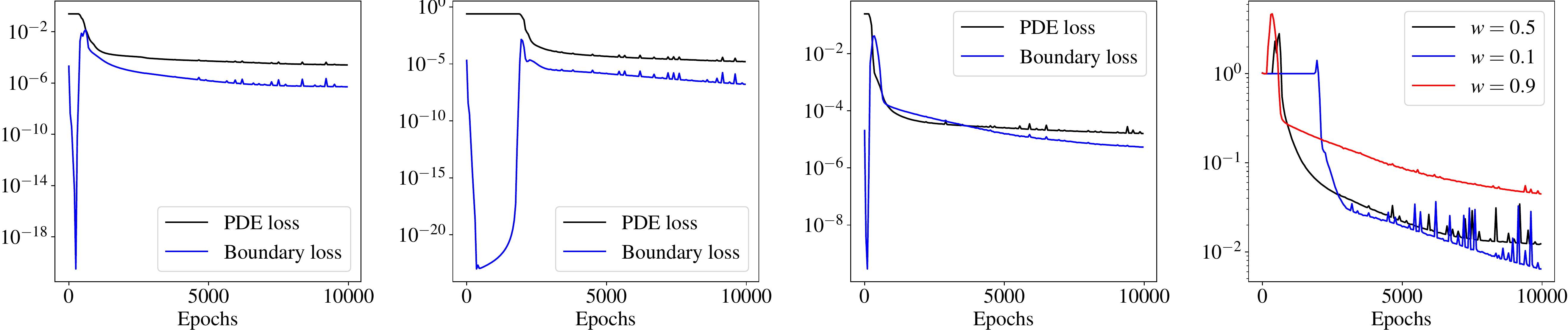}}; 
\node at (-2.32in,-0.85in) { (a) $w=0.5$ };
\node at (-0.72in,-0.85in) { (b) $w=0.1$};
\node at ( 0.90in,-0.85in) { (c) $w=0.9$};
\node at (2.48in,-0.85in) { (d) Error};
\end{tikzpicture}
\caption{Loss function for $\unnxt$ in Example 1 is a convex combination of
PDE loss and boundary loss terms ($k = 1$).
Evolution of loss function for 
(a) $w = 0.5$, (b) $w = 0.1$ and
(c) $w = 0.9$. (d) Evolution of normalized absolute error
during training 
for $w = 0.5,\,0.1,\,0.9$.} \label{fig:HeatEq-LossMetric}
\end{figure}

Figures~\ref{fig:HeatEq-LossMetric}a--\ref{fig:HeatEq-LossMetric}c show the evolution of the PDE loss and the
boundary loss
for different values of $w$ ($k = 1$). The value $w = 0.1$ weighs the boundary loss term 9 times more than the PDE loss term and it achieves smaller error than both $ w = 0.5, \, 0.9$ (\fref{fig:HeatEq-LossMetric}d). This suggests that there likely
is a sweet spot for $w$ that results in very low errors. However, the exact regime for this solution is likely dependent on the problem under consideration and on the boundary conditions, the determination of which may be impossible to ascertain in problems where the exact solution is unknown. A distinguishing attribute of our approach is that no such tuning of relative weights is needed. 

\subsubsection{Example 2}\label{subsubsec:heat:differentebcs}
We consider the Laplace equation ($f = 0$) over
the unit square, $ \Omega = (0,1)^2$, with boundary conditions
\begin{equation}\label{eq:Example2:EBC}
u(\vx) = 0 \ \ \textrm{on } \Gamma_1, \ \Gamma_2, \ \Gamma_3, \quad
u(\vx) = g_4(\vx) = \sin \pi x \ \ \textrm{on } \Gamma_4,
\end{equation}
where $\Gamma_1 = \{ (x,y) : x = 0, \ 0 \le y \le 1 \}$,  
$\Gamma_2 = \{ (x,y) : 0 \le x \le 1, \ y = 0 \}$,
$\Gamma_3 = \{ (x,y) : x = 1, \ 0 \le y \le 1\}$, and
$\Gamma_4 = \{ (x,y) : 0 \le x \le 1, \ y = 1\}$ are the
boundary edges. The exact solution for this problem is:
\begin{equation*}
u(\vx) = \frac{ ( e^{-\pi y} + e^{\pi y} ) \sin \pi x }
              {e^{-\pi} + e^{\pi}} .
\end{equation*}

We chose this problem to demonstrate how to exactly satisfy
nonzero essential boundary conditions on different subsets of
the boundary through the use of transfinite interpolation. To construct a trial solution that satisfies the boundary 
conditions, we first create a composite approximate distance function, $\phi(\vx)$, to $\Gamma=\Gamma_1\cup\Gamma_2\cup\Gamma_3\cup\Gamma_4$. This ADF can either be formed by the joining operation via R-equivalence,
or directly via~\eqref{eq:phi_mvc} that uses mean value potential 
fields on polygons. The resultant $\phi(\vx)$ is similar to the $\phi(\vx)$ used in Example 1
(see~\fref{fig:Heat-Laplacian}). We combine the Dirichlet boundary data into one function $g(\vx)$ by using the 
transfinite interpolant in~\eqref{eq:ti}.

\begin{figure}[htp]
\centering
\begin{tikzpicture}
\node at (0,0) {\includegraphics[width=0.88\textwidth]{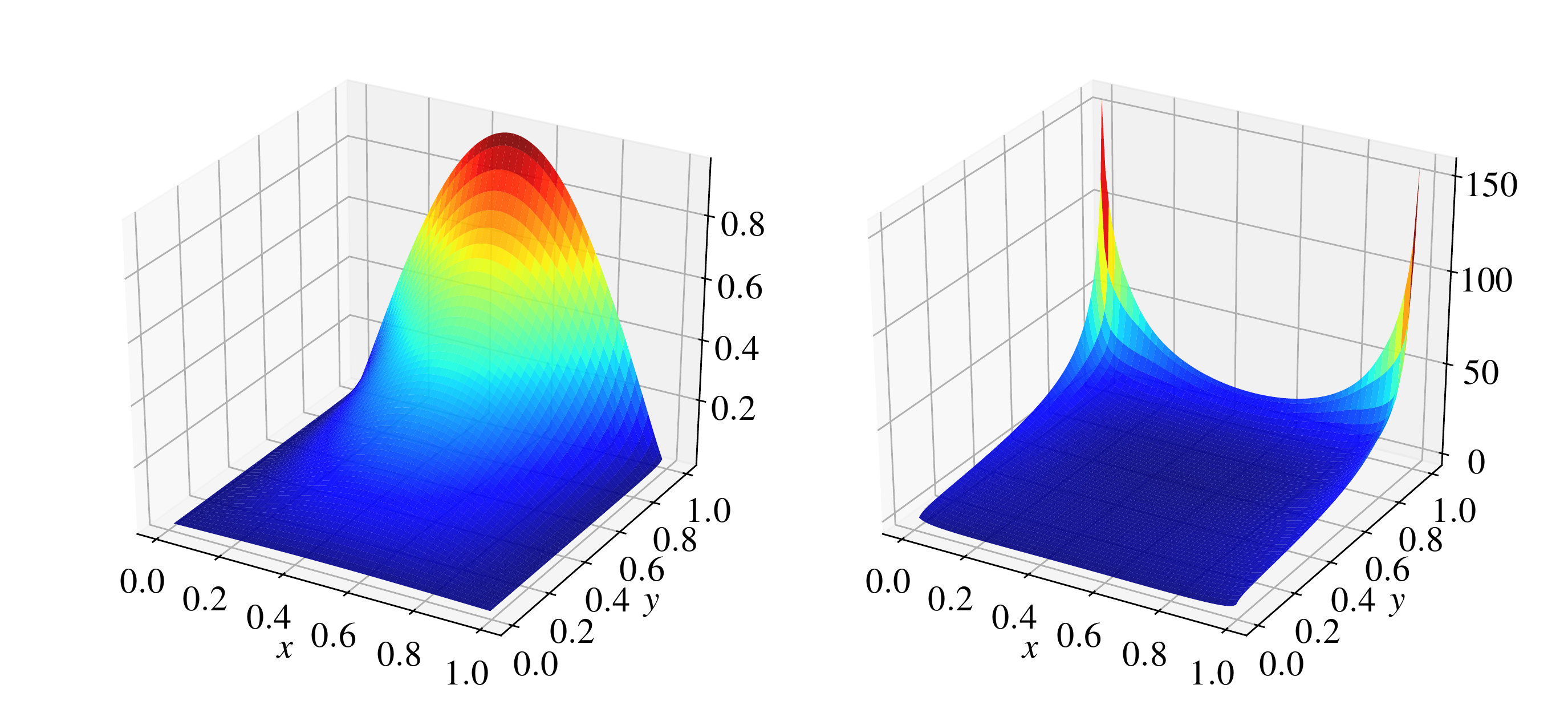}}; 
\node at (-1.2in,-1.45in) { (a) $g(\vx)$ };
\node at (1.40in,-1.45in) { (b) $\nabla^2 g(\vx)$};
\end{tikzpicture}
\caption{Plots of (a) $g(\vx)$ and (b) $\nabla^2 g(\vx)$  over the unit square.} 
\label{fig:Heat-Laplacian-g}
\end{figure}

In~\fref{fig:Heat-Laplacian-g}, the function $g(\vx)$ and 
its Laplacian are plotted over the unit square. We observe that $g(\vx)$ is zero on $\Gamma_\alpha$ ($\alpha = 1,2,3$), and it is equal to $\sin\pi x$ on $\Gamma_4$. Referring to~\eqref{eq:trial_EBC}, the
trial function for PINN is:
\begin{equation*}
\unnbcxt = g(\vx) + \phi(\vx) \, \unnRxt. 
\end{equation*}
Since $g(\vx)$ satisfies the boundary conditions 
in~\eqref{eq:Example2:EBC} and $\phi = 0$ on $\Gamma$,
$\unnbcxt$ satisfies the Dirichlet boundary conditions on all edges of the boundary. From~\fref{fig:Heat-Laplacian-g}b, we observe that the Laplacian of $g(\vx)$ is singular at two
of the four vertices on the boundary and again this is handled by performing collocation over the smaller square $[0.01,\, 0.99]^2$.

\begin{figure}[!htb]
\centering
\begin{tikzpicture}
\node at (0,0) {\includegraphics[width=0.98\textwidth]{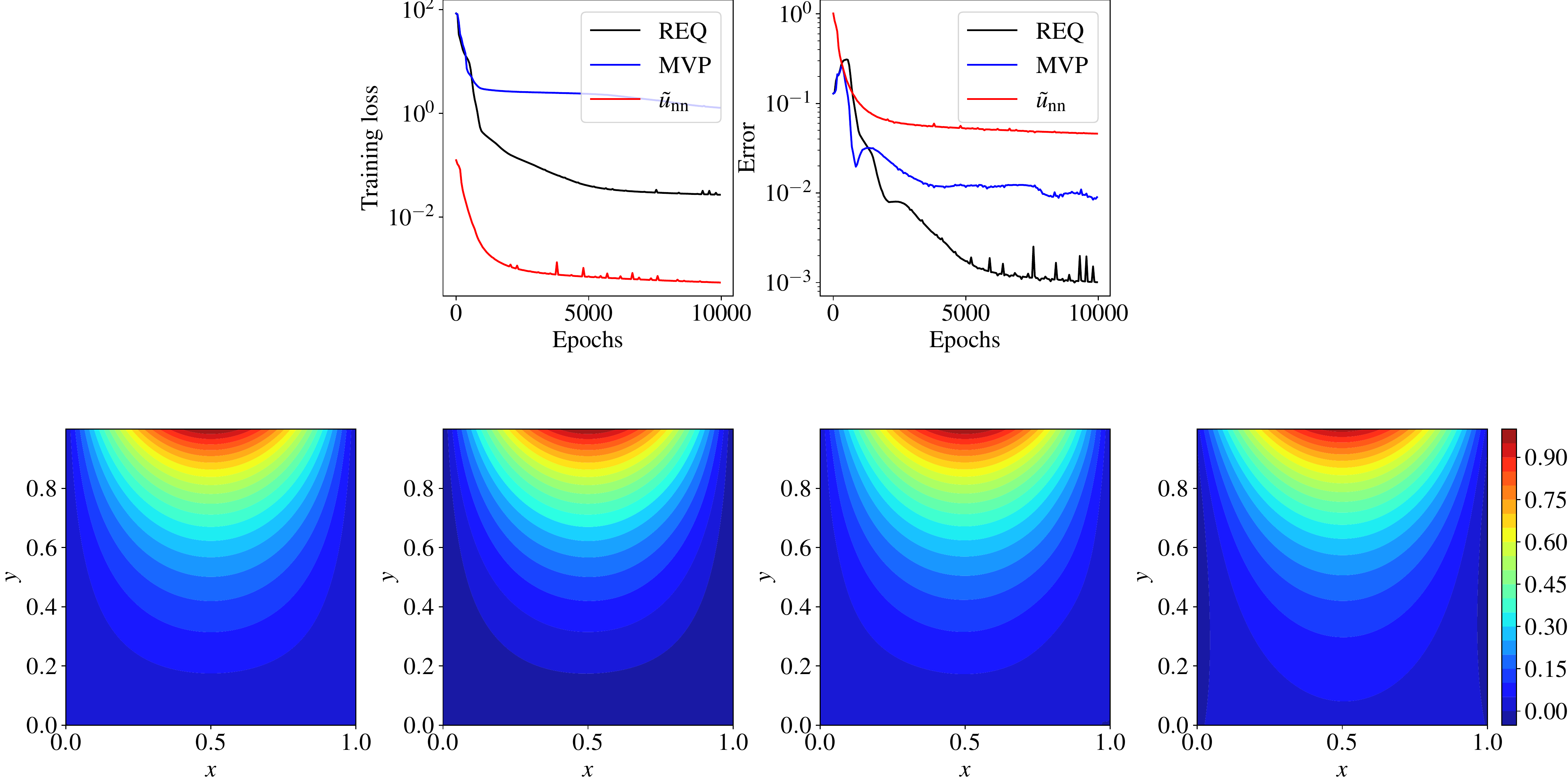}}; 
\node at (-0.80in,0.02in) {(a) };
\node at (0.75in,0.02in) {(b)};
\node at (-2.32in,-1.71in) {(c) };
\node at (-0.80in,-1.71in) {(d)};
\node at ( 0.75in,-1.71in) {(e)};
\node at (2.28in,-1.71in) {(f)};
\end{tikzpicture}
\caption{Numerical results for the Laplace problem on the
unit square with nonzero essential boundary conditions.
(a), (b) Training loss and normalized absolute errors for $\unnbcxt$ (REQ and MVP) and $\unnxt$. Contour plots over the unit square of the
(c) exact solution, 
(d) $\unnbcx$ (REQ),
(e) $\unnbcx$ (MVP), and
(f) $\unnx$.
} 
\label{fig:Heat-NonZeroDirichlet}
\end{figure}

We determine the parameters of this network by minimizing $\Lnnbct$ as described earlier (with $f$ set to zero). For standard PINN, the parameters of $\unnxt$ are determined by minimizing the following 
loss function:
\begin{equation*}
\Lnnt = ||\nabla^2\unnxt||^2_{\Omega, N_I} 
+ \sum_{\alpha=1}^{4}||\unnxt - g_\alpha(\vx) ||^2_{\Gamma_\alpha,N_B^\alpha},
\end{equation*}
where $N_B^\alpha$ is the number of collocation points on
$\Gamma_\alpha$ ($\alpha=1,2,3,4$) and $g_\alpha(\vx) =0$ ($\alpha= 1,2,3,4$) in this problem. 
In the computations, we select a total of 400 boundary collocation points
and 5,000 interior collocation within $[.01,.99]^2$. The 
numerical results are shown in~\fref{fig:Heat-NonZeroDirichlet}, where the training loss, the normalized absolute error, and exact and approximate solutions are presented for $\unnbcxt$ and $\unnxt$. As in the previous example, 
here also we notice that the training loss for $\unnxt$ is orders of magnitude smaller than both REQ and MVP 
(\fref{fig:Heat-NonZeroDirichlet}a). The losses for $\unnbcxt$ start at relatively high values during the initial stages of the training, 
which is due to the large contributions from the Laplacian in the vicinity of the vertices. However,
the network is quickly able to optimize and bring the losses down by almost two orders of magnitude in a couple of thousand training epochs. Still, the losses for $\unnbcxt$ remain several orders of magnitude larger than $\unnxt$ at the end of the training. However, this example also reveals that the 
absolute value of the loss in itself is not very meaningful. The normalized absolute errors of the three schemes are presented 
in~\fref{fig:Heat-NonZeroDirichlet}b as a function of the training epochs. It is evident from this plot that the errors in $\unnbcxt$ are orders of magnitude smaller than the error in $\unnxt$. 
The error achieved by $\unnbcxt$ (REQ) is almost an order of magnitude smaller than $\unnbcxt$ (MVP),
and almost two orders of magnitude better than $\unnxt$. 
Contour plots of the exact and approximate solutions appear in Figs.~\ref{fig:Heat-NonZeroDirichlet}c--\ref{fig:Heat-NonZeroDirichlet}f.

\begin{figure}[!htb]
\centering
\begin{tikzpicture}
\node at (0,0) {\includegraphics[width=0.98\textwidth]{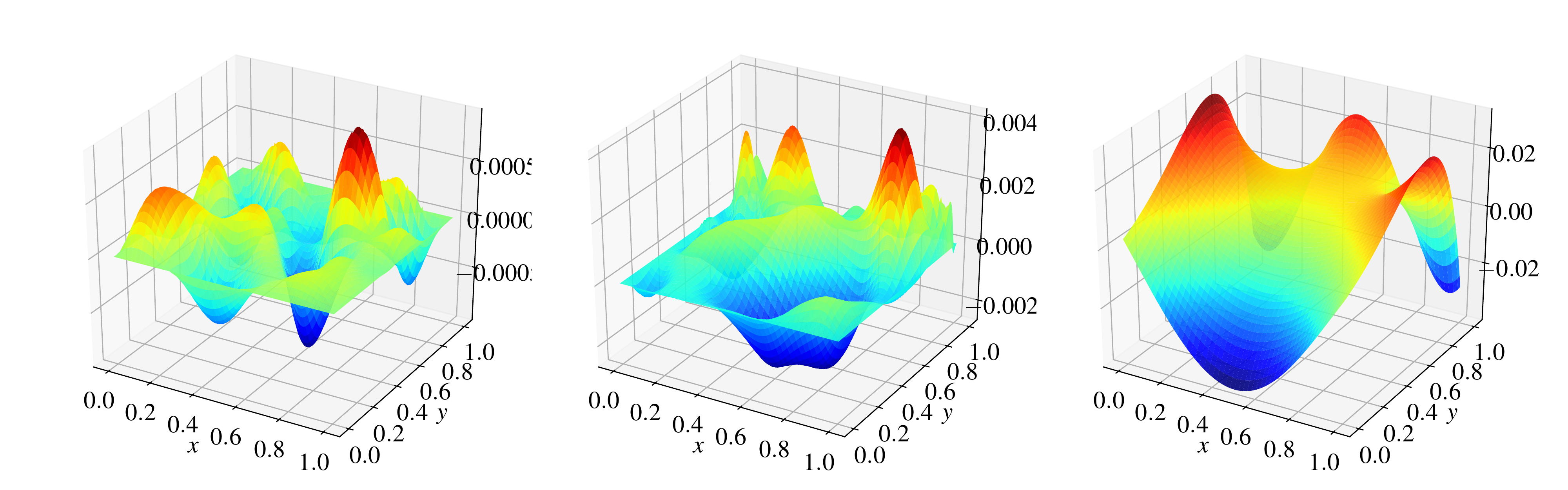}};
\node at (-2.0in,-1.15in) { (a)};
\node at (0.1in,-1.15in) { (b)};
\node at ( 2.1in,-1.15in) { (c)};
\end{tikzpicture}
\caption{Surface plots of the errors in the numerical solutions 
for the Laplace
problem with nonzero essential boundary condition.
(a) $\unnbcx$ (REQ), (b) $\unnbcx$ (MVP), and (c) $\unnx$.} 
\label{fig:Heat-NonZeroDirichlet-errors}
\end{figure}
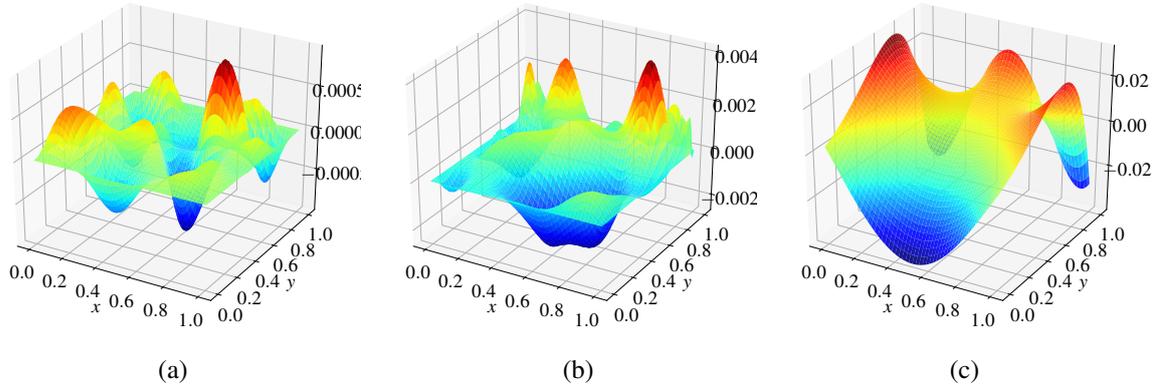

The predicted errors of the three schemes are displayed 
in~\fref{fig:Heat-NonZeroDirichlet-errors}.  It can be seen that the boundary errors in both REQ and MVP are precisely zero. This is expected, since $\unnbcxt$ has been designed to satisfy the boundary conditions. On the other hand, $\unnxt$ has large errors on the boundary of the domain. The errors from $\unnxt$ are roughly an order of magnitude larger than MVP and two orders of magnitude larger than REQ. As in the previous example, one can improve the results of $\unnxt$ by weighing the boundary loss more than the PDE loss in $\Lnnt$.

\begin{figure}[htp]
\centering
\begin{tikzpicture}
\node at (0,0) {\includegraphics[width=0.78\textwidth]{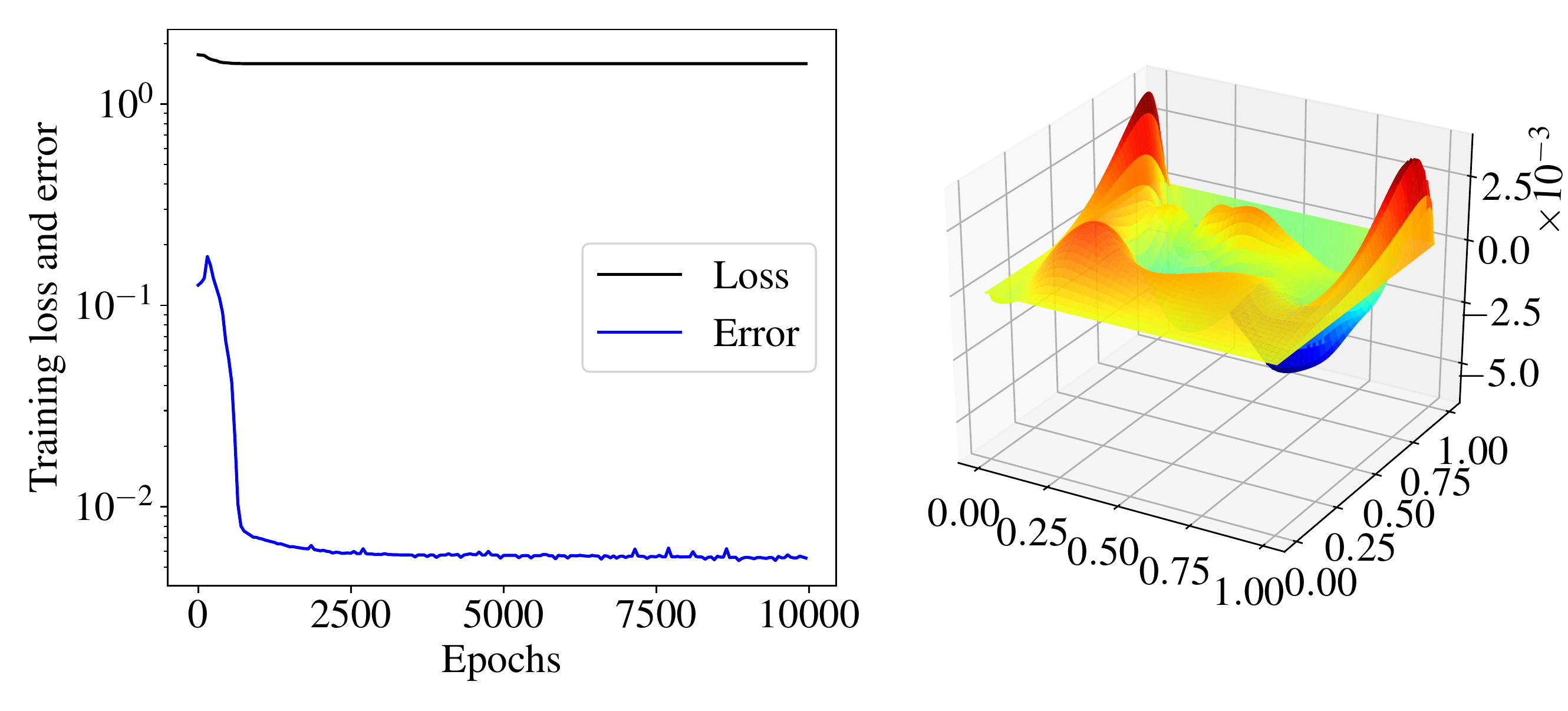}};
\node at (-0.9in,-1.2in) {(a)};
\node at ( 1.43in,-1.2in) {(b)};
\end{tikzpicture}
\caption{Ritz solution using $\unnbcxt$ (REQ) for the Laplace
problem in Example 2. (a) Evolution of training loss and 
normalized absolute error, and (b) Prediction errors over the
domain after training.} 
\label{fig:Heat-NonZeroDirichlet-Ritz}
\end{figure}

This problem can also be solved using a Ritz scheme, which is
appealing since only first-order derivatives are required in
the loss function. We form a trial function that
satisfies the essential conditions using~\eqref{eq:trial_EBC}:
\begin{equation*}
\unnbcxt = g(\vx) + \phi(\vx) \, \unnRxt.
\end{equation*}
The parameters of the network can now be found by minimizing the loss 
in~\eqref{eq:poisson_var_Lt}. To numerically evaluate the integral, we divide the square into a uniform grid with $N_I$ number of interior points. Since the Ritz loss does not involve second derivatives of $\phi(\vx)$ or $g(\vx)$, all terms in the loss are well-defined and bounded even arbitrarily close to the boundaries of the domain. We
select 5,000 interior points on the square $[0.0001,\, 0.9999]^2$.
For a Dirichlet problem, it is especially important to sample
close to the boundaries, because in the absence of doing so, the loss may be trivially minimized by a $u(\vx)$ that is a
constant. Sampling close to the boundaries informs the algorithm that a constant $u(\vx)$ leads to large errors near the 
boundaries, which are manifested in the loss term as large gradients. Numerical results for the Ritz method using 
$\unnxt$ with REQ are presented in~\fref{fig:Heat-NonZeroDirichlet-Ritz}.
The training loss and normalized absolute errors as a function of epochs 
are shown in~\fref{fig:Heat-NonZeroDirichlet-Ritz}a, which
reveal that the error reduces to ${\cal O}(10^{-2})$ in less than 2,000 epochs. 
In~\fref{fig:Heat-NonZeroDirichlet-Ritz}b, the
prediction errors over the square are displayed. Compared 
to~\fref{fig:Heat-NonZeroDirichlet-errors}, we find that 
the errors in the Ritz scheme are smaller than MVP-based collocation but larger than REQ-based collocation. 

\subsubsection{Curved domain}\label{subsubsec:heatM}
For a problem with a curved domain, we consider the Laplace
equation, $\nabla^2 u = 0$, in an annulus that is bounded between circles
of radii, $R_1 = 1$ (boundary $\Gamma_1$)
and $R_2 = 1/4$ (boundary $\Gamma_2$)~\cite{Tsukanov:2011:HME}. 
The essential boundary conditions are:
\begin{equation*}
u = 1 \ \ \textrm{on } \Gamma_1, \quad
u = 2 \ \ \textrm{on } \Gamma_2.
\end{equation*}
The exact solution for this problem is:
\begin{equation}\label{eq:exact_Laplace_annulus}
u(\vx) = 1 - \frac{\ln \sqrt{x^2+y^2} }{\ln 4}.
\end{equation}

To impose the boundary conditions, we need a distance function to both $\Gamma_1$ and $\Gamma_2$, and also a
composite boundary data function $g(\vx)$. In this case, on using
the exact distance functions to the two circles, we have
$\phi_1(\vx)= 1 - \sqrt{x^2+y^2} $ (positive in the interior of the larger disk), and
$\phi_2(\vx)= \sqrt{x^2+y^2} - 1/4$ (positive outside the smaller disk). Since the origin is not part of the computational
domain, here we can use the exact distance functions to form
$\phi_1$ and $\phi_2$. Now, on combining these two ADFs using the R-equivalence 
operation ($m = 1$) in~\eqref{eq:phi_eq}, we obtain
$\phi(\vx) = \phi_1 \sim \phi_2$ (positive in the annulus).
Finally, we use the transfinite formula~\eqref{eq:ti} to 
construct the composite boundary data as 
\begin{equation*}
    g(\vx) = \frac{ 2\phi_1+\phi_2}{\phi_1+\phi_2}
           = \frac{7 - 4 \sqrt{x^2 + y^2}}{3}.
\end{equation*}

\begin{figure}[!htb]
\centering
\begin{tikzpicture}
\node at (0,0) {\includegraphics[width=0.78\textwidth]{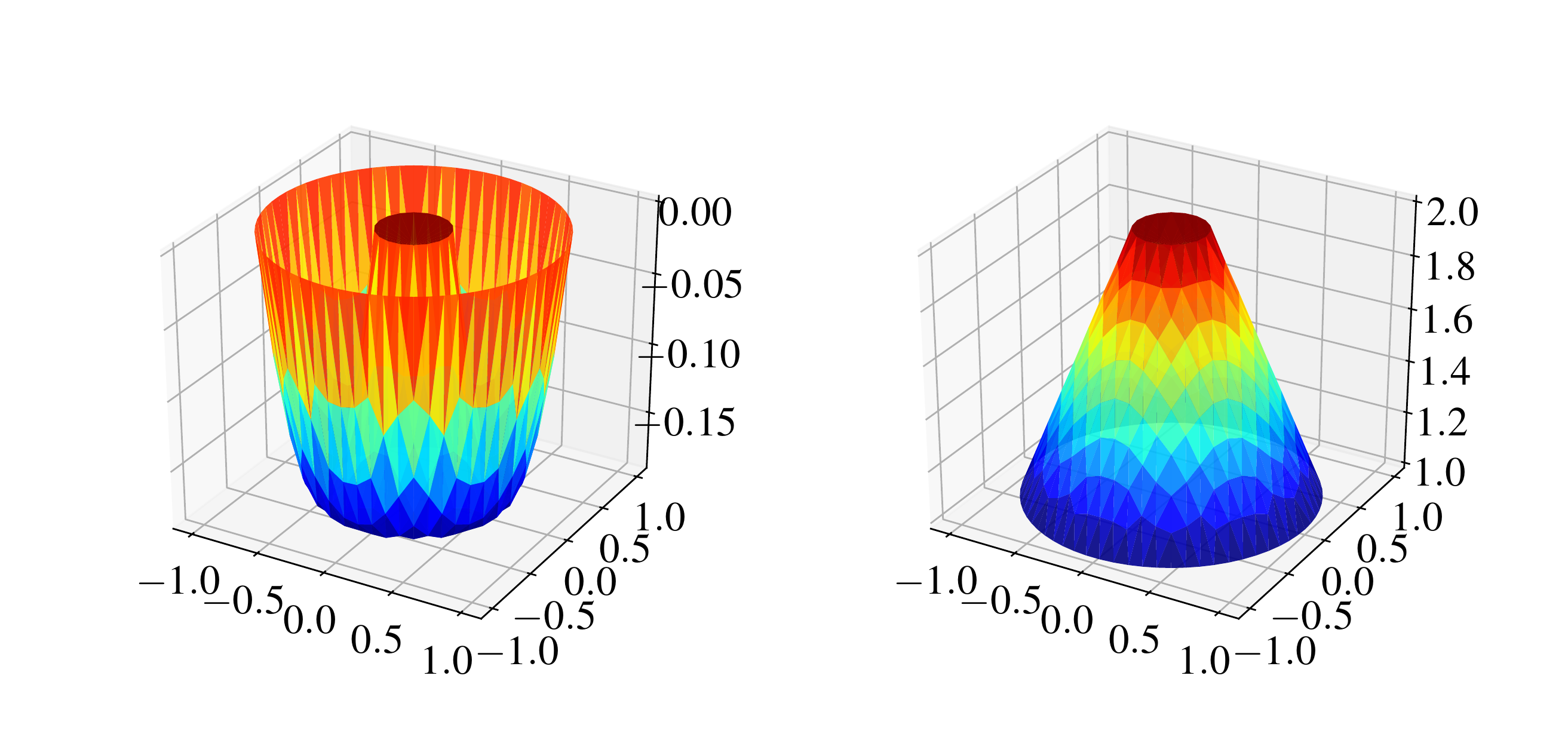}};
\node at (-1.15in,-1.22in) {(a)};
\node at ( 1.29in,-1.22in) {(b)};
\end{tikzpicture}
\caption{Approximate distance function to the boundaries of the 
annulus using R-equivalence composition of the exact distance functions
to $\Gamma_1$ and $\Gamma_2$. 
Plots of (a) $-\phi(\vx)$ and (b) $g(\vx)$ that
interpolate essential boundary data on $\Gamma_1$ and $\Gamma_2$.} \label{fig:Heat-Circular-phi-g}
\end{figure}

In~\fref{fig:Heat-Circular-phi-g}, $\phi(\vx)$ and
$g(\vx)$ are plotted over the annulus. To clearly see that  $\phi(\vx)$ is zero on $\Gamma_1$ and $\Gamma_2$, we show
$-\phi(\vx)$ in~\fref{fig:Heat-Circular-phi-g}a.  From the plot 
in~\fref{fig:Heat-Circular-phi-g}b, we observe that $g(\vx)$ matches the imposed boundary data on $\Gamma_1$
and $\Gamma_2$. An ansatz that exactly satisfies the boundary conditions is: $\unnbcxt = g(\vx) + \phi(\vx) \, \unnRxt$. 
We compare the performance of $\unnbcxt$ with a standard PINN trial function, $\unnxt$, where the boundary conditions need to be enforced through the loss function. The loss function for the two cases are similar to those discussed in previous examples, the only difference being that now we sample the boundary data at $N_B$ points on the curved boundary $\Gamma = \Gamma_1 \cup \Gamma_2$ and the interior collocation data at $N_I$ points in $\Omega_1 \cap \Omega_2$. We use \texttt{dmsh}~\cite{Schlomer:2020} to triangulate the annulus and choose the centroid of the triangles as interior collocation points and the center of the 
edges on the boundary as the boundary collocation points.  For this problem, we pick
$N_I = 612$ points in the interior of the domain,
$66$ points on $\Gamma_1$ and $30$ points on 
$\Gamma_2$ for a total of $N_B = 96$ points on 
$\Gamma_1 \cup \Gamma_2$. 
A representative mesh that display the interior collocation
points over the annulus is
shown in~\fref{fig:FEM-mesh}.
We use the \texttt{Adam} 
optimizer with a
learning rate (step size) of $10^{-3}$ for training both $\unnbcxt$ and $\unnxt$, and the training is stopped at 10,000 epochs for both networks in order to perform 
a fair comparison.

\begin{figure}[!htb]
\centering
\begin{tikzpicture}
\node at (0,0) {\includegraphics[width=0.98\textwidth]{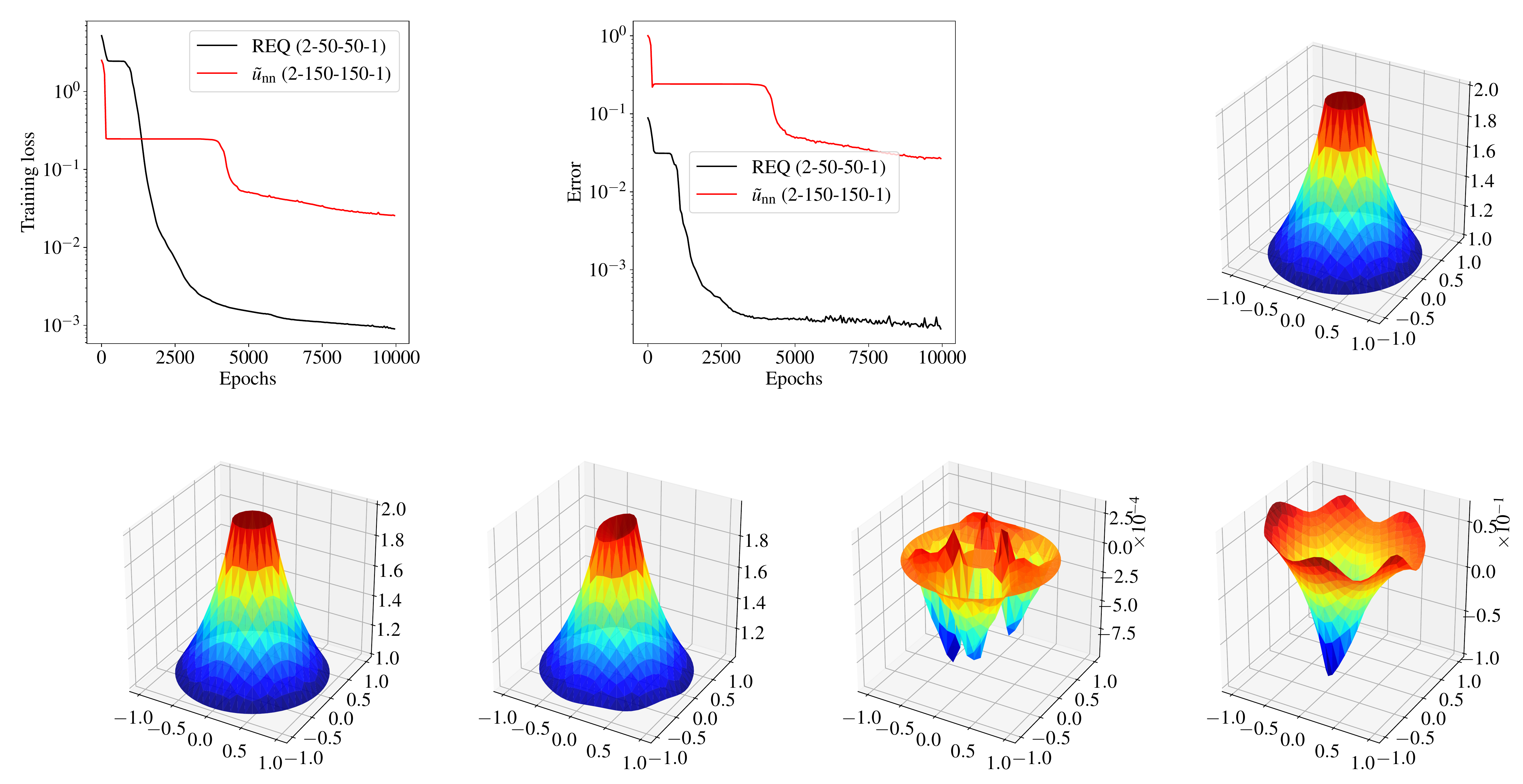}}; 
\node at (-2.11in,-0.12in) {(a)};
\node at (0.13in,-0.12in) {(b)};
\node at (2.5in,-0.12in) {(c)};
\node at (-2.11in,-1.75in) {(d) };
\node at (-0.58in,-1.75in) {(e)};
\node at ( 0.95in,-1.75in) { (f)};
\node at (2.46in,-1.75in) { (g)};
\end{tikzpicture}
\caption{Numerical solutions for the Laplace problem on
         an annulus with Dirichlet boundary conditions.
         (a), (b) Training loss and normalized absolute
         errors of $\unnbcxt$ (REQ) and $\unnxt$. Surface plots of the (c) exact solution,
         (d) $\unnbcxt$ using 2--50--50--1 network, and (e) $\unnx$ using 2--150--150--1 network. Surface plots of the error for the
         numerical solutions
         (f) $\unnbcx$ and (g) $\unnx$.
         } 
\label{fig:Heat-Circular}
\end{figure}

In~\fref{fig:Heat-Circular}, the results of the training as well as the approximate solutions produced by $\unnbcxt$ (REQ) and
$\unnxt$ are 
presented. It is pertinent to mention here that $\unnxt$ required a much larger network architecture (2--150--150--1) compared to $\unnbcxt$, which only required a 2--50--50--1 network in order to converge to acceptable results. However, even with a much larger network, as observable in Figs.~\ref{fig:Heat-Circular}a 
and~\ref{fig:Heat-Circular}b, the error in $\unnxt$ by the end of the training 
is two orders of magnitude larger than $\unnbcxt$. It can be seen 
from~\fref{fig:Heat-Circular}e that $\unnxt$ does not satisfy the boundary conditions. The 
approximation errors of $\unnbcx$ (REQ) and $\unnx$ appear
in Figs.~\ref{fig:Heat-Circular}f and~\ref{fig:Heat-Circular}g. The errors in $\unnx$ are especially large on the boundaries of the domain.

We now use the Laplace problem over the annulus to also demonstrate how to solve the problem using 
mixed boundary conditions. To 
this end, we retain the essential boundary condition on
$\Gamma_1$ and convert the boundary condition on 
$\Gamma_2$ to a Robin boundary condition. 
For this problem, we find that $\partial u/ \partial n = 4 /\ln 4$ on the
inner boundary. 
We use the following mixed boundary conditions:
\begin{equation}\label{eq:annulus_mixed_BCs}
u = 1 \ \ \textrm{on } \Gamma_1, \quad
\frac{\partial u}{\partial n} + u = 2 + \frac{4}{\ln 4} =: h\ \ \textrm{on }
\Gamma_2.  
\end{equation}

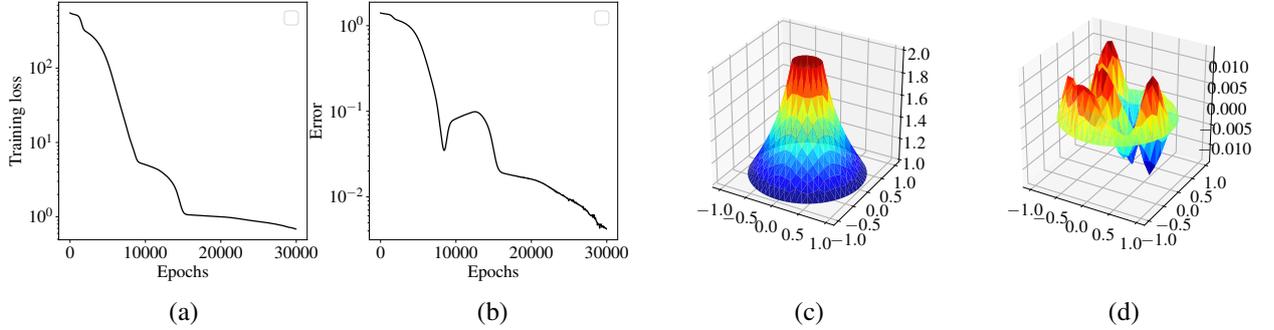
\begin{figure}[!htb]
\centering
\begin{tikzpicture}
\node at (0,0) {\includegraphics[width=\textwidth]{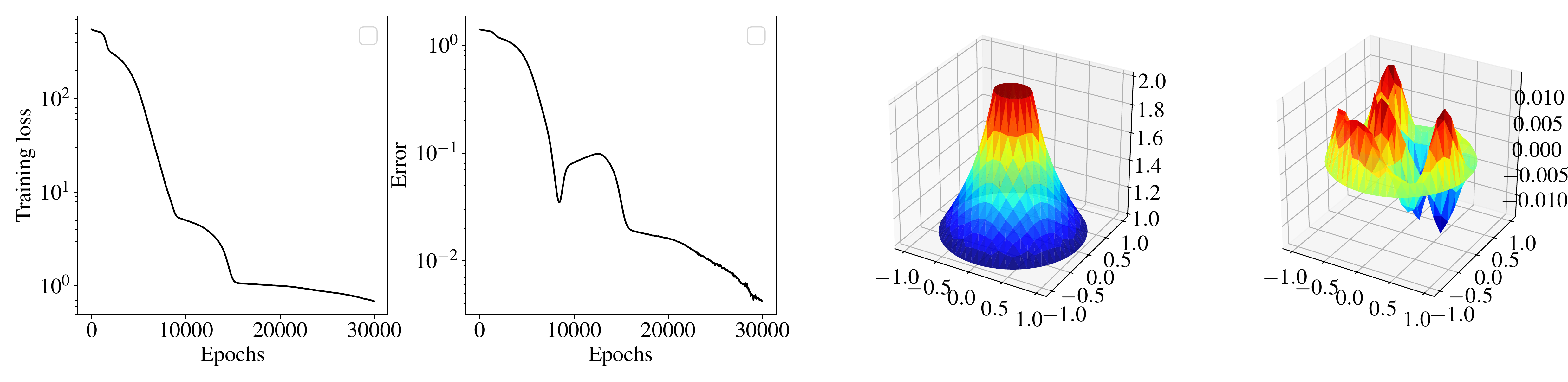}}; 
\node at (-2.28in,-0.9in) {(a)};
\node at (-0.68in,-0.9in) {(b)};
\node at ( 0.96in,-0.9in) {(c)};
\node at (2.59in,-0.9in) { (d)};
\end{tikzpicture}
\caption{Numerical solution using $\unnbcxt$ (REQ) for the
Laplace problem in an annulus with mixed boundary conditions. 
(a), (b) Training loss and normalized absolute errors as a function of
epochs. Surface plots of (c) $\unnbcx$ and (d) error in
$\unnbcx$.}
\label{fig:Heat-Circular-Mixed}
\end{figure}

The exact solution remains unchanged and is given 
in~\eqref{eq:exact_Laplace_annulus}. To create an ansatz for this problem with mixed boundary conditions, we follow the formulation in~\sref{subsec:EBCRobin}.  We form
$\phi_1$, $\phi_2$, which remain unchanged from the previous
case when essential boundary conditions are imposed 
on $\Gamma_1$ and $\Gamma_2$. Since $g = 1$, referring to~\eqref{eq:trial_mixedBC-II}, we can write
\begin{equation}
u_1(\vx) = g(\vx) = 1, \quad
u_2(\vx) = \left[ 1 + \phi_2 (1 + D_1^{\phi_2} ) \right] ( \unnRxt ) - \phi_2 h ,
\end{equation}
where $h$ is given in~\eqref{eq:annulus_mixed_BCs},
and then on using transfinite interpolation given 
in~\eqref{eq:trial_mixedBC-II-a} and~\eqref{eq:trial_mixedBC-II-b}, we form the trial function $\unnbcxt$. 
So we now have
an ansatz that satisfies both the
Dirichlet and Robin boundary conditions. The training loss
and normalized absolute errors for the numerical solution using $\unnbcxt$ are presented in Figs~\ref{fig:Heat-Circular-Mixed}a 
and~\ref{fig:Heat-Circular-Mixed}b. 
On comparing $\unnbcx$ in~\fref{fig:Heat-Circular-Mixed}c to the exact solution 
shown in~\fref{fig:Heat-Circular}, we see that the boundary conditions have been exactly satisfied. The error plot 
in~\fref{fig:Heat-Circular-Mixed}d
reveals that the numerical solution is within 
1 percent of the exact solution over the domain.

\subsection{Generalized barycentric coordinates over polygons}\label{subsec:GBCs}
Consider a planar polygon with $n$ vertices (nodal coordinate
$\{\vx_i\}_{i=1}^n$) that are in counterclockwise orientation. On an $n$-gon, harmonic coordinates~\cite{Joshi:2007:HCF} are
one of the instances of generalized barycentric coordinates~\cite{Anisimov:2017:BCP}. 
Each coordinate (shape function), 
$\varphi_i := \varphi_i(\vx)$, is  associated with vertex $i$ and is obtained by solving the Laplace equation with piecewise affine Dirichlet boundary conditions.  The boundary-value problem for harmonic coordinates is: find $\phi_i \ge 0$ $(i=1,2,\dots,n)$ that solves
\begin{subequations}
\begin{align}
\nabla^2 \varphi_i &= 0 \ \ \textrm{in } \Omega, \\
\varphi_i &= g_i \ \ \textrm{on } \partial \Omega,
\end{align}
\end{subequations}
where $g_i := g_i(\vx)$ is a piecewise affine (hat)
function 
that is unity at $\vx_i$ and is zero at all other vertices, i.e., $g_i(\vx_j) = \delta_{ij}$, where $\delta_{ij}$
is the Kronecker-delta. By virtue of the maximum principle for the Laplace equations, $\phi_i > 0$ in the interior of the polygon.

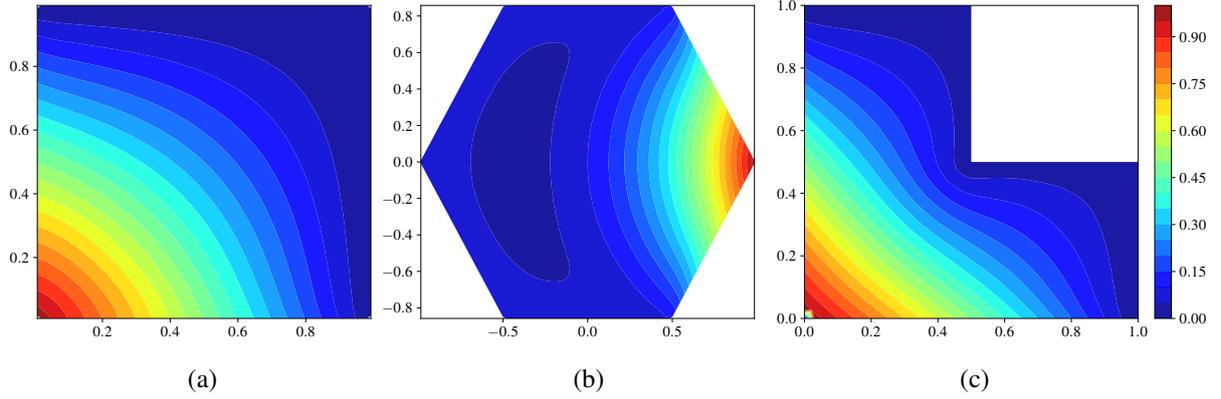
\begin{figure}[!htb]
\centering
\begin{tikzpicture}
\node at (0,0) {\includegraphics[width=0.98\textwidth]{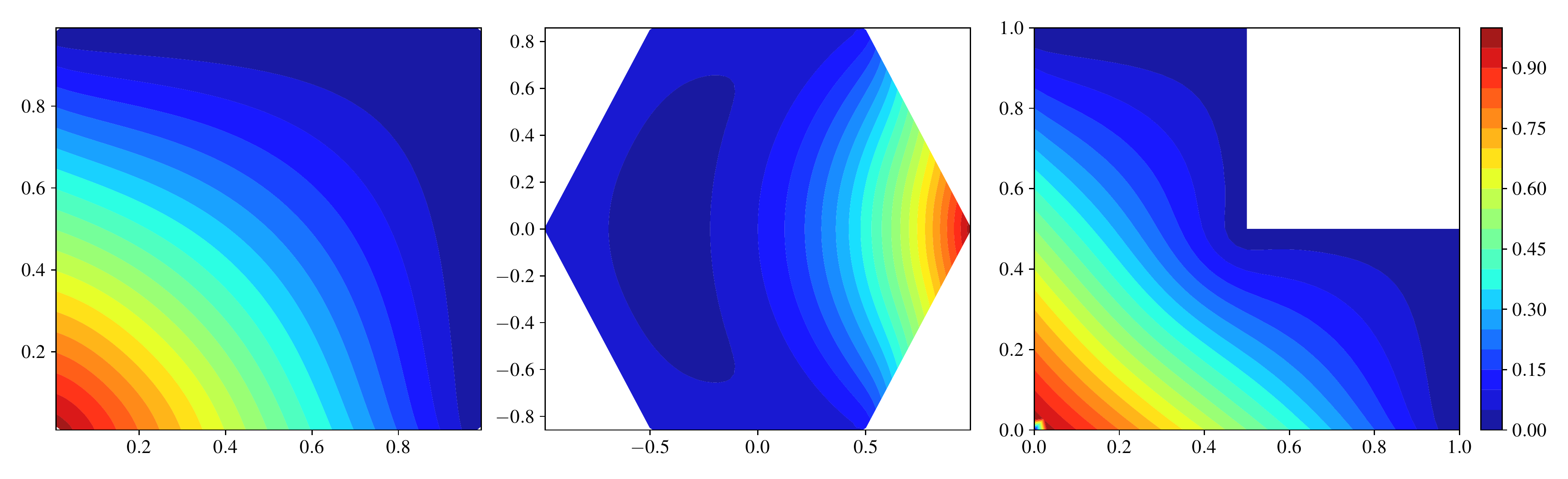}};
\node at (-2.1in,-1.1in) { (a)};
\node at (-0.1in,-1.1in) { (b)};
\node at ( 1.9in,-1.1in) { (c)};
\end{tikzpicture}
\caption{Computation of harmonic coordinate for a vertex in a polygon. Contour plots of $g(\vx)$ over a
(a) square (vertex at the origin), (b) 
regular hexagon (rightmost vertex), and (c) L-shaped polygon 
(vertex at the origin).
}
\label{fig:Barycentric-g}
\end{figure}

Here we will solve the harmonic coordinate problem on two representative polygons: a square and an L-shaped (nonconvex)
polygon. For both examples, the network architecture 2--50--50--1 is used.
In both examples, an important step is to assemble the boundary data into a function $g(\vx)$ through transfinite interpolation. Figure~\ref{fig:Barycentric-g} shows the function $g(\vx)$ for specific choices of the vertex $i$. If the vertices of the square are numbered 1--2--3--4 (counterclockwise), starting
at vertex 1 that is at $(0,0)$, then the harmonic coordinate $u(\vx)$ must satisfy $u(0,0) = 1$, 
$u(1,0) = u(1,1) = u(0,1) = 0$. So the
boundary conditions are affine along edges 1--2 and 4--1 and
identically zero along edges 2--3 and 3--4. All
these boundary conditions are 
simultaneously captured in the $g(\vx)$ function. 
It can be seen from the colormaps that $g(\vx)$ appropriately interpolates the boundary data in all cases.
\subsubsection{Harmonic coordinates on a square}
We first consider the case of computing the harmonic coordinate over a square. On a square, there exists an exact solution for this problem---harmonic coordinates coincide with bilinear
finite element shape functions. So the solution that is
associated with vertex 1 is: $u(x,y) = (1-x)(1-y)$. 
To compute the harmonic coordinates, we adopt the Ritz method
to determine the approximate solutions. Since harmonic coordinates minimize the Dirichlet energy, use of the Ritz formulation is
natural. As done earlier, we consider the approximations $\unnbcxt$ and $\unnxt$. For $\unnbcxt$, the loss to be minimized is given
in~\eqref{eq:poisson_var_Lt},
whereas for $\unnxt$, it is supplemented with an additional term $\sum_{\alpha=1}^{n}||\unnxt - g_\alpha(\vx) ||^2_{\partial\Omega_\alpha,N_B^\alpha}$ to impose the essential
boundary conditions. Here, $n$ is the number of boundary segments over which the Dirichlet boundary conditions are specified.

\begin{figure}[!htb]
\centering
\begin{tikzpicture}
\node at (0,0) {\includegraphics[width=0.98\textwidth]{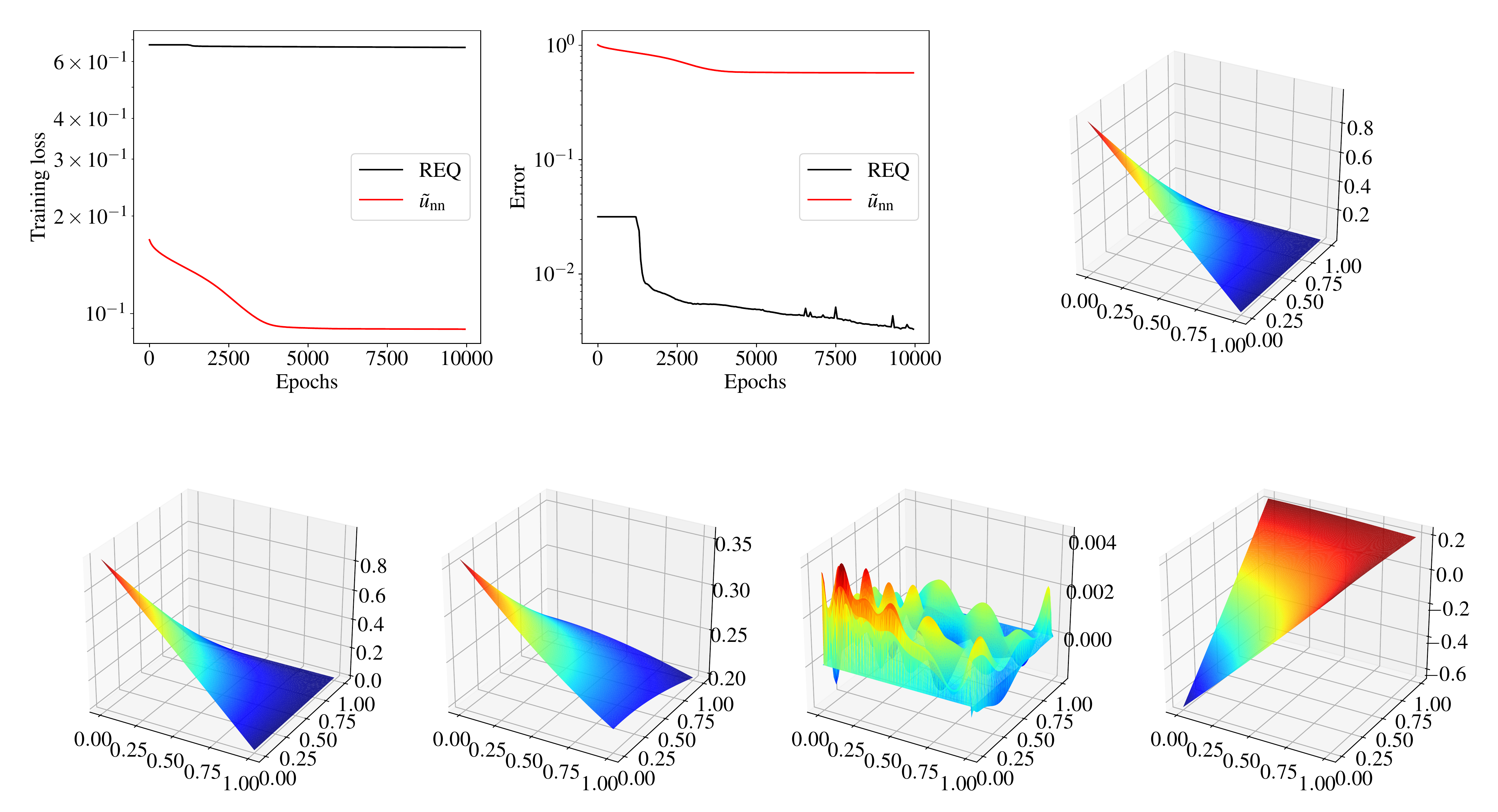}}; 
\node at (-1.89in,-0.07in) {(a)};
\node at (0.04in,-0.07in) {(b)};
\node at (2in,-0.07in) {(c)};
\node at (-2.23in,-1.87in) {(d) };
\node at (-0.68in,-1.87in) {(e)};
\node at ( 0.85in,-1.87in) {(f)};
\node at (2.38in,-1.87in) {(g)};
\end{tikzpicture}
\caption{Computation of harmonic coordinates on the unit square
for the vertex at the origin.
(a), (b) 
Training loss and errors for $\unnbcxt$ (REQ) and
$\unnxt$. Surface plots of the (c) exact solution,
(d) $\unnbcx$, and (e) $\unnx$. Surface plot of the error for the
numerical solutions
(f) $\unnbcx$ and (g) $\unnx$.}
\label{fig:barycentric-square}
\end{figure}

In~\fref{fig:barycentric-square}, numerical results using
$\unnbcxt$ (REQ) and $\unnxt$ are presented. In 
Figs.~\ref{fig:barycentric-square}a and~\ref{fig:barycentric-square}b, the training loss and normalized absolute error as
a function of epochs are presented. The exact solution along
with the numerical solutions are displayed
in Figs.~\ref{fig:barycentric-square}c--\ref{fig:barycentric-square}e. Surface plots of the errors in $\unnbcx$ and $\unnx$
appear in Figs.~\ref{fig:barycentric-square}f and~\ref{fig:barycentric-square}g. 
Once again, and
consistent with prior findings, $\unnbcxt$ exactly satisfies the 
essential boundary conditions and has far smaller errors than $\unnxt$. The large errors in $\unnxt$ are particularly 
noticeable in~\fref{fig:barycentric-square}g. The issue again has to do with the relative scaling of the PDE loss versus the boundary loss. As done
earlier, the numerical results from $\unnxt$ can be improved 
by considering a loss of the form:
\begin{equation*}
\Lnnt = w  \left\{  \frac{1}{N_I} \sum_{k=1}^{N_I}
\left[ 
\frac{1}{2} \, a \, \Bigl( \unn (\vx_k; \vm{\theta}) , 
\unn (\vx_k; \vm{\theta}) \Bigr) 
- \ell \, \Bigl( \unn (\vx_k ; \vm{\theta} ) \Bigr)
\right]  \right\}
+ (1 - w)  \left\{ 
\sum_{\alpha=1}^{n}||\unnxt - g_\alpha(\vx) ||^2_{\partial\Omega_\alpha,N_B^\alpha}\right\} ,
\end{equation*}
and then tuning the weight $w \in [0,1]$. This tunes the relative importance of the boundary loss term with respect to the PDE loss term. In fact, if one considers $w = 10^{-3}$, $\unnxt$ achieves errors that are comparable to $\unnbcxt$. In other words, the boundary loss in $\Lnnt$ has to be weighed a thousand times more than the PDE loss in order for $\unnxt$ to produce results comparable to 
$\unnbcxt$.

\subsubsection{Harmonic coordinates on an L-shaped polygon}
We repeat the computations for the square over an L-shaped polygon.  Here we consider $\phi(\vx)$ that is formed using 
REQ and MVP. The plot of $\phi(\vx)$ for REQ and MVP over the
L-shaped polygon are shown in Figs.~\ref{fig:phi_polygons-d}
and~\ref{fig:mvc_ADF-b}. The function
$g(\vx)$ over the L-shaped polygon is shown in~\fref{fig:Barycentric-g}c. We only present numerical
simulation results for REQ and MVP. As shown for the case
of the square, the loss terms have to be weighed judiciously to obtain acceptable accuracy for $\unnx$.  Harmonic coordinates
associated with the vertex at the origin and for the vertex at the reentrant corner are computed.  Since an exact analytical solution for this problem is not available, we compute an accurate finite
element solution that we use as the reference solution. This finite
element
solution is used to compute the errors in $\unnbcx$. A Delaunay triangular mesh is created using the mesh generation package \texttt{Triangle}~\cite{Shewchuk:1996:TRI}: mesh has 13,952 elements with very small elements in the vicinity of the reentrant corner and larger elements near other vertices. The mesh size, $h = 10^{-3}$, is used near the reentrant corner to
capture the weakly singular behavior of the Laplace equation at the reentrant corner.

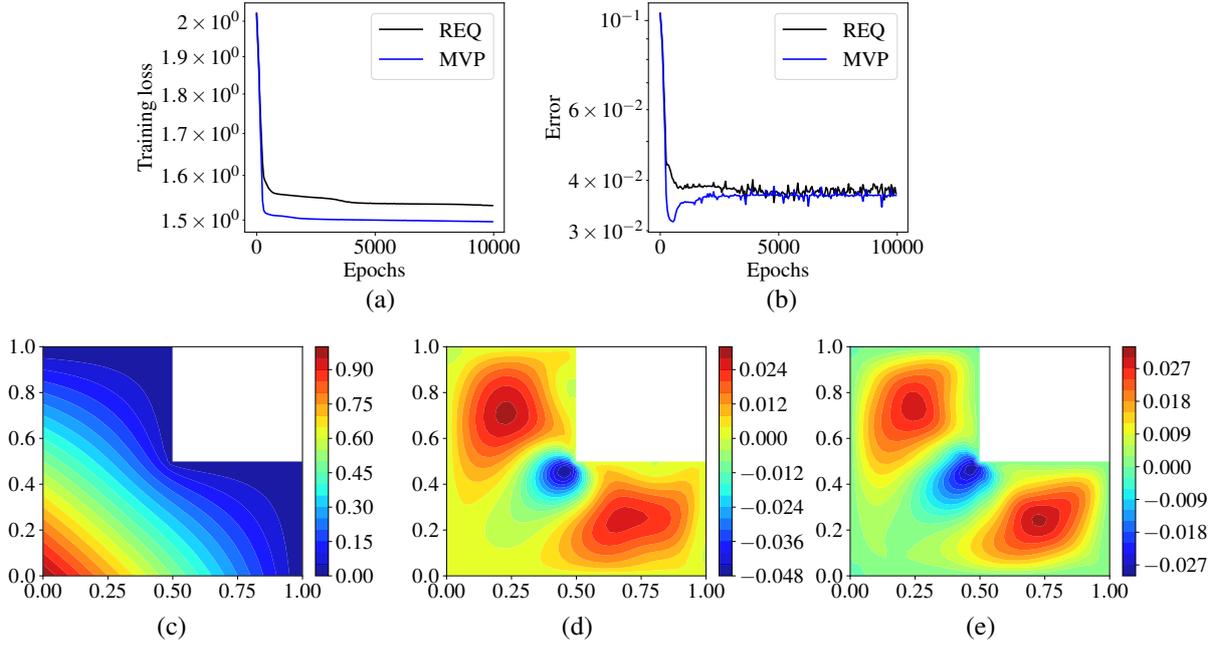
\begin{figure}[!htb]
\centering
\begin{tikzpicture}
\node at (0,0) {\includegraphics[width=0.98\textwidth]{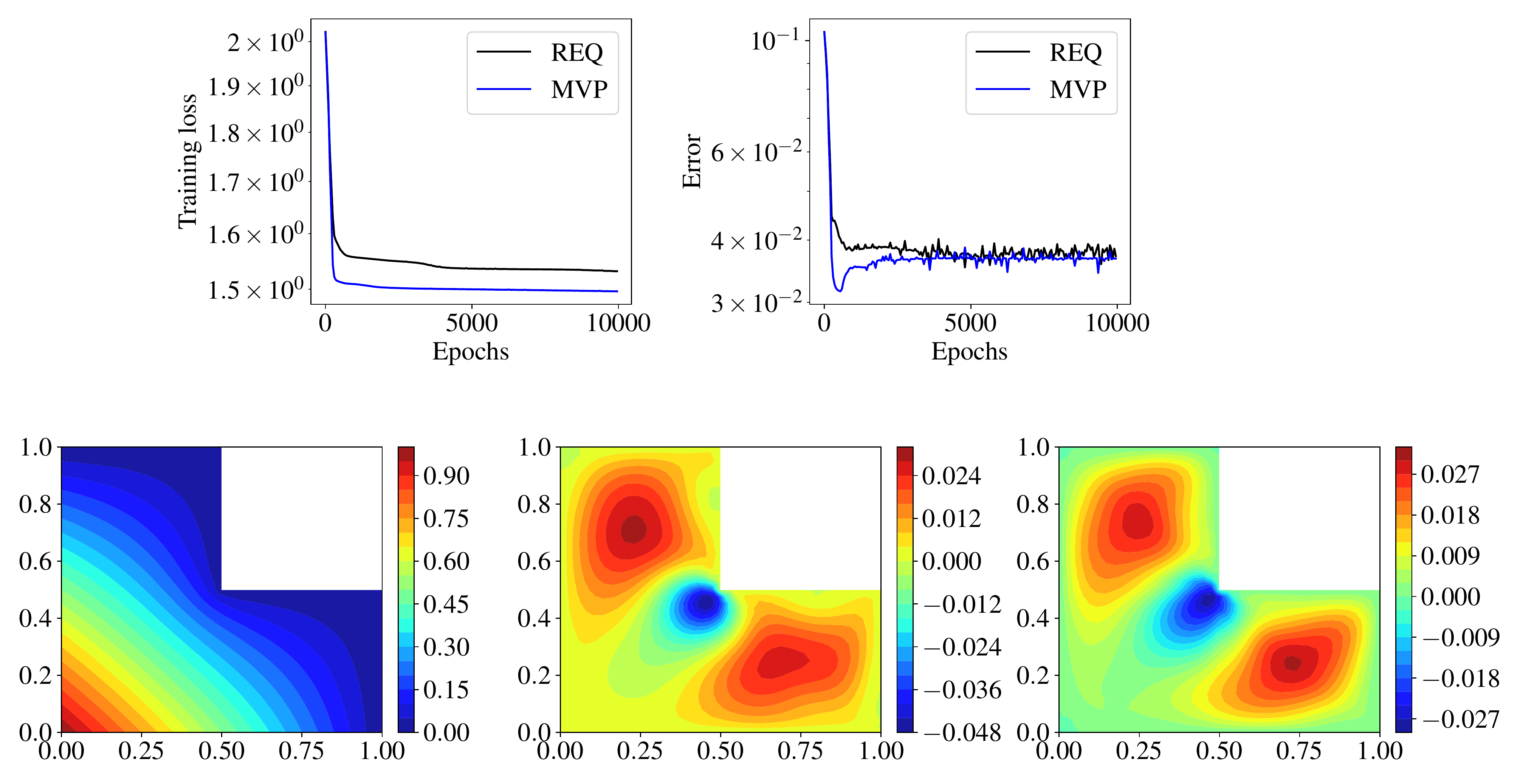}}; 
\node at (-1.18in,0.01in) {(a) };
\node at (0.90in,0.01in) {(b)};
\node at (-2.26in,-1.7in) { (c)};
\node at (-0.16in,-1.7in) { (d)};
\node at ( 1.93in,-1.7in) {(e)};
\end{tikzpicture}
\caption{Computation of harmonic coordinates (vertex
         at $(0,0)$) on an L-shaped
        polygon. (a), (b) Training loss and
        normalized absolute error for $\unnbcxt$ (REQ and MVP).
        (c) Reference finite element solution.
        Contour plots
         of the error for 
         (d) $\unnbcx$ (REQ) and
         (e) $\unnbcx$ (MVP).
} 
\label{fig:barycentric-L1}
\end{figure}

\begin{figure}[!htb]
\centering
\begin{tikzpicture}
\node at (0,0) {\includegraphics[width=0.98\textwidth]{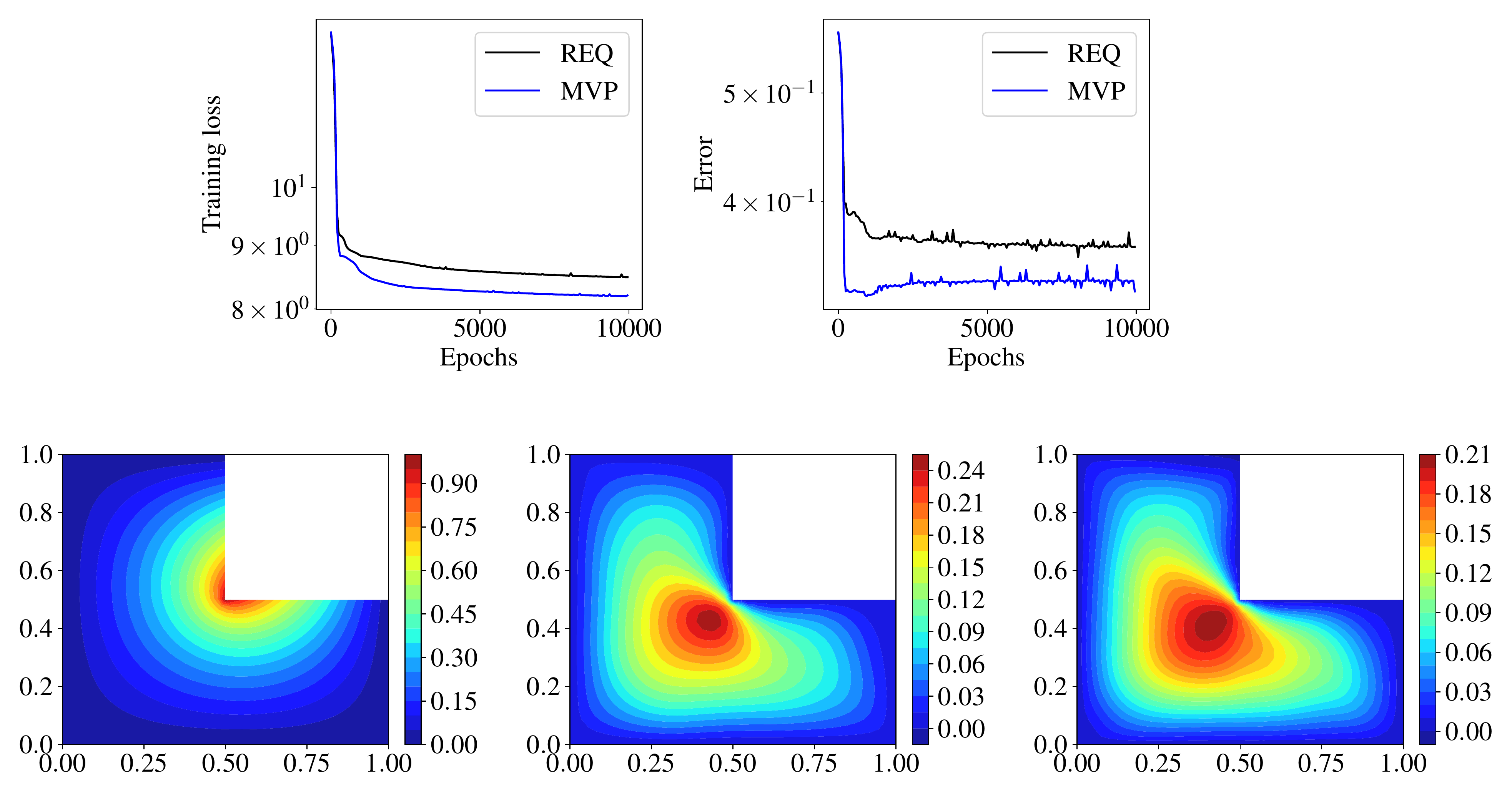}}; 
\node at (-1.14in,0.01in) {(a) };
\node at (0.98in,0.01in) {(b)};
\node at (-2.22in,-1.72in) { (c)};
\node at (-0.09in,-1.72in) { (d)};
\node at ( 2.04in,-1.72in) {(e)};
\end{tikzpicture}
\caption{Computation of harmonic coordinates (vertex
         at $(1/2,1/2)$) on an L-shaped
        polygon. (a), (b) Training loss and
        normalized absolute error for $\unnbcxt$ (REQ and MVP).
        (c) Reference finite element solution.
        Contour plots of the error for 
        (d) $\unnbcx$ (REQ) and
        (e) $\unnbcx$ (MVP).
} 
\label{fig:barycentric-L2}
\end{figure}
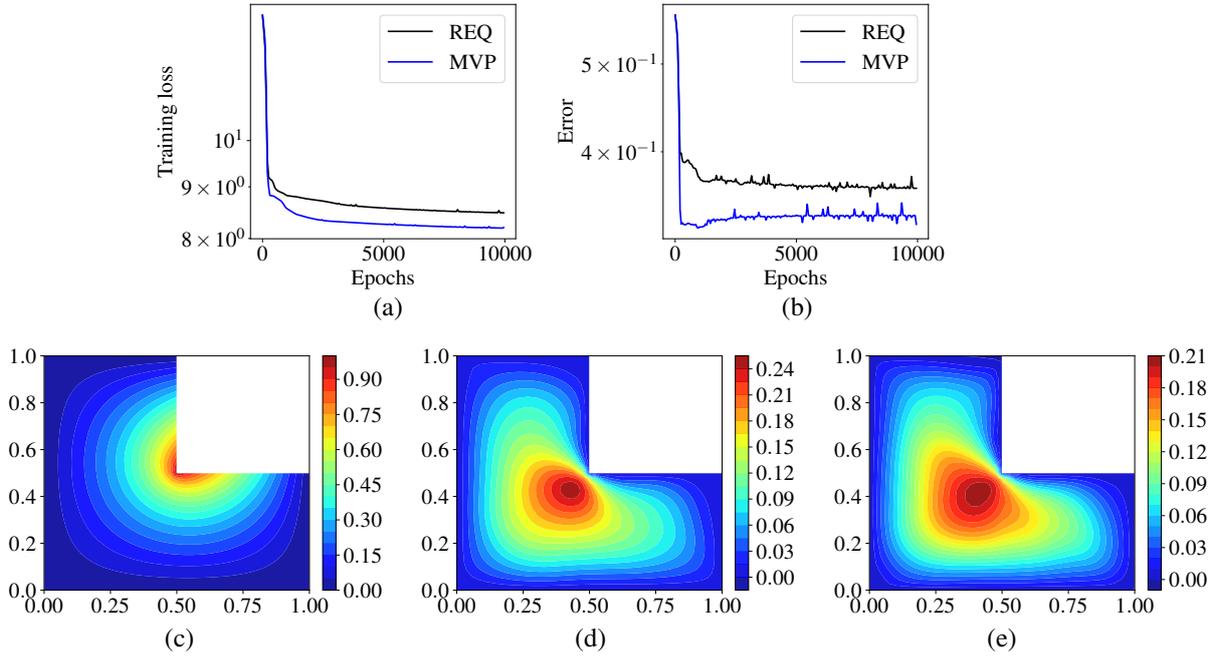

A representative mesh that displays the interior collocation
points in the L-shaped polygon
is shown in~\fref{fig:FEM-mesh}. 
In Figs.~\ref{fig:barycentric-L1}
and~\ref{fig:barycentric-L2},  numerical solutions
obtained from $\unnbcxt$ (REQ and MVP) are presented for the
computation of $u$ that is associated with vertices at
$(0,0)$ and $(1/2,1/2)$, respectively.
In~\fref{fig:barycentric-L1}a, we observe that
the training
loss stabilizes within a few thousand epochs to about 1.5 for
REQ and MVP, and this correspond to a normalized absolute
error (see~\fref{fig:barycentric-L1}b) of ${\cal O}(10^{-2})$. The reference finite element solution is presented
in~\fref{fig:barycentric-L1}c, and the error in
$\unnbcx$ (REQ and MVP) are displayed in Figs.~\ref{fig:barycentric-L1}d and~\ref{fig:barycentric-L1}e.
Numerical solutions using 
REQ and MVP have maximum errors of about 3 percent.
From~\fref{fig:barycentric-L2}a, we observe that
the training
loss stabilizes within a few thousand epochs to about 8 for
REQ and MVP, and this correspond to a normalized absolute
error (see~\fref{fig:barycentric-L2}b) of ${\cal O}(10^{-1})$. The reference finite element solution is depicted 
in~\fref{fig:barycentric-L2}c, which display sharp gradients
near the rentrant corner. The error in
$\unnbcx$ (REQ and MVP) are shown in
Figs.~\ref{fig:barycentric-L2}d and~\ref{fig:barycentric-L2}e,
with maximum errors near the singularity on the order of
20 percent. Compared to the errors in Figs.~\ref{fig:barycentric-L1}d
and~\ref{fig:barycentric-L1}e,
this is a 10-fold increase in the maximum error (due to the presence of the derivative singularity at the reentrant corner).
\subsection{Clamped circular Kirchhoff plate}\label{subsec:plate}
We consider the boundary-value problem for a clamped plate that is given in~\eqref{eq:plate_bvp}. For a clamped
circular plate of unit radius and transverse load $f = 1$,
the boundary-value problem is: 
\begin{subequations}\label{eq:Kirchhoff_bvp}
\begin{align}
\label{eq:Kirchhoff_bvp-a}
\nabla^4 u &= 1 \ \ \textrm{in } 
\Omega = \{ (x,y): x^2 + y^2 < 1 \},
\\
\label{eq:Kirchhoff_bvp-b}
u &= 0 \ \ \textrm{on } \partial \Omega,  \quad 
u,n := \frac{\partial u}{\partial n} = 
0 \ \ \textrm{on } 
\partial \Omega .
\end{align}
\end{subequations}
The exact solution for this problem in
polar coordinates is~\cite{Timoshenko:1959:TPS}:
\begin{equation}\label{eq:exact_axisymmplate}
   u(r) = \frac{(1 - r^2)^2}{64}.
\end{equation}
Given that both $u$ and $\partial u/\partial n$ are specified on
the boundary, this problem illustrates the use of a different solution structure than the ones considered until now. To impose the boundary conditions, we first create a distance function $\phi$ to the circular boundary. In this case, the distance function can be exactly determined as $\phi(\vx)= 1 - \sqrt{x^2+y^2}$; however, this exact distance function has derivative singularities at the origin. Therefore, 
we use an approximate distance function to a unit circle that is 
given in~\eqref{eq:phi_circle}, which is reproduced below:
\begin{equation*}
    \phi(\vx) = \frac{1 - \vx \cdot \vx}{2} ,
\end{equation*}
which is a bivariate polynomial. Now we construct
an ansatz that satisfies both essential boundary conditions 
using~\eqref{eq:trial_plate}: $\unnbcxt=\phi^2\tilde{u}(\vx;\vm{\theta})$. As with previous problems, this problem can be solved using either the collocation approach (see Guo et al.~\cite{Guo:2021:DCM}) or the Ritz method.
For collocation, we minimize the loss function
\begin{equation*}
\Lnnbct = || \nabla^4\unnbcxt - 1 ||_{\Omega,N_I}^2, 
\end{equation*}
whereas for the Ritz approach, we minimize the loss function
that is presented in~\sref{subsubsec:4th_order}.

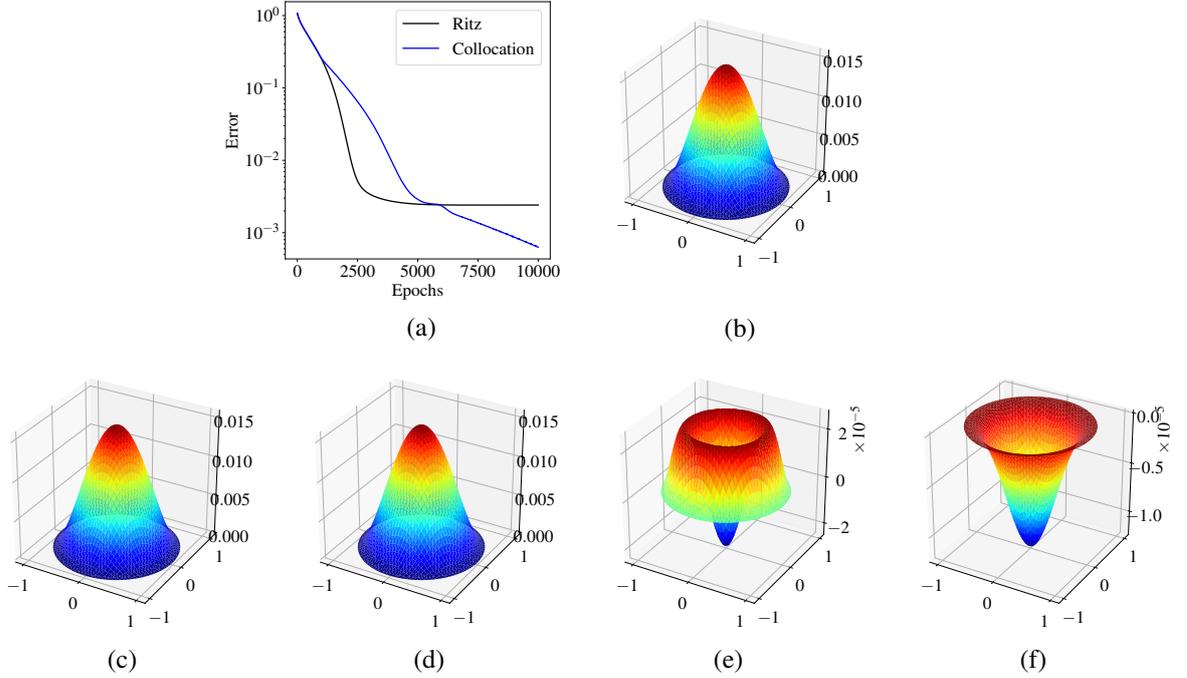
\begin{figure}[!htb]
\centering
\begin{tikzpicture}
\node at (0,0) {\includegraphics[width=0.98\textwidth]{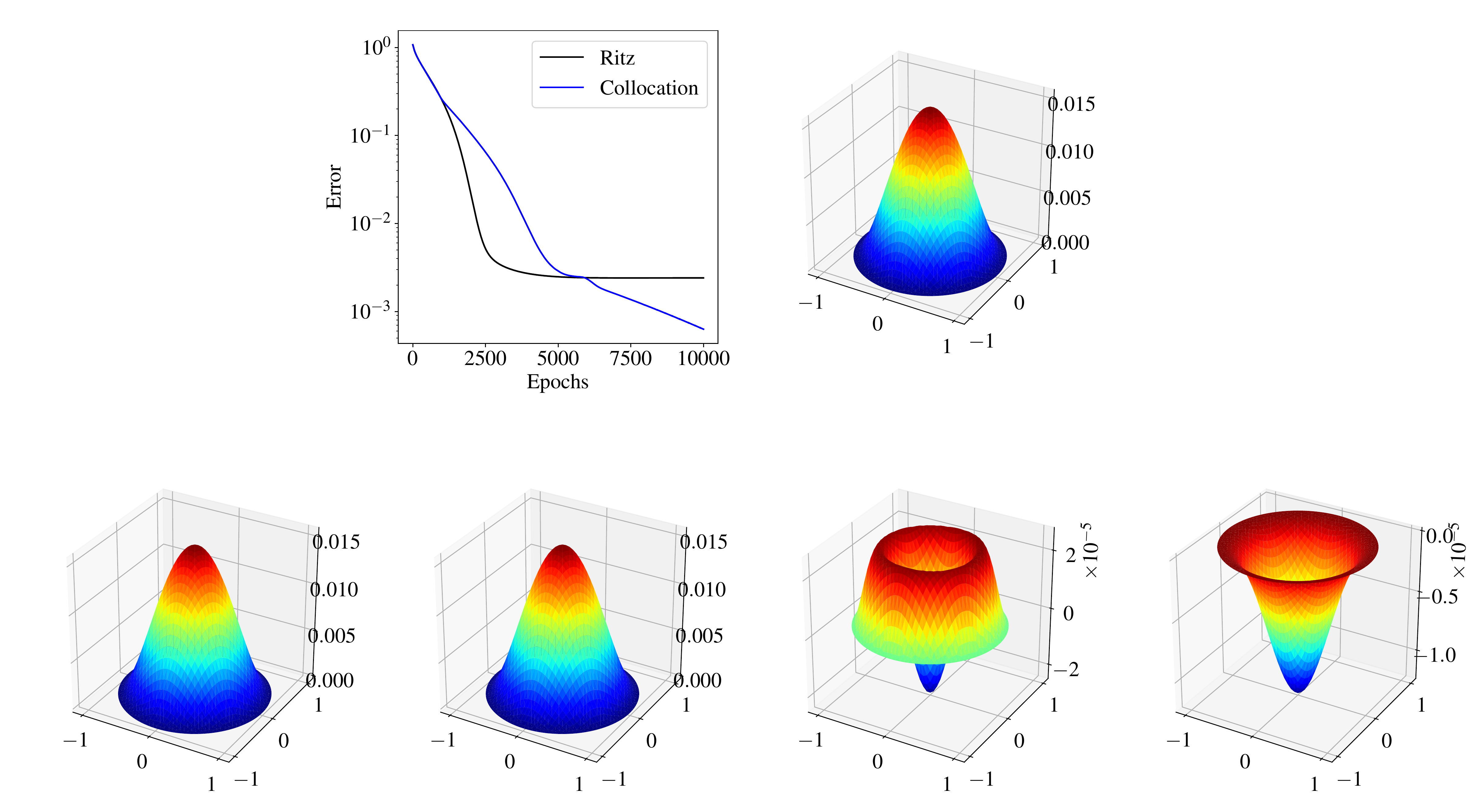}};
\node at (-0.77in,-0.1in) { (a)};
\node at (0.88in,-0.1in) { (b)};
\node at (-2.32in,-1.83in) { (c)};
\node at (-0.73in,-1.83in) { (d)};
\node at ( 0.82in,-1.83in) { (e)};
\node at (2.4in,-1.83in) { (f)};
\end{tikzpicture}
\caption{Numerical solutions for clamped circular Kirchhoff 
plate bending problem using $\unnbcxt$ (Ritz and collocation). 
Network architecture used is 2--50--50--1.
(a) Normalized absolute error in training. Surface plots of
of the
(b) exact solution, and (c), (d) Ritz and collocation solutions. 
Surface plots of (e) error in Ritz solution
and (f) error in collocation solution.
} 
\label{fig:Kirchhoff-plate-bending}
\end{figure}

Numerical results for collocation and Ritz using $\unnbcxt$ are presented in~\fref{fig:Kirchhoff-plate-bending}.  In the
computations, 2800 interior points are used for both methods. For both collocation and Ritz we use the cubic ReLU activation
function.  Note that for collocation 
with standard PINN one cannot use the
cubic ReLU since it is a biharmonic function. The trial function
$\unnbcxt$ has other terms that are present in it, and hence
in general it is not biharmonic ($\nabla^4 \unnbc \neq 0$). A
representative mesh that display the interior points over a disk
is shown in~\fref{fig:FEM-mesh}.
The network architecture is 2--50--50--1. The training loss as a function of epochs for the collocation and Ritz methods are shown
in~\fref{fig:Kirchhoff-plate-bending}a. 
The exact solution
is shown in~\fref{fig:Kirchhoff-plate-bending}b, and the
collocation and Ritz solution for $\unnbcx$ after 10,000 epochs
are plotted 
in Figs~\ref{fig:Kirchhoff-plate-bending}c 
and~\ref{fig:Kirchhoff-plate-bending}d.  The surface plots of 
the error in the Ritz
and collocation solutions are presented
in Figs~\ref{fig:Kirchhoff-plate-bending}e 
and~\ref{fig:Kirchhoff-plate-bending}f. We observe that the
unnormalized error in the Ritz and collocation methods are 
${\cal O}(10^{-5})$; the latter is consistent with the accuracy reported in Guo et al.~\cite{Guo:2021:DCM}.

\subsection{Eikonal equation}\label{subsec:Eikonal}
We consider the Eikonal equation, which is a first-order nonlinear hyperbolic PDE.
The boundary-value problem of the Eikonal equation is:
\begin{subequations}\label{eq:eikonal}
\begin{align}
|| \nabla u || &= \frac{1}{f} \ \  \textrm{in } \Omega \subset \Re^2, \\
u & = 0 \ \ \textrm{on } \Gamma, 
\end{align}
\end{subequations}
where $\Gamma$ is an interface in two dimensions and
$f(\vx) > 0$ is the speed on the interface. If $f(\vx) = 1$, then $u(\vx)$ is the shortest (signed) 
distance from $\vx$ to the boundary $\partial \Omega$. 
If the zero level curve
of $u(\vx)$ represents the initial location of the interface, then
$u^{-1}(t)$ yields the location of the interface at time $t$. Hence,
$u(\vx)$ represents the shortest time (arrival time) that is required to 
travel from the 
boundary $\Gamma$ to $\vx$. For monotonically advancing fronts ($f> 0$),
the fast marching method~\cite{Sethian:1999:LSM} and the fast sweeping method~\cite{Zhao:2005:FSM} are highly efficient
methods to solve~\eqref{eq:eikonal}. When upwind finite-differences are used to 
solve~\eqref{eq:eikonal}, it implies a causality: value of $u(\vx)$ only
depends on values of $u(\vm{y})$ for which $u(\vx) > u(\vm{y})$. No such restriction is used in PINN during training---we use a collocation 
method, where the exact satisfaction of the Dirichlet condition on $\Gamma$ is met by constructing ADFs that use
R-equivalence from~\sref{subsec:Requiv} and the generalized mean value potential from~\sref{subsec:mvp_polygon}.

\begin{figure}[!htp]
\centering
\begin{tikzpicture}
\node at (0,0) {\includegraphics[width=0.98\textwidth]{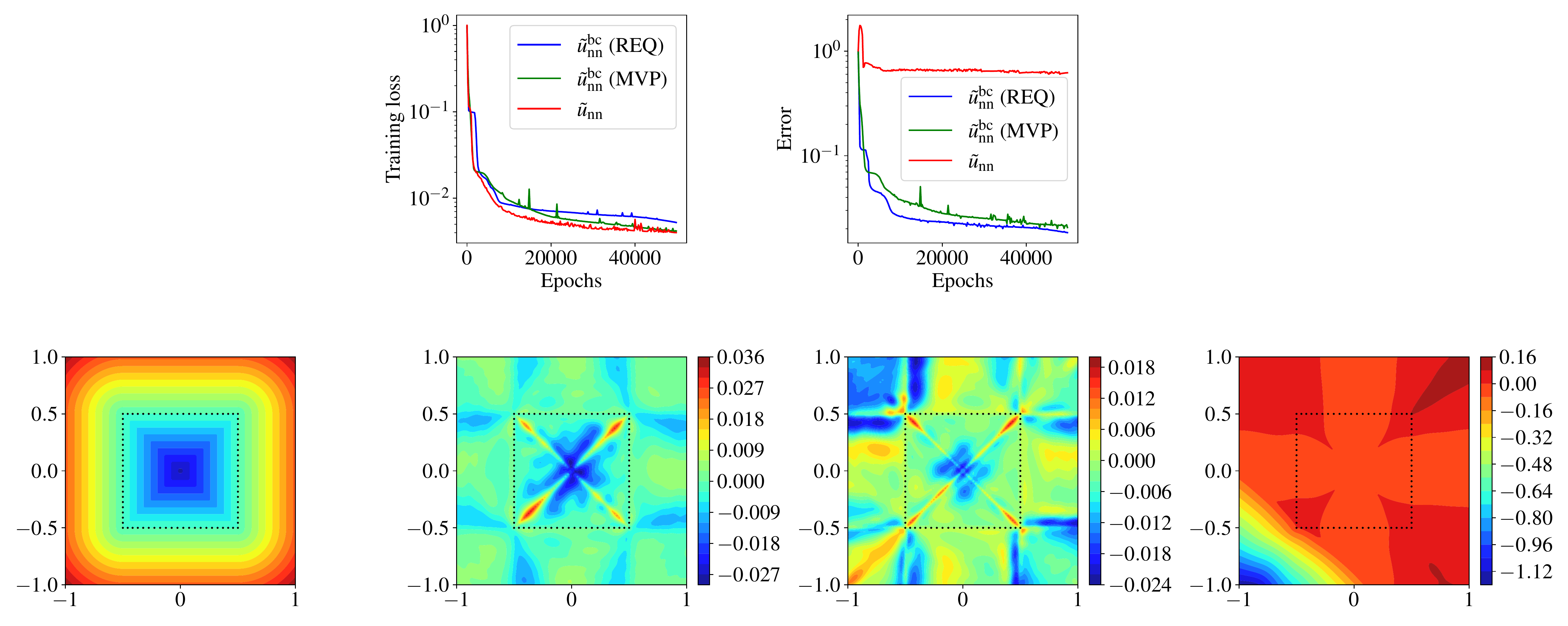} }; 
\node at (-0.85in,-0.03in) {(a)};
\node at (0.75in,-0.03in) {(b) };
\node at (-2.45in,-1.35in) {(c)  };
\node at (-0.86in,-1.35in) {(d)};
\node at ( 0.73in,-1.35in) {(e)};
\node at (2.32in,-1.35in) {(f)};
\end{tikzpicture}
\caption{Solving the Eikonal equation using $\unnbcxt$ and
         $\unnxt$ to compute 
         the signed distance to the boundary of a square. 
         The network architecture used is 2--30--30--30--1 for
         $\unnbc$ and 2--70--70--1 for $\unn$.
         (a), (b) Training loss and normalized absolute error as a function
         of epochs for $\unnbcxt$ (REQ and 
         MVP) and $\unnxt$. 
         (c) Exact signed distance function.  Contour plots
         of the error for 
         (d) $\unnbcx$ (REQ),
         (e) $\unnbcx$ (MVP), and (f) $\unnx$.
}
\label{fig:Eikonal-square}
\end{figure}

We solve~\eqref{eq:eikonal} with $f = 1$ to compute the signed distance function using PINN. The cubic ReLU activation function, which
is a $C^2$ function,
is used for all problems.
The closed interface $\Gamma$ is embedded within the biunit square, $\Omega_0 = (-1,1)^2$, and as benchmark problems
we consider affine (polygonal) and curved interfaces for $\Gamma$. The first problem that we consider is the computation
of the signed distance function to 
the boundary of a smaller square,
$\Omega = (-1/2,1/2)^2$.  The Dirichlet
boundary condition $u = 0$ is imposed on $\Gamma = \partial
\Omega$. The network architecture used is 2--30--30--30--1 for $\unnbc$ and
2--70--70--1 for $\unn$. For collocation,
10,000 points are used in the interior of the biunit square, and 
\suku{400 points} on the boundary $\Gamma$ for $\unnxt$. 
Numerical results are presented in~\fref{fig:Eikonal-square}.
The training loss and normalized absolute errors as a function of epochs for $\unnbcxt$ (REQ and MVP) and $\unnxt$ are shown 
in Figs.~\ref{fig:Eikonal-square}a and~\ref{fig:Eikonal-square}b.
Once again, we notice that 
while $\unnxt$ attains the lowest loss among the three schemes, 
it has the highest error. 
This problem has an exact solution, which is shown as a contour plot in~\fref{fig:Eikonal-square}c. The exact distance, $u(\vx)$, achieves its maximum value of $1/\sqrt{2}$ at the corners of the biunit square and its minimum value of $-1/2$ at the center. The error in the numerical solutions are plotted in Figs.~\ref{fig:Eikonal-square}d--\ref{fig:Eikonal-square}f.
While $\unnbcx$ (REQ and MVP) satisfy the boundary condition exactly, MVP results in a slightly more accurate solution in the entire domain. The $L_\infty$ norm of the error using REQ, MVP
and $\unnx$ are 0.03, 0.026, and 1.34, respectively. The 
standard PINN,
$\unnxt$, does poorly on
satisfying the boundary condition, which in turn leads to
larger pointwise errors over the whole domain.

\begin{figure}[!htp]
\centering
\begin{tikzpicture}
\node at (0,0) {\includegraphics[width=0.98\textwidth]{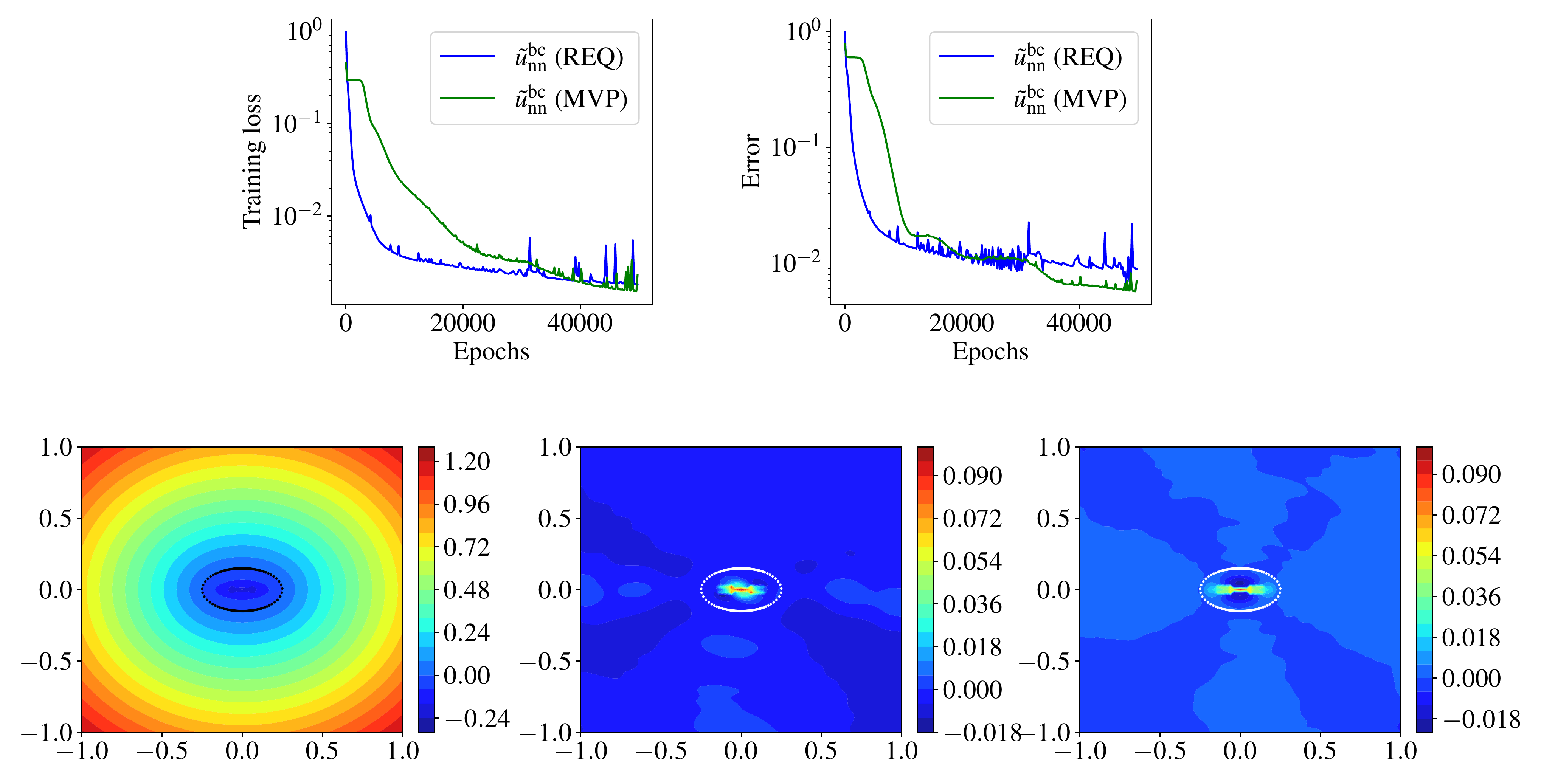} }; 
\node at (-1.14in,-0.02in) {(a)};
\node at (0.92in,-0.02in) {(b) };
\node at (-2.18in,-1.725in) {(c)  };
\node at (-0.12in,-1.725in) {(d)};
\node at (1.94in,-1.725in) {(e)};
\end{tikzpicture}
\caption{Solving the Eikonal equation using $\unnbcxt$ to compute
         the signed distance to the boundary of an ellipse. 
         The network architecture used in 2--50--50--1.
         (a), (b) Training loss and normalized absolute error in training as a function
         of epochs for $\unnbcxt$ (REQ and 
         MVP) .
         (c) Exact signed distance function.
         Contour plots
         of the error for 
         (d) $\unnbcx$ (REQ) and
         (e) $\unnbcx$ (MVP). 
}
\label{fig:Eikonal-ellipse}
\end{figure}
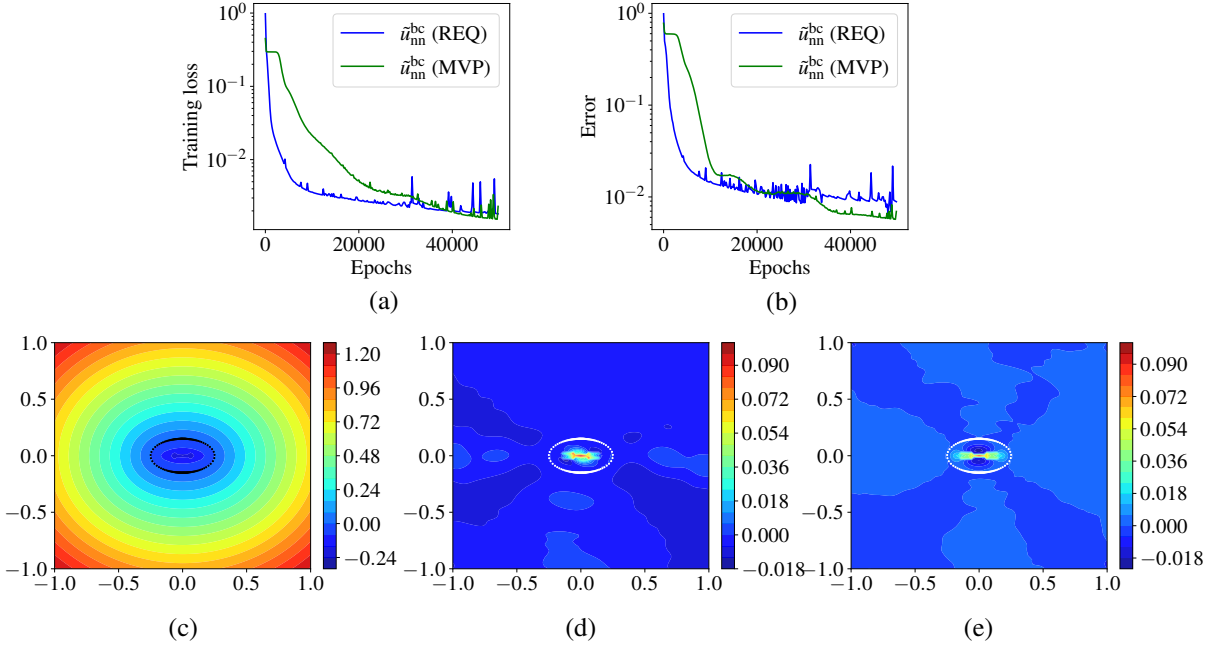

As the next problem, we consider the signed distance function
to the boundary $\Gamma$ of an ellipse that is centered
at the origin and with semi-major and semi-minor axes of
0.25 and 0.15, respectively.  The approximate distance function
to the ellipse, $\phi(\vx)$, is computed using~\eqref{eq:phi_ellipse}
(see~\fref{fig:phi_circle_ellipse-b}) for REQ and
using~\eqref{eq:phi_tmvi} for MVP (see~\fref{fig:tmvi_ADF-a}) for MVP.
The network architecture used is 2--50--50--1. The results that
$\unnx$ produced are very poor, and hence are not included. Numerical results using   
 $\unnbcxt$ (REQ and MVP) are presented in~\fref{fig:Eikonal-ellipse}. The training loss and the normalized error as a function of epochs
 is shown in Figs.~\ref{fig:Eikonal-ellipse}a and~\ref{fig:Eikonal-ellipse}b. The exact distance function
 (computed numerically) is presented 
 in~\fref{fig:Eikonal-ellipse}c, and the error plots for the numerical solutions
 $\unnbcx$ (REQ and MVP) are shown in Figs.~\ref{fig:Eikonal-ellipse}d and~\ref{fig:Eikonal-ellipse}e, respectively. Larger errors are concentrated in the region that is close to the center of the ellipse; away from the center the errors are less than 1 percent.  The exact distance function is $C^0$ (derivative
 discontinuities at the center of the ellipse), whereas the
 numerical solution is $C^2$ smooth.

\begin{figure}[!htp]
\centering
\begin{tikzpicture}
\node at (0,0) {\includegraphics[width=0.98\textwidth]{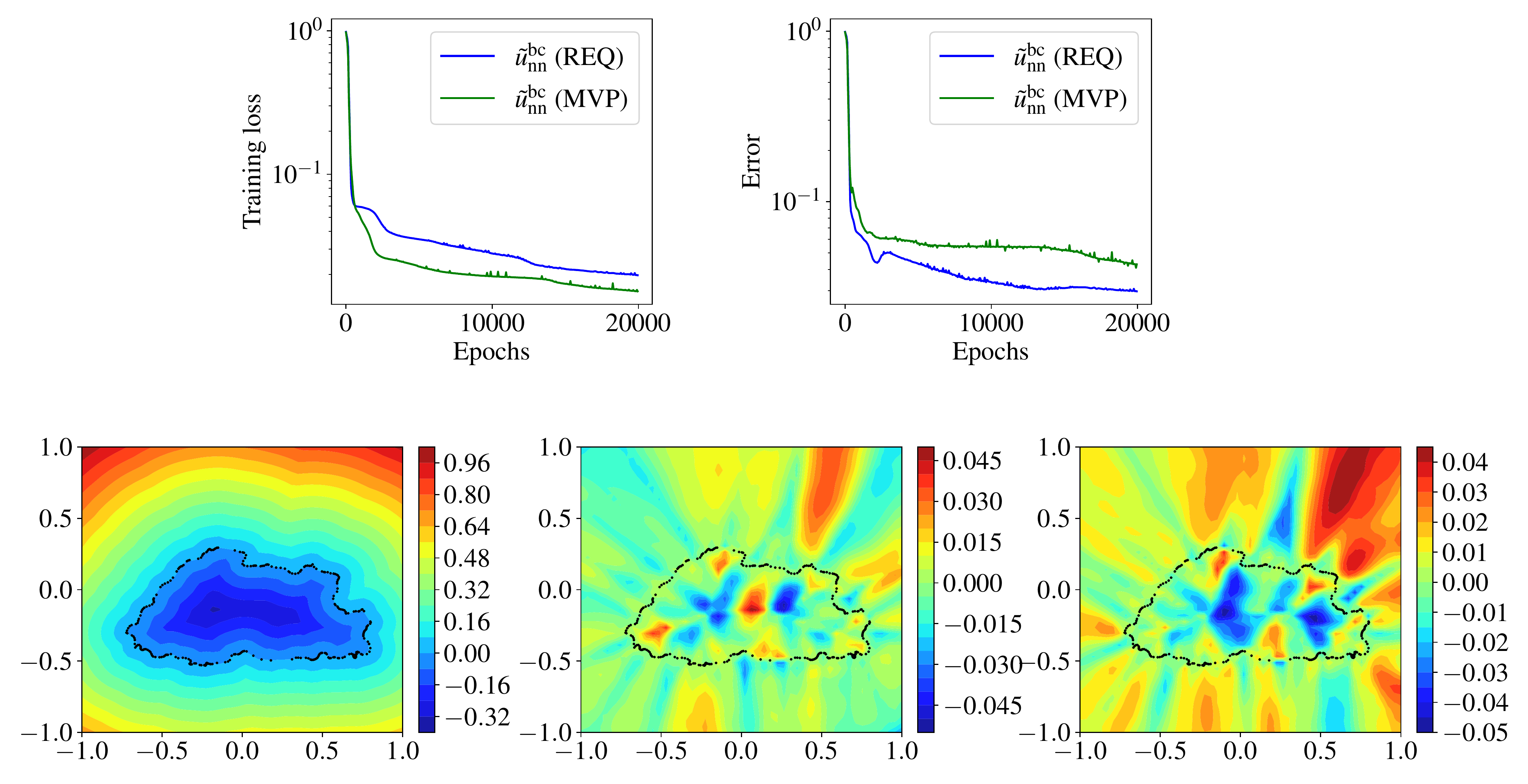}};
\node at (-1.12in,0in) { (a) };
\node at (0.96in,0in) { (b) };
\node at (-2.18in,-1.7in) { (c) };
\node at ( -0.1in,-1.7in) { (d)};
\node at (1.98in,-1.7in) { (e)};
\end{tikzpicture}
\caption{Solving the Eikonal equation using $\unnbcxt$ to compute
         the signed distance to the boundary of the
         polygonal map of Bhutan.
         The network architecture is \suku{2--50--50--1}.
         (a), (b) Training loss and normalized absolute error in training as a function
         of epochs for $\unnbcxt$ (REQ and 
         MVP).
         (c) Exact signed distance function.
         Contour plots
         of the error for 
         (d) $\unnbcx$ (REQ) and
         (e) $\unnbcx$ (MVP).
}
\label{fig:Bhutan}
\end{figure}
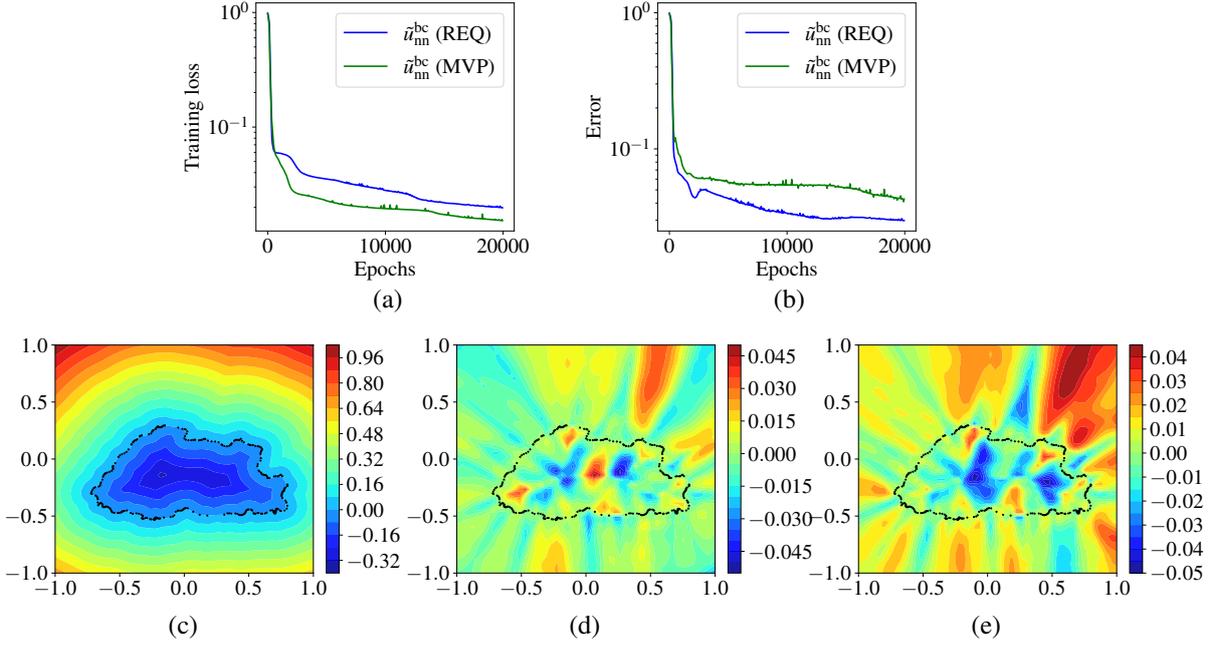

Lastly, we consider the signed distance function
to the boundary $\Gamma$ of the polygonalized map of Bhutan;
see plots of the approximate distance functions to $\Gamma$ 
using REQ and MVP that are presented in Figs.~\ref{fig:phi_bhutan_map} and~\ref{fig:mvc_bhutan_map}, respectively.  The 
network architecture used is 2--50--50--1. Here too the
results of $\unnx$ are not included since they are very poor.
Numerical results using $\unnbcxt$ (REQ and MVP) are presented in~\fref{fig:Bhutan}. 
The training loss and the normalized absolute error as a function of epochs are shown in Figs~\ref{fig:Bhutan}a and~\ref{fig:Bhutan}b. 
At 20,000 epochs, the absolute normalized error is ${\cal O}(10^{-1})$.
The exact distance function
is presented 
in~\fref{fig:Bhutan}c, and the error plots for the numerical solutions
$\unnbcx$ (REQ and MVP) are presented in Figs.~\ref{fig:Bhutan}d and~\ref{fig:Bhutan}e, respectively. The maximum error is about 4 percent.

\section{Poisson Problem over the Four-Dimensional Hypercube} \label{sec:4D}
As the last problem, we consider the model isotropic steady-state heat conduction (Poisson) problem over the 4-dimensional hypercube to show that essential boundary conditions can be readily imposed in higher dimensions as well. Consider the following Poisson problem
with homogeneous Dirichlet boundary conditions:
\begin{subequations}\label{eq:heat_bvp4d}
\begin{align}
-\nabla^2 u &= f \ \ \textrm{in } \Omega = (-1,1)^4 \\
u &= 0 \ \  \textrm{on }  \partial\Omega ,
\end{align}
\end{subequations}
where $u(\vx) : \Re^4 \to \Re$ is sought and
$f(\vx)$ is the forcing function.  Let
$\vx := (x_1,x_2,x_3,x_4) \in \Omega$ denote 
a point in the hypercube. We
choose $f(\vx)=\prod_{i=1}^4\sin(\pi x_i)$, so that the exact solution is:
\begin{equation*}
u(\vx) = \frac{ \prod_{i=1}^4\sin(\pi x_i) } {4\pi^2}. 
\end{equation*}
As noted in prior high-dimensional studies using PINN~\cite{E:2018:DRM,Sirignano:2018:DGM}, numerical solutions for problems in high-dimensions are challenging since there is no easy way to mesh the domain and they are also subject to the curse of dimensionality. Among meshfree methods, since construction of 
radial basis functions is dimension-independent, RBF-based meshfree methods have had success in solving high-dimensional problems~\cite{Cecil:2004:NMH}.
It is in these problems that the power and potential
of a meshfree method such as PINN becomes most apparent. 
For this problem, we only consider trial functions, $\unnbcxt$, which exactly enforce the homogeneous Dirichlet boundary condition on $\partial \Omega$. One way of enforcing the boundary conditions is to assume $\unnbc$ to be of the form
\begin{equation*}
    \unnbcxt= \left( \prod_{i=1}^4(1 - x_i^2) \right) \unnRxt ,
\end{equation*}
which we refer to as the `product method' and note that while this is an obvious approach for the present problem, it leads to very small numbers inside the domain and away from the boundaries. In this case, the multiplicative factor scales as ${\cal O}(x^8)$ inside the biunit hypercube, and therefore the network parameters have to compensate for this highly nonlinear behavior during training. It is preferable to have a multiplicative factor that is better behaved in order to aid the training. To this end, we construct $\phi(\vx)$ using 
R-equivalence in~\eqref{eq:phin_eq}, which
seamlessly extends to higher dimensions. For this choice, the trial
function $\unnbcxt$ is of the form:
\begin{equation}
   \unnbcxt = \phi(\vx) \, \unnRxt  ,
\end{equation}
where $\phi(\vx)$
consists of R-equivalence (REQ) operations on $\phi_i(\vx)$, where 
$\phi_i(\vx)$ is the R-function for the region (strip) bounded by the hyperplanes $1 - x_i$ and $x_i - 1$. We form $\phi_i(\vx) = (1-x_i^2)/2$, which is an ADF that is normalized to order 1. On using the REQ composition in~\eqref{eq:phin_eq}, we write
\begin{equation}\label{eq:REQ_4D}
    \phi (\vx) = \phi_1(\vx) \sim \phi_2(\vx) \sim \phi_3(\vx)
    \sim \phi_4(\vx) ,
\end{equation}
which generalizes to the hypercube in $\Re^d$.
Note that $\phi(\vx)$ only scales as $x_i^2$ in each coordinate
direction, and therefore is much better behaved. For this
problem, another choice for $\phi(\vx)$ is to define two
$\phi_i$'s in each dimension, i.e., $\phi_{2i-1} = 1 + x_i$
and $\phi_{2i} = 1 - x_i$ and then define
$\phi(\vx) = \phi_1 \sim \phi_2 \dots \phi_8$ for the
4-dimensional hypercube. 
\suku{For $m = 1$, this construction
coincides with the expression for $\phi$ in~\eqref{eq:REQ_4D}.
}

\begin{figure}[htp]
\centering
\begin{subfigure}{0.30\textwidth}
\includegraphics[width=\textwidth]{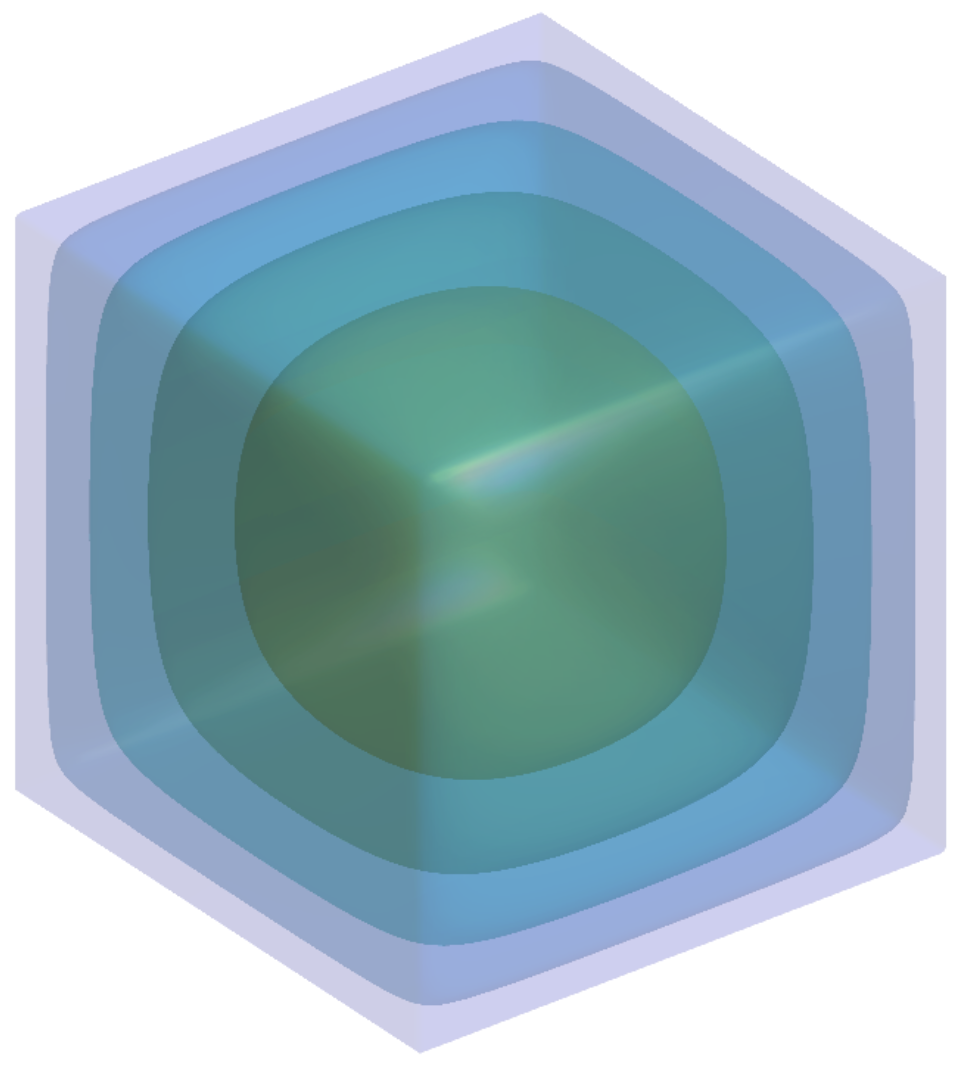} 
\caption{}
\label{fig:HeatEq4D-a}
\end{subfigure}
\begin{subfigure}{0.68\textwidth}
\includegraphics[width=\textwidth]{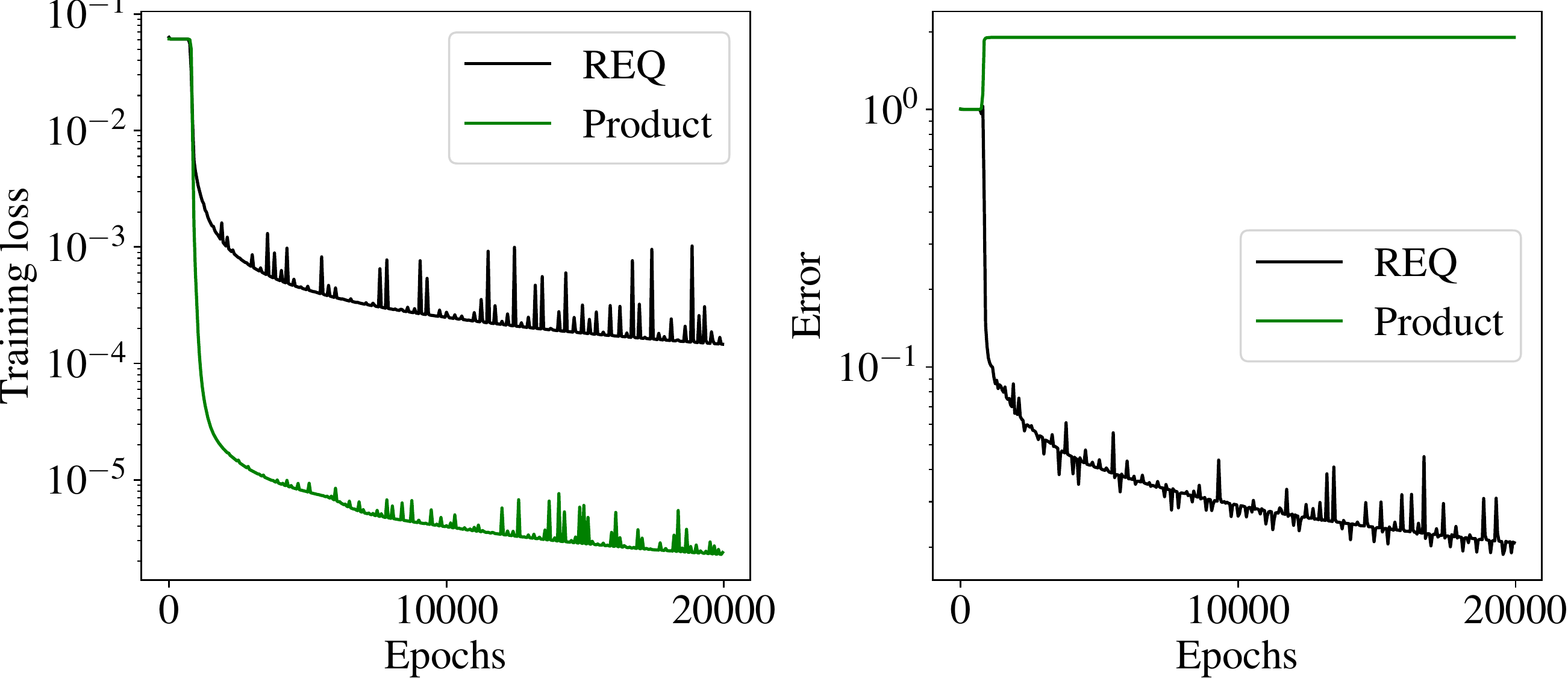} 
\caption{}
\label{fig:HeatEq4D-b}
\end{subfigure}
\caption{Numerical solution of the Poisson problem over the four-dimensional hypercube. The network architecture is 2--100--100--1. The trial function, $\unnbcx$, is constructed using
R-equivalence (REQ) and the product method. (a) Isosurface plot of the approximate distance function using REQ, $\phi(\vx)$, over the 3-dimensional biunit cube ($x_4 = 0$), with $\phi = 0$ being satisfied on the boundary of the cube.
(b) Training loss and normalized absolute error for REQ 
and the product method as a function of epochs. 
}
\label{fig:HeatEq4D}
\end{figure}

In~\fref{fig:HeatEq4D}, we present the numerical solutions using
the product and REQ ($m = 1$) forms of the trial function. 
For both choices, we consider 5,000 randomly generated interior points in $\Omega$ for training and compute the normalized error at a separate set of 5,000 interior points. The network architecture for both choices is 4--100--100--1. 
The
isosurface plot of $\phi(\vx)$ ($x_4 = 0$ plane) using REQ
is shown in~\fref{fig:HeatEq4D-a}. In~\fref{fig:HeatEq4D-b}, the
evolution of the training loss and normalized absolute error is presented
for the product
and REQ trial functions.
We observe that while the REQ method is able to reach error levels of about 1 percent, the product method does not converge.
The $\unnbcx$ solution with REQ yields
${\cal O}(10^{-2})$ error, whereas the product method
has ${\cal O}(1)$ error even \suku{though} it has a much smaller
PDE loss of ${\cal O}(10^{-5})$.

\section{Conclusions}\label{sec:conclusions}
Starting from the seminal works of Lagaris et al.~\cite{Lagaris:1998:ANN,Lagaris:1997:ANN,Lagaris:2000:NNM} and the recent extensions and major advancements by
Raissi et al.~\cite{Raissi:2019:PIN} and
E and Yu~\cite{E:2018:DRM}, there has been a surge in the development and application of physics-informed neural networks to solve partial different equations. In this paper, we have introduced a new approach based on distance fields to
construct {\em geometry-aware} approximations in physics-informed neural networks by ensuring that the necessary boundary conditions are met a priori: all boundary conditions in a collocation method, and the essential boundary conditions (kinematic admissibility) in a Ritz method. Our approach relied on the theory of R-functions~\cite{Rvachev:1995:RBV,Rvachev:2000:OCR} to construct approximate distance fields and their use within a meshfree method to exactly impose boundary conditions to solve PDEs~\cite{Shapiro:1999:MSD,Shapiro:2007:SAG}. Apart from R-functions, we also showed that mean value potential fields~\cite{Floater:2003:MVC,Dyken:2009:TMV,Belyaev:2013:SLD} can be
used to construct suitable distance field to solve PDEs over domains with affine as well as curved boundaries. 

We presented several numerical examples to reveal the benefits of exactly imposing the boundary conditions versus the current state-of-the-art in deep collocation and deep Ritz methods for physics-informed neural networks.  Notably, requiring only the interior residual error contribution \suku{in} the loss function simplifies the training of the network and leads to more accurate numerical solutions.  This was shown through several verification tests on benchmark one- and two-dimensional boundary-value problems---and consistently revealed the pitfalls of being guided by the magnitude of the loss function in standard 
PINN-based collocation approaches~\cite{Raissi:2019:PIN}. In PINN methods, when there
are multiple terms (PDE loss and boundary losses) that are present in the loss function \suku{and the loss weights are fixed a priori}, the magnitude of the loss function at the end of the training does not provide a measure of the accuracy of the approximation. There can be a many-fold difference in the two error measures, which is revealed in our simulations.  This is not surprising since the weights associated with each loss term is not known a priori since it depends on the PDE, boundary conditions, and the training. There is no clear rationale way to set these weights \suku{in order to ensure that} the approach is robust and guaranteed to lead to reliable results. This study has reinforced that it is important that the PDE loss stands on its own, and boundary conditions are enforced via the ansatz. One approach is to construct 
a separate neural
network to meet the boundary conditions as is done by 
\suku{Berg and Nystr{\"o}m~\cite{Berg:2018:UDA}}, but the inherent inaccuracy in the satisfaction of the boundary conditions then propagates when training for the PDE is conducted. Moreover, when the geometry is complicated and different types of boundary conditions are imposed on different subset of the boundary, then this approach may soon become impractical. Our approach ensures exact satisfaction of \suku{all} necessary boundary condition, which makes it appealing---so that training of the network (training loss and PDE loss coincide) is more efficient and accurate solutions can be realized. 

This study has provided a method to perform meshfree analysis---solving PDEs without domain discretization---on complex two-dimensional geometries using physics-informed neural networks. \suku{This was made possible on using approximate distance functions, which were based on
R-functions (R-equivalence composition) and generalized mean value potentials, in conjunction with transfinite inverse-distance based interpolation to exactly satisfy the boundary conditions.
For the problems that we considered, both ADFs delivered the same order of accuracy with the constant being smaller for R-equivalence in most cases (Eikonal equation was the exception). If $u = 0$ is imposed on the boundary of a
polygon with many (tens to hundreds) edges, then the ADF using
mean value potential is more efficient since it results in an ADF by construction, whereas R-equivalence requires first constructing an ADF to each edge and then the joining operation to be performed. However,
if different boundary conditions are imposed on distinct boundary segments then only R-equivalence is viable since 
this cannot be realized using the ADF based on mean value potential fields.
While we have used R-functions and the mean value potential to construct smooth approximate distance functions, fields such as constructive geometric modeling~\cite{Ricci:1973:CGC} with implicit
functions~\cite{Bloomenthal:1997:BEC,Sherstyuk:1999:KFC,Barthe:2004:CBO,Gourmel:2013:GBI}, PDE-based solutions for distance computations~\cite{Belyaev:2015:VPB,Crane:2017:HMD,Belyaev:2019:VMA},
and deep learning in computer vision~\cite{Park:2019:SDF} 
are rapidly advancing and may offer other attractive alternatives to construct smooth distance fields
for use in deep neural networks to solve PDEs.
A separate and more in depth investigation is needed to
explore if there are other network architectures and optimizers for network training that are better-suited for the PINN ansatz with approximate distance functions, and to also quantify the
accuracy that is obtained using different ADFs.
Lastly, the formulation permitted exact modeling of 
affine and curved boundaries, thereby providing a pathway to conducting simulations on
the exact geometry (isogeometric analysis)~\cite{Hughes:2005:IAC}.}

The ideas herein can be extended for higher dimensional problems, since R-equivalence composition for the approximate distance function is additive and does not suffer from the curse of dimensionality. This was demonstrated in~\sref{sec:4D}, where we obtained an accurate PINN solution
for a Poisson problem over the 4-dimensional hypercube. 
Extending our formulation to complex geometries in 3D, and
the development of a deep Petrov-Galerkin domain-decomposition method are topics that we plan to pursue.

\section{Acknowledgments}
AS acknowledges support from the NSF CAREER grant \#1554033 to the Illinois Institute of
Technology.
NS thanks Anand Reddy, Eric Chin and
Kai Hormann for many helpful discussions.


\begin{thebibliography}{100}
\expandafter\ifx\csname url\endcsname\relax
  \def\url#1{\texttt{#1}}\fi
\expandafter\ifx\csname urlprefix\endcsname\relax\def\urlprefix{URL }\fi
\expandafter\ifx\csname href\endcsname\relax
  \def\href#1#2{#2} \def\path#1{#1}\fi

\bibitem{Lagaris:1998:ANN}
I.~E. Lagaris, A.~Likas, D.~I. Fotiadis, Artifical neural networks for solving
  ordinary and partial differential equations, IEEE Transactions on Neural
  Networks 9~(5) (1998) 987--1000.

\bibitem{Lagaris:1997:ANN}
I.~E. Lagaris, A.~Likas, D.~I. Fotiadis, Artifical neural network methods in
  quantum mechanics, Computer Physics Communications 104 (1997) 1--14.

\bibitem{Lagaris:2000:NNM}
I.~E. Lagaris, A.~C. Likas, D.~G. Papageorgiou, Neural-network methods for
  boundary value problems with irregular boundaries, IEEE Transactions on
  Neural Networks 11~(5) (2000) 1041--1049.

\bibitem{McFall:2006:ANN}
K.~S. McFall, An artificial neural network method for solving boundary value
  problems with arbitrary irregular boundaries, Ph.D. thesis, Georgia Institute
  of Technology, Atlanta, GA, USA (2006).

\bibitem{McFall:2009:ANN}
K.~S. McFall, J.~R. Mahan, Artificial neural network method for solution of
  boundary value problems with exact satisfaction of arbitrary boundary
  conditions, IEEE Transactions on Neural Networks 20~(8) (2009) 1221--1233.

\bibitem{Raissi:2019:PIN}
M.~Raissi, P.~Perdikaris, G.~E. Karniadakis, Physics-informed neural networks:
  {A} deep learning framework for forward and inverse problems involving
  nonlinear partial differential equations, Journal of Computational Physics
  378 (2019) 686--707.

\bibitem{Berg:2018:UDA}
J.~Berg, K.~Nystr{\"o}m, A unified deep artificial neural network approach to
  partial differential equations in complex geometries, Neuralcomputing 317
  (2018) 28--41.

\bibitem{Sirignano:2018:DGM}
J.~Sirignano, K.~Spiliopoulos, {DGM}: {A} deep learning algorithm for solving
  partial differential equations, Journal of Computational Physics 375 (2018)
  1339--1364.

\bibitem{E:2018:DRM}
W.~E, B.~Yu, The deep {Ritz} method: a deep learning-based numerical algorithm
  for solving variational problems, Communications in Mathematics and
  Statistics 6~(1) (2018) 1--12.

\bibitem{Han:2018:SHD}
J.~Han, A.~Jentzen, W.~E, Solving high-dimensional partial differential
  equations using deep learning, Proceedings of the National Academy of
  Sciences 115~(34) (2018) 8505--8510.

\bibitem{Kharazmi:2019:VPI}
E.~Kharazmi, Z.~Zhang, G.~E. Karniadakis, Variational physics-informed neural
  networks for solving partial differential equations (2019).
\newblock \href {http://arxiv.org/abs/1912.00873} {\path{arXiv:1912.00873}}.

\bibitem{Kharazmi:2020:HPV}
E.~Kharazmi, Z.~Zhang, G.~E. Karniadakis, {hp-VPINNs}: {Variational}
  physics-informed neural networks with domain decomposition, Computer Methods
  in Applied Mechanics and Engineering 374 (2020) 113547.

\bibitem{Lu:2020:DLL}
L.~Lu, X.~Meng, Z.~Mao, G.~E. Karniadakis, {DeepXDE}: {A} deep learning library
  for solving differential equations, SIAM Review 63~(1) (2021) 208--228.

\bibitem{Abadi:2016:TEF}
M.~Abadi, P.~Barham, J.~Chen, Z.~Chen, A.~Davis, J.~Dean, M.~Devin,
  S.~Ghemawat, G.~Irving, M.~Isard, et~al., Tensorflow: A system for
  large-scale machine learning, in: 12th $\{$USENIX$\}$ symposium on operating
  systems design and implementation ($\{$OSDI$\}$ 16), 2016, pp. 265--283.

\bibitem{Paszke:2016:PYT}
A.~Paszke, S.~Gross, F.~Massa, A.~Lerer, J.~Bradbury, G.~Chanan, T.~Killeen,
  Z.~Lin, N.~Gimelshein, L.~Antiga, A.~Desmaison, A.~Kopf, E.~Yang, Z.~DeVito,
  M.~Raison, A.~Tejani, S.~Chilamkurthy, B.~Steiner, L.~Fang, J.~Bai,
  S.~Chintala, Pytorch: An imperative style, high-performance deep learning
  library, in: H.~Wallach, H.~Larochelle, A.~Beygelzimer, F.~d'Alch\'{e} Buc,
  E.~Fox, R.~Garnett (Eds.), Advances in Neural Information Processing Systems,
  Vol.~32, Curran Associates, Inc., 2019.

\bibitem{Wang:2020:UMG}
S.~Wang, Y.~Teng, P.~Perdikaris, Understanding and mitigating gradient
  pathologies in physics-informed neural networks (2020).
\newblock \href {http://arxiv.org/abs/2001.04536} {\path{arXiv:2001.04536}}.

\bibitem{Chen:2020:CSD}
J.~Chen, R.~Du, K.~Wu, A comparison study of deep {Galerkin} method and deep
  {Ritz} method for elliptic problems with different boundary conditions,
  Communications in Mathematical Research 36~(3) (2020) 354–376.

\bibitem{Lyu:2020:EEB}
L.~Lyu, K.~Wu, R.~Du, J.~Chen, Enforcing exact boundary and initial conditions
  in the deep mixed residual method (2020).
\newblock \href {http://arxiv.org/abs/2008.01491} {\path{arXiv:2008.01491}}.

\bibitem{Babuska:2003:SMG}
I.~Babu\v{s}ka, U.~Banerjee, J.~E. Osborn, Survey of meshless and generalized
  finite element methods: a unified approach, Acta Numerica 12 (2003) 1--125.

\bibitem{Huerta:2018:MM}
A.~Huerta, T.~Belytschko, S.~Fern\'andez-M\'endez, T.~Rabczuk, X.~Zhuang,
  M.~Arroyo, Meshfree methods, 2nd Edition, Vol.~2 of Encyclopedia of
  Computational Mechanics, Wiley, 2017, Ch.~3, pp. 1--38.

\bibitem{Kantorovich:1958:AMH}
L.~V. Kantorovich, V.~I. Krylov, Approximate Methods of Higher Analysis,
  Interscience, New York, NY, USA, 1958.

\bibitem{Rvachev:1982:TRF}
V.~L. Rvachev, Theory of R-functions and Some Applications, Naukova Dumka,
  Kiev. In Russian, 1982.

\bibitem{Rvachev:1995:RBV}
V.~L. Rvachev, T.~I. Sheiko, R-functions in boundary value problems in
  mechanics, Applied Mechanics Reviews 48~(4) (1995) 151--188.

\bibitem{Rvachev:2000:OCR}
V.~L. Rvachev, T.~I. Sheiko, V.~Shapiro, I.~Tsukanov, On completeness of {RFM}
  solution structures, Computational Mechanics 25 (2000) 305--3163.

\bibitem{Rvachev:2001:TII}
V.~L. Rvachev, T.~I. Sheiko, V.~Shapiro, I.~Tsukanov, Transfinite interpolation
  over implicitly defined sets, Computer Aided Geometric Design 18 (2001)
  195--220.

\bibitem{Shapiro:1991:TRF}
V.~Shapiro, Theory of {R-functions} and applications: {A} primer, Tech. Rep.
  CPA88-3, Cornell Programmable Automation, Sibley School of Mechanical
  Engineering, Ithaca, NY 14853, USA (1991).

\bibitem{Shapiro:1999:MSD}
V.~Shapiro, I.~Tsukanov, Meshfree simulation of deforming domains,
  Computer-Aided Design 31~(7) (1999) 459--471.

\bibitem{Shapiro:2002:ASA}
V.~Shapiro, I.~Tsukanov, The architecture of {SAGE}--a meshfree system based on
  {RFM}, Engineering with Computers 18~(4) (2002) 295--311.

\bibitem{Biswas:2004:ADF}
A.~Biswas, V.~Shapiro, Approximate distance fields with non-vanishing
  gradients, Graphical Models 66~(3) (2004) 133--159.

\bibitem{Shapiro:2007:SAG}
V.~Shapiro, Semi-analytic geometry with {R}-functions, Acta Numerica 16 (2007)
  239--303.

\bibitem{Freytag:2011:FEA}
M.~Freytag, V.~Shapiro, I.~Tsukanov, Finite element analysis in situ, Finite
  Elements in Analysis and Design 47~(9) (2011) 957--972.

\bibitem{Hollig:2001:WEB}
K.~H\"{o}llig, U.~Reif, J.~Wipper, Weighted extended {B}-spline approximation
  of {D}irichlet problems, SIAM Journal on Numerical Analysis 39~(2) (2001)
  442--462.

\bibitem{Millan:2015:CBM}
D.~Mill{\'a}n, N.~Sukumar, M.~Arroyo, Cell-based maximum-entropy approximants,
  Computer Methods in Applied Mechanics and Engineering 284 (2015) 712--731.

\bibitem{Floater:2003:MVC}
M.~S. Floater, Mean value coordinates, Computer Aided Geometric Design 20~(1)
  (2003) 19--27.

\bibitem{Dyken:2009:TMV}
C.~Dyken, M.~S. Floater, Transfinite mean value interpolation, Computer Aided
  Geometric Design 26~(1) (2009) 117--134.

\bibitem{Belyaev:2013:SLD}
A.~Belyaev, P.-A. Fayolle, A.~Pasko, Signed {$L_p$}-distance fields,
  Computer-Aided Design 45~(2) (2013) 523--528.

\bibitem{Hughes:2005:IAC}
T.~J.~R. Hughes, J.~A. Cottrell, Y.~Bazilevs, Isogeometric analysis: {CAD},
  finite elements, {NURBS}, exact geometry and mesh refinement, Computer
  Methods in Applied Mechanics and Engineering 194~(39--41) (2005) 4135--4195.

\bibitem{Thoutireddy:2004:VRA}
P.~Thoutireddy, M.~Ortiz, A variational $r$-adaption and shape-optimization
  method for finite-deformation elasticity, International Journal for Numerical
  Methods in Engineering 61~(1) (2004) 1--21.

\bibitem{He:2018:RDN}
J.~He, L.~Li, J.~Xu, C.~Zheng, {ReLU} deep neural networks and linear finite
  elements (2018).
\newblock \href {http://arxiv.org/abs/1807.03973} {\path{arXiv:1807.03973}}.

\bibitem{Grinspun:2003:BRM}
E.~Grinspun, The basis refinement method, Ph.D. thesis, California Institute of
  Technology, Pasadena, CA, USA (2003).

\bibitem{Cyr:2020:RTI}
E.~C. Cyr, M.~A. Gulian, R.~G. Patel, M.~Perego, N.~A. Trask, Robust training
  and initialization of deep neural networks: {An} adaptive basis viewpoint,
  in: Mathematical and Scientific Machine Learning, PMLR, 2020, pp. 512--536.

\bibitem{Opschoor:2020:DRN}
J.~A.~A. Opschoor, P.~C. Petersen, C.~Schwab, Deep {ReLU} networks and
  high-order finite element methods, Analysis and Applications 18~(05) (2020)
  715--770.

\bibitem{Kansa:1990a:MUQ}
E.~J. Kansa, Multiquadrics---{A} scattered data approximation scheme for
  applications to computational fluid-dynamics. 1. {S}urface approximations and
  partial derivative estimates, Computers \& Mathematics with Applications
  19~(8/9) (1990) 127--145.

\bibitem{Kansa:1990b:MUQ}
E.~J. Kansa, Multiquadrics---{A} scattered data approximation scheme for
  applications to computational fluid-dynamics. 2. {S}olutions to parabolic,
  hyperboloc and elliptic partial-differential equations, Computers \&
  Mathematics with Applications 19~(8/9) (1990) 147--161.

\bibitem{Buhmann:2003:RBF}
M.~D. Buhmann, Radial basis functions: theory and implementations, Cambridge
  University Press, Cambridge, UK, 2003.

\bibitem{Fasshauer:2007:MAM}
G.~Fasshauer, Meshfree Approximation Methods in {MATLAB}, Interdisciplinary
  Mathematical Sciences -- Vol. 6, World Scientific Publishers, Singapore,
  2007.

\bibitem{Schaback:2006:KTM}
R.~Schaback, H.~Wendland, Kernel techniques: from machine learning to meshless
  methods, Acta Numerica 15 (2006) 543.

\bibitem{Arroyo:2006:LME}
M.~Arroyo, M.~Ortiz, Local \emph{maximum-entropy} approximation schemes: a
  seamless bridge between finite elements and meshfree methods, International
  Journal for Numerical Methods in Engineering 65~(13) (2006) 2167--2202.

\bibitem{Babuska:1997:PUM}
I.~Babu\v{s}ka, J.~M. Melenk, The partition of unity method, International
  Journal for Numerical Methods in Engineering 40 (1997) 727--758.

\bibitem{Rajan:1994:ODT}
V.~Rajan, Optimality of the {D}elaunay triangulation in {$\mathbb{R}^d$},
  Discrete \& Computational Geometry 12~(1) (1994) 189--202.

\bibitem{Sukumar:2004:COP}
N.~Sukumar, Construction of polygonal interpolants: a maximum entropy approach,
  International Journal for Numerical Methods in Engineering 61~(12) (2004)
  2159--2181.

\bibitem{Sukumar:2005:MEA}
N.~Sukumar, Maximum entropy approximation, AIP Conference Proceedings 803~(1)
  (2005) 337--344.

\bibitem{Arroyo:2007:LME}
M.~Arroyo, M.~Ortiz, Local maximum-entropy approximation schemes, in:
  M.~Griebel, M.~A. Schweitzer (Eds.), Meshfree Methods for Partial
  Differential Equations {III}, Vol.~57 of Lecture Notes in Computational
  Science and Engineering, Springer, Berlin, Germany, 2007, pp. 1--16.

\bibitem{Sukumar:2007:OAC}
N.~Sukumar, R.~W. Wright, Overview and construction of meshfree basis
  functions: {From} moving least squares to entropy approximants, International
  Journal for Numerical Methods in Engineering 70~(2) (2007) 181--205.

\bibitem{Rosolen:2010:OSS}
A.~Rosolen, D.~Mill{\'a}n, M.~Arroyo, On the optimum support size in meshfree
  methods: a variational adaptivity approach with maximum entropy approximants,
  International Journal for Numerical Methods in Engineering 82~(7) (2010)
  868--895.

\bibitem{Park:1991:UAU}
J.~Park, I.~W. Sandberg, Universal approximation using radial-basis-function
  networks, Neural computation 3 (1991) 246--257.

\bibitem{Mhaskar:1996:NNO}
H.~N. Mhaskar, Neural networks for optimal approximation of smooth and analytic
  functions, Neural computation 8~(1) (1996) 164--177.

\bibitem{Lee:2020:PUN}
K.~Lee, N.~A. Trask, R.~G. Patel, M.~A. Gulian, E.~C. Cyr, Partition of unity
  networks: deep hp-approximation (2021).
\newblock \href {http://arxiv.org/abs/2101.11256} {\path{arXiv:2101.11256}}.

\bibitem{Ramabathiran:2021b:SPI}
A.~A. Ramabathiran, P.~Ramachandran, {SPINN: Sparse, physics-based, and
  partially interpretable neural networks for PDEs}, Journal of Computational
  Physics 445 (2021) 110600.

\bibitem{Greco:2020:HOM}
F.~Greco, M.~Arroyo, High-order maximum-entropy collocation methods, Computer
  Methods in Applied Mechanics and Engineering 367 (2020) 113115.

\bibitem{Sheng:2021:PFN}
H.~Sheng, C.~Yang, {PFNN}: A penalty-free neural network method for solving a
  class of second-order boundary-value problems on complex geometries, Journal
  of Computational Physics 428 (2021) 110085.

\bibitem{Dwivedi:2020:PIE}
V.~Dwivedi, B.~Srinivasan, Physics informed extreme learning machine
  {(PIELM)}--{A} rapid method for the numerical solution of partial
  differential equations, Neuralcomputing 391 (2020) 96--118.

\bibitem{Dwivedi:2020:SBE}
V.~Dwivedi, B.~Srinivasan, Solution of biharmonic equation in complicated
  geometries with physics informed extreme learning machine, Journal of
  Computing and Information Science in Engineering 20~(6) (2020).

\bibitem{Liao:2019:DNM}
Y.~Liao, P.~Ming, Deep {Nitsche} method: {Deep Ritz} method with essential
  boundary conditions (2019).
\newblock \href {http://arxiv.org/abs/1912.01309} {\path{arXiv:1912.01309}}.

\bibitem{Li:2004:MPP}
S.~Li, W.~K. Liu, Meshfree Particle Methods, Springer-Verlag, New York, NY,
  USA, 2004.

\bibitem{Hornik:1989:MFN}
K.~Hornik, M.~Stinchcombe, H.~White, Multilayer feedforward networks are
  universal approximators, Neural Networks 2 (1989) 359--366.

\bibitem{Hornik:1991:ACM}
K.~Hornik, Approximation capabilities of multilayer perceptrons, Neural
  Networks 4 (1991) 251--257.

\bibitem{Strang:1973:AFE}
G.~Strang, G.~J. Fix, An Analysis of the Finite Element Method, Prentice--Hall,
  New York, NY, USA, 1973.

\bibitem{Rohrhofer:2021:PFP}
F.~M. Rohrhofer, S.~Posch, B.~C. Geiger, On the pareto front of
  physics-informed neural networks (2021).
\newblock \href {http://arxiv.org/abs/2105.00862} {\path{arXiv:2105.00862}}.

\bibitem{Hennigh:2021:NVI}
O.~Hennigh, S.~Narasimhan, M.~A. Nabian, A.~Subramaniam, K.~Tangsali, Z.~Fang,
  M.~Rietmann, W.~Byeon, S.~Choudhry, {NVIDIA SimNet™}: An {AI}-accelerated
  multi-physics simulation framework, in: International Conference on
  Computational Science, Springer, 2021, pp. 447--461.

\bibitem{Tsukanov:2011:HME}
I.~Tsukanov, S.~R. Posireddy, Hybrid method of engineering analysis:
  {Combining} meshfree method with distance fields and collocation technique,
  Journal of Computing and Information Science in Engineering 11~(3) (2011).

\bibitem{Sethian:1999:LSM}
J.~A. Sethian, {Level Set Methods and Fast Marching Methods: Evolving
  Interfaces in Computational Geometry, Fluid Mechanics, Computer Vision, and
  Materials Science}, Cambridge University Press, Cambridge, U.K., 1999.

\bibitem{Bloomenthal:1997:BEC}
J.~Bloomenthal, Bulge elimination in convolution surfaces, Computer Graphics
  Forum 16~(1) (1997) 31--41.

\bibitem{Shapiro:1999:IFG}
V.~Shapiro, I.~Tsukanov, Implicit functions with guaranteed differential
  properties, in: Proceedings of the Fifth {ACM} Symposium on Solid Modeling
  and Applications, 1999, pp. 258--269.

\bibitem{Upreti:2014:ADE}
K.~Upreti, T.~Song, A.~Tambat, G.~Subbarayan, Algebraic distance estimations
  for enriched isogeometric analysis, Computer Methods in Applied Mechanics and
  Engineering 280 (2014) 28--56.

\bibitem{Chin:2019:MCI}
E.~B. Chin, N.~Sukumar, Modeling curved interfaces without element-partitioning
  in the extended finite element method, International Journal for Numerical
  Methods in Engineering 120~(5) (2019) 607--649.

\bibitem{Belyaev:2017:TBC}
A.~G. Belyaev, P.-A. Fayolle, Transfinite barycentric coordinates, in: Hormann
  and Sukumar  \cite{Hormann:2017:GBC}, pp. 43--62.

\bibitem{Floater:2013:GBC}
M.~S. Floater, Generalized barycentric coordinates and applications, Acta
  Numerica 24 (2015) 161--214.

\bibitem{Anisimov:2017:BCP}
D.~Anisimov, Barycentric coordinates and their properties, in: Hormann and
  Sukumar  \cite{Hormann:2017:GBC}, pp. 3--22.

\bibitem{Hormann:2017:GBC}
K.~Hormann, N.~Sukumar (Eds.), Generalized Barycentric Coordinates in Computer
  Graphics and Computational Mechanics, CRC Press, New York, NY, 2017.

\bibitem{Hormann:2006:MVC}
K.~Hormann, M.~S. Floater, Mean value coordinates for arbitrary planar
  polygons, ACM Transactions on Graphics 25~(4) (2006) 1424--1441.

\bibitem{Bruvoll:2010:TMV}
S.~Bruvoll, M.~S. Floater, Transfinite mean value interpolation in general
  dimension, Journal of Computational and Applied Mathematics 233~(7) (2010)
  1631--1639.

\bibitem{Ju:2005:MVC}
T.~Ju, S.~Schaefer, J.~Warren, Mean value coordinates for closed triangular
  meshes, ACM Transactions on Graphics 24~(3) (2005) 561--566.

\bibitem{Chin:2021:SBC}
E.~B. Chin, N.~Sukumar, Scaled boundary cubature scheme for numerical
  integration over planar regions with affine and curved boundaries, Computer
  Methods in Applied Mechanics and Engineering 380 (2021) 113796.

\bibitem{Shepard:1968:ATD}
D.~Shepard, A two-dimensional interpolation function for irregularly-spaced
  data, in: Proceedings of the 23rd ACM national conference, Association for
  Computing Machinery, New York, New York, 1968, pp. 517--524.

\bibitem{rosenblatt1958perceptron}
F.~Rosenblatt, {The perceptron: A probabilistic model for information storage
  and organization in the brain.}, Psychological Review 65~(6) (1958) 386.

\bibitem{lecun1998mnist}
Y.~LeCun, {The MNIST database of handwritten digits},
  http://yann.lecun.com/exdb/mnist/ (1998).

\bibitem{finol2019deep}
D.~Finol, Y.~Lu, V.~Mahadevan, A.~Srivastava, Deep convolutional neural
  networks for eigenvalue problems in mechanics, International Journal for
  Numerical Methods in Engineering 118~(5) (2019) 258--275.

\bibitem{Hornik:1993:SNR}
K.~Hornik, Some new results on neural network approximation, Neural Networks 6
  (1993) 1069--1072.

\bibitem{LeNail:2019:NNS}
A.~{LeNail}, Nn-svg: Publication-ready neural network architecture schematics,
  Journal of Open Source Software 4~(33) (2019) 747.

\bibitem{kingma2014adam}
D.~P. Kingma, J.~Ba, Adam: {A} method for stochastic optimization (2014).
\newblock \href {http://arxiv.org/abs/1412.6980} {\path{arXiv:1412.6980}}.

\bibitem{Samaniego:2020:EAS}
E.~Samaniego, C.~Anitescu, S.~Goswami, V.~M. Nguyen-Thanh, H.~Guo, K.~Hamdia,
  X.~Zhuang, T.~Rabczuk, An energy approach to the solution of partial
  differential equations in computational mechanics via machine learning:
  Concepts, implementation and applications, Computer Methods in Applied
  Mechanics and Engineering 362 (2020) 112790.

\bibitem{jax2018github}
J.~Bradbury, R.~Frostig, P.~Hawkins, M.~J. Johnson, C.~Leary, D.~Maclaurin,
  G.~Necula, A.~Paszke, J.~Vander{P}las, S.~Wanderman-{M}ilne, Q.~Zhang, {JAX}:
  composable transformations of {P}ython+{N}um{P}y programs\ {URL}
  \url{http://github.com/google/jax} (2018).

\bibitem{bisong2019google}
E.~Bisong, Google colaboratory, in: Building Machine Learning and Deep Learning
  Models on Google Cloud Platform, Springer, 2019, pp. 59--64.

\bibitem{Schlomer:2020}
N.~Schl{\"o}mer, J.~Hariharan, dmsh, (2020). {A}vailable at
  \url{https://github.com/nschloe/dmsh}. Accessed on April 1, 2021.

\bibitem{Persson:2004:SMG}
P.-O. Persson, G.~Strang, A simple mesh generator in {MATLAB}, SIAM Review
  46~(2) (2004) 329--345.

\bibitem{Rahaman:2019:SBN}
N.~Rahaman, A.~Baratin, D.~Arpit, F.~Draxler, M.~Lin, F.~Hamprecht, Y.~Bengio,
  A.~Courville, On the spectral bias of neural networks, in: International
  Conference on Machine Learning, PMLR, 2019, pp. 5301--5310.

\bibitem{Wang:2021:EBF}
S.~Wang, H.~Wang, P.~Perdikaris, On the eigenvector bias of {Fourier} feature
  networks: From regression to solving multi-scale {PDEs} with physics-informed
  neural networks, Computer Methods in Applied Mechanics and Engineering 384
  (2021) 113938.

\bibitem{Joshi:2007:HCF}
J.~Pushkar, M.~Meyer, T.~{DeRose}, B.~Green, T.~Sanocki, Harmonic coordinates
  for character articulation, ACM Transactions on Graphics 26~(3) (2007)
  Article 71.

\bibitem{Shewchuk:1996:TRI}
J.~R. Shewchuk, Triangle: {E}ngineering a {2D} {Q}uality {M}esh {G}enerator and
  {D}elaunay {T}riangulator, in: M.~C. Lin, D.~Manocha (Eds.), Applied
  Computational Geometry: Towards Geometric Engineering, Vol. 1148 of Lecture
  Notes in Computer Science, Springer-Verlag, 1996, pp. 203--222.

\bibitem{Timoshenko:1959:TPS}
S.~P. Timoshenko, S.~Woinowsky-Krieger, Theory of Plates and Shells, 2nd
  Edition, McGraw Hill, New York, NY, USA, 1959.

\bibitem{Guo:2021:DCM}
H.~Guo, X.~Zhuang, T.~Rabczuk, A deep collocation method for the bending
  analysis of {Kirchhoff} plate (2021).
\newblock \href {http://arxiv.org/abs/2102.02617} {\path{arXiv:2102.02617}}.

\bibitem{Zhao:2005:FSM}
H.~Zhao, A fast sweeping method for {Eikonal} equations, Mathematics of
  Computation 74~(250) (2005) 603--627.

\bibitem{Cecil:2004:NMH}
T.~Cecil, J.~Qian, S.~Osher, Numerical methods for high dimensional
  {Hamilton--Jacobi} equations using radial basis functions, Journal of
  Computational Physics 196~(1) (2004) 327--347.

\bibitem{Ricci:1973:CGC}
A.~Ricci, A constructive geometry for computer graphics, The Computer Journal
  16~(2) (1973) 157--160.

\bibitem{Sherstyuk:1999:KFC}
A.~Sherstyuk, Kernel functions in convolution surfaces: a comparative analysis,
  The Visual Computer 15~(4) (1999) 171--182.

\bibitem{Barthe:2004:CBO}
L.~Barthe, B.~Wyvill, E.~{De Groot}, Controllable binary {CSG} operators for
  {``}soft objects{''}, International Journal of Shape Modelling 10~(02) (2004)
  135--154.

\bibitem{Gourmel:2013:GBI}
O.~Gourmel, L.~Barthe, M.-P. Cani, B.~Wyvill, A.~Bernhardt, M.~Paulin,
  H.~Grasberger, A gradient-based implicit blend, ACM Transactions on Graphics
  32~(2) (2013) 12:1--12:12.

\bibitem{Belyaev:2015:VPB}
A.~Belyaev, P.-A. Fayolle, On variational and {PDE}-based distance function
  approximations, Computer Graphics Forum 34~(8) (2015) 104--118.

\bibitem{Crane:2017:HMD}
K.~Crane, C.~Weischedel, M.~Wardetzky, The heat method for distance
  computation, Communications of the ACM 60~(11) (2017) 90--99.

\bibitem{Belyaev:2019:VMA}
A.~G. Belyaev, P.-A. Fayolle, A variational method for accurate distance
  function estimation, in: Numerical Geometry, Grid Generation and Scientific
  Computing, Springer, Cham, 2019, pp. 175--181.

\bibitem{Park:2019:SDF}
J.~J. Park, P.~Florence, J.~Straub, R.~Newcombe, S.~Lovegrove, Deepsdf:
  Learning continuous signed distance functions for shape representation, in:
  Proceedings of the IEEE/CVF Conference on Computer Vision and Pattern
  Recognition, 2019, pp. 165--174.

\end{thebibliography}
\end{document}